\documentclass[a4paper, 10pt, openany]{book}
\usepackage{graphicx,graphics}
\usepackage[cp1251]{inputenc}
\usepackage[russian,english]{babel}
\usepackage{calc}

\usepackage[cmtip,arrow,all,2cell]{xy}
\usepackage{pb-diagram,pb-xy}
\usepackage{amsmath,amsthm,amssymb,tabularx}
\usepackage[mathscr]{eucal}

\usepackage{makeidx}

\newcommand{\norm}[1]{\left\Vert#1\right\Vert}
\newcommand{\abs}[1]{\left\vert#1\right\vert}

\allowdisplaybreaks

\theoremstyle{plain}

\newtheorem{tm}{Theorem}[section]

\newtheorem{prop}{Proposition}[section]

\newtheorem{lm}{Lemma}[section]

\newtheorem{cor}{Corollary}[section]

\theoremstyle{definition}

\newtheorem{df}{Definition}[section]

\newtheorem{rem}{Remark}[section]

\newtheorem{ex}{Example}[section]

\newcommand{\beq}{\begin{equation}}
\newcommand{\eeq}{\end{equation}}
\newcommand{\bit}{\begin{itemize}}
\newcommand{\eit}{\end{itemize}}
\newcommand{\btm}{\begin{tm}}
\newcommand{\etm}{\end{tm}}
\newcommand{\blm}{\begin{lm}}
\newcommand{\bprop}{\begin{prop}}
\newcommand{\eprop}{\end{prop}}
\newcommand{\elm}{\end{lm}}
\newcommand{\bcor}{\begin{cor}}
\newcommand{\ecor}{\end{cor}}
\newcommand{\bex}{\begin{ex}}
\newcommand{\eex}{\end{ex}}
\newcommand{\bcx}{\begin{cex}}
\newcommand{\ecx}{\end{cex}}
\newcommand{\bers}{\begin{ers}}
\newcommand{\eers}{\end{ers}}
\newcommand{\bdf}{\begin{df}}
\newcommand{\edf}{\end{df}}
\newcommand{\brem}{\begin{rem}}
\newcommand{\erem}{\end{rem}}
\newcommand{\bpr}{\begin{proof}}
\newcommand{\epr}{\end{proof}}

\newcommand{\bextm}{\begin{extm}}
\newcommand{\eextm}{\end{extm}}

\def\le{\leqslant}
\def\ge{\geqslant}

\def \N {\mathbb{N}}
\def \Z {\mathbb{Z}}
\def \R {\mathbb{R}}
\def \T {\mathbb{T}}
\def \C {\mathbb{C}}

\def\Ker{\operatorname{\sf Ker}}
\def\Aut{\operatorname{\sf Aut}}
\def\id{\operatorname{\sf id}}
\def\Int{\operatorname{\sf Int}}
\def\d{\operatorname{\sf d}}
\def\M{\operatorname{\sf M}}
\def\GL{\operatorname{\sf GL}}
\def\supp{\operatorname{\sf supp}}

\def\cabsconv{\operatorname{\overline{\sf absconv}}}

\def\exp{\operatorname{\sf exp}}
\def\card{\operatorname{\sf card}}

\def\B{\operatorname{\sf B}}
\def\a{\operatorname{\sf{a}}}
\def\ll{\operatorname{\sf{l}}}
\def\rr{\operatorname{\sf{r}}}
\def\cc{\operatorname{\sf{c}}}
\def\d{\operatorname{\sf{d}}}
\def\ee{\operatorname{\sf{e}}}
\def\i{\operatorname{\sf{i}}}
\def\h{\operatorname{\sf{h}}}

\def\M{\operatorname{\sf{M}}}
\def\G{\operatorname{\sf{G}}}

\def \e {\varepsilon}
\def \ph {\varphi}

\def\BSQ{\resizebox{4pt}{5pt}{$\blacksquare$}}

\def\SQ{\resizebox{4pt}{5pt}{$\square$}}

\def\BLZ{\resizebox{5pt}{5pt}{$\blacklozenge$}}

\def\LZ{\resizebox{5pt}{5pt}{$\lozenge$}}

\def\boxedast{\raise.8pt\hbox{\rlap{$\mkern2.5mu *$}}\square}

\numberwithin{equation}{section}

\makeindex

\addtocounter{section}{-1}

\begin{document}

\title{\bf HOLOMORPHIC FUNCTIONS OF EXPONENTIAL TYPE AND DUALITY FOR STEIN GROUPS
WITH ALGEBRAIC CONNECTED COMPONENT OF IDENTITY}

\author{\bf S.~S.~Akbarov}

\maketitle

\def\theequation{\Alph{equation}}

\eject\addcontentsline{toc}{section}{Introduction}
\section*{Introduction}

Since in 1930s L.~S.~Pontryagin published his famous duality theorem for
Abelian locally compact groups \cite{Pontryagin}, the following problem engages
the imagination of specialists in harmonic analysis from time to time: {\it is
it possible to generalize Pontryagin duality to non-Abelian locally compact
groups in such a way that the dual object has the same nature as the initial
one?}

As is known, the first attempts to generalize Pontryagin duality did not meet
this requirement: in the M.~G.~Krein theory, for instance, the dual object
$\widehat{G}$ for a group $G$ is a block-algebra \cite{Krein} (but not a group,
unlike Pontryagin theory). Apparently, a deep peculiarity in human psychology
manifests itself here, but such a harmless trait like asymmetry between $G$ and
$\widehat{G}$ in the theory of representations -- a trait that can be compared
with difference between the left and the right in anatomy -- leads to numerous
and, because of the changing with time understanding of what the notion of
group should mean, continuing attempts to build a duality theory, where, on the
one hand, ``all'' groups are covered, and, on the other, the Pontryagin
symmetry between initial objects and their duals is preserved.

In the category theory language, the unique one fitting for such speculative
aspirations, one can formulate this task correctly using the following two
definitions.

\bit

\item[1.] Let us call a contravariant functor $A\mapsto A^*:\mathfrak{K}\to
\mathfrak{K}$ on a given category $\mathfrak{K}$ a {\it duality functor
in}\index{duality functor} $\mathfrak{K}$, if its square, i.e. the covariant
functor $A\mapsto (A^*)^*:\mathfrak{K}\to \mathfrak{K}$, is isomorphic to the
identity functor $\id_{\mathfrak{K}}:\mathfrak{K}\to \mathfrak{K}$.
 \beq\label{functor-dvoistvennosti}
 \xymatrix @R=1.pc @C=1.pc
 {
 & \mathfrak{K} \ar[rd]^{*} & \\
 \mathfrak{K} \ar[ru]^{*} \ar[rr]_{\id_{\mathfrak{K}}} & & \mathfrak{K}
 }
 \eeq
 \eit
The passage to the dual group $G\mapsto G^\bullet$ in Pontryagin theory is an
example of duality functor: the natural isomorphism between
$G^{\bullet\bullet}$ and $G$ here is the mapping
$$
i_G:G\to G^{\bullet\bullet}\quad \Big| \quad i_G(x)(\chi)=\chi(x),\qquad x\in
G,\ \chi\in G^\bullet
$$
On the contrary, say, in the category of Banach spaces the passage to the dual
Banach space $X\mapsto X^*$ is not a duality functor (because there are
non-reflexive Banach spaces).
 \bit
\item[2.] Suppose we have:
 \bit

\item[(a)] three categories $\mathfrak{K}$, $\mathfrak{L}$, $\mathfrak{M}$ with
two full and faithful covariant functors $A:\mathfrak{K}\to \mathfrak{L}$ and
$B:\mathfrak{L}\to \mathfrak{M}$ defining a chain of embeddings:
$$
\mathfrak{K}\subset \mathfrak{L}\subset \mathfrak{M},
$$

\item[(b)] two duality functors $K\mapsto K^\bullet:\mathfrak{K}\to
\mathfrak{K}$ and $M\mapsto M^*:\mathfrak{M}\to \mathfrak{M}$ such that the
functors $K\mapsto B(A(K^\bullet))$ and $K\mapsto \Big(B(A(K))\Big)^*$ are
isomorphic:
 \beq\label{diagr-K-L-M}
 \xymatrix @R=1.pc @C=1.pc
 {
 \mathfrak{M} \ar[r]^{*} & \mathfrak{M} \\
 \mathfrak{L} \ar[u]^B & \mathfrak{L} \ar[u]_B \\
 \mathfrak{K} \ar[u]^A \ar[r]^{\bullet} & \mathfrak{K} \ar[u]_A
 }
 \eeq
 \eit

We shall call this construction a {\it generalization of the duality $\bullet$
from the category $\mathfrak{K}$ to the category
$\mathfrak{L}$}\index{generalization of duality}.
 \eit

On adding these terms into our armoury, we can formulate the task we are
discussing as follows: {\it are there any generalizations of Pontryagin duality
from the category of Abelian locally compact groups to the category of
arbitrary locally compact groups, and if yes, then which namely?} The category
diagram \eqref{diagr-K-L-M} in this formulation becomes as follows:
 {\sf
 \beq
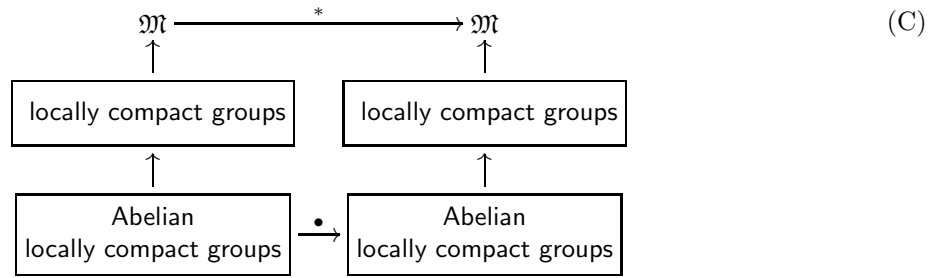
\label{diagramma-obobsheniya-dvoistv-pontryagina}
 \xymatrix @R=1.pc @C=1.pc
 {
 \mathfrak{M} \ar[r]^{*} & \mathfrak{M} \\
 \boxed{\phantom{\Big|} \text{locally compact groups}}
 \ar[u] &
 \boxed{\phantom{\Big|} \text{locally compact groups}}
 \ar[u] \\
 \boxed{\begin{matrix}
 \text{Abelian}\\
 \text{locally compact groups}\end{matrix}} \ar[u] \ar[r]^{\bullet} &
  \boxed{\begin{matrix}
 \text{Abelian}\\
 \text{locally compact groups}
 \end{matrix}}\ar[u]
 }
 \eeq}\noindent
-- and the key example here is the duality theory for finite groups, which can
be regarded as a generalization of Pontryagin duality from the category of
Abelian finite groups to the category of all finite groups:
 {\sf
 \beq
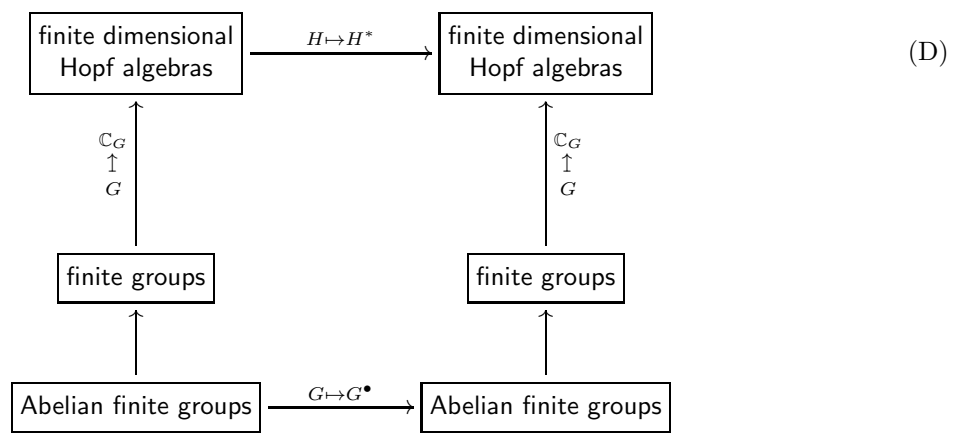
\label{diagramma-kategorij-dlya-konechnyh-grupp}
 \xymatrix
 {
 \boxed{\begin{matrix}
 \text{finite dimensional}\\
 \text{Hopf algebras}
 \end{matrix}}
 \ar[rr]^{H\mapsto H^*} & &
\boxed{\begin{matrix}
 \text{finite dimensional}\\
 \text{Hopf algebras}
 \end{matrix}}
 \\ & & \\
 \boxed{\begin{matrix}
  \text{finite groups}
 \end{matrix}} \ar[uu]^{\scriptsize\begin{matrix} \C_G\\
 \text{\rotatebox{90}{$\mapsto$}} \\ G\end{matrix}} & &
 \boxed{\begin{matrix}
  \text{finite groups}
 \end{matrix}} \ar[uu]_{\scriptsize\begin{matrix} \C_G \\
 \text{\rotatebox{90}{$\mapsto$}} \\ G\end{matrix}} \\
 \boxed{\begin{matrix}
 \text{Abelian finite groups}
 \end{matrix}} \ar[u] \ar[rr]^{G\mapsto G^\bullet} & &
  \boxed{\begin{matrix}
 \text{Abelian finite groups}
  \end{matrix}}\ar[u]
 }
 \eeq }\noindent
(here $G\mapsto \C_G$ is the passage to group algebra, and $H\mapsto H^*$ the
passage to dual Hopf algebra). This example, apart from everything else,
illustrates another guiding idea of the ``general duality theory'': if we want
to reduce the representations of groups to representations of algebras, and
therefore claim that the category $\mathfrak{M}$ consists of associative
algebras, then these algebras $H$ must have some supplementary structure, which
allows us to endow dual objects $H^*$ with a natural structure of associative
algebra. Natural objects of that kind in general algebra are Hopf algebras. As
a corollary, the constructions in ``duality theory'' usually resemble Hopf
algebras, although as a rule differ from them, except for trivial cases when,
for instance, the algebra has finite dimension.

A generalization of Pontryagin duality as it is presented in diagram
\eqref{diagramma-obobsheniya-dvoistv-pontryagina} was suggested in 1973
independently by L.~I.~Vai\-ner\-man and G.~I.~Kac from one side
(\cite{Vainerman}, \cite{Vainerman-Kac-1}, \cite{Vainerman-Kac-2}), and by
M.~Enock and J.-M.~Schwartz from another (\cite{Enock-Schwartz-1},
\cite{Enock-Schwartz-2}, \cite{Enock-Schwartz-3}). The theory of Kac algebras
they developed summarized a series of attempts made by different
mathematicians, history of which, as well as the theory of Kac algebras itself,
one can learn from monograph \cite{Enock-Schwartz} by Enock and Schwartz. After
1973 the work in this direction did not cease, because on the one hand, some
improvements were added into the theory (again see details in
\cite{Enock-Schwartz}), and on the other, after the discovery of quantum groups
in 1980s, mathematicians began to generalize the Pontryagin duality to this
class as well. Moreover, the latter work is not finished by now -- the theory of
locally compact quantum groups which appeared on this wave is being actively
developed by J.~Kustermans, W.~Pusz, P.~M.~So\l tan, S.~Vaes, L.~Vainerman,
A.~Van Daele, S.~L.~Woronowicz (an impression of this topic one can find in
collective monograph \cite{QSNG}). The idea of multiplier Hopf algebra
suggested by A.~Van Daele in 1990-s \cite{VanDaele-1,VanDaele-2}, seems to be
strategic in these investigations.

Despite the active work and impressive enthusiasm demonstrated by
mathematicians engaged in this theme, the theories they suggest have a serious
shortcoming: {\it all the enveloping categories $\mathfrak{M}$ in these
theories consist of objects, which are formally not Hopf algebras.} The Kac
algebras, for instance, although being chosen as a subclass among the objects
called Hopf-von~Neumann algebras, in fact are not Hopf algebras, because from
the point of view of category theory in the definition of Hopf-von~Neumann
algebras, unlike the pure algebraic situation, two different tensor products
are used simultaneously -- the projective tensor product for multiplication,
and the tensor product of von Neumann algebras for comultiplication. In the
theory of locally compact quantum groups the situation on this point seems to give even
less hopes, because here the claim that the comultiplication must indispensably
act into tensor product is rejected: in accordance with the above mentioned Van
Daele's idea of multiplier Hopf algebras the comultiplication here is defined
as an operator from $A$ into the algebra $M(A\otimes A)$ of multipliers on
tensor product $A\otimes A$ (which is chosen here as minimal tensor product of
$C^*$-algebras, see \cite{Kustermans-Vaes,QSNG}).

In this paper we suggest another approach to the generalization of Pontryagin
duality, where this unpleasant effect does not occur: in our theory the
enveloping category $\mathfrak{M}$ consists of ``true'' Hopf algebras, of
course in the categorical sense, i.e. the Hopf algebras defined in the same way as
the usual Hopf algebras, but after replacing the category of vector spaces by a
given symmetrical (in the more general case -- braided) monoidal category (such
Hopf algebras are sometimes called Hopf monoids, see
\cite{Rivano,Schauenburg,Street,Zhang}).

The mission we set to ourselves is not quite a generalization of Pontryagin
duality to the class of locally compact groups, as it is presented in diagram
\eqref{diagramma-obobsheniya-dvoistv-pontryagina}, but a solution of the same
problem in the class of complex Lie groups. The fact is that among the four
main branches of mathematics, where the idea of invariant integration manifests
itself --
 \bit
\item[--] general topology (where the groups with invariant integral are
exactly locally compact groups),

\item[--] differential geometry (where this role belongs to Lie groups),

\item[--] complex analysis (where the reductive complex groups can be regarded
as the groups with integral), and

\item[--] algebraic geometry (again with reductive complex groups),
 \eit
-- at least in the three first disciplines the generalization of Pontryagin
duality has sense. In topology this problem takes form of diagram
\eqref{diagramma-obobsheniya-dvoistv-pontryagina}, in differential geometry one
can consider the problem of generalization of Pontryagin duality from Abelian
compactly generated Lie groups, say, to all compactly generated Lie groups, and in
the complex analysis one can consider the problem of generalization of
Pontryagin duality from the class of Abelian compactly generated Stein groups
to the class of reductive complex groups (following \cite{Vinberg}, by a
reductive group we mean complexification of a compact real Lie group).

We give a solution for the third of these problems. We begin our considerations
with the algebra ${\mathcal O}(G)$ of holomorphic functions on a complex Lie
group $G$, since this is a natural functional algebra in complex analysis (in
the other three disciplines this role belongs to the algebra ${\mathcal C}(G)$ of
continuous functions, the algebra ${\mathcal E}(G)$ of smooth functions, and
the algebra ${\mathcal R}(G)$ of polynomials). The idea we suggest as a heuristic
hypothesis in this paper -- and its justification we see in our work -- is as
follows. It seems likely, that to each of the three first disciplines in our
list -- general topology, differential geometry and complex analysis -- from
the point of view of the topological algebras used in these disciplines,
corresponds some class of seminorms, intrinsically connected to this class of
algebras. What those seminorms should be in general topology and in
differential geometry we cannot say right now, but the class of
submultiplicative seminorms, i.e. seminorms defined by the inequality
$$
p(x\cdot y)\le p(x)\cdot p(y),
$$
is intrinsically connected with complex analysis (this idea is inspired to us
by A.~Yu.~Pirkovskii's results on the Arens-Michael envelopes of topological
algebras -- see \cite{Pirkovskii}).

To give exact meaning to the words ``intrinsically connected'', we need to
introduce the following two functors and one category of Hopf algebras:
 \bit
\item[1)] {\it The Arens-Michael envelope $ A\mapsto A^\heartsuit $.} The
algebras whose topology is generated by submultiplicative seminorms and
satisfies supplementary condition of completeness, are called Arens-Michael
algebras (we discuss them in \ref{SEC:Arens-Michael}). Each topological algebra
$A$ has a nearest ``from outside'' Arens-Michael algebra -- we denote it by
$A^\heartsuit$ -- called Arens-Michael envelope of $A$. This algebra is a
completion of $A$ with respect to the system of all continuous
submultiplicative seminorms on it.

\item[2)] {\it Passage to dual stereotype Hopf algebra $H\mapsto H^\star $.} In
the work \cite{Akbarov} we discussed in detail the symmetrical monoidal
categories $(\mathfrak{Ste},\odot)$ and $(\mathfrak{Ste},\circledast)$ of
stereotype spaces (with injective $\odot$ and projective $\circledast$ tensor
products). The Hopf algebras in these categories are called respectively
injective and projective stereotype Hopf algebras. The passage to the dual
stereotype space $H\mapsto H^\star$ establishes an antiequivalence between
these categories of Hopf algebras.

\item[3)] {\it Category of holomorphically reflexive Hopf algebras.} The
functors $\heartsuit$ and $\star$ allow us to consider the category of
projective and at the same time injective stereotype Hopf algebras $H$, for
which the successive application of the operations $\heartsuit$ and $\star$ always
leads to injective and at the same time projective Hopf algebras and, if we
begin with $\heartsuit$, then at the fourth step this chain leads back to the
initial Hopf algebra (of course up to an isomorphism). We depict this closed
chain by the following {\it reflexivity diagram},
 \beq\label{int:diagr-refl}
 \xymatrix @R=2.pc @C=4.pc @M=14pt
 {
 H
 \ar@{|->}[r]^{\heartsuit} &
 H^{\heartsuit}
 \ar@{|->}[d]^{\star}
 \\
 H^\star\ar@{|->}[u]^{\star}
 &
 \ar@{|->}[l]_{\heartsuit}
 (H^{\heartsuit})^\star
 }
 \eeq
and call such Hopf algebras $H$ {\it holomorphically reflexive} (see accurate
definition in \ref{SEC:Arens-Michael}\ref{holom-reflexivnost}). The duality
functor in the category of such Hopf algebras is, certainly, the operation
$H\mapsto (H^\heartsuit)^\star$.
 \eit

The main result of our work, which clarifies our hints about ``seminorms,
intrinsically connected with complex analysis'', is that the algebra ${\mathcal
O}^\star(G)$ of analytical functionals on a compactly generated Stein group $G$
with algebraic connected component of identity is a holomorphically reflexive
Hopf algebra. To prove this we introduce in this paper the algebra ${\mathcal
O}_{\exp}(G)$ of holomorphic functions of exponential type on group $G$, which
is a subalgebra in the algebra ${\mathcal O}(G)$ of all holomorphic functions
on $G$. The reflexivity diagram \eqref{int:diagr-refl} for ${\mathcal
O}^\star(G)$ is as follows:
 \beq\label{int:chetyrehugolnik-O-O*}
 \xymatrix @R=2.pc @C=4.pc @M=14pt
 {
 {\mathcal O}^\star(G)
  \ar@{|->}[r]^{\heartsuit} &
 {\mathcal O}_{\exp}^\star(G)
 \ar@{|->}[d]^{\star}
 \\
 {\mathcal O}(G) \ar@{|->}[u]^{\star}
 &
 {\mathcal O}_{\exp}(G) \ar@{|->}[l]_{\heartsuit}
 }
 \eeq
For the case of Abelian groups the operation $\heartsuit$ becomes naturally
isomorphic to the usual Fourier transform, so this diagram takes the form:
 \beq\label{int:chetyrehugolnik-O-O*-Abel}
 \xymatrix @R=2.pc @C=4.pc @M=14pt
 {
 {\mathcal O}^\star(G)
 \ar@{|->}[r]^{\heartsuit}_{\scriptsize\begin{matrix}\text{Fourier}\\ \text{transform}\end{matrix}}
 &
 {\mathcal O}(G^\bullet)
 \ar@{|->}[d]^{\star}
 \\
 {\mathcal O}(G)
 \ar@{|->}[u]^{\star}
 &
 \ar@{|->}[l]_{\heartsuit}^{\scriptsize\begin{matrix}\text{Fourier}\\ \text{transform}\end{matrix}}
  {\mathcal O}^\star(G^\bullet)
 }
 \eeq
(the dual group $G^\bullet$ for an Abelian compactly generated Stein group $G$
is defined as the group of all homomorphisms from $G$ into the multiplicative
group $\C^\times:=\C\setminus\{0\}$ of non-zero complex numbers). This, in
particular, implies the isomorphism of functors
$$
{\mathcal O}^\star(G^\bullet)\cong \left(\Big({\mathcal
O}^\star(G)\Big)^\heartsuit\right)^\star,
$$
which gives a generalization of Pontryagin duality in the complex case, and the
diagram \eqref{diagr-K-L-M} here takes the form:
 {\sf
 \beq\label{diagramma-kategorij-dlya-grupp-steina}
 \xymatrix @R=1.pc @C=1.pc
 {
 \boxed{\begin{matrix}
 \text{holomorphically reflexive}\\
 \text{Hopf algebras}
 \end{matrix}}
 \ar[rr]^{H\mapsto (H^\heartsuit)^\star} & &
 \boxed{\begin{matrix}
 \text{holomorphically reflexive}\\
 \text{Hopf algebras}
 \end{matrix}}
 \\ & & \\
 \boxed{\begin{matrix}
  \text{compactly generated} \\ \text{Stein groups} \\ \text{with algebraic} \\ \text{component of identity}
 \end{matrix}} \ar[uu]^{\scriptsize\begin{matrix} {\mathcal O}^\star(G)\\
 \text{\rotatebox{90}{$\mapsto$}} \\ G\end{matrix}} & &
 \boxed{\begin{matrix}
  \text{compactly generated} \\ \text{Stein groups} \\ \text{with algebraic} \\ \text{component of identity}
 \end{matrix}} \ar[uu]_{\scriptsize\begin{matrix} {\mathcal O}^\star(G)\\
 \text{\rotatebox{90}{$\mapsto$}} \\ G\end{matrix}} \\
 \boxed{\begin{matrix}
  \text{Abelian compactly generated} \\ \text{Stein groups}
 \end{matrix}} \ar[u] \ar[rr]^{G\mapsto G^\bullet} & &
  \boxed{\begin{matrix}
  \text{Abelian compactly generated} \\ \text{Stein groups}
  \end{matrix}}\ar[u]
 }
 \eeq
 }\noindent
Since every reductive group is algebraic, this indeed will be a solution of the
third problem in our list.

We show in addition that the holomorphic duality we introduce here does not
limit itself to the class of compactly generated Stein groups with algebraic
connected component of identity, but extends to quantum groups. As an example
we consider the quantum group `$az+b$' of quantum affine automorphisms of
complex plane (see \cite{Woronowicz,VanDaele,Wang,QSNG}). We prove that
`$az+b$' is a holomorphically reflexive Hopf algebra in the sense of our
definition.

The author thanks sincerely D.~N.~Akhieser, O.~Yu.~Aristov, P.~Gaucher,
A.~Ya.~Helemskii, A.~Huckleberry, E.~B.~Katsov, Yu.~N.~Kuznetsova, T.~Maszczyk,
S.~Yu.~Nemirovskii, A.~Yu.~Pirkovskii, V.~L.~Popov, P.~So\l tan, A.~Van Daele
for innumerous consultations and help during the work on this paper. Besides
that the idea of proof of Propositions
\ref{PROP-stroenie-poluharakterov-na-GL_n}, \ref{PROP:p-na-R_q(C-times-C^x)}
and Lemma \ref{LM-nepr-subm-polunorm-na-R(C)-psi-odot-O^star(C^times)} belongs
to Yu.~N.~Kuznetsova.

\def\theequation{\arabic{section}.\arabic{equation}}

\section{Stereotype spaces}\label{SEC-ster-spaces}

Stereotype spaces we are speaking about in this section were studied by the
author in detail in \cite{Akbarov} (see also \cite{Akbarov-2,Akbarov-3}).

\subsection{Definition and typical examples}

Let $X$ be a locally convex space over $\Bbb C$. Denote by $X^\star$ the space
of all linear continuous functionals $f:X\to \Bbb C$, endowed with the topology
of uniform convergence on totally bounded sets in $X$. The space $X$ is called
{\it stereotype}\index{stereotype!space}, if the natural mapping
$$
\i_X:X\to (X^\star)^\star \quad | \quad \i_X(x)(f)=f(x), \quad x\in X, f\in
X^\star
$$
is an isomorphism of locally convex spaces. Clearly the following theorem
holds:

\btm\label{TH:X-ste=>X*-ste} If $X$ is a stereotype space, then $X^\star$ is
also a stereotype space. \etm

It turns out that stereotype spaces form a very wide class, what can be illustrated by the following diagram:
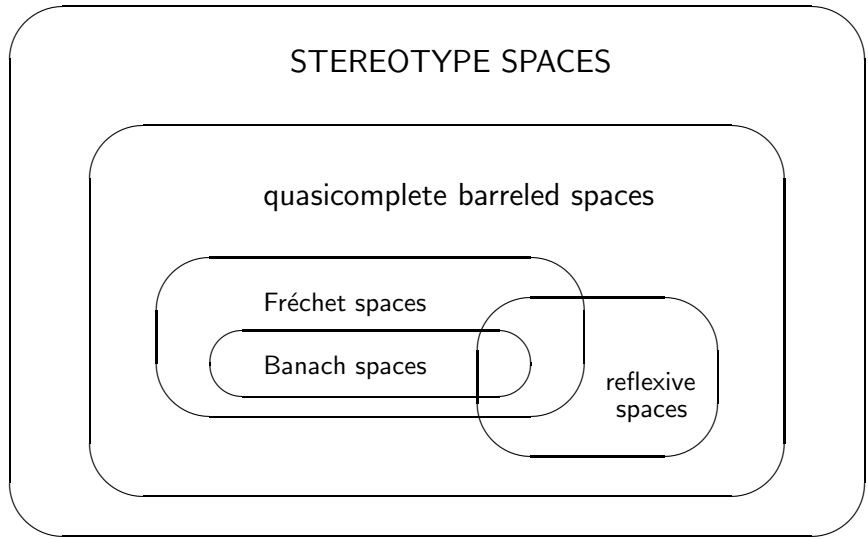
\vglue-20pt
$$
 \begin{picture}(380,220)
\put(195,105){\oval(320,200)} \put(140,180){\text{\sf\large STEREOTYPE SPACES}}
\put(195,90){\oval(260,140)} \put(130,130){\text{\sf\large quasicomplete
barreled spaces}}
 \put(170,80){\oval(160,60)} \put(130,90){\text{\sf Fr\'echet spaces}}
 \put(170,70){\oval(120,25)} \put(130,66){\text{\sf Banach spaces}}
 \put(255,65){\oval(90,60)} \put(258,60){\text{\sf reflexive}}
 \put(262,50){\text{\sf spaces}}
\end{picture}
$$
Fr\'echet spaces and Banach spaces will be of special interest for us in this picture, so we shall consider them in detail.

\bex[\bf Fr\'echet spaces and Brauner spaces] Every Fr\'echet space $X$ is
stereotype \cite{Brauner}. Its dual space $Y=X^\star$ is also stereotype by theorem
\ref{TH:X-ste=>X*-ste}. If $\{U_n\}$ is a countable local base in $X$, then the
polars $K_n=U^\circ_n$ form a countable {\it fundamental system of compact
sets}\index{fundamental system!of compact sets} in $Y$: every compact set
$T\subseteq Y$ is contained in some compact set $K_n$ (this means by the way
that $Y$ cannot be Fr\'echet space, if $X$ infinite dimensional). The spaces
$Y$ dual to Fr\'echet spaces $X$ (in the sense of our definition) were
originally considered by K.~Brauner in \cite{Brauner}, and we call them {\it
Brauner spaces}\index{Brauner space}. Their characteristic properties are
listed in the following proposition. \eex

\bprop\label{Brauner-crit} For a locally convex space $Y$ the following
conditions are equivalent:
\begin{itemize}
\item[(i)] $Y$ is a Brauner space;

\item[(ii)] $Y$ is a complete Kelley space (i.e. every set $M\subseteq Y$ that has
a closed intersection $M\cap K$ with any compact set $K\subseteq Y$, is closed
in $Y$) and has a countable fundamental system of compact sets $K_n$: for each
compact set $T\subset Y$ there exists $n\in\N$ \footnote{Everywhere in our
paper $\N$ means the set of non-negative integers: $\N:=\{0,1,2,3,...\}$.} such
that $T\subseteq K_n$;

\item[(iii)] $Y$ is a stereotype space and has a countable fundamental system
of compact sets $K_n$: for each compact set $T\subset Y$ there exists $n\in\N$
such that $T\subseteq K_n$;

\item[(iv)] $Y$ is a stereotype space and has a countable exhausting system of
compact sets $K_n$: $\bigcup_{n=1}^{\infty} K_n=Y$.
\end{itemize}
\eprop
 \bpr
The countable fundamental system of compact sets in $Y$ is the system of polars
$K_n={^\circ U_n}$ of a local base $U_n$ in the dual Fr\'echet space $Y^\star$.
Modulo this remark all the statements in proposition \ref{Brauner-crit} are
obvious, except one -- that a Brauner space $Y$ is always a Kelley space. This
result belongs to K.~Brauner and is deduced in his paper \cite{Brauner} from the Banach-Dieudonn\'e theorem (see \cite{Jarchow}).\epr

\bcor\label{nepr-na-Brauner} If $Y$ is a Brauner space with the fundamental
system of compact sets $K_n$, then any linear mapping $\ph:Y\to Z$ into a
locally convex space $Z$ is continuous if and only if it is continuous on each
compact set $K_n$. \ecor

\bex[\bf Banach spaces and Smith spaces] These are special cases of Fr\'echet
spaces and Brauner spaces. If $X$ is a Banach space, then $X$ and $Y=X^\star$
are stereotype spaces \cite{Smith}. The polar $K=B^\circ$ of a ball $B$ in $Y$ is a {\it
universal compact set}\index{universal compact set} in $Y$, i.e. a compact set
that swallows any other compact set $T$ in $Y$. The spaces $Y=X^\star$ dual to
Banach spaces $X$ (in the sense of our definition) were originally considered
by M.~F.~Smith in \cite{Smith} -- that is why we call them {\it Smith
spaces}\index{Smith space}. Their characteristic properties are listed in the
following proposition: \eex

\bprop\label{Smith-crit} For a locally convex space $Y$ the following
conditions are equivalent:
\begin{itemize}
\item[(i)] $Y$ is a Smith space;

\item[(ii)] $Y$ is a complete Kelley space and has a universal compact set $K$:
for any compact set $T\subset X$ there exists $\lambda\in \Bbb C$ such that
$T\subseteq \lambda K$;

\item[(iii)] $Y$ is a stereotype space with a universal compact set;

\item[(iv)] $Y$ is a stereotype space with a compact barrel.
\end{itemize}
\eprop

\bcor\label{nepr-na-Smith} If $Y$ is a Smith space with a universal compact set
$K$, then a linear mapping $\ph:Y\to Z$ into a locally compact space $Z$ is
continuous if and only if it is continuous on $K$. \ecor

The connections between the spaces of Fr\'echet, Brauner, Banach and Smith are
illustrated in the following diagram (where turnover corresponds to the passage to
the dual class):
$$
 
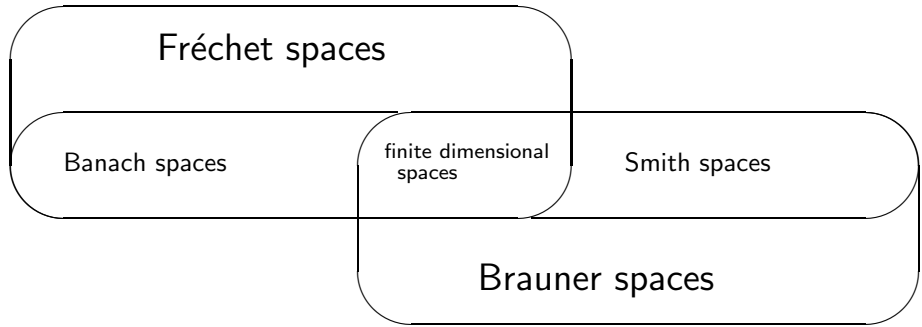
\begin{picture}(400,140)
 \put(130,90){\oval(210,80)} \put(80,110){\text{\sf\Large Fr\'echet spaces}}
 \put(170,70){\oval(290,40)[l]} \put(45,68){\text{\sf Banach spaces}}
 \put(260,50){\oval(210,80)} \put(200,23){\text{\sf\Large Brauner spaces}}
 \put(220,70){\oval(290,40)[r]} \put(255,68){\text{\sf Smith spaces}}
 \put(165,73){\text{\sf\footnotesize finite dimensional}}
 \put(170,65){\text{\sf\footnotesize spaces}}
\end{picture}
$$

\subsection{Smith space generated by a compact set}

Let $X$ be a stereotype space and $K$ an absolutely convex compact set in $X$.
We denote by $\C K$ the linear subspace in $X$, generated by the set $K$:
 \beq
\C K=\bigcup_{\lambda>0}\lambda K
 \eeq
Endow $\C K$ by the {\it Kelley topology\index{Kelley topology} generated by
compact sets $\lambda K$}: a set $M\subseteq \C K$ is considered closed in $\C
K$, if its intersection $M\cap \lambda K$ with any compact set $\lambda K$ is
closed in $X$ (or, equivalently, in $\lambda K$).

\btm\label{CK-smith} The Kelley topology on $\C K$, generated by compact sets
$\lambda K$, is a unique topology on $\C K$, which turns $\C K$ into a Smith
space with the universal compact set $K$.
 \etm
\bpr Denote by $(\C K)'$ the set of all {\it linear} functionals on $\C K$,
continuous on the compact set $K$:
$$
f\in(\C K)'\quad\Longleftrightarrow\quad f:\C K\to\C\quad\&\quad f|_K\in C(K)
$$
Clearly, $(\C K)'$ is a Banach space with respect to the norm
$$
||f||=\max_{x\in K}|f(x)|
$$
(formally this turns $(\C K)'$ into a closed subspace in $C(K)$). Note that
functionals $f\in (\C K)'$ separate points in $K$, because those of them who
are restrictions on $\C K$ of functionals from $g\in X^\star$ already possess
this property:
$$
\forall x,y\in K\qquad \Big(x\ne y\quad\Longrightarrow\quad \exists g\in
X^\star\quad g(x)\ne g(y)\Big)
$$
This means that the weak topology $\sigma$ on $K$, generated by functionals
$f\in (\C K)'$, coincides with the initial topology $\tau$ of this compact set
(because $\sigma$ is Hausdorff and is majorized by $\tau$). This implies in its
turn that the topology of the space $(\C K)'$ is the topology of uniform
convergence on $(\C K)'$-weak absolutely convex compact sets of the form
$\lambda K$ (these sets form a saturated system in $\C K$). Hence by the
Mackey-Arens theorem \cite{Schaefer}, the system $((\C K)')^\star$ of linear
continuous functionals on $(\C K)'$ coincides with $\C K$:
 $$
\C K=((\C K)')^\star
 $$

We see that $\C K$ can be identified with the space of linear continuous
functionals on $(\C K)'$. As a corollary, $\C K$ can be endowed with the
topology of dual space (in stereotype sense) to the Banach space $(\C K)'$:
 \beq
\C K\cong ((\C K)')^\star
 \eeq
This topology turns $\C K$ into a Smith space, and, by Proposition
\ref{Smith-crit}, it coincides with the Kelley topology on $\C K$, generated by
compact sets $\lambda K$.

Let us denote this topology on $\C K$ by $\varkappa$, and show that it is a
unique topology under which $\C K$ is a Smith space with the universal compact
set $K$. Indeed, if $\rho$ is another topology on $\C K$ with the same
property, then the identity mapping
$$
(\C K)_\rho\to (\C K)_\varkappa
$$
is continuous on $K$ (since it preserves the topology on $K$), hence by
Corollary \ref{nepr-na-Smith}, it is a continuous mapping of Smith spaces. In
the same way, the inverse mapping
$$
(\C K)_\varkappa\to (\C K)_\rho
$$
is continuous, and this means that the topologies $\varkappa$ and $\rho$
coincide.
 \epr

\bcor The topology of the space $\C K$ can be equivalently described as the
topology of uniform convergence on sequences of functionals $\{f_k\}\subset (\C
K)'$ tending to zero, i.e. as the topology generated by seminorms of the form:
$$
p_{\{f_k\}}(x)=\sup_{k\in\N}|f_n(x)|
$$
where $f_k$ is a sequence of linear functionals on $\C K$, continuous on $K$,
and such that $\max_{t\in K}|f_k(t)|\underset{k\to\infty}{\longrightarrow}
0$.\ecor

\bprop\label{vlozh-Smith} For any two absolutely convex compact sets
$K,L\subseteq X$ such that
$$
K\subseteq L,
$$
the natural imbedding of the Smith spaces which they generate
$$
\iota_K^L: \C K\to \C L
$$
is a continuous mapping.\eprop
 \bpr
The mapping $\iota_K^L$ is continuous on the compact set $K$, hence, by
Corollary \ref{nepr-na-Smith}, it is continuous on all $\C K$.\epr

\subsection{Brauner spaces generated by an expanding sequence of compact sets}

A sequence of absolutely convex compact sets $K_n$ in a stereotype space $X$
will be called {\it expanding}\index{expanding system of sets}, if
$$
\forall n\in\N\qquad  K_n+K_n\subseteq K_{n+1}
$$
For every such sequence the set
$$
\bigcup_{n=1}^\infty \C K_n=\bigcup_{n\in \N, \lambda\in\C} \lambda
K_n=\bigcup_{n\in \N} K_n
$$
is a subspace in the vector space $X$. We endow it with the {\it Kelley
topology generated by compact sets $K_n$}: a set $M\subseteq
\bigcup_{n=1}^\infty K_n$ is considered closed in $\bigcup_{n=1}^\infty K_n$,
if its intersection $M\cap K_n$ with any compact set $K_n$ is closed in $X$
(and equivalently, in $K_n$). The following proposition is proved similarly
with Theorem \ref{CK-smith}:

\btm The Kelley topology on $\bigcup_{n=1}^\infty K_n$, generated by compact
sets $K_n$, is a unique topology on $\bigcup_{n=1}^\infty K_n$, under which
this space is a Brauner space with the fundamental sequence of compact sets
$\{K_n\}$. \etm

\bcor The topology of the space $\bigcup_{n=1}^\infty K_n$ can be equivalently
described as the topology of uniform convergence on sequences of functionals
$\{f_k\}\subset \left(\bigcup_{n=1}^\infty K_n\right)'$ tending to zero, i.e.
as the topology generated by seminorms of the form:
$$
p_{\{f_k\}}(x)=\sup_{k\in\N}|f_k(x)|
$$
where $f_k$ is an arbitrary sequence of linear functionals on
$\bigcup_{n=1}^\infty K_n$, continuous on each $K_n$, such that $\forall n$
$\max_{t\in K_n}|f_k(t)|\underset{k\to\infty}{\longrightarrow} 0$.\ecor

\subsection{Projective Banach systems and injective Smith systems}

Let $X$ be a locally convex space. A standard construction in the theory of
topological vector spaces assigns to each absolutely convex neighborhood of
zero $U$ in $X$ a Banach space, which it is convenient to denote by $X/U$ and
to call a {\it quotient space of $X$ over the neighborhood of zero
$U$}\index{quotient space!over a neighborhood of zero}. It is defined as
follows. First, we consider the set
$$
\Ker U=\bigcap_{\varepsilon>0}\varepsilon\cdot U,
$$
which is called {\it kernel of the neighborhood of zero}\index{kernel of a
neighborhood of zero} $U$ -- this is a closed subspace in $X$, since $U$ is
absolutely convex. Then we construct the quotient space $X/\Ker U$, and endow
it with the topology of normed space with $U+\Ker U$ as unit ball (this
topology in general is weaker than the usual topology of quotient space on
$X/\Ker U$). This normed space $X/\Ker U$ is usually not complete. Its
completion is declared the final result:
 \beq\label{X/U}
X/U:=(X/\Ker U)^\blacktriangledown
 \eeq
(here $\blacktriangledown$ means completion).

\btm\label{TH-(C K)^star=X^star-K^circ} Let $X$ be a stereotype space. For any
absolutely convex compact set $K\subseteq X$ the Smith space it generates, $\C
K$, is connected with the quotient space $X^\star/K^\circ$ of the dual space
$X^\star$ with respect to the neighborhood of zero $K^\circ$ through the
formula
 \beq\label{(C K)^star=X^star-K^circ}
(\C K)^\star\cong X^\star/K^\circ
 \eeq
\etm

Let $X$ be a stereotype space and let $\mathcal K$ be a {\it expanding} system
of absolutely convex compact set in $X$, i.e. a system satisfying the following
condition:
$$
\forall K,L\in {\mathcal K}\qquad \exists M\in {\mathcal K}\qquad K\cup
L\subseteq M
$$
By Proposition \ref{vlozh-Smith}, for any compact sets $K,L\in{\mathcal K}$
such that $K\subseteq L$, the Smith spaces they generate are connected through
a natural linear continuous mapping $\iota_K^L: \C K\to \C L$. Since,
obviously, for any three compact sets $K\subseteq L\subseteq M$ the
corresponding mappings are connected through the equality
$$
\iota_L^M\circ\iota_K^L=\iota_K^M,
$$
the arising system of mappings $\{\iota_K^L; \; K,L\in{\mathcal K}:\,
K\subseteq L\}$ is an injective system in the category $\mathfrak{Ste}$ of
stereotype spaces. Like any other injective system in $\mathfrak{Ste}$, it has
a limit -- this is the pseudocompletion of its locally convex injective limit
\cite[Theorem 4.21]{Akbarov}:
$$
\mathfrak{Ste}\text{-}\kern-5pt\underset{K\to\infty}{\underset{\longrightarrow}{\lim}}
\C K =
\left(\mathfrak{LCS}\text{-}\kern-5pt\underset{K\to\infty}{\underset{\longrightarrow}{\lim}}
\C K\right)^\triangledown
$$

The dual construction is often used in the theory of topological vector spaces.
Let $X$ be a stereotype space and suppose $\mathcal U$ is a {\it
decreasing}\index{decreasing system of sets} system of absolutely convex
neighborhoods of zero in $X$, i.e. a system satisfying the following condition:
$$
\forall U,V\in {\mathcal U}\qquad \exists W\in {\mathcal U}\qquad W\subseteq
U\cap V
$$
For any two neighborhoods of zero $U,V\in{\mathcal U}$ such that $V\subseteq U$
the Banach spaces they generate $X/\Ker U$ and $X/\Ker V$ are connected with
each other through a natural linear continuous mapping $\pi_U^V: X/\Ker V\to
X/\Ker U$. If we consider three neighborhoods of zero $W\subseteq V\subseteq U$
the corresponding mappings are connected through the equality
$$
\pi_U^V\circ\pi_V^W=\pi_U^W,
$$
This means that the system of mappings $\{\pi_U^V; \; U,V\in{\mathcal U}:\,
V\subseteq U\}$ is a projective system in the category $\mathfrak{Ste}$ of
stereotype spaces. Its limit is the pseudosaturation of its locally convex
projective limit \cite[Theorem 4.21]{Akbarov}:
$$
\mathfrak{Ste}\text{-}\kern-2pt\underset{0\gets
U}{\underset{\longleftarrow}{\lim}} X/ U
=\left(\mathfrak{LCS}\text{-}\kern-2pt\underset{0\gets
U}{\underset{\longleftarrow}{\lim}} X/ U\right)^\vartriangle
$$
From \eqref{(C K)^star=X^star-K^circ} we have:

\btm If $\mathcal K$ is an expanding system of absolutely convex compact sets
in a stereotype space $X$, then the system of polars ${\mathcal U}=\{K^\circ;
\; K\in{\mathcal K}\}$ is a decreasing system of absolutely convex
neighborhoods of zero in the dual space $X^\star$. The limits of these systems
are dual to each other:
 \beq
\left(\mathfrak{Ste}\text{-}\kern-5pt\underset{K\to\infty}{\underset{\longrightarrow}{\lim}}
\C K \right)^\star= \mathfrak{Ste}\text{-}\kern-2pt\underset{0\gets
U}{\underset{\longleftarrow}{\lim}} X^\star/ K^\circ
 \eeq
\etm

\bex Let ${\mathcal K}={\mathcal K}(X)$ be a system of all absolutely convex
compact sets in $X$. Then the limit of the injective system $\{\iota_K^L; \;
K,L\in{\mathcal K}:\, K\subseteq L\}$ in the category of stereotype spaces is
the saturation\footnote{Saturation $X^\blacktriangle$ of a locally convex space
$X$ was defined in \cite[1.2]{Akbarov}.} of the space $X$:
$$
\mathfrak{Ste}\text{-}\kern-5pt\underset{K\to\infty}{\underset{\longrightarrow}{\lim}}\C
K=X^\blacktriangle
$$
 \eex

\bex Dually, if ${\mathcal U}={\mathcal U}(X)$ is the set of all absolutely
convex neighborhoods of zero in $X$, then the limit of the projective system
$\{\pi_U^V; \; U,V\in{\mathcal U}:\, V\subseteq U\}$ in the category of
stereotype spaces is the completion of the space $X$:
$$
\mathfrak{Ste}\text{-}\kern-2pt\underset{0\gets
U}{\underset{\longleftarrow}{\lim}}X/U=X^\blacktriangledown
$$
 \eex

\btm\label{Brauner-K_n} If $K_n$ is an expanding sequence of absolutely convex
compact sets in a stereotype space $X$, then the limit of the injective system
$\C K_n$ in the category of stereotype spaces coincides with locally convex
limit of this system and with the Brauner space generated by the sequence
$K_n$:
 \beq
\mathfrak{Ste}\text{-}\kern-3pt\underset{n\to\infty}{\underset{\longrightarrow}{\lim}}\C
K_n=\mathfrak{LCS}\text{-}\kern-3pt\underset{n\to\infty}{\underset{\longrightarrow}{\lim}}\C
K_n=\bigcup_{n=1}^\infty \C K_n
 \eeq
\etm

\subsection{Banach representation of a Smith space}

If $X$ is a Banach space  (with the norm $||\cdot||_X$), then let us denote by
$X^*$ its {\it dual Banach space} in the usual sense, i.e. the space of linear
continuous functionals on $X$ with the norm
$$
||f||_{X^*}=\sup_{||x||_X\le 1}|f(x)|
$$
If $\ph:X\to Y$ is a continuous linear mapping of Banach spaces, then the
symbol $\ph^*:Y^*\to X^*$ denotes the dual mapping:
$$
\ph^*(f)=f\circ\ph,\qquad f\in Y^*
$$
The natural mapping from $X$ into $X^{**}$ will be denoted by $s_X$:
$$
s_X:X\to X^{**}.
$$

Let $Y$ be a Smith space with the universal compact set $T$. Denote by $Y^{\B}$
the normed space with $Y$ as support and $T$ as unit ball.

\btm\label{TH-Y^B-cong-Y^star^*} For any Smith space $Y$
 \beq\label{Y^B-cong-Y^star^*}
Y^{\B}\cong (Y^\star)^*.
 \eeq
hence, $Y^{\B}$ is a (complete and so a) Banach space.
 \etm\bpr
Consider the space $X=Y^\star$. Then $Y$ becomes a space of linear continuous
functionals on a Banach space $X$ with the topology of uniform convergence on
compact sets in $X$. The universal compact $T$ in $Y$ becomes a polar of the
unit ball $B$ in $X$. If we endow $Y$ with the topology, where $T$ is a unit
ball, this is the same as if we endow $Y$ with the topology of the normed
space, dual to $X$: $Y^{\B}=X^*=(Y^\star)^*$. Since the (Banach) dual to a
Banach space is always a Banach space, $Y^{\B}$ must be a Banach space as well.
\epr

We call the space $Y^{\B}$ {\it Banach representation}\index{Banach
representation!of a Smith space} of the Smith space $Y$. Note that the natural
mapping
$$
\iota_Y:Y^{\B}\to Y,\quad \iota_Y(y)=y
$$
is universal in the following sense: for any Banach space $Z$ and for any
linear continuous mapping $\ph:Z\to Y$ there exists a unique linear continuous
mapping $\chi:Z\to Y^{\B}$ such that the following diagram is commutative:
$$
\begin{diagram}
\node{Y^{\B}} \arrow[2]{e,t}{\iota_Y} \node[2]{Y}
\\
\node[2]{Z}\arrow{nw,b,--}{\chi}\arrow{ne,b}{\ph}
\end{diagram}
$$
If $\ph:X\to Y$ is a linear continuous mapping of Smith spaces, then it turns
the universal compact set $S$ in $X$ into a compact set $\ph(S)$ in $Y$, which,
like any other compact set, is contained in some homothety of the universal
compact set $T$ in $Y$:
$$
\ph(S)\subseteq \lambda T,\qquad \lambda>0
$$
This means that the mapping $\ph$, being considered as the mapping between the
corresponding Banach representations $X^{\B}\to Y^{\B}$, is also continuous. We
shall denote this mapping by the symbol $\ph^{\B}$
$$
\ph^{\B}:X^{\B}\to Y^{\B}
$$
and call it {\it Banach representation of the mapping}\index{Banach
representation!of a mapping} $\ph$. Obviously,
 \beq\label{ph^B}
\ph^{\B}\cong (\ph^\star)^*
 \eeq
and the following diagram is commutative
$$
\begin{diagram}
\node{X} \arrow[2]{e,t}{\ph} \node[2]{Y}
\\
\node{X^{\B}}\arrow{n,b,}{\iota_X}\arrow[2]{e,b}{\ph^{\B}}
\node[2]{Y^{\B}}\arrow{n,b,}{\iota_Y}
\end{diagram}
$$

\subsection{Injective systems of Banach spaces, generated by compact sets}

The following construction is often used in the theory of topological vector
spaces. If $B$ is a bounded absolutely convex closed set in a locally convex
space $X$, then $X_B$ denotes the space $\bigcup_{\lambda>0}\lambda B$, endowed
with the topology of normed space with the unit ball $B$. If $X_B$ turns out to
be complete (i.e. a Banach) space, then the set $B$ is called a {\it Banach
disk}\index{Banach disk}. From theorem \ref{TH-Y^B-cong-Y^star^*} we have

\bprop\label{banach-disk} Every absolutely convex compact set $K$ in a
stereotype space $X$ is a Banach disk, and the Banach space $X_K$ generated by
$K$ is the Banach representation of the Smith space $\C K$:
$$
X_K=(\C K)^{\B}
$$
\eprop

If $\mathcal K$ is an expanding system of absolutely convex compact sets in a
stereotype space $X$, then, as we told above, $\mathcal K$ generates an
injective system of Smith spaces
$$
\{\C K\}_{K\in{\mathcal K}},\qquad \iota_K^L:\C K\to\C L\qquad (K,L\in{\mathcal
K},\quad K\subseteq L).
$$
Specialists in topological vector spaces used to replace this system of Smith
spaces with the system of the corresponding Banach representations:
$$
\{(\C K)^{\B}\}_{K\in{\mathcal K}},\qquad (\iota_K^L)^{\B}:(\C K)^{\B}\to (\C
L)^{\B}\qquad (K,L\in{\mathcal K},\quad K\subseteq L)
$$
The following result shows that the limits of those systems coincide in the
case when the injections $\iota_K^L:\C K\to\C L$ are compact mappings.

\btm\label{CK<->CK^B} Let $\mathcal K$ be an expanding system of absolutely
convex compact sets in a stereotype space $X$. Then the following conditions
are equivalent:
 \bit
\item[(i)] for any compact set $K\in{\mathcal K}$ there is a compact set
$L\in{\mathcal K}$ such that $K\subseteq L$ and the mapping of Smith spaces
$\iota_K^L:\C K\to\C L$ is compact,

\item[(ii)] for any compact set $K\in{\mathcal K}$ there is a compact set
$L\in{\mathcal K}$ such that $K\subseteq L$ and the mapping of Banach spaces
$(\iota_K^L)^{\B}:(\C K)^{\B}\to (\C L)^{\B}$ is compact.
 \eit
If these conditions hold then the locally convex injective limits of the
systems $\{\C K\}_{K\in{\mathcal K}}$ and $\{(\C K)^{\B}\}_{K\in{\mathcal K}}$
coincide,
 \beq\label{LCS-lim-X=LCS-lim-X^B}
\mathfrak{LCS}\text{-}\kern-5pt\underset{K\to\infty}{\underset{\longrightarrow}{\lim}}
\C K =
\mathfrak{LCS}\text{-}\kern-5pt\underset{K\to\infty}{\underset{\longrightarrow}{\lim}}
(\C K)^{\B}
 \eeq
and the same is true for their stereotype injective limits:
 \beq\label{Ste-lim-X=Ste-lim-X^B}
\mathfrak{Ste}\text{-}\kern-5pt\underset{K\to\infty}{\underset{\longrightarrow}{\lim}}
\C K =
\mathfrak{Ste}\text{-}\kern-5pt\underset{K\to\infty}{\underset{\longrightarrow}{\lim}}
(\C K)^{\B}
 \eeq
If in addition to {\rm (i)-(ii)} the system $\mathcal K$ is countable (or
contains a countable cofinal subsystem), then those four limits coincide,
 \beq\label{Ste-lim-X=Ste-lim-X^B-2}
\mathfrak{Ste}\text{-}\kern-5pt\underset{K\to\infty}{\underset{\longrightarrow}{\lim}}
\C K =
\mathfrak{LCS}\text{-}\kern-5pt\underset{K\to\infty}{\underset{\longrightarrow}{\lim}}
\C K =
\mathfrak{LCS}\text{-}\kern-5pt\underset{K\to\infty}{\underset{\longrightarrow}{\lim}}
(\C K)^{\B} =
\mathfrak{Ste}\text{-}\kern-5pt\underset{K\to\infty}{\underset{\longrightarrow}{\lim}}
(\C K)^{\B}
 \eeq
and define a Brauner space.
 \etm

We shall premise the proof of this theorem by several auxiliary propositions.

First we note that by {\it compact mapping}\index{compact mapping} of
stereotype spaces we mean what is usually meant, i.e. a linear continuous
mapping $\ph:X\to Y$ such that
$$
\ph(U)\subseteq T
$$
for some neighborhood of zero $U\subseteq X$ and some compact set $T\subseteq
Y$.

\bprop\label{comp-maps} Let $X$ and $Y$ be Smith spaces and $\ph:X\to Y$ a
linear continuous mapping. The following conditions are equivalent:
 \bit
\item[(i)] $\ph:X\to Y$ is a compact mapping;

\item[(ii)] $\ph^\star:Y^\star\to X^\star$ is a compact mapping;

\item[(iii)] $\ph^{\B}:X^{\B}\to Y^{\B}$ is a compact mapping.

 \eit
\eprop \bpr The equivalence (i)$\Leftrightarrow$(ii) is obvious, and
(ii)$\Leftrightarrow$(iii) follows from \eqref{ph^B} and a classical result on
compact mappings \cite[Theorem 4.19]{Rudin}. \epr

\bprop\label{kompakt^**} If $\ph:X\to Y$ is a compact mapping of Banach spaces,
then its (Banach) second dual mapping $\ph^{**}:X^{**}\to Y^{**}$ turns
$X^{**}$ into $Y$:
 \beq\label{ph^**-X^**-subseteq-Y}
\ph^{**}(X^{**})\subseteq Y
  \eeq
\eprop
 \bpr Let $B$ be unit ball in $X$ and $T$ a compact set in $Y$ such that
$$
\ph(B)\subseteq T
$$
By the bipolar theorem, $B$ is $X^*$-weakly dense in the unit ball
$B^{\circ\circ}$ of the space $X^{**}$. Hence for each $z\in B^{\circ\circ}$
one can choose a net $z_i\in B$ tending to $z$ $X^*$-weakly
 \beq\label{z_i-X^*-slabo->-z}
B\owns z_i\overset{\text{$X^*$-weakly}}{\underset{i\to\infty}{\longrightarrow}}
z\in B^{\circ\circ}
 \eeq
Since $\ph$, like any weakly compact operator, turns $X^*$-weak Cauchy nets
into `strong' Cauchy nets, we obtain that $\ph(z_i)$ is a Cauchy net in the
compact set $T$, hence it tends to some $y\in T$.
$$
\ph(z_i)\overset{Y}{\underset{i\to\infty}{\longrightarrow}} y
$$
From this we have:
 \beq\label{ph^**-z_i-g-1}
\forall g\in Y^*\qquad
\ph^{**}(z_i)(g)=\ph*(g)(z_i)=g(\ph(z_i))\underset{i\to\infty}{\longrightarrow}
g(y)=\i_Y(y)(g)
 \eeq
(where $\i_Y:Y\to Y^{**}$ is a natural embedding).

On the other hand from \eqref{z_i-X^*-slabo->-z} we also obtain that
$$
\forall f\in X^*\qquad f(z_i)\underset{i\to\infty}{\longrightarrow} f(z)
$$
and, in particular, this must be true for functionals $f=\ph*(g)=g\circ\ph$,
where $g\in Y^*$:
 \beq\label{ph^**-z_i-g-2}
\forall g\in Y^*\qquad
\ph^{**}(z_i)(g)=\ph*(g)(z_i)\underset{i\to\infty}{\longrightarrow}
\ph*(g)(z)=\ph^{**}(z)(g)
 \eeq
From \eqref{ph^**-z_i-g-1} and \eqref{ph^**-z_i-g-2} we have
$$
\ph^{**}(z)=y
$$
I.e., $\ph^{**}(z)\in Y$. This holds for each point $z\in B^{\circ\circ}$, so
\eqref{ph^**-X^**-subseteq-Y} must be true.
 \epr

\bprop Let $\ph:X\to Y$ be a linear continuous mapping of stereotype spaces.
For any absolutely convex closed set $V\subseteq Y$ the following formula
holds:
 \beq\label{ph^-1(V)}
\ph^{-1}(V)={^\circ}\Big(\ph^\star(V^\circ)\Big)
 \eeq
\eprop
 \bpr
 \begin{multline*}
x\in {^\circ}\Big(\ph^\star(V^\circ)\Big)\quad\Longleftrightarrow\quad \forall
g\in V^\circ\quad |\ph^\star(g)(x)|=|g(\ph(x))|\le 1 \quad\Longleftrightarrow
\\ \Longleftrightarrow\quad \ph(x)\in {^\circ}(V^\circ)=V\quad\Longleftrightarrow\quad x\in
\ph^{-1}(V)
 \end{multline*}
 \epr

\bprop\label{sigma-1(V)} Let $\sigma:X\to Y$ be a compact mapping of Smith
spaces, and $V$ an absolutely convex closed neighborhood of zero in the Banach
representation $Y^{\B}$ of the space $Y$. Then its preimage
$(\sigma^{\B})^{-1}(V)=\sigma^{-1}(V)$ is a neighborhood of zero in the space
$X$ (and not only in its Banach representation $X^{\B}$).
$$
\begin{diagram}
\node{X^{\B}} \arrow{s,t}{\iota_X}\arrow[2]{e,t}{\sigma^{\B}}
\node[2]{Y^{\B}}\arrow{s,b}{\iota_Y}
\\
\node{X}\arrow[2]{e,t}{\sigma} \node[2]{Y}
\end{diagram}
$$
\eprop
 \bpr
By formula \eqref{ph^B}, we can consider $\sigma^{\B}$ as a mapping of dual
Banach spaces for $X^\star$ and $Y^\star$:
$$
\sigma^{\B}=(\sigma^\star)^*:(X^\star)^*\to (Y^\star)^*
$$
Let us consider the Banach dual mapping
$(\sigma^{\B})^*=(\sigma^\star)^{**}:(Y^\star)^{**}\to (X^\star)^{**}$. Since
$\sigma^\star$ is compact, from Proposition \ref{kompakt^**} it follows that
$$
(\sigma^\star)^{**}\Big((Y^\star)^{**}\Big)\subseteq X^\star
$$
This can be illustrated by the diagram
$$
\begin{diagram}
\node{(X^\star)^{**}}
\node[2]{(Y^\star)^{**}}\arrow[2]{w,t}{(\sigma^\star)^{**}} \arrow{sww,t,--}{}
\\
\node{X^\star}\arrow{n,t}{\i_X}
\node[2]{Y^\star}\arrow[2]{w,t}{\sigma^\star}\arrow{n,t}{\i_Y}
\end{diagram}
$$
Now we have:
$$
\text{$V$ is a neighborhood of zero in $Y^{\B}=(Y^\star)^*$}
$$
$$
\Downarrow
$$
$$
\text{$V^\circ$ is a bounded set in $(Y^\star)^{**}$}
$$
$$
\Downarrow
$$
 \begin{multline*}
\text{$(\sigma^\star)^{**}(V^\circ)$ is a totally bounded set in
$(X^\star)^{**}$} \\
\text{(since $(\sigma^\star)^{**}$ is a compact mapping)}\quad \&\quad
(\sigma^\star)^{**}(V^\circ)\subseteq X^\star
 \end{multline*}
$$
\Downarrow
$$
$$
\text{$(\sigma^\star)^{**}(V^\circ)$ is a totally bounded set in $X^\star$}
$$
$$
\Downarrow
$$
$$
\text{$(\sigma^{\B})^{-1}(V)=\eqref{ph^-1(V)}={^\circ}\Big((\sigma^{\B})^*(V^\circ)\Big)
={^\circ}\Big((\sigma^\star)^{**}(V^\circ)\Big)$ is a neighborhood of zero in
$X$}
$$
 \epr

\bpr[Proof of Theorem \ref{CK<->CK^B}] The equivalence of (i) and (ii)
immediately follows from Propositions \ref{comp-maps}. After that formula
\eqref{LCS-lim-X=LCS-lim-X^B} follows from Proposition \ref{sigma-1(V)}: if $U$
is an absolutely convex neighborhood of zero in
$\mathfrak{LCS}\text{-}\kern-5pt\underset{K\to\infty}{\underset{\longrightarrow}{\lim}}
(\C K)^{\B}$, i.e. $U=\{U_K;\; K\in{\mathcal K}\}$ is a system of neighborhoods
of zeroes in the spaces $(\C K)^{\B}$, satisfying the condition
$$
\forall K,L\in{\mathcal K}\qquad \Big(K\subseteq L\quad\Longrightarrow\quad
U_K=(\iota_K^L)^{-1}(U_L)\Big)
$$
then from Proposition \ref{sigma-1(V)} it follows that  all those neighborhoods
are neighborhoods of zero in the spaces $\C K$. This proves the continuity of
the mapping
$$
\mathfrak{LCS}\text{-}\kern-5pt\underset{K\to\infty}{\underset{\longrightarrow}{\lim}}
\C K \to
\mathfrak{LCS}\text{-}\kern-5pt\underset{K\to\infty}{\underset{\longrightarrow}{\lim}}
(\C K)^{\B}
$$
The continuity of the inverse mapping is obvious. Thus, the topologies in those
spaces coincide. In other words, those limits coincide, and we obtain
\eqref{LCS-lim-X=LCS-lim-X^B}. This implies in its turn
\eqref{Ste-lim-X=Ste-lim-X^B}. Finally if $\mathcal K$ is countable, then by
Theorem \ref{Brauner-K_n} the locally convex limits coincide with stereotype
limits, and we have \eqref{Ste-lim-X=Ste-lim-X^B-2} . \epr

\subsection{Nuclear stereotype spaces}

A mapping of stereotype spaces $\ph:X\to Y$ we shall call {\it
nuclear}\index{nuclear!mapping}, if it satisfies the following four equivalent
conditions:
 \bit
\item[(i)] there exist two sequences
$$
f_n\overset{X^\star}{\underset{n\to\infty}{\longrightarrow}} 0,\qquad
b_n\overset{Y}{\underset{n\to\infty}{\longrightarrow}} 0,\qquad \lambda_n\ge
0,\qquad \sum_{n=1}^\infty \lambda_n<\infty,
$$
such that
 \beq\label{nucl-oper}
\ph(x)=\sum_{n=1}^\infty \lambda_n\cdot f_n(x)\cdot b_n,\qquad x\in X;
 \eeq

\item[(ii)] there exist two totally bounded sequences $\{f_n\}$ in $X^\star$
and $\{b_n\}$ in $Y$, and a number sequence $\lambda_n\ge 0$,
$\sum_{n=1}^\infty \lambda_n<\infty$ such that \eqref{nucl-oper} holds;

\item[(iii)] there exist two bounded sequences $\{f_n\}$ in $X^\star$, and
$\{b_n\}$ in $Y$, and a number sequence $\lambda_n\ge 0$, $\sum_{n=1}^\infty
\lambda_n<\infty$ such that \eqref{nucl-oper} holds;

\item[(iv)] there exist a sequence of functionals $\{f_n\}\subseteq X^\star$,
equicontinuous on $X$, and a sequence of vectors $\{b_n\}$ in some Banach disk
$B\subset Y$, and a number sequence $\lambda_n\ge 0$, $\sum_{n=1}^\infty
\lambda_n<\infty$ such that \eqref{nucl-oper} holds.
 \eit
Clearly, a nuclear mapping is always continuous.
 \bpr
The implications (i) $\Rightarrow$ (ii) $\Rightarrow$ (iii) are obvious. Let us
prove the implication (iii) $\Rightarrow$ (i). Suppose (iii) is true, i.e. the
sequences $\{f_n\}\subset X^\star$ and $\{b_n\}\subset Y$ in the formula
\eqref{nucl-oper} are just bounded. Let us choose a number sequences
$\sigma_n>0$ such that
$$
\sigma_n\underset{n\to\infty}{\longrightarrow}+\infty,\qquad \sum_{n=1}^\infty
\lambda_n\cdot\sigma_n <\infty
$$
(we can denote $\sum_{n=1}^\infty \lambda_n=C$ for this, and then divide the
sequence $\lambda_n$ into blocks $\sum_{n_k<n\le n_{k+1}}
\lambda_n\le\frac{C}{2^k}$, $k\ge 0$, $n_0=0$, and put $\sigma_n=k+1$ for
$n_k<n\le n_{k+1}$). Then we can replace $f_n$, $b_n$, $\lambda_n$ with
$$
\widetilde{f_n}:=\frac{f_n}{\sqrt{\sigma_n}},\qquad
\widetilde{b_n}:=\frac{b_n}{\sqrt{\sigma_n}},\qquad
\widetilde{\lambda_n}:=\lambda_n\cdot \sigma_n,
$$
because
$$
\sum_{n=1}^\infty \widetilde{\lambda_n}\cdot \widetilde{f_n}(x)\cdot
\widetilde{b_n}=\sum_{n=1}^\infty \lambda_n\cdot \sigma_n\cdot
\frac{f_n(x)}{\sqrt{\sigma_n}}\cdot \frac{b_n}{\sqrt{\sigma_n}}=
\sum_{n=1}^\infty \lambda_n\cdot f_n(x)\cdot b_n=\ph(x)
$$
And since $f_n$ and $b_n$ are bounded, then $\widetilde{f_n}$ and
$\widetilde{b_n}$ tend to zero
$$
\widetilde{f_n}\underset{n\to\infty}{\longrightarrow}0,\qquad
\widetilde{b_n}\underset{n\to\infty}{\longrightarrow}0,\qquad \sum_{n=1}^\infty
\widetilde{\lambda_n} <\infty
$$
i.e. (i) holds.

It remains now to show that the conditions (i)-(iii) are equivalent to (iv),
i.e. the standard definition of nuclear mapping \cite{Schaefer,Jarchow}. The
implication (iv) $\Rightarrow$ (iii) is obvious. On the other hand, the
implication (ii)$\Rightarrow$(iv) holds: from (ii) it follows that the sequence
of functionals $\{f_n\}$, being totally bounded in $X^\star$, is equicontinuous
on $X$ by \cite[Theorem 2.5]{Akbarov}, and the sequence $\{b_n\}$ belong to
Banach disk $T=\cabsconv\{b_n\}$ by theorem \ref{banach-disk}. This is exactly
(iv). \epr

\btm\label{nucl->nucl^star} A linear continuous mapping between stereotype
spaces $\ph:X\to Y$ is nuclear if and only if its dual mapping
$\ph^\star:Y^\star\to X^\star$ is nuclear. \etm

Let us call a mapping of stereotype spaces $\ph:X\to Y$ {\it
quasinuclear}\index{quasinuclear mapping}, if for any compact set $T\subset
Y^\star$ there exist a sequence of functionals $\{f_n\}\subseteq X^\star$
equicontinuous on $X$ and a number sequence $\lambda_n\ge 0$,
$\sum_{n=1}^\infty \lambda_n<\infty$ such that
 \beq
\max_{g\in T}|g\Big(\ph(x)\Big)|\le \sum_{n=1}^\infty \lambda_n |f_n(x)|
 \eeq
If $X$ and $Y$ are Banach spaces, then quasinuclearity of $\ph:X\to Y$ means
inequality
 \beq
||\ph(x)||_Y\le \sum_{n=1}^\infty \lambda_n |f_n(x)|,
 \eeq
where $||\cdot||_Y$ is the norm in $Y$, $f_n$ a bounded sequence of functionals
on $X$, $\lambda_n$ a summable number sequence.

\blm\label{nucl->nucl^*} Let $\ph:X\to Y$ be a linear continuous mapping of
Banach spaces, and let $\ph^*:Y^*\to X^*$ be its Banach dual mapping. Then
 \bit
\item[--] if $\ph$ is nuclear, then $\ph^*$ is nuclear as well,

\item[--] if $\ph^*$ is nuclear, then $\ph$ is quasinuclear.
 \eit
\elm
 \bpr
The first proposition here is obvious (and well-known, see
\cite[3.1.8]{Pietsch}) and implies the second one: if $\ph^*:Y^*\to X^*$ is
nuclear, then $\ph^{**}:X^{**}\to Y^{**}$ is also nuclear, and as a corollary,
quasinuclear. From this we have that $\ph$ is also quasinuclear, since it is a
restriction (and a corestriction) of a quasinuclear mapping $\ph^{**}$ on the
subspace $X\subseteq X^{**}$ (and subspace $Y\subseteq Y^{**}$).
 \epr

As usual, we call a stereotype space $X$ {\it nuclear}\index{nuclear!space}
\cite{Pietsch}, if every its continuous linear mapping into an arbitrary Banach
space $X\to Y$ is nuclear.

\btm[Brauner, \cite{Brauner}]\label{TH-X-nucl<->X*-nucl} A Brauner space $X$ is
nuclear if and only if its dual Fr\'echet space $X^\star$ is nuclear. \etm

\btm\label{CK<->CK^B-nucl} Let $\mathcal K$ be an expanding system of
absolutely convex compact sets in a stereotype space $X$. Then the following
conditions are equivalent:
 \bit
\item[(i)] for any compact set $K\in{\mathcal K}$ there is a compact set
$L\in{\mathcal K}$ such that $K\subseteq L$ and the mapping of Smith spaces
$\iota_K^L:\C K\to\C L$ is nuclear,

\item[(ii)] for any compact set $K\in{\mathcal K}$ there is a compact set
$L\in{\mathcal K}$ such that $K\subseteq L$ and the mapping of Banach spaces
$(\iota_K^L)^{\B}:(\C K)^{\B}\to (\C L)^{\B}$ is nuclear,
 \eit
If these conditions hold, then the locally convex injective limits of the
systems $\{\C K\}_{K\in{\mathcal K}}$ and $\{(\C K)^{\B}\}_{K\in{\mathcal K}}$
coincide,
 \beq
\mathfrak{LCS}\text{-}\kern-5pt\underset{K\to\infty}{\underset{\longrightarrow}{\lim}}
\C K =
\mathfrak{LCS}\text{-}\kern-5pt\underset{K\to\infty}{\underset{\longrightarrow}{\lim}}
(\C K)^{\B}
 \eeq
and the same is true for their stereotype injective limits:
 \beq
\mathfrak{Ste}\text{-}\kern-5pt\underset{K\to\infty}{\underset{\longrightarrow}{\lim}}
\C K =
\mathfrak{Ste}\text{-}\kern-5pt\underset{K\to\infty}{\underset{\longrightarrow}{\lim}}
(\C K)^{\B}.
 \eeq
If in addition to {\rm (i)-(ii)} the system $\mathcal K$ is countable (or
contains a countable cofinal subsystem), then all those four limits are equal
 \beq
\mathfrak{Ste}\text{-}\kern-5pt\underset{K\to\infty}{\underset{\longrightarrow}{\lim}}
\C K =
\mathfrak{LCS}\text{-}\kern-5pt\underset{K\to\infty}{\underset{\longrightarrow}{\lim}}
\C K =
\mathfrak{LCS}\text{-}\kern-5pt\underset{K\to\infty}{\underset{\longrightarrow}{\lim}}
(\C K)^{\B} =
\mathfrak{Ste}\text{-}\kern-5pt\underset{K\to\infty}{\underset{\longrightarrow}{\lim}}
(\C K)^{\B}
 \eeq
and define a nuclear Brauner space. \etm

\bpr All statements here follow from Theorem \ref{CK<->CK^B}, except the
equivalence between (i) and (ii).

(i)$\Longrightarrow$(ii). If $\iota_K^L$ is nuclear, then $(\iota_K^L)^\star$
is nuclear by Theorem \ref{nucl->nucl^star}. Hence
$(\iota_K^L)^{\B}=((\iota_K^L)^\star)^*$ is nuclear by Lemma
\ref{nucl->nucl^*}.

(i)$\Longrightarrow$(ii). For a given compact set $K$ let us choose a compact
set $L\supseteq K$ such that $(\iota_K^L)^{\B}$ is nuclear. Similarly, take a
compact set $M\supseteq L$ such that $(\iota_L^M)^{\B}$ is nuclear. Then from
nuclearity of $(\iota_K^L)^{\B}=((\iota_K^L)^\star)^*$ and
$(\iota_L^M)^{\B}=((\iota_L^M)^\star)^*$ by Lemma \ref{nucl->nucl^*} we have
that $(\iota_K^L)^\star$ and $(\iota_L^M)^\star$ are quasinuclear. from this we
have that $(\iota_K^L)^\star\circ (\iota_L^M)^\star$ is nuclear as a
composition of quasinuclear mappings \cite[3.3.2]{Pietsch}. Now by Theorem
\ref{nucl->nucl^star}, $\iota_L^M\circ \iota_K^L=\iota_K^M$ becomes nuclear, as
a stereotype dual mapping.
 \epr

\subsection{Spaces $\C^M$ and $\C_M$}\label{C^M-i-C_M}

The spaces $\C^M$ and $\C_M$ we are talking about here are usually mentioned in
textbooks on topological vector spaces as objects for exercises (\cite[Chapter
IV,\S 1, Exercise 11,13]{Bourbaki-tvs}, \cite[Chapter IV, Exerscise
6]{Schaefer}). We list here some of their properties for the further
references.

\paragraph{Space of functions $\C^M$.}
Let $M$ be an arbitrary set. Denote by $\C^M$ the locally convex space of all
complex-valued functions on $M$,
$$
u\in\C^M\quad\Longleftrightarrow\quad u:M\to\C
$$
with the pointwise operations, and the topology of pointwise convergence on
$M$, i.e. the topology generated by seminorms
 \beq\label{polun-na-C^M}
\abs{u}_N=\sup_{x\in N} \abs{u(x)}
 \eeq
where $N$ is an arbitrary finite set in $M$. We call $\C^M$ {\it the space of
functions} on the set $M$. Note that $\C^M$ is isomorphic to the direct product
of $\card M$ copies of the field $\C$:
$$
\C^M\cong \C^{\card M}
$$
As a corollary, $\C^M$ is always nuclear (since nuclearity is inherited by
direct products \cite[Theorem 7.4]{Schaefer}).

\btm\label{TH:Martineau} For a locally convex space $X$ over $\C$ the following
conditions are equivalent:
 \bit
\item[(a)] $X\cong\C^M$ for some set $M$;

\item[(b)] $X$ is a complete space with the weak topology\footnote{A locally
convex space $X$ is called {\it a space with the weak topology}, if its
topology is generated by the seminorms of the form $|x|_f=|f(x)|$, where $f$
are linear continuous functionals on $X$.};

\item[(c)] $X$ is a space of minimal type\footnote{A locally convex space $X$
is called {\it a space of minimal type}\index{space!of minimal type}, if there
is no weaker Hausdorff locally convex topology on $X$.}.
 \eit
 \etm

\btm $\C^M$ is a Fr\'echet space if and only if the set $M$ is at most
countable. \etm

\paragraph{The space of point charges $\C_M$.}
Let again $M$ be an arbitrary set. Denote by $\C_M$ the set of all number
families $\{\alpha_x;\; x\in M\}$ indexed by elements of $M$ and satisfying the
following finiteness condition: all the numbers $\alpha_x$, but a finite
subfamily, vanish
 \beq\label{def:C_M}
\alpha\in \C_M\quad \Longleftrightarrow\quad \alpha=\{\alpha_x;\; x\in
M\},\quad \alpha_x\in\C,\quad \card\{x\in M:\alpha_x\ne 0\}<\infty
 \eeq
(clearly the families $\{\alpha_x;\; x\in M\}$ can be considered as functions
$\alpha:M\to\C$ with finite support). The set $\C_M$ is endowed with pointwise
algebraic operations (sum and multiplication by a scalar) and a topology
generated by seminorms
 \beq\label{polun-na-C_M}
\abs{\alpha}_r=\sup_{|u|\le r} \abs{\langle u,\alpha \rangle}=\sum_{x\in
G}r(x)\cdot |\alpha_x|
 \eeq
where $r:G\to\R_+$ is an arbitrary nonnegative function on $M$. We call $\C_M$
{\it the space of point charges}\index{space!of point charges} on the set $M$
(and its elements -- point charges on $M$). We can note that $\C_M$ is
isomorphic to a locally convex direct sum of $\card M$ copies of the field
$\C$:
$$
\C_M\cong \C_{\card M}
$$
As a corollary, $\C_M$ is nuclear if and only if $M$ is at most countable
(nuclearity is inherited by countable direct sums \cite[Theorem 7.4]{Schaefer},
and if $M$ is not countable, then the imbedding $\C_M\to \ell_1(M)$ is not a
nuclear mapping).

\btm\label{TH-C_M} For a locally convex space $X$ over the field $\C$ the
following conditions are equivalent:
 \bit
\item[(a)] $X\cong\C_M$ for some set $M$;

\item[(b)] $X$ is a cocomplete\footnote{We say that a locally convex space $X$ is {\it
co-complete}\index{co-completed space}, if every linear functional $f:X\to\C$
which is continuous on every compact set $K\subset X$, is continuous on $X$.}
Mackey space, where every compact set is finite-dimensional;

\item[(c)] topology of $X$ is maximal in the class of locally convex topologies
on $X$ (i.e. there is no stronger locally convex topology on $X$).
 \eit
\etm

The following theorem in essential part belongs to S.~Kakutani and V.~Klee
\cite{Kakutani-Klee}:

\btm\label{TH:C_M-Brauner} $\C_M$ is a Brauner space if and only if the set $M$
is at most countable. In this (and only in this) case the topology of $\C_M$
coincides with the so called {\rm finite topology} on $\C_M$ (a set $A$ is said
to be closed in finite topology, if its intersection $A\cap F$ with any
finite-dimensional subspace $F\subseteq\C_M$ is closed in $F$). \etm

\paragraph{Duality between $\C^M$ and $\C_M$.}

\btm The bilinear form
 \beq\label{bil-forma-C^M-C_M}
\langle u,\alpha \rangle=\sum_{x\in M}u(x)\cdot \alpha_x, \qquad
u\in\C^M,\;\alpha\in\C_M
 \eeq
turns $\C^M $ and $\C_M$ into a dual pair $\langle\C^M,\C_M\rangle$ of
stereotype spaces:
 \bit
\item[(i)] every point charge $\alpha\in\C_M$ generates a linear continuous
functional $f$ on $\C^M$ by the formula
$$
f(u)=\langle u,\alpha \rangle,\qquad u\in\C^M,
$$
and the mapping $\alpha\mapsto f$ is an isomorphism of locally convex spaces
$$
\C_M\cong (\C^M)^\star;
$$

\item[(ii)] on the contrary, every function $u\in\C^M$ generates a linear
continuous functional $f$ on $\C_M$ by the formula
$$
f(\alpha)=\langle u,\alpha \rangle,\qquad \alpha\in\C_M,
$$
and the mapping $u\mapsto f$ is an isomorphism of locally convex spaces
$$
\C^M\cong (\C_M)^\star.
$$
 \eit
\etm

\paragraph{Bases in $\C^M$ and $\C_M$.}

A {\it basis}\index{basis} in a topological vector space $X$ over $\C$ is a
family of vectors $\{e_i;\; i\in I\}$ in $X$ such that every vector $x\in X$
can be uniquely represented as a sum of a converging series in $X$
 \beq\label{bazis}
x=\sum_{i\in I} \lambda_i\cdot e_i,
 \eeq
with coefficients $\lambda_i=\lambda_i(x)$ continuously depending on $x\in X$.
The summability of series \eqref{bazis} is understood in the sense of Bourbaki
\cite{Bourbaki}: for each neighborhood of zero $U$ in $X$ there exists a finite
set $J\subseteq I$ such that for any its finite superset $K$, $J\subseteq
K\subseteq I$
$$
x-\sum_{i\in I} \lambda_i\cdot e_i\in U
$$
(the summability of series \eqref{bazis} does not mean that this series must
have finite or countable number of nonzero terms).

\btm The characteristic functions $\{1_x;\; x\in G\}$ of singletons
$\{x\}\subseteq M$:
 \beq\label{1_x-infty}
1_x(y)=\begin{cases}1,& x=y \\ 0& x\ne y \end{cases},\qquad y\in M
 \eeq
form a basis in the topological vector space $\C^M$: every function $u\in\C^M$
is a sum of a series
 \beq\label{u=sum_x-in-M-u(x)-1_x-infty}
u=\sum_{x\in M}u(x)\cdot 1_x
 \eeq
with coefficients $u(x)\in\C$ continuously depending on $u\in\C^M$. \etm

\btm The characteristic functions of singletons $\{x\}\subseteq M$, which as
elements of $\C_M$ we denote by the Kronecker symbols $\{\delta^x;\; x\in G\}$
 \beq\label{delta^x-infty}
\delta^x(y)=\begin{cases}1,& x=y \\ 0& x\ne y \end{cases},\qquad y\in M
 \eeq
form a basis in the topological vector space $\C_M$: every point charge
$\alpha\in\C_M$ is a sum of a series
 \beq\label{u=sum_x-in-M-u(x)-1_x-infty}
\alpha=\sum_{x\in M}\alpha_x\cdot \delta^x
 \eeq
with coefficients $\alpha_x\in\C$ continuously depending on $\alpha\in\C_M$.
\etm

\btm The bases $\{1_x;\; x\in M\}$ (in $\C^M$) and $\{\delta^x;\; x\in M\}$ (in
$\C_M$) are dual:
$$
\langle 1_x,\delta^y\rangle=\begin{cases}1, & x=y \\ 0, & x\ne y \end{cases}
$$
\etm

\btm\label{isomorph-bazisov} In the spaces $\C^M$ and $\C_M$ any two bases can
be transformed into each other through some automorphism (i.e. a linear
homeomorphism of the space into itself).\etm

\section{Stereotype Hopf algebras}\label{SEC:ster-algebry-Hopfa}

\subsection{Tensor products and the structure of monoidal category on $\mathfrak{Ste}$}

If $X$ and $Y$ are stereotype spaces, then by $Y:X$ we denote a set of all
linear continuous mappings $\ph:X\to Y$ endowed with the topology of uniform
convergence on totally bounded sets in $X$. The symbol $Y\oslash X$ means
pseudosaturation of this space:
$$
Y\oslash X=(Y:X)^\vartriangle
$$
(the operation of pseudosaturation $\vartriangle$ means some special
strengthening of the topology of the initial space -- see \cite[\S
1]{Akbarov}). The space $Y\oslash X$ is always stereotype (if $X$ and $Y$ are
stereotype).

In the category $\mathfrak{Ste}$ of stereotype spaces there are two natural
tensor products:
 \bit

\item[--] a projective tensor product is defined by the equality
$$
X\circledast Y=(X^\star\oslash Y)^\star;
$$
the corresponding elementary tensor $x\circledast y\in X\circledast Y$ ($x\in
X$, $y\in Y$) is defined by the formula
 \beq\label{x-circledast-y}
x\circledast y(\ph)=\ph(y)(x),\qquad \ph\in X^\star\oslash Y;
 \eeq

\item[--] an injective tensor product is defined by the equality
$$
X\odot Y=X\oslash Y^\star;
$$
the corresponding elementary tensor $x\odot y\in X\odot Y$ ($x\in X$, $y\in Y$)
is defined by the formula
 \beq\label{x-odot-y}
x\odot y(f)=f(y)\cdot x,\qquad f\in Y^\star.
 \eeq
 \eit

These two operations are connected with each other by two
isomorphisms of functors:
$$
\d:(X\circledast Y)^\star\cong X^\star\odot Y^\star, \qquad \ee:(X\odot
Y)^\star\cong X^\star\circledast Y^\star
$$

\btm\label{TH-Ste-circledast} The identities
 \begin{align}
& \a^\circledast_{X,Y,Z}\Big((x\circledast y)\circledast
z\Big)=x\circledast (y\circledast z) & & x\in X,\quad y\in Y,\quad z\in Z \\
& \ll^\circledast_{X}(\lambda\circledast x)=\lambda\cdot
x=\rr^\circledast_{X}(x\circledast\lambda)
 & & \lambda\in \C,\quad x\in X \\
& \cc^\circledast_{X}(x\circledast y)=y\circledast x & & x\in X,\quad y\in Y
 \end{align}
correctly define natural isomorphisms of functors in the category
$\mathfrak{Ste}$
 \begin{align*}
&\a^\circledast_{X,Y,Z}:(X\circledast Y)\circledast Z\to X\circledast
(Y\circledast Z) \\
& \ll^\circledast_{X}:\C\circledast X\to X \\
& \rr^\circledast_{X}:X\circledast \C\to X \\
& \cc^\circledast_{X,Y}:X\circledast Y\to Y\circledast X
 \end{align*}
These isomorphisms in their turn define a structure of symmetrical monoidal
category on $\mathfrak{Ste}$ \cite{General-algebra,MacLane} with respect to the
bifunctor $\circledast$. \etm

Before formulating the next theorem let us agree to denote by $\h$ the mapping
from $\C$ into $\C^\star$, which every number $\lambda\in\C$ turns into the
linear functional on $\C$, acting as multiplication by $\lambda$:
 \beq
\h:\C\to \C^\star\quad\big|\quad\h(\lambda)(\mu)=\lambda\cdot\mu\quad
(\lambda,\mu\in\C)
 \eeq
Clearly, $\h:\C\to \C^\star$ is an isomorphism of (finite-dimensional)
stereotype spaces.

\btm\label{TH-Ste-odot} The formulas
 \begin{align*}
&
\a^{\odot}_{X,Y,Z}=\Big(\i_X\odot(\i_Y\odot\i_Z)\circ1_{X^{\star\star}}\odot\d^{-1}_{Y^\star,Z^\star}\circ
\d^{-1}_{X^\star,Y^\star\circledast
Z^\star}\circ(\a^\circledast_{X^\star,Y^\star,Z^\star})^\star\circ
\\ &\kern100pt \circ\d_{X^\star\circledast Y^\star,Z^\star}\circ\d_{X^\star,Y^\star}\odot
1_{Z^{\star\star}}\circ(\i^{-1}_{X^\star}\odot \i^{-1}_{Y^\star})\odot
\i^{-1}_{Z^\star}\Big)^\star
\\
& \ll^\odot_{X}=
\left(\h^{-1}\odot\i^{-1}_X\circ\d_{\C,X^\star}\circ(\ll^\circledast_{X^\star})^\star\circ\i_X\right)^{-1}
 \\
& \rr^\odot_{X}=
\left(\i^{-1}_X\odot\h^{-1}\circ\d_{X^\star,\C}\circ(\rr^\circledast_{X^\star})^\star\circ\i_X\right)^{-1}
 \\
& \cc^\odot_{X,Y}=\left(\i^{-1}_X\odot\i^{-1}_Y\circ\d_{X^\star,Y^\star}\circ
(\cc^\circledast_{X^\star,Y^\star})^\star\circ(\d_{Y^\star,X^\star})^{-1}\circ\i_Y\odot\i_X\right)^{-1}
 \end{align*}
define isomorphisms of functors in the category $\mathfrak{Ste}$
 \begin{align*}
&\a^\odot_{X,Y,Z}:(X\odot Y)\odot Z\to X\odot
(Y\odot Z) \\
& \ll^\odot_{X}:\C\odot X\to X \\
& \rr^\odot_{X}:X\odot \C\to X \\
& \cc^\odot_{X,Y}:X\odot Y\to Y\odot X
 \end{align*}
such that
 \begin{align}
& \a^\odot_{X,Y,Z}\Big((x\odot y)\odot z\Big)=x\odot (y\odot z) & & x\in
X,\quad y\in Y,\quad z\in Z
\label{a-odot}\\
& \ll^\odot_{X}(\lambda\odot x)=\lambda\cdot x=\rr^\odot_{X}(x\odot\lambda) & &
\lambda\in \C,\quad x\in X
\label{l-r-odot}\\
& \cc^\odot_{X}(x\odot y)=y\odot x & & x\in X,\quad y\in Y \label{c-odot}
 \end{align}
These isomorphisms define a structure of symmetric monoidal category on
$\mathfrak{Ste}$ \cite{General-algebra} with respect to the bifunctor $\odot$.
\etm

The following fact was noted in \cite[Theorem 7.9]{Akbarov}: the identity
 \beq
@_{X,Y}(x\circledast y)=x\odot y,\qquad x\in X,\quad y\in Y
 \eeq
defines a natural transformation of the bifunctor $\circledast$ into the
bifunctor $\odot$
$$
@_{X,Y}:X\circledast Y\to X\odot Y,
$$
called {\it Grothendieck transformation}\index{Grothendieck transformation}.

\btm If $X$ and $Y$ are Fr\'echet spaces (respectively, Brauner spaces), then
 \bit
\item[(i)] the stereotype tensor products $X\circledast Y$ and $X\odot Y$ are
also Fr\'echet spaces (respectively, Brauner spaces),

\item[(ii)] if at least one of the spaces $X$ and $Y$ is nuclear, then the
Grothendieck mapping $@_{X,Y}:X\circledast Y\to X\odot Y$ is an isomorphism of
stereotype spaces
$$
X\circledast Y\stackrel{@_{X,Y}}{\cong} X\odot Y,
$$

\item[(iii)] if both of the spaces $X$ and $Y$ are nuclear, then the space
$X\circledast Y\cong X\odot Y$ is also nuclear.
 \eit
 \etm
 \bpr
The fact that the spaces $X\circledast Y$ and $X\odot Y$ are Fr\'echet spaces
(respectively, Brauner spaces), was noted in \cite[Theorems 7.22,
7.23]{Akbarov}). If $X$ or $Y$ is nuclear, then the isomorphism $X\circledast
Y\cong X\odot Y$ is a corollary of the fact that for two Fr\'echet spaces, of
which one possesses the classical approximation property, the tensor products
$\circledast$ and $\odot$ coincide with the usual projective $\hat{\otimes}$
and injective $\check{\otimes}$ tensor products \cite[Theorems 7.17,
7.21]{Akbarov}. Finally, if both $X$ and $Y$ are nuclear, then in the case when
both $X$ and $Y$ are Fr\'echet spaces, their tensor products $X\circledast
Y\cong X\odot Y\cong X\hat{\otimes} Y$ are nuclear, since nuclearity is
inherited by projective tensor product \cite[5.4.2]{Pietsch}. If $X$ and $Y$
are Brauner spaces, then by what we have already proved, the Fr\'echet space
$X^\star\odot Y^\star$ is nuclear, hence, by the Brauner theorem
\ref{TH-X-nucl<->X*-nucl}, the space $X\circledast Y\cong (X^\star\odot
Y^\star)^\star$ is also nuclear.
 \epr

As a corollary, we have

\btm The categories $\mathfrak{NFre}$ of nuclear Fr\'echet spaces and
$\mathfrak{NBra}$ of nuclear Brauner spaces are symmetrical monoidal categories
with respect to bifunctors $\circledast$ and $\odot$ (which coincide on each of
those categories).
 \etm

\subsection{Stereotype Hopf algebras}\label{ster-alg-Hopf}

\paragraph{Algebras, coalgebras and Hopf algebras in a symmetric monoidal category.}

Recall \cite{MacLane}, that an {\it algebra}\index{algebra in a symmetrical
monoidal category} or a {\it monoid} in a symmetrical monoidal category
$(\mathfrak{K},\otimes,I,\a,\rr,\ll,\c)$ is a triple $(A,\mu,\iota)$, where $A$
is an object in $\mathfrak{K}$, and $\mu:A\otimes A\to A$ (multiplication) and
$\iota:I\to A$ (identity) are morphisms, satisfying the following axioms of
associativity and identity:
  \beq\label{monoid}
\begin{diagram}
\node{(A\otimes A)\otimes A} \arrow{s,l}{\mu\otimes 1}
\arrow[2]{e,l}{\a_{A,A,A}} \node[2]{A\otimes (A\otimes A)}
\arrow{s,r}{1\otimes\mu}
\\
\node{A\otimes A} \arrow{e,t}{\mu} \node{A} \node{A\otimes A} \arrow{w,t}{\mu}
\end{diagram}
\qquad
\begin{diagram}
\node{I\otimes A} \arrow{se,b}{\iota\otimes 1} \arrow{e,l}{\ll_A} \node{A}
 \node{A\otimes I} \arrow{sw,b}{1\otimes\iota}\arrow{w,l}{\rr_A}
\\
\node[2]{A\otimes A} \arrow{n,r}{\mu}
\end{diagram}.
 \eeq
To any two monoids $(A,\mu_A,\iota_A)$ and $(B,\mu_B,\iota_B)$ one can assign a
monoid $(A\otimes B,\mu_{A\otimes B},\iota_{A\otimes B})$ in which the
structure morphisms are defined by formulas
 \beq\label{def-mu-iota-v-A-otimes-B}
\mu_{A\otimes B}:=\mu_A\otimes\mu_B\circ\theta_{A,B,A,B},\qquad \iota_{A\otimes
B}:=\iota_A\otimes\iota_B\circ\ll_I^{-1}
 \eeq
 $$
 \begin{diagram}
  \node{(A\otimes B)\otimes (A\otimes B)} \arrow{se,r,--}{\mu_{A\otimes B}}
  \arrow[2]{e,r}{\theta_{A,B,A,B}}
  \node[2]{(A\otimes A)\otimes (B\otimes B)}
  \arrow{sw,r}{\mu_A\otimes\mu_B}\\
  \node[2]{A\otimes B}
  \end{diagram}
  \qquad
   \begin{diagram}
  \node{I} \arrow[2]{e,l}{\ll^{-1}_I}
  \arrow{se,r,--}{\iota_{A\otimes B}}
  \node[2]{I\otimes I}
  \arrow{sw,r}{\iota_A\otimes\iota_B}\\
  \node[2]{A\otimes B}
  \end{diagram}
$$
and $\theta$ is the isomorphism of functors
 \beq\label{ABCD->ACBD}
 \Big((A,B,C,D)\mapsto (A\otimes B)\otimes (C\otimes D)\Big)
 \stackrel{\theta}{\rightarrowtail}
 \Big((A,B,C,D)\mapsto (A\otimes C)\otimes (B\otimes D)\Big)
 \eeq
coming in as a combination of structure isomorphisms $\a$, $\ll$, $\rr$, $\cc$,
$\a^{-1}$, $\ll^{-1}$, $\rr^{-1}$, $\cc^{-1}$ in the tensor category
$\mathfrak{K}$ (by the coherence theorem \cite{MacLane} this morphism is unique
in this class).

The notion of {\it coalgebra}\index{coalgebra in a symmetrical monoidal
category} or {\it comonoid} in a symmetric monoidal category $\mathfrak{K}$ is
defined dually as an arbitrary triple $(A,\varkappa,\e)$, where $A$ is an
object of $\mathfrak{K}$, and $\varkappa:A\to A\otimes A$ (comultiplication)
and $\e:A\to I$ (counit) are morphisms satisfying the dual conditions of
coassociativity and counit:
 \beq\label{komonoid}
\begin{diagram}
\node{(A\otimes A)\otimes A} \arrow[2]{e,l}{\a_{A,A,A}} \node[2]{A\otimes
(A\otimes A)}
\\
\node{A\otimes A}\arrow{n,l}{\varkappa\otimes 1_A}
\node{A}\arrow{e,t}{\varkappa}\arrow{w,t}{\varkappa} \node{A\otimes A}
\arrow{n,r}{1_A\otimes\varkappa}
\end{diagram} \qquad
\begin{diagram}
\node{I\otimes A}   \node{A}\arrow{w,l}{\ll_A^{-1}}\arrow{e,l}{\rr_A^{-1}}
 \arrow{s,r}{\varkappa}
 \node{A\otimes I}
\\
\node[2]{A\otimes A}\arrow{nw,b}{\e\otimes 1_A}\arrow{ne,b}{1_A\otimes\e}
\end{diagram}.
 \eeq
Like in the case of algebras, the tensor product $A\otimes B$ of two coalgebras
$(A,\varkappa_A,\e_A)$ and $(B,\varkappa_B,\e_B)$ in $\mathfrak{K}$ possesses a
natural structure of comonoid in the category $\mathfrak{K}$ with the structure
morphisms
$$
\varkappa_{A\otimes
B}:=\theta_{A,A,B,B}\circ\varkappa_A\otimes\varkappa_B,\qquad \e_{A\otimes
B}:=\lambda_I\circ\e_A\otimes\e_B
$$

{\it Hopf algebra}\index{Hopf algebra in a symmetrical monoidal category}
(another term -- {\it Hopf monoid}) in the symmetrical monoidal category
$\mathfrak{K}$ is a quintuple $(H,\mu,\iota,\varkappa,\e,\sigma)$, where $H$ is
an object in $\mathfrak{K}$, and morphisms
 \begin{align*}
 &\mu:H\otimes H\to H &\text{(multiplication)}, & \\
 & \iota:I\to H & \text{(unit)}, & \\
 &\varkappa:H\to H\otimes H & \text{(comultiplication)},& \\
 &\e:H\to I & \text{(counit)}, &\\
 &\sigma:H\to H & \text{(antipode)} &
 \end{align*}
satisfy the following conditions:
 \bit{
\item[1)] the triple $(H,\mu,\iota)$ is an algebra in $\mathfrak{K}$,

\item[2)] the triple $(H,\varkappa,\e)$ is a coalgebra in $\mathfrak{K}$,

\item[3)] the following diagrams are commutative:
 \beq\label{5-angles}
\dgARROWLENGTH=-2em
\begin{diagram}
  \node[3]{H}\arrow{see,r}{\varkappa}
 \\
 \node{H\otimes H}
 \arrow{sse,t,3}{\varkappa\otimes \varkappa}\arrow{nee,r}{\mu}
 \node[4]{H\otimes H} \\ \\
 \node[2]{(H\otimes H)\otimes (H\otimes H)}\arrow[2]{e,r}{\theta_{H,H,H,H}}
 \node[2]{(H\otimes H)\otimes (H\otimes H)}\arrow{nne,r,3}{\mu\otimes \mu}
\end{diagram}
 \eeq
 \beq\label{lambda}
\begin{diagram}
\node{I} \arrow{s,l}{\iota} \arrow{e,l}{\ll_I^{-1}}
 \node{I\otimes I} \arrow{s,r}{\iota\otimes \iota}
\\
\node{H}\arrow{e,r}{\varkappa} \node{H\otimes H}
\end{diagram};
\qquad
\begin{diagram}
\node{H\otimes H} \arrow{s,l}{\e\otimes\e} \arrow{e,l}{\mu}
 \node{H} \arrow{s,r}{\e}
\\
\node{I\otimes I}\arrow{e,r}{\ll_I} \node{I}
\end{diagram};
 \eeq
 \beq\label{i-e}
\begin{diagram}
\node[2]{H} \arrow{se,l}{\e}
\\
\node{I}\arrow{ne,l}{\iota}\arrow[2]{e,r}{1_I} \node[2]{I}
\end{diagram};
 \eeq
-- they mean that morphisms $\varkappa:H\to H\otimes H$ and $\e:H\to I$ are
homomorphisms of algebras, and morphisms $\mu:H\otimes H\to H$ and $\iota:I\to
H$ are homomorphisms of coalgebras in the category $\mathfrak{K}$ (here again
$\theta$ is the transformation \eqref{ABCD->ACBD});

\item[4)] the following diagram called {\it antipode axiom} is commutative:

\beq\label{AX-antipode}
 \xymatrix  @R=2.0pc @C=.8pc
 {
 & H\otimes H\ar[rr]^{\sigma\otimes 1_H} & & H\otimes H \ar[dr]^{\mu}& \\
 H\ar[ur]^{\varkappa}\ar[dr]_{\varkappa}\ar[rr]^{\e} &  & I\ar[rr]^{\iota} & & H \\
 & H\otimes H\ar[rr]^{1_H\otimes\sigma} & & H\otimes H \ar[ur]_{\mu}&
 }
\eeq
 }\eit
If only conditions 1)-3) are fulfilled, then the quadruple
$(H,\mu,\iota,\varkappa,\e)$ is called a {\it bialgebra}\index{bialgebra in a
symmetrical monoidal category} in the category $\mathfrak{K}$.

\paragraph{Projective and injective stereotype algebras.}

In accordance with the general definition, a {\it projective} (respectively,
{\it injective}) {\it stereotype algebra}\index{stereotype!algebra} is an
algebra in the symmetric monoidal category of stereotype spaces
$(\mathfrak{Ste},\circledast)$ (respectively, $(\mathfrak{Ste},\odot)$).

For the case of projective algebras this definition admits a simple
reformulation \cite[\S 10]{Akbarov}:

\bprop\label{PROP:proj-algebra} The structure of projective stereotype algebra
on a stereotype space $A$ is equivalent to the structure of associative
(unital) algebra on $A$, where the multiplication $(x,y)\mapsto x\cdot y$
satisfies the following two equivalent conditions of continuity:
 \bit{
\item[(i)] for every compact set $K$ in $A$ and for every neighborhood of zero
$U$ in $A$ there exists a neighborhood of zero $V$ in $A$ such that
$$
K\cdot V\subseteq U\quad \& \quad V\cdot K\subseteq U
$$

\item[(ii)] for every compact set $K$ in $A$ and for any net $a_i$ in $A$,
tending to zero, $a_i\underset{i\to\infty}{\longrightarrow} 0$, the nets
$x\cdot a_i$ and $a_i\cdot x$ tend to zero in $A$ uniformly in $x\in K$.
 }\eit
\eprop

\bex Banach algebras and Fr\'echet algebras are examples of projective
stereotype algebras. Another example is the stereotype algebra ${\mathcal
L}(X)=X\oslash X$ of operators on an arbitrary stereotype space $X$. \eex

\bex The standard functional algebras ${\mathcal C}(M)$, ${\mathcal E}(M)$,
${\mathcal O}(M)$, ${\mathcal R}(M)$ of continuous, smooth, holomorphic
functions and polynomials are examples of injective stereotype algebras (see
details in \cite[\S 10]{Akbarov}). \eex

\paragraph{Stereotype Hopf algebras.}

Again, following the general definition, we call
 \bit

\item[--] {\it a projective stereotype Hopf algebra}\index{stereotype!Hopf
algebra} a Hopf algebra in the symmetric monoidal category of stereotype spaces
with the projective tensor product $(\mathfrak{Ste},\circledast)$;

\item[--] {\it an injective stereotype Hopf algebra} a Hopf algebra in the
symmetric monoidal category of stereotype spaces with the injective tensor
product  $(\mathfrak{Ste},\odot)$;

\item[--] {\it a nuclear Hopf-Fr\'echet algebra} a Hopf algebra in the
symmetric monoidal category of nuclear Fr\'echet spaces $\mathfrak{NFre}$;

\item[--] {\it a nuclear Hopf-Brauner algebra} a Hopf algebra in the symmetric
monoidal category of nuclear Brauner spaces $\mathfrak{NBra}$.

 \eit

Suppose in addition that a stereotype space $H$ is such that the Grothendieck
transformation for the pair $(H;H)$, for the triple $(H;H;H)$, and for the quadruple
$(H;H;H;H)$ give isomorphisms of stereotype spaces:
 \begin{align*}\label{@_H,H-@_H,H,H}
& @_{H,H}:H\circledast H\cong H\odot H,\\ & @_{H,H,H}:H\circledast H\circledast
H\cong H\odot H\odot H\\ & @_{H,H,H,H}:H\circledast H\circledast H\circledast
H\cong H\odot H\odot H\odot H
 \end{align*}
(this is always the case if $H$ is a nuclear Fr\'echet space or a nuclear
Brauner space). Then, obviously, every structure of projective stereotype Hopf
algebra on $H$ is equivalent to some structure of injective stereotype Hopf
algebra on $H$ and vice versa: the structure elements of Hopf algebra in
$(\mathfrak{Ste},\circledast)$ and $(\mathfrak{Ste},\odot)$ (we differ them by
indices $\circledast$ and $\odot$) either coincide
$$
\iota_{\circledast}=\iota_{\odot},\qquad \e_{\circledast}=\e_{\odot},\qquad
\sigma_{\circledast}=\sigma_{\odot}
$$
or are connected by diagrams
$$
\begin{diagram}
\node{H\circledast H} \arrow{se,r}{\mu_{\circledast}} \arrow[2]{e,l}{@_{H,H}}
\node[2]{H\odot H} \arrow{sw,r}{\mu_{\odot}}
\\
\node[2]{H}
\end{diagram}
\qquad
\begin{diagram}
\node{H\circledast H}  \arrow[2]{e,l}{@_{H,H}} \node[2]{H\odot H}
\\
\node[2]{H}\arrow{nw,r}{\varkappa_{\circledast}}
\arrow{ne,r}{\varkappa_{\odot}}
\end{diagram}
$$
We call those (projective and at the same time injective) Hopf algebras {\it
rigid stereotype Hopf algebras}\index{rigid Hopf algebra}.

\bex Clearly, every nuclear Hopf-Fr\'echet algebra and every nuclear
Hopf-Brauner algebra are rigid stereotype Hopf algebras. \eex

\paragraph{Duality for stereotype Hopf algebras.}

\btm[\bf on duality for stereotype Hopf algebras]
\label{TH-duality-in-st-Hopf-algebras} A structure of injective (projective,
rigid) Hopf algebra on a stereotype space $H$ automatically defines a structure
of projective (injective, rigid) Hopf algebra on the dual stereotype space
$H^\star$ -- the structure elements of Hopf algebra on $H^\star$ are defined as
the dual morphisms for the structure elements of Hopf algebra on $H$:
 \begin{align*}
& \mu_{H^\star}=(\varkappa_H)^\star,&& \iota_{H^\star}=(\e_H)^\star,\\
& \varkappa_{H^\star}=(\mu_H)^\star,&& \e_{H^\star}=(\iota_H)^\star,\\
& \sigma_{H^\star}=(\sigma_H)^\star.&&
 \end{align*}
 \etm

\bex Every Hopf algebra in usual sense $H$ (i.e. a Hopf algebra in the category
of vector spaces with the usual algebraic tensor product $\otimes$) becomes a
rigid Hopf algebra, being endowed with the strongest locally convex topology.
The dual space $H^\star$ with respect to this topology is a space of minimal
type (in the sense of Theorem \ref{TH:Martineau}) and a rigid Hopf algebra, as
well as $H$. This, by the way, illustrates one of the advantages of stereotype
theory: here we do not need to narrow the space of linear functionals to make a
Hopf algebra from them, like it is usually done (see e.g.
\cite[1.5]{Dascalescu} or \cite[4.1.D]{Chari-Pressley}) -- the space $H^\star$,
being a space of all linear functionals (automatically continuous, by the
choice of the topology in $H$), is a ``true'' Hopf algebra, but to see this we
have to take one of the stereotype tensor products $\circledast$ or $\odot$
instead of the algebraic tensor product $\otimes$. \eex

\paragraph{Dual pairs.}

Let $H$ be an injective, and $M$ a projective stereotype Hopf algebras. Let us
say that $H$ and $M$ form a {\it dual pair of Hopf algebras}\index{dual!pair of
Hopf algebras}, if there is a non-degenerate continuous (in the sense of
\cite[Section 5.6]{Akbarov}) bilinear form $\langle\cdot,\cdot\rangle:H\times
M\to\C$, which turns algebraic operations in $H$ into dual operations in $M$:
 \begin{align*}
&\langle \mu_H(U),\alpha\rangle=\langle U,\varkappa_M(\alpha)\rangle,&& \langle 1_H,\alpha\rangle=\e_A(\alpha),\\
&\langle \varkappa_H(u), A\rangle=\langle u,\mu_M(A)\rangle,&& \e_H(u)=\langle
u,1_M\rangle,\\
&\langle \sigma_H(u),\alpha\rangle=\langle u,\sigma_M(\alpha)\rangle & &
 \end{align*}
($u\in H$, $U\in H\odot H$, $\alpha\in M$, $A\in M\circledast M$). If $\langle
H, M\rangle$ is a dual pair of stereotype Hopf algebras, where $H$ is an
injective, $M$ a projective Hopf algebra, then on default we shall denote by
dot $\cdot$ the multiplication in $H$, by snowflake $*$ the multiplication in
$M$:
$$
\mu_H(u\odot v)=u\cdot v,\qquad \mu_M(\alpha\circledast\beta)=\alpha*\beta
$$

\bex Certainly, the pair $\langle H,H^\star\rangle$, where $H$ is an injective
stereptype Hopf algebra, and $\langle \cdot ,\cdot\rangle$  the canonical
bilinear form,
$$
\langle a,\alpha\rangle:=\alpha(a),\qquad a\in H,\ \alpha\in H^\star
$$
is an example of dual pair of Hopf algebras. Below we use the notation $\langle
\cdot ,\cdot\rangle$ for this form without supplementary explanations.\eex

\subsection{Key example: Hopf algebras $\C^G$ and $\C_G$}\label{C^G-i-C_G}

Recall that in \ref{SEC-ster-spaces}\ref{C^M-i-C_M} we defined the spaces
$\C^M$ and $\C_M$ of functions and of point charges on a set $M$. If $M$ is a
group, then the spaces $\C^M$ and $\C_M$ naturally turn into a dual pair of
Hopf algebras.

\paragraph{Algebra $\C^G$ of functions on $G$.}

Let $G$ be an arbitrary group (not necessarily finite) and let $\C^G$ be the
space of all (complex-valued) functions on $G$, defined in
\ref{SEC-ster-spaces}\ref{C^M-i-C_M},
$$
u\in\C^G\quad\Longleftrightarrow\quad u:G\to\C
$$
with the topology of pointwise convergence (generated by seminorms
\eqref{polun-na-C^M}). We endow $\C^G$ with the supplementary structure of
algebra with the pointwise algebraic operations:
 \beq\label{umnozhenie-v-C^G-infty}
(u\cdot v)(x)=u(x)\cdot v(x),\qquad 1_{\C^G}(x)=1\qquad (x\in G).
 \eeq
Recall that by formula \eqref{1_x-infty} above we defined the characteristic
functions $1_x$ of singletons in $G$. The multiplication in $\C^G$ can be
written in the decomposition by the basis $\{1_x;\; x\in G\}$ by formula
 \beq\label{umnozhenie-v-C^G-infty-po-bazisu}
u\cdot v=\left(\sum_{x\in G} u(x)\cdot 1_x\right)\cdot\left(\sum_{x\in G}
v(x)\cdot 1_x\right)=\sum_{x\in G} u(x)\cdot v(x) \cdot 1_x
 \eeq
and on the elements of this basis looks as follows:
 \beq\label{mult-in-C^G-infty}
1_x\cdot 1_y= \begin{cases}1_x,& x=y \\ 0& x\ne y \end{cases}
 \eeq

\paragraph{The algebra $\C_G$ of point charges on $G$.} Again let $G$ be an arbitrary group
and $\C_G$ the space of point charges on $G$, defined in
\ref{SEC-ster-spaces}\ref{C^M-i-C_M},
 $$
\alpha\in \C_G\quad \Longleftrightarrow\quad \alpha=\{\alpha_x;\; x\in
G\},\quad \alpha_x\in\C,\quad \card\{x\in G:\alpha_x\ne 0\}<\infty
 $$
We endow $\C_G$ with the topology generated by seminorms \eqref{polun-na-C_M}
(by theorem \ref{TH-C_M}, this is equivalent to the strongest locally convex
topology on $\C_G$). Besides this $\C_G$ is endowed with the structure of
algebra under the multiplication
$$
(\alpha*\beta)_y=\sum_{x\in G}\alpha_x\cdot\beta_{x^{-1}\cdot y}
$$
By formula \eqref{delta^x-infty} we defined the characteristic function of
singletons in $G$, which, being considered as elements in $\C_G$ were denoted
by the symbols of delta-functionals: $\delta^x$. The unit in $\C_G$ is the
characteristic function $\delta^e$ supported in the unit $e$ of the group $G$:
$$
\delta^e_x=\begin{cases}1,& x=e \\ 0& x\ne e \end{cases}
$$
The multiplication in the algebra $\C_G$ is written in the decomposition by
elements of the basis $\{\delta^x;\; x\in G\}$ by formula
 \beq\label{umnozhenie-v-C_G-infty}
\alpha * \beta=\left(\sum_{x\in G}\alpha_x\cdot\delta^x\right) *
\left(\sum_{y\in G}\beta_y\cdot\delta^y\right)=\sum_{x,y\in
G}\alpha_x\cdot\beta_y \cdot\delta^{x\cdot y}=\sum_{z\in G}\left(\sum_{x\in
G}\alpha_x\cdot\beta_{x^{-1}\cdot z}\right) \cdot\delta^z
 \eeq
and on the elements of this basis looks as follows:
 \beq\label{umnozhenie-na-bazise-v-C_G-infty}
\delta^x*\delta^y=\delta^{x\cdot y}.
 \eeq

\paragraph{$\C^G$ and $\C_G$ as stereotype Hopf algebras.}

For two arbitrary sets $S$ and $T$ and for two functions $u:S\to\C$ and
$v:T\to\C$ let the symbol $u\boxdot v$ denote the function on the Cartesian
product $S\times T$ defined by the identity
 \beq\label{g-h-na-S-T}
(u\boxdot v)(s,t):=u(s)\cdot v(t),\qquad s\in S,\quad t\in T
 \eeq
The definition of the Hopf algebra in the infinite-dimensional algebras $\C^G$
and $\C_G$ is based on the following observation:

\btm\label{TH-C^S-times-T-cong-C^S-odot-C^T} The formula
 \beq\label{u-boxdot-v->u-odot-v}
\rho_{S,T}(u\boxdot v)=u\odot v
 \eeq
defines an isomorphism of topological vector spaces
$$
\rho_{S,T}:\C^{S\times T}\cong \C^S\odot \C^T
$$
This isomorphism is an isomorphism of functors,
$$
\Big((S;T)\mapsto \C^{S\times T}\Big)\stackrel{\rho_{S,T}}{\rightarrowtail}
\Big((S;T)\mapsto \C^S\odot \C^T\Big),
$$
since for any mappings $\pi:S\to S'$ and $\sigma:T\to T'$ the following diagram
is commutative:
$$
\begin{diagram}
\node{\C^{S\times T}}\arrow{e,t}{\rho_{S,T}} \node{\C^S\odot\C^T}
\\
\node{\C^{S'\times T'}}\arrow{n,l}{\id_{\C}\oslash(\pi\times\sigma)}
\arrow{e,t}{\rho_{S',T'}}\node{\C^{S'}\odot\C^{T'}}
\arrow{n,r}{(\id_{\C}\oslash\pi)\odot(\id_{\C}\oslash\sigma)}
\end{diagram}
$$
-- here the mappings $\id_{\C}\oslash ?$ are defined by formula
$$
\id_{\C}\oslash\pi:\C^{S'}\to\C^S,\qquad\id_{\C}\oslash\pi(v)=v\circ\pi
$$
\etm

This theorem allows to define the structure elements of Hopf algebra on $\C^G$
with respect to $\odot$: the mapping $\mu:\C^G\odot\C^G\to\C^G$ is initially
defined on the set of functions on the Cartesian product $G\times G$,
$$
\widetilde{\mu}:\C^{G\times G}\to\C^G\quad\Big|\quad
\widetilde{\mu}(v)(t)=v(t,t)
$$
and then is extended to tensor square by the isomorphism $\rho_{G,G}$:
$$
\mu=\widetilde{\mu}\circ\rho_{G,G}.
$$
Similarly the comultiplication $\varkappa: \C^G\to\C^G\odot\C^G$ is initially
defined with the values in the space of functions on the Cartesian square
$$
\widetilde{\varkappa}:\C^G\to\C^{G\times G}\quad\Big|\quad
\widetilde{\varkappa}(u)(s,t)=u(s\cdot t)
$$
and after that is extended to tensor square by the isomorphism $\rho_{G,G}$:
$$
\varkappa=\rho_{G,G}\circ\widetilde{\varkappa}.
$$
The other structure elements of Hopf algebra on $\C^G$ are obvious:
 \begin{align*}
\text{unit:} & & & \iota:\C\to\C^G,& &
\iota(\lambda)(t)=\lambda \\
\text{counit:} & & &  \e:\C^G\to\C,& &
\e(u)=u(1_G) \\
\text{antipode:} & & & \sigma:\C^G\to\C^G,& & \sigma(u)(t)=u(t^{-1})
 \end{align*}
The following picture illustrates these definitions:
 \beq\label{str-morph-v-C^G}
 \begin{diagram}
 \node[3]{\C^{G\times G}}
 \arrow{se,l}{\widetilde{\mu}}
 \arrow[2]{s,r}{\rho_{G,G}}
 \\
 \node{\C}\arrow{e,t}{\iota} \node{\C^G}
 \arrow{ne,l}{\widetilde{\varkappa}}
 \arrow{se,r}{\varkappa}
 \node[2]{\C^G}\arrow{e,t}{\e}
 \node{\C}
 \\
 \node[3]{\C^G\odot\C^G}
 \arrow{ne,r}{\mu}
 \end{diagram}
 \eeq
On the space $\C_G$ of point charges on $G$ the structure of Hopf algebra with
respect to tensor product $\circledast$ is defined dually by Theorem
\ref{TH-duality-in-st-Hopf-algebras}, as on the dual space to $\C^G$ in the
sense of bilinear form \eqref{bil-forma-C^M-C_M}. The following theorem shows
that these definitions indeed define a structure of Hopf algebra on $\C^G$ and
$\C_G$:

\btm\label{TH-C^G-algebra-Hopfa} For any group $G$
 \bit

\item[--] the space $\C^G$ of functions on $G$ is an injective (and moreover, a
rigid) stereotype Hopf algebra; if in addition $G$ is countable, then $\C^G$ is
a nuclear Hopf-Fr\'echet algebra;

\item[--] the space $\C_G$ of point charges on $G$ is a projective (and
moreover, a rigid) stereotype Hopf algebra; if in addition $G$ is countable,
then $\C_G$ is a nuclear Hopf-Brauner algebra.
 \eit
The Hopf algebras $\C^G$ and $\C_G$ form a dual pair with respect to bilinear
form \eqref{bil-forma-C^M-C_M}, and the algebraic operations on them act on
bases $\{1_x\}$ and $\{\delta^x\}$ by formulas:
 \begin{align}
\label{antipode-in-C_G} & \C^G: && 1_x\cdot 1_y=\begin{cases}1_x, & x=y \\ 0, & x\ne y \end{cases} && 1_{\C^G}=\sum_{x\in G} 1_x && \sigma(1_x)=1_{x^{-1}}\\
\label{koumn-v-C^G} & && \varkappa(1_x)=\sum_{y\in G}1_y\odot 1_{x\cdot y^{-1}} && \e(1_x)=\begin{cases}1, & x=e \\ 0, & x\ne e \end{cases}&& \\
\nonumber & && && &&\\
& \C_G: && \delta^x*\delta^y=\delta^{x\cdot y} && 1_{\C_G}=\delta^e && \sigma(\delta^x)=\delta^{x^{-1}} \\
\label{koumn-v-C_G}& && \varkappa(\delta^x)=\delta^x\circledast\delta^x &&
\e(\delta^x)=1 &&
 \end{align}
 \etm
\bpr Formulas \eqref{antipode-in-C_G}-\eqref{koumn-v-C_G} are verified by
direct calculation. The rigidity follows from \cite[(7.37)]{Akbarov}:
$$
\C^I\circledast\C^J\cong \C^{I\times J}\cong \C^I\odot\C^J,\qquad
\C_I\circledast\C_J\cong \C_{I\times J}\cong \C_I\odot\C_J.
$$
The structure of Hopf algebra on $\C_G$ is generated by the structure of Hopf
algebra on $\C^G$ due to Theorem \ref{TH-duality-in-st-Hopf-algebras}. Thus, we
need only to prove that $\C^G$ is a Hopf algebra with respect to $\odot$.

1. Let us check diagram \eqref{5-angles}. Replace everywhere $H$ by $\C^G$ and
$\otimes$ by $\odot$, and after that let us overbuild this diagram to the
following prism:

 $$\dgARROWLENGTH=-6em
\begin{diagram}
  \node[3]{\C^G }\arrow[2]{se,l}{\widetilde{\varkappa}}\arrow[4]{s,-}
 \\
 \\
 \node{\C^{G\times G}}\arrow[8]{s,t}{\rho_{G,G}}
 \arrow{sse,l}{\widetilde{\widetilde{\varkappa}}}
 \arrow[2]{ne,l}{\widetilde{\mu}}
 \node[4]{\C^{G\times G}}\arrow[8]{s,r}{\rho_{G,G}} \\
  \\
 \node[2]{\C^{(G\times G)\times (G\times G)}}
 \arrow[2]{e,l,3}{\widetilde{\theta}}
 \arrow[3]{s,r}{\rho_{G\times G,G\times G}}
 \node{}\arrow[4]{s,r}{\id_{\C^G}}
 \node{\C^{(G\times G)\times (G\times G)}}
 \arrow{nne,l}{\widetilde{\widetilde{\mu}}}
 \arrow[3]{s,r}{\rho_{G\times G,G\times G}}
 \\ \\ \\
 \node[2]{\C^{G\times G}\odot \C^{G\times G}}\arrow[5]{s,r,3}{\rho_{G,G}\odot\rho_{G,G}}
 \node[2]{\C^{G\times G}\odot \C^{G\times G}}\arrow[5]{s,t,3}{\rho_{G,G}\odot\rho_{G,G}}
 \arrow{ssse,l}{\widetilde{\mu}\odot\widetilde{\mu}}
 \\
  \node[3]{\C^G}\arrow{se,-}
 \\
 \node[2]{}\arrow{ne,r,3}{\mu}
 \node[2]{}\arrow{se,r,3}{\varkappa}
 \\
 \node{\C^G\odot \C^G}
 \arrow{nnne,l}{\widetilde{\varkappa}\odot\widetilde{\varkappa}}
 \arrow{sse,l}{\varkappa\odot\varkappa}
 \arrow{ne,-}
 \node[4]{\C^G\odot \C^G}
  \\ \\
 \node[2]{\big(\C^G\odot\C^G\big)\odot \big(\C^G \odot\C^G\big)}
 \arrow[2]{e,l}{\theta}
 \node[2]{\big(\C^G\odot\C^G\big)\odot \big(\C^G \odot\C^G\big)}
 \arrow{nne,l}{\mu\odot\mu}
\end{diagram}
 $$
Here $\theta=\theta_{\C^G,\C^G,\C^G,\C^G}$ is the isomorphism of functors from
\eqref{ABCD->ACBD}, and the other morphisms are defined as follows:
 $$
\widetilde{\theta} w(a,b,c,d)=w(a,c,b,d), \quad
\widetilde{\widetilde{\varkappa}}v(a,b,c,d)=v(a\cdot b,c\cdot d), \quad
\widetilde{\widetilde{\mu}}v(a,b)=v(a,a,b,b)
 $$
To prove that the base of the prism is commutative, it is sufficient to verify
that all the other faces are commutative. The remote lateral faces
 $$
\begin{diagram}
 \node{\C^{G\times G}}
 \arrow{s,l}{\rho_{G,G}}
 \arrow{e,l}{\widetilde{\mu}}
 \node{\C^G}
 \arrow{s,r}{\id_{\C^G}}
 \\
 \node{\C^G\odot\C^G}
 \arrow{e,l}{\mu}
 \node{\C^G}
 \end{diagram}
\qquad\qquad
\begin{diagram}
 \node{\C^G}
 \arrow{s,l}{\id_{\C^G}}
 \arrow{e,l}{\widetilde{\varkappa}}
 \node{\C^{G\times G}}
 \arrow{s,r}{\rho_{G,G}}
 \\
 \node{\C^G}
 \arrow{e,l}{\varkappa}
 \node{\C^G\odot\C^G}
 \end{diagram}
 $$
-- are just distorted triangles from diagram \eqref{str-morph-v-C^G}.

In the left nearby face
 $$
\begin{diagram}
 \node{\C^{G\times G}}
 \arrow[2]{s,l}{\rho_{G,G}}
 \arrow{e,l}{\widetilde{\widetilde{\varkappa}}}
 \node{\C^{(G\times G)\times (G\times G)}}
 \arrow{s,r}{\rho_{G\times G,G\times G}}
 \\
 \node[2]{\C^{G\times G}\odot \C^{G\times G}}\arrow{s,r}{\rho_{G,G}\odot\rho_{G,G}}
 \\
 \node{\C^G\odot\C^G}
 \arrow{ne,l}{\widetilde{\varkappa}\odot\widetilde{\varkappa}}
 \arrow{e,r}{\varkappa\odot\varkappa}
 \node{\big(\C^G\odot\C^G\big)\odot\big(\C^G \odot\C^G\big)}
 \end{diagram}
 $$
-- the lower triangle is just left triangle in \eqref{str-morph-v-C^G}
multiplied by itself via the operation $\odot$, and the commutativity of the
inner quadrangle is verified by the substituting the function $u\boxdot v\in
\C^{G\times G}$, $u,v\in \C^G$, as an argument: if we move down, and then right
and up, we turn this function into $\widetilde{\varkappa}(u)\odot
\widetilde{\varkappa}(v)$,
 \begin{flalign*}
&u\boxdot v\in \C^{G\times G} &\\
 \mapsto\qquad & u\odot v \in \C^G\odot \C^G&
\\
 \mapsto\qquad &(\widetilde{\varkappa}\odot\widetilde{\varkappa})(u\odot v)=
\widetilde{\varkappa}(u)\odot \widetilde{\varkappa}(v)
 \in \C^{G\times G}\odot \C^{G\times G}&
 \end{flalign*}
-- and if we move right and then down, we obtain the same result:
 \begin{flalign*}
&u\boxdot v\in \C^{G\times G} &\\
\mapsto\qquad &\widetilde{\widetilde{\varkappa}}(u\boxdot v)=
\widetilde{\varkappa}(u)\boxdot \widetilde{\varkappa}(v) \in \C^{(G\times G)\times (G\times G)}: \\
& \widetilde{\widetilde{\varkappa}}(u\boxdot v)(a,b,c,d)=(u\boxdot v)(a\cdot
b,c\cdot d)=u(a\cdot b)\cdot v(c\cdot d) =\widetilde{\varkappa}(u)(a,b)\cdot
\widetilde{\varkappa}(v)(c,d)=& \\ &= \big(\widetilde{\varkappa}(u)\boxdot
\widetilde{\varkappa}(v)\big)(a,b,c,d)&
 \\
 \mapsto\qquad &
 \rho_{G\times G,G\times G}\left(\widetilde{\varkappa}(u)\boxdot
\widetilde{\varkappa}(v)\right)=\eqref{u-boxdot-v->u-odot-v}=
 \widetilde{\varkappa}(u))\odot \widetilde{\varkappa}(v) \in \C^{G\times G}\odot \C^{G\times G}.
 &
 \end{flalign*}

The commutativity of the central nearby face:
 $$
\begin{diagram}
 \node{\C^{(G\times G)\times (G\times G)}}
 \arrow{s,l}{\rho_{G\times G,G\times G}}
 \arrow{e,l}{\widetilde{\theta}}
 \node{\C^{(G\times G)\times (G\times G)}}
 \arrow{s,r}{\rho_{G\times G,G\times G}}
 \\
 \node{\C^{G\times G}\odot \C^{G\times G}}\arrow{s,l}{\rho_{G,G}\odot\rho_{G,G}}
 \node{\C^{G\times G}\odot \C^{G\times G}}\arrow{s,r}{\rho_{G,G}\odot\rho_{G,G}}
 \\
 \node{\big(\C^G\odot\C^G\big)\odot\big(\C^G \odot\C^G\big)}
 \arrow{e,l}{\theta}
 \node{\big(\C^G\odot\C^G\big)\odot\big(\C^G \odot\C^G\big)}
 \end{diagram}
 $$
-- is verified by substituting a function $(u\boxdot v)\boxdot (p\boxdot q)\in
\C^{G\times G}$ as an argument, its motion due to \eqref{u-boxdot-v->u-odot-v}
will be as follows:
 $$
 \xymatrix
 {
 (u\boxdot v)\boxdot (p\boxdot q)
 \ar@{|->}[d]_{\rho_{G\times G,G\times G}}
  \ar@{|->}[r]^{\widetilde{\theta}}
 &
 (u\boxdot p)\boxdot (v\boxdot q)
  \ar@{|->}[d]^{\rho_{G\times G,G\times G}}
 \\
 (u\boxdot v)\odot (p\boxdot q)
 \ar@{|->}[d]_{\rho_{G,G}\odot\rho_{G,G}}
 &
 (u\boxdot p)\odot (v\boxdot q)
  \ar@{|->}[d]^{\rho_{G,G}\odot\rho_{G,G}}
 \\
 (u\odot v)\odot (p\odot q)
  \ar@{|->}[r]^{\theta}
 &
 (u\odot p)\odot (v\odot q)
 }
 $$

In the right nearby face
 $$
\begin{diagram}
 \node{\C^{(G\times G)\times (G\times G)}}
 \arrow{s,l}{\rho_{G\times G,G\times G}}
 \arrow{e,l}{\widetilde{\widetilde{\mu}}}
 \node{\C^{G\times G}}
 \arrow[2]{s,r}{\rho_{G,G}}
 \\
 \node{\C^{G\times G}\odot \C^{G\times G}}
 \arrow{s,l}{\rho_{G,G}\odot\rho_{G,G}}
 \arrow{se,l}{\widetilde{\mu}\odot\widetilde{\mu}}
 \\
 \node{\big(\C^G\odot\C^G\big)\odot\big(\C^G \odot\C^G\big)}
 \arrow{e,r}{\mu\odot\mu}
 \node{\C^G\odot\C^G}
 \end{diagram}
 $$
-- the lower triangle is just the right triangle in \eqref{str-morph-v-C^G}
multiplied by itself via the operation $\odot$, and the commutativity of the
inner quadrangle is verified by taking as an argument the function $u\boxdot
v\in \C^{(G\times G)\times (G\times G)}$, $u,v\in \C^{G\times G}$: if we move
down, and then right and down, we obtain
 \begin{flalign*}
&u\boxdot v\in \C^{(G\times G)\times (G\times G)} &\\
 \mapsto\qquad & \rho_{G\times G,G\times G}(u\boxdot v)=u\odot v \in \C^{G\times G}\odot \C^{G\times G}&
\\
 \mapsto\qquad &(\widetilde{\mu}\odot\widetilde{\mu})(u\odot v)=
\widetilde{\mu}(u)\odot \widetilde{\mu}(v)
 \in \C^G\odot \C^G&
 \end{flalign*}
-- and if we move right and then down, we obtain the same result:
 \begin{flalign*}
& u\boxdot v\in \C^{(G\times G)\times (G\times G)} &\\
\mapsto\qquad &\widetilde{\widetilde{\mu}}(u\boxdot v)=
\widetilde{\mu}(u)\boxdot \widetilde{\mu}(v)\in \C^{G\times G}: \\
& \widetilde{\widetilde{\mu}}(u\boxdot v)(a,b)=(u\boxdot
v)(a,a,b,b)=u(a,a)\cdot v(b\cdot b)= \widetilde{\mu}(u)(a)\cdot
\widetilde{\mu}(v)(b)= \big(\widetilde{\mu}(u)\boxdot
\widetilde{\mu}(v)\big)(a,b)
& \\
 \mapsto\qquad &
 \rho_{G,G}\left(\widetilde{\mu}(u)\boxdot
\widetilde{\mu}(v)\right)=\eqref{u-boxdot-v->u-odot-v}=
 \widetilde{\mu}(u))\odot \widetilde{\mu}(v) \in \C^G\odot\C^G.
 &
 \end{flalign*}

It remains to check the commutativity of the upper base of the prism.
 $$\dgARROWLENGTH=-2em
\begin{diagram}
  \node[3]{\C^G }\arrow[2]{se,l}{\widetilde{\varkappa}}
 \\
 \\
 \node{\C^{G\times G}}
 \arrow{sse,l}{\widetilde{\widetilde{\varkappa}}}
 \arrow[2]{ne,l}{\widetilde{\mu}}
 \node[4]{\C^{G\times G}}
  \\
  \\
 \node[2]{\C^{(G\times G)\times (G\times G)}}
 \arrow[2]{e,l,3}{\widetilde{\theta}}
 \node[2]{\C^{(G\times G)\times (G\times G)}}
 \arrow{nne,l}{\widetilde{\widetilde{\mu}}}
\end{diagram}
$$
A function $v\in \C^{G\times G}$, being moved through the upper two edges of
the pentagon, undergoes the following transmutations:
 \begin{flalign*}
&v\in \C^{G\times G} &\\
 \mapsto\qquad &\widetilde{\mu}v \in \C^G,\qquad
\widetilde{\mu}v(a)=v(a,a) &
\\
 \mapsto\qquad &\widetilde{\varkappa}\left(\widetilde{\mu}v\right) \in \C^{G\times G},\qquad \widetilde{\varkappa}\left(\widetilde{\mu}v\right)(a,b)=
\widetilde{\mu}v(a\cdot b)=v(a\cdot b,a\cdot b) &
 \end{flalign*}
-- and the result is the same as if we move it through the three lower edges:
 \begin{flalign*}
&v\in \C^{G\times G}\\
\mapsto\qquad &\widetilde{\widetilde{\varkappa}}v \in \C^{(G\times G)\times
(G\times G)},\qquad \widetilde{\widetilde{\varkappa}}v(a,b,c,d)=v(a\cdot
b,c\cdot d)
& \\
 \mapsto\qquad &\widetilde{\theta}\left(\widetilde{\widetilde{\varkappa}}v\right) \in \C^{(G\times G)\times (G\times G)},\qquad
\widetilde{\theta}\left(\widetilde{\widetilde{\varkappa}}v\right)(a,b,c,d)=
\widetilde{\widetilde{\varkappa}}v(a,c,b,d)=v(a\cdot c,b\cdot d)
& \\
\mapsto\qquad &
\widetilde{\widetilde{\mu}}\left(\widetilde{\theta}\left(\widetilde{\widetilde{\varkappa}}v\right)\right)
 \in \C^{G\times G},\qquad
\widetilde{\widetilde{\mu}}\left(\widetilde{\theta}\left(\widetilde{\widetilde{\varkappa}}v\right)\right)
(a,b)=\widetilde{\theta}\left(\widetilde{\widetilde{\varkappa}}v\right)(a,a,b,b)=
\widetilde{\widetilde{\varkappa}}v(a,b,a,b)=v(a\cdot b,a\cdot b) &
 \end{flalign*}

2. After that we verify diagrams \eqref{lambda}, which for the algebra $\C^G$
come to the following form:
 $$
 \xymatrix
 {
 \C \ar[d]_{\iota} \ar[r]^{\ll_\C^{-1}}
 &
 \C\odot \C \ar[d]^{\iota\odot \iota}
 \\
 \C^G\ar[r]^(0.4){\varkappa}
 &
 \C^G\odot \C^G
 }
\qquad
 \xymatrix
 {
 \C^G\odot \C^G \ar[d]_{\e\odot\e} \ar[r]^{\mu}
 &
 \C^G \ar[d]^{\e}
 \\
 \C\odot \C
 \ar[r]^{\ll_\C}
 &
 \C
 }
 $$
Let us input an arbitrary number $\zeta\in\C$ as a argument into the first
diagram, and an elementary tensor $u\odot v\in \C^G\odot \C^G$ into the second
diagram, and apply identities \eqref{l-r-odot}:
 $$
 \xymatrix
 {
 \zeta=\zeta\cdot 1_\C \ar@{|->}[d]_{\iota} \ar@{|->}[r]^{\ll_\C^{-1}}
 &
 \zeta\cdot 1_\C\odot 1_\C \ar@{|->}[d]^{\iota\odot \iota}
 \\
 \zeta\cdot 1_{\C^G} \ar@{|->}[r]^(0.4){\varkappa}
 &
 \zeta\cdot 1_{\C^G}\odot 1_{\C^G}
 }
\qquad
 \xymatrix
 {
 u\odot v \ar@{|->}[d]_{\e\odot\e} \ar@{|->}[r]^{\mu}
 &
 u\cdot v \ar@{|->}[d]^{\e}
 \\
 u(1_G)\odot v(1_G)
 \ar@{|->}[r]^{\ll_\C}
 &
 u(1_G)\cdot v(1_G)
 }
 $$

3. It remains to verify diagram \eqref{AX-antipode}, which for $\C^G$ comes to
the form:
 \beq\label{AX-antipode-O(G)}
 \xymatrix  @R=3.0pc @C=.8pc
 {
 & \C^G\odot \C^G\ar[rr]^{\sigma\odot 1_{\C^G}} &
 & \C^G\odot \C^G \ar[dr]^{\mu}& \\
 \C^G
 \ar[ur]^{\varkappa}\ar[dr]_{\varkappa}\ar[rr]^{\e} &  & \C\ar[rr]^{\iota} & & \C^G \\
 & \C^G\odot \C^G\ar[rr]^{1_{\C^G}\odot\sigma} & &
 \C^G\odot \C^G \ar[ur]_{\mu}&
 }
 \eeq
We overbuild it to the diagram
$$
 \xymatrix  @R=3.0pc @C=.8pc
 {
 & \C^{G\times G}\ar[d]^{\rho_{G,G}}
 \ar[rr]^{\sigma_1} & & \C^{G\times G}\ar[d]_{\rho_{G,G}}
 \ar@/^3ex/[ddr]^{\widetilde{\mu}}
 & \\
 & \C^G\odot \C^G\ar[rr]^{\sigma\odot 1_{\C^G}} & & \C^G\odot \C^G \ar[dr]^{\mu}& \\
 \C^G
 \ar@/^3ex/[uur]^{\widetilde{\varkappa}}
 \ar@/_3ex/[ddr]_{\widetilde{\varkappa}}
 \ar[ur]^{\varkappa}\ar[dr]_{\varkappa}\ar[rr]^{\e} &  & \C\ar[rr]^{\iota} & & \C^G \\
 & \C^G\odot \C^G\ar[rr]^{1_{\C^G}\odot\sigma} & &
 \C^G\odot \C^G \ar[ur]_{\mu}& \\
 & \C^{G\times G}\ar[u]_{\rho_{G,G}}\ar[rr]^{\sigma_2} & & \C^{G\times G}\ar[u]^{\rho_{G,G}}
 \ar@/_3ex/[uur]_{\widetilde{\mu}}
 &
 }
$$
where the mappings $\sigma_1$ and $\sigma_2$ are defined by the identities
$$
\sigma_1(v)(s,t):=v(s^{-1},t),\qquad \sigma_2(v)(s,t):=v(s,t^{-1})
$$
Obviously, all the triangles and quadrangles siding with the borders of this
picture, are commutative here. Hence, to prove the commutativity of the two
inner pentagons (i.e. the commutativity of \eqref{AX-antipode-O(G)}) it is
sufficient to check the commutativity of the diagram arising after throwing out
the vertexes $\C^G\odot \C^G$:
$$
 \xymatrix  @R=3.0pc @C=.8pc
 {
 & \C^{G\times G}
 \ar[rr]^{\sigma_1} & & \C^{G\times G}
 \ar@/^1ex/[dr]^{\widetilde{\mu}}
 & \\
 \C^G
 \ar@/^1ex/[ur]^{\widetilde{\varkappa}}
 \ar@/_1ex/[dr]_{\widetilde{\varkappa}}
 \ar[rr]^{\e} &  & \C\ar[rr]^{\iota} & & \C^G \\
 & \C^{G\times G}\ar[rr]^{\sigma_2} & & \C^{G\times G}
 \ar@/_1ex/[ur]_{\widetilde{\mu}}
 &
 }
$$
This is done by the direct calculation: for any function $u\in\C^G$ its final
image after motion through the diagram is the function $u(1_G)\cdot
1_{\C^G}\in\C^G$. Indeed,
 \bit{

\item[--] if we move through the upper arrows, we obtain the following:
  \begin{flalign*}
&u\in \C^G &\\
\mapsto\qquad &\widetilde{\varkappa}(u) \in \C^{G\times G},\qquad
\widetilde{\varkappa}(u)(s,t)=u(s\cdot t) &\\
\mapsto\qquad &\sigma_1(\widetilde{\varkappa}(u)) \in \C^{G\times G},\qquad
\sigma_1(\widetilde{\varkappa}(u))(s,t)=
\widetilde{\varkappa}(u)(s^{-1},t)=u(s^{-1}\cdot t) &\\
\mapsto\qquad &\widetilde{\mu}\big(\sigma_1(\widetilde{\varkappa}(u))\Big) \in
\C^G,\qquad \widetilde{\mu}\big(\sigma_1(\widetilde{\varkappa}(u))\Big)(s)=
\sigma_1(\widetilde{\varkappa}(u))(s,s)=
\widetilde{\varkappa}(u)(s^{-1},s)=u(s^{-1}\cdot s)=u(1_G) &
 \end{flalign*}

\item[--] if we move aflat through the center of the pentagon, we obtain the
same:
  \begin{flalign*}
&u\in \C^G &\\
\mapsto\qquad &\e(u) \in \C,\qquad
\e(u)=u(1_G)\cdot 1_\C &\\
\mapsto\qquad &\iota(\e(u)) \in \C^G,\qquad \iota(\e(u))=u(1_G)\cdot 1_{\C^G} &
 \end{flalign*}

\item[--] and if we move through the lower arrows, we come to the same result:
  \begin{flalign*}
&u\in \C^G &\\
\mapsto\qquad &\widetilde{\varkappa}(u) \in \C^{G\times G},\qquad
\widetilde{\varkappa}(u)(s,t)=u(s\cdot t) &\\
\mapsto\qquad &\sigma_2(\widetilde{\varkappa}(u)) \in \C^{G\times G},\qquad
\sigma_2(\widetilde{\varkappa}(u))(s,t)=
\widetilde{\varkappa}(u)(s,t^{-1})=u(s\cdot t^{-1}) &\\
\mapsto\qquad &\widetilde{\mu}\big(\sigma_2(\widetilde{\varkappa}(u))\Big) \in
\C^G,\qquad \widetilde{\mu}\big(\sigma_2(\widetilde{\varkappa}(u))\Big)(s)=
\sigma_1(\widetilde{\varkappa}(u))(s,s)=
\widetilde{\varkappa}(u)(s,s^{-1})=u(s\cdot s^{-1})=u(1_G) &
 \end{flalign*}
 }\eit
 \epr

\subsection{Sweedler's notations and the stereotype approximation property}

A useful instrument for proving results in the theory of Hopf algebra are
Sweedler's notations \cite{Kassel}. This technique can be applied in the
stereotype theory as well, at least in situations, where a given stereotype
Hopf algebra $H$, being considered as a stereotype space possesses the
stereotype approximation property (see \cite{Akbarov}).

The following result explains this:

\btm\label{TH:ster-sweedler} If $H$ is an injective (respectively, a
projective) stereotype coalgebra with the stereotype approximation property,
then for any $x\in H$
 \bit
\item[(i)] the comultiplication $\varkappa(x)$ can be approximated in the
topology of $H$ by the finite sums of the form
$$
\sum_{i=1}^n x_i'\odot x_i''\qquad \left(\sum_{i=1}^n x_i'\circledast
x_i''\right)
$$

\item[(ii)] the identity
 $$
\langle\varkappa(x),\alpha\circledast\beta\rangle=0\qquad
\Big(\langle\varkappa(x),\alpha\odot\beta\rangle=0\Big),\qquad \alpha,\beta\in
H^\star
 $$
is equivalent to the identity
$$
\varkappa(x)=0
$$
 \eit
 \etm
\bpr This follows from the definition of the stereotype approximation \cite[\S
9]{Akbarov}. \epr

If now $H$, say, is an injective stereotype coalgebra with the stereotype
approximation property, then for each element $x\in H$ the symbol
$\sum_{(x)}x'\odot x''$ denotes the class of nets of the form
$\sum_{i=1}^{n_{\nu}} x_{\nu,i}'\odot x_{\nu,i}''$ tending to $\varkappa(x)$:
 \beq\label{def-sweedler-1}
\varkappa(x)\underset{\infty\gets\nu}{\longleftarrow}\sum_{i=1}^{n_{\nu}}
x_{\nu,i}'\odot x_{\nu,i}''
 \eeq
and the record
 \beq\label{def-sweedler}
\varkappa(x)=\sum_{(x)} x'\odot x''
 \eeq
should be understood as follows: the right side denotes one of those nets (but
without indices), and the arrow is replaced by the equality.

The formulas for comultiplication and antipode, and the others, like
$$
\varkappa(\sigma(x))=\sum_{(x)} \sigma(x'')\odot \sigma(x')
$$
are interpreted as follows: for every net $\sum_{i=1}^{n_{\nu}} x_{\nu,i}'\odot
x_{\nu,i}''$ the condition \eqref{def-sweedler-1} automatically implies the
condition
 $$
\varkappa(\sigma(x))\underset{\infty\gets\nu}{\longleftarrow}\sum_{i=1}^{n_{\nu}}
\sigma(x_{\nu,i}'')\odot \sigma(x_{\nu,i}')
 $$
The proof can also be conducted with the help of the record
\eqref{def-sweedler}:
 \begin{multline*}
\langle \varkappa(\sigma(x)),\alpha\circledast\beta\rangle=\langle
\sigma(x),\alpha*\beta\rangle=\langle
x,\sigma^\star(\alpha*\beta)\rangle=\langle
x,\sigma^\star(\beta)*\sigma^\star(\alpha)\rangle=\langle
\varkappa(x),\sigma^\star(\beta)\circledast\sigma^\star(\alpha)\rangle=\\=
\left\langle \sum_{(x)} x'\odot x'',
\sigma^\star(\beta)\circledast\sigma^\star(\alpha)\right\rangle=\sum_{(x)}
\langle x'\odot x'', \sigma^\star(\beta)\circledast\sigma^\star(\alpha)\rangle=
\sum_{(x)} \langle x', \sigma^\star(\beta)\rangle\cdot \langle
x'',\sigma^\star(\alpha)\rangle=\\= \sum_{(x)} \langle \sigma(x'),
\beta\rangle\cdot \langle \sigma(x''),\alpha\rangle= \sum_{(x)} \langle
\sigma(x'')\odot\sigma(x'), \alpha\circledast\beta\rangle= \left\langle
\sum_{(x)} \sigma(x'')\odot\sigma(x'), \alpha\circledast\beta\right\rangle
 \end{multline*}

As an example of the application of these ``generalized'' Sweedler's notations
let us consider the following situation. In the theory of quantum groups the
verification of the diagram for antipode \eqref{AX-antipode}, i.e. the identity
 \beq\label{antipod-x-y}
\mu\Big((\sigma\otimes 1)(\varkappa(x))\Big)= \e(x)\cdot 1_H=
\mu\Big((1\otimes\sigma)(\varkappa(x))\Big),\qquad x\in H
 \eeq
often leads to some bulky computations. In those cases one remark, made by
A.~Van Daele in \cite{VanDaele-Bulletin} is useful. Being applied to stereotype
algebras it looks as follows:

\blm\label{LM:antipod-x-y} Suppose $H$ is a stereotype bialgebra (no matter,
projective, or injective) with the stereotype approximation property, $\sigma$
is its (continuous) antihomomorphism and the equalities \eqref{antipod-x-y} are
true for two elements $x\in H$ and $y\in H$. Then they are true for their
multiplication $x\cdot y$. \elm

\bpr Both those equalities for $x\cdot y$ are proved by direct computations,
for instance, the left one is obtained as follows (here the tensor product
$\otimes$ means $\circledast$ or $\odot$):
\begin{multline*}
\mu\Big((\sigma\otimes 1)(\varkappa(x\cdot y))\Big)= \mu\Big((\sigma\otimes
1)(\varkappa(x)\cdot\varkappa(y))\Big)= \mu\left((\sigma\otimes
1)\left(\sum_{(x)} x'\otimes x''\cdot \sum_{(y)} y'\otimes
y''\right)\right)=\\= \sum_{(x),(y)}\mu\Big((\sigma\otimes 1)\left( x'\cdot
y'\otimes x''\cdot y''\right)\Big)=\sum_{(x),(y)}\mu\Big(\sigma(x'\cdot
y')\otimes x''\cdot y''\Big)= \sum_{(x),(y)}\mu\Big(\sigma(y')\cdot
\sigma(x')\otimes x''\cdot y''\Big)=\\= \sum_{(x),(y)}\sigma(y')\cdot
\sigma(x')\cdot x''\cdot y''= \sum_{(y)}\sigma(y')\cdot\left(\sum_{(x)}
\sigma(x')\cdot x''\right)\cdot y''= \sum_{(y)}\sigma(y')\cdot \e(x)\cdot
1_H\cdot y''=\\=\e(x)\cdot 1_H\cdot\sum_{(y)}\sigma(y')\cdot  y''=\e(x)\cdot
1_H\cdot \e(y)\cdot 1_H\cdot=\e(x\cdot y)\cdot 1_H\cdot
\end{multline*}
 \epr

\subsection{Grouplike elements}

\bprop\label{PROP:grup-elementy} For an element $a\in H$ in a stereotype Hopf
algebra $H$ with the stereotype approximation property the following conditions
are equivalent:
 \bit
\item[(i)] $\varkappa(a)=a\odot a\qquad (\varkappa(a)=a\circledast a)$;

\item[(ii)] the functional $\langle a,\cdot\rangle:H^\star\to\C$ is
multiplicative:
 \beq\label{mult-grupp-elem}
\langle a,\alpha*\beta\rangle=\langle a,\alpha\rangle\cdot \langle a,
\beta\rangle
 \eeq

\item[(iii)] the operator $\M_a^\star:H^\star\to H^\star$ dual to the operator
of multiplication by element $a$,
$$
\M_a(x):=a\cdot x
$$
is a homomorphism of stereotype algebra $H^\star$.
 \eit
 \eprop
\bpr Obviously, (i) and (ii) are equivalent. Let us show that (i)
$\Longleftrightarrow$ (iii). If $a$ satisfies (i), then
\begin{multline*}
\langle x, \M_a^\star(\alpha*\beta)\rangle=\langle \M_a(x),
\alpha*\beta\rangle= \langle a\cdot x, \alpha*\beta\rangle=\langle
\varkappa(a\cdot x), \alpha\circledast\beta\rangle=\langle \varkappa(a)\cdot
\varkappa(x), \alpha\circledast\beta\rangle=\\=\left\langle a\odot a\cdot
\sum_{(x)} x'\odot x'', \alpha\circledast\beta\right\rangle= \left\langle
\sum_{(x)} (a\cdot x')\odot (a\cdot x''), \alpha\circledast\beta\right\rangle=
\sum_{(x)} \langle a\cdot x',\alpha\rangle \cdot \langle a\cdot
x'',\beta\rangle=\\= \sum_{(x)} \langle x',\M_a^\star(\alpha)\rangle \cdot
\langle x'',\M_a^\star(\beta)\rangle= \left\langle \sum_{(x)} x'\odot x'',
\M_a^\star(\alpha)\circledast \M_a^\star(\beta)\right\rangle= \langle
\varkappa(x), \M_a^\star(\alpha)\circledast \M_a^\star(\beta)\rangle=\\=
\langle x, \M_a^\star(\alpha)* \M_a^\star(\beta)\rangle
\end{multline*}
This is true for every $x\in H$, therefore
 \beq\label{M_z^star(alpha*beta)}
\M_a^\star(\alpha*\beta)=\M_a^\star(\alpha)* \M_a^\star(\beta)
 \eeq
On the contrary, if \eqref{M_z^star(alpha*beta)} is true, then
\begin{multline*}
\langle\varkappa(a),\alpha\circledast\beta\rangle= \langle
a,\alpha*\beta\rangle=\langle 1,\M_a^\star(\alpha*\beta)\rangle= \langle
1,\M_a^\star(\alpha) * \M_a^\star(\beta)\rangle=\langle
\varkappa(1),\M_a^\star(\alpha) \circledast \M_a^\star(\beta)\rangle=\\=
\langle 1\odot 1,\M_a^\star(\alpha) \circledast \M_a^\star(\beta)\rangle=
\langle 1,\M_a^\star(\alpha)\rangle\cdot \langle 1, \M_a^\star(\beta)\rangle=
\langle a,\alpha\rangle\cdot \langle a, \beta\rangle=\langle a\odot
a,\alpha\circledast \beta\rangle
\end{multline*}
and since this is true for each $\alpha$ and $\beta$, we obtain (i).
 \epr

An element $a$ of an injective (resp., projective) stereotype Hopf algebra $H$
is called {\it grouplike element}\index{grouplike element}, if $a\ne 0$ and $a$
satisfies the conditions (i)-(ii)-(iii) of Proposition
\ref{PROP:grup-elementy}. The set of grouplike elements in $H$ is denoted by
$\G(H)$. As in the pure algebraic situation (see \cite{Sweedler}), $\G(H)$ is a
group with respect to the multiplication in $H$, since it possesses the
following properties:
 {\it
 \bit
\item[$1^\circ$.] $\forall a\in\G(H)\quad \e(a)=1$.

\item[$2^\circ$.] $\forall a\in\G(H)\quad a^{-1}=\sigma(a)\in\G(H)$.

\item[$3^\circ$.] $\forall a,b\in\G(H)\quad a\cdot b\in\G(H)$.
 \eit }
\bpr $1^\circ$ is proved by applying the comultiplication axiom:
$$
1\odot a=\ll_H^{-1}(a)=\eqref{komonoid}=
(\e\odot\id_H)(\varkappa(a))=(\e\odot\id_H)(a\odot a )=\e(a)\odot a
\quad\Longrightarrow\quad \e(a)=1
$$
In $2^\circ$ we need to apply the fact that $\sigma^\star$ is an
antihomomorphism: on the one hand,
$\varkappa(\sigma(a))=\sigma(a)\odot\sigma(a)$, since
 \begin{multline*}
\langle\varkappa(\sigma(a)),\alpha\circledast\beta\rangle= \langle
\sigma(a),\alpha*\beta\rangle=\langle a,\sigma^\star(\alpha*\beta)\rangle=
\langle a,\sigma^\star(\beta)*\sigma^\star(\alpha)\rangle= \langle
\varkappa(a),\sigma^\star(\beta)\circledast\sigma^\star(\alpha)\rangle=\\=
\langle a\odot a,\sigma^\star(\beta)\circledast\sigma^\star(\alpha)\rangle=
\langle a,\sigma^\star(\beta)\rangle\cdot\langle a,\sigma^\star(\alpha)\rangle=
\langle \sigma(a),\beta\rangle\cdot\langle \sigma(a),\alpha\rangle= \langle
\sigma(a)\odot\sigma(a),\alpha\circledast\beta\rangle
 \end{multline*}
And, on the other hand, $\e(\sigma(a))=\langle\sigma(a),1_H\rangle=\langle
a,\sigma^\star(1_H)\rangle= \langle a, 1_H\rangle=\e(a)=1$. Together these
conditions mean that $\sigma(a)\in\G(H)$. Apart from that,
$$
\sigma(a)\cdot a=\mu(\sigma(a)\odot a)=\mu\Big((\sigma\odot \id_H)(a\odot
a)\Big)=\mu\Big((\sigma\odot
\id_H)(\varkappa(a))\Big)=\eqref{AX-antipode}=\e(a)\cdot 1_H=(\text{property
$1^\circ$})=1_H
$$
and, similarly, $a\cdot \sigma(a)=1_H$. Hence, $\sigma(a)=a^{-1}$.

Finally, $3^\circ$: if $a,b\in\G(H)$, then, first, $\varkappa(a\cdot
b)=\varkappa(a)\cdot\varkappa(b)=a\odot a\cdot b\odot b=(a\cdot b)\odot(a\cdot
b)$, and, second, $\e(a\cdot b)=\e(a)\cdot\e(b)=1$. Therefore, $a\cdot
b\in\G(H)$.
 \epr

Recall that an element $a$ of an algebra $A$ is called {\it
central}\index{central element}, if it commutes with all other elements of $A$:
 $$
\forall x\in A\quad a\cdot x=x\cdot a
 $$

\bprop If $a$ is a grouplike and in addition a central element in a Hopf
algebra $H$, then
 \bit
\item[1)] the following identities hold:
 \begin{align}
 \label{M_z-s-M_z=s}
& \M_a\circ\; \sigma\circ \M_a= \sigma, &&
\sigma\circ\M_a=\M_{a^{-1}}\circ\sigma
 \\
 \label{M_z*-s-M_z*=s}
& \M_a^\star\circ\; \sigma^\star\circ \M_a^\star= \sigma^\star, &&
\sigma^\star\circ\M_a^\star=\M_{a^{-1}}^\star\circ\sigma^\star
 \end{align}

\item[2)] if $H$ has the stereotype approximation property, then the following
identities hold:
 \beq\label{(M_z*)^k-1}
\varkappa\Big(\M_a^\star(\alpha)\Big)=\sum_{(\alpha)}\M_a^\star(\alpha')\circledast\alpha''=
\sum_{(\alpha)}\alpha'\circledast\M_a^\star(\alpha''),\qquad \alpha\in H^\star
 \eeq
 \beq\label{(M_z*)^k}
\varkappa\Big((\M_a^\star)^{i+j}(\alpha)\Big)=\sum_{(\alpha)}(\M_a^\star)^i(\alpha')\circledast
(\M_a^\star)^j(\alpha''),\qquad \alpha\in H^\star,\quad i,j\in\N
 \eeq
 \eit
 \eprop
\bpr 1. For each $x\in H$ we have
$$
(\M_a\circ \sigma\circ \M_a)(x)=a\cdot\sigma(a\cdot
x)=a\cdot\sigma(x)\cdot\sigma(a)=a\cdot\sigma(x)\cdot a^{-1}=\sigma(x)\cdot
a\cdot a^{-1}=\sigma(x)
$$
2. For any $u,v\in H$, $\alpha\in H^\star$
 \begin{multline*}
\langle u\odot v, \varkappa(\M_a^\star(\alpha))\rangle= \langle u\cdot v,
\M_a^\star(\alpha)\rangle= \langle a\cdot u\cdot v, \alpha\rangle= \langle
(a\cdot u)\odot v, \varkappa(\alpha)\rangle=\left\langle (a\cdot u)\odot v,
\sum_{(\alpha)} \alpha'\circledast\alpha''\right\rangle=\\= \sum_{(\alpha)}
\langle (a\cdot u)\odot v, \alpha'\circledast\alpha''\rangle= \sum_{(\alpha)}
\langle (a\cdot u), \alpha'\rangle\cdot\langle v,\alpha''\rangle=
\sum_{(\alpha)} \langle u, \M_a^\star(\alpha')\rangle\cdot\langle
v,\alpha''\rangle= \sum_{(\alpha)} \langle u\odot v,
\M_a^\star(\alpha')\circledast \alpha''\rangle=\\= \left\langle u\odot v,
\sum_{(\alpha)} \M_a^\star(\alpha')\circledast \alpha''\right\rangle
 \end{multline*}
By Theorem \ref{TH:ster-sweedler} this means that the first equality in
\eqref{(M_z*)^k-1} holds. The rest equalities are proved by analogy. \epr

\section{Stein manifolds: rectangles in $\mathcal O(M)$ and rhombuses in
${\mathcal O}^\star(M)$}\label{rectangles-buses}

Here we discuss some special properties of the space of holomorphic functions
on a complex manifold. For illustration purposes it is convenient for us to use
the condition of holomorphic separability, so we formulate our results only for
Stein manifolds. We use terminology from \cite{Shabat,Taylor,Grauert-Remmert}.

\subsection{Stein manifolds}

Let $M$ be a complex manifold. Symbol $\mathcal O(M)$ denotes the algebra of
all holomorphic functions on $M$ (with the usual pointwise algebraic operations
and topology of uniform convergence on compact sets in $M$). It is well-known (
see \cite{Taylor}), that as a topological vector space, $\mathcal O(M)$ is a
Montel space.

A manifold $M$ is called a {\it Stein manifold}\index{Stein!manifold}
\cite{Shabat}, if the following three conditions are fulfilled:
 \bit{
\item[1)] {\it holomorphic separability:}\label{holom-separability} for any two
points $x,y\in M$, $x\ne y$, there exists a function $u\in\mathcal O(M)$ such
that
$$
u(x)\ne u(y)
$$
\item[2)] {\it holomorphic uniformization:} for any point $x\in M$ there exist
functions $u_1,...,u_n\in\mathcal O(M)$, forming local coordinates of the
manifold $M$ in a neighborhood of $x$;

\item[3)] {\it holomorphic convexity:} for any compact set $K\subseteq M$ its
{\it holomorphically convex hull}, i.e. a the set
$$
\hat{K}=\{x\in M: \ \forall u\in\mathcal O(M)\ |u(x)|\le\max_{y\in K}|u(y)| \}
$$
is a compact set in $M$.
 }\eit

The complex space $\C^n$ and various domains of holomorphy in $\C^n$ are
examples of Stein manifolds. A simplest example of a complex manifold which is
not a Stein manifold is {\it complex torus}\index{complex!torus}, i.e. the
quotient group of the complex plane $\C$ over the lattice $\Z+i\Z$.

\subsection{Outer envelopes on $M$ and rectangles in $\mathcal O(M)$}

\paragraph{Operations $\protect\BSQ$ and $\protect\SQ$.}

Here we shall spend some time on studying real functions $f$ on a manifold $M$,
bounded by 1 from below,
$$
f\ge 1,
$$
i.e. $f$ having range in the interval $[1,+\infty)$. Certainly, we shall use
record $f:M\to[1;+\infty)$ for such functions. As usual, we call a function
$f:M\to[1;+\infty)$ {\it locally bounded}, if for each point $x\in M$ one can
find a neighborhood $U\owns x$ such that
$$
\sup_{y\in U} |f(y)|<\infty
$$
Since $f$ is bounded by 1 from below, this condition is equivalent to the
condition
$$
\sup_{y\in U} f(y)<\infty
$$

 \bprop
For each locally bounded function $f:M\to[1;+\infty)$ the formula
 \beq\label{f^blacksquare}
f^\text{\BSQ} :=\{u\in \mathcal O(M): \quad \forall x\in M\quad |u(x)|\le
f(x)\}
 \eeq
defines an absolutely convex set of functions $f^\text{\BSQ}\subseteq\mathcal
O(M)$, containing the identity function:
$$
1\in f^\text{\BSQ}.
$$
\eprop \bpr The set $f^\text{\BSQ}$ is compact since it is closed and bounded
in the Montel space $\mathcal O(M)$. \epr

 \bprop\label{PROP-nepr-vnesh-ogib}
For any bounded set of functions $D\subseteq \mathcal O(M)$ containing the
identity function,
$$
1\in D,
$$
the formula
 \beq\label{fi^square}
 D^\text{\SQ} (x):=\sup_{u\in D} |u(x)|, \qquad x\in M
 \eeq
defines a continuous real function $D^\text{\SQ}:M\to\R$ bounded by 1 from
below:
$$
 D^\text{\SQ}\ge 1.
$$
\eprop \bpr Note from the very beginning that $D$ can be considered as compact.
Consider for this its closure $\overline{D}$. Since $D$ is bounded in a Montel
space $\mathcal O(M)$, $\overline{D}$ is compact in $\mathcal O(M)$. The
mapping $u\mapsto \delta^x(u)=u(x)$ is continuous, so the image of the closure
$\delta^x(\overline{D})$ is contained in the closure of image
$\overline{\delta^x( D)}$
$$
\delta^x(\overline{D})\subseteq \overline{\delta^x( D)},
$$
As a corollary, we obtain the following chain of inequalities:
$$
\sup_{\lambda\in \delta^x( D)} |\lambda|\le \sup_{\lambda\in
\delta^x(\overline{D})} |\lambda|\le \sup_{\lambda\in \overline{\delta^x( D)}}
|\lambda|=\sup_{\lambda\in \delta^x( D)} |\lambda|
$$
Therefore,
$$
\sup_{\lambda\in \delta^x(\overline{D})} |\lambda|=\sup_{\lambda\in \delta^x(
D)} |\lambda|
$$
and the functions $\overline{D}^\text{\SQ}$ and $D^\text{\SQ}$ coincide:
$$
\overline{D}^\text{\SQ} (x)=\sup_{u\in \overline{D}}
|\delta^x(u)|=\sup_{\lambda\in \delta^x(\overline{D})}
|\lambda|=\sup_{\lambda\in \delta^x( D)} |\lambda|=\sup_{u\in D} |\delta^x(u)|=
D^\text{\SQ} (x)
$$

Thus it is sufficient to consider the case when $D$ is compact. Let us take an
arbitrary compact set $K\subseteq M$ and consider the space $C(K)$ of
continuous functions on $K$ (with the usual topology of uniform convergence on
$K$). The restriction mapping $u\in \mathcal O(M)\mapsto u|_K\in C(K)$ is
continuous from $\mathcal O(M)$ into $C(K)$, so the image $D|_K$ of a compact
$D$ in $\mathcal O(M)$ must be compact in $C(K)$. Hence, by the Arcela theorem
$D|_K$, is pointwise bounded and equicontinuous on $K$. Therefore the function
$$
 D^\text{\SQ} (x):=\sup_{u\in D} |u(x)|, \qquad x\in K
$$
is continuous on $K$. This is true for every compact set $K$ in $M$, so we
obtain that $D^\text{\SQ}$ is continuous on $M$. \epr

\bigskip

\centerline{\bf Properties of operations $\text{\BSQ}$ and $\text{\SQ}$:}
 \begin{align}\label{sv-0}
f&\le g\quad\Longrightarrow\quad f^\text{\BSQ}\subseteq g^\text{\BSQ}, &
 D &\subseteq E\quad\Longrightarrow\quad D^\text{\SQ}\le
 E^\text{\SQ}
 \end{align}
 \begin{align}\label{sv-1}
(f^\text{\BSQ})^\text{\SQ} &\le f, & D &\subseteq ( D^\text{\SQ})^\text{\BSQ}
 \end{align}
 \begin{align}\label{sv-2}
((f^\text{\BSQ})^\text{\SQ})^\text{\BSQ} &=f^\text{\BSQ}, & ((
D^\text{\SQ})^\text{\BSQ})^\text{\SQ} &= D^\text{\SQ}
 \end{align}
\bpr Properties \eqref{sv-0} and \eqref{sv-1} are evident, and \eqref{sv-2}
follows from them:
$$
\left\{\begin{matrix}(f^\text{\BSQ})^\text{\SQ} \le f
\quad\Longrightarrow\quad\text{(apply the operation
$\text{\BSQ}$)}\quad\Longrightarrow\quad
((f^\text{\BSQ})^\text{\SQ})^\text{\BSQ} \subseteq f^\text{\BSQ} \\ D \subseteq
( D^\text{\SQ})^\text{\BSQ} \quad\Longrightarrow\quad \text{(substitution:
$D=f^\text{\BSQ}$)}\quad\Longrightarrow\quad
f^\text{\BSQ}\subseteq((f^\text{\BSQ})^\text{\SQ})^\text{\BSQ}\end{matrix}\right\}
\quad\Longrightarrow\quad ((f^\text{\BSQ})^\text{\SQ})^\text{\BSQ} =
f^\text{\BSQ}
$$
$$
\left\{\begin{matrix} D \subseteq ( D^\text{\SQ})^\text{\BSQ}
\quad\Longrightarrow\quad\text{(apply the operation
$\text{\SQ}$)}\quad\Longrightarrow\quad D^\text{\SQ}\le
(( D^\text{\SQ})^\text{\BSQ})^\text{\SQ} \\
(f^\text{\BSQ})^\text{\SQ} \le f\quad\Longrightarrow\quad \text{(substitution:
$f= D^\text{\SQ}$)}\quad\Longrightarrow\quad ((
D^\text{\SQ})^\text{\BSQ})^\text{\SQ} \le
 D^\text{\SQ}
\end{matrix}\right\}\quad\Longrightarrow\quad
(( D^\text{\SQ})^\text{\BSQ})^\text{\SQ} = D^\text{\SQ}
$$
\epr

\paragraph{Outer envelopes on $M$.}

Let us introduce the following notations:
 \begin{align}\label{sv-3}
f^{\text{\BSQ}\text{\SQ}} &:=(f^\text{\BSQ})^\text{\SQ} &
 D^{\text{\SQ}\text{\BSQ}} &:=( D^\text{\SQ})^\text{\BSQ}
 \end{align}
Then  \eqref{sv-0}, \eqref{sv-1}, and \eqref{sv-2} imply
 \begin{align}\label{sv-4}
f^{\text{\BSQ}\text{\SQ}} &\le f, & D &\subseteq D^{\text{\SQ}\text{\BSQ}}
 \end{align}
 \begin{align}\label{sv-0-0}
f&\le g\quad\Longrightarrow\quad f^{\text{\BSQ}\text{\SQ}}\le
g^{\text{\BSQ}\text{\SQ}}, & D &\subseteq E\quad\Longrightarrow\quad
 D^{\text{\SQ}\text{\BSQ}}\subseteq E^{\text{\SQ}\text{\BSQ}}
 \end{align}
 \begin{align}\label{sv-6}
(f^{\text{\BSQ}\text{\SQ}})^{\text{\BSQ}\text{\SQ}} &=f^{\text{\BSQ}\text{\SQ}}
& ( D^{\text{\SQ}\text{\BSQ}})^{\text{\SQ}\text{\BSQ}} &=
D^{\text{\SQ}\text{\BSQ}}
 \end{align}

Let us call a locally bounded function $g:M\to[1,+\infty)$,
 \bit
\item[---] an {\it outer envelope for a bounded set}\index{outer envelope}
$D\subseteq \mathcal O(M)$, $1\in D$, if
$$
g= D^\text{\SQ}
$$

\item[---] an {\it outer envelope for a locally bounded function}
$f:M\to[1;+\infty)$, if
$$
g=f^{\text{\BSQ}\text{\SQ}}
$$

\item[---] an {\it outer envelope on $M$}, if it satisfies the following
equivalent conditions:
 \bit
\item[(i)] $g$ is an outer envelope for some bounded set $D\subseteq \mathcal
O(M)$, $1\in D$,
$$
g= D^\text{\SQ}
$$

\item[(ii)] $g$ is an outer envelope for some locally bounded function
$f:M\to[1;+\infty)$,
$$
g=f^{\text{\BSQ}\text{\SQ}}
$$

\item[(iii)] $g$ is an outer envelope for some for itself:
$$
g^{\text{\BSQ}\text{\SQ}} = g
$$
 \eit
 \eit
\bpr The equivalence of conditions (i), (ii), (iii) requires some comments.

$(i)\Longrightarrow (ii)$. If $g$ is an outer envelope for some bounded set
$D$, i.e. $g= D^\text{\SQ}$, then $g= D^\text{\SQ}=\eqref{sv-2}=((
D^\text{\SQ})^\text{\BSQ})^\text{\SQ}=( D^\text{\SQ})^{\text{\BSQ}\text{\SQ}}$,
i.e. $g$ is an outer envelope for the function $f= D^\text{\SQ}$.

$(ii)\Longrightarrow (iii)$. If  $g$  is an outer envelope for some function
$f$, i.e. $g=f^{\text{\BSQ}\text{\SQ}}$, then
$g^{\text{\BSQ}\text{\SQ}}=(f^{\text{\BSQ}\text{\SQ}})^{\text{\BSQ}\text{\SQ}}=\eqref{sv-6}=f^{\text{\BSQ}\text{\SQ}}=g$,
i.e. $g$  is an outer envelope for itself.

$(iii)\Longrightarrow (i)$. If $g$  is an outer envelope for itself, i.e.
$g=g^{\text{\BSQ}\text{\SQ}}=(g^\text{\BSQ})^\text{\SQ}$, then we put
$D=g^\text{\BSQ}$, and after that $g$ becomes an outer envelope for $D$:  $g=
D^\text{\SQ}$. \epr

\bigskip

\centerline{\bf Properties of outer envelopes:}
 {\it
 \bit
\item[(i)] Every outer envelope $g$ on $M$ is a continuous (positive) function
on $M$.

\item[(ii)] For any locally bounded function $f:M\to[1;+\infty)$ its outer
envelope $f^{\text{\BSQ}\text{\SQ}}$ is the greatest outer envelope on $M$,
majorized by $f$:
 \bit
\item[(a)] $f^{\text{\BSQ}\text{\SQ}}$ is an outer envelope on $M$, majorized
by $f$:
$$
f^{\text{\BSQ}\text{\SQ}} \le f
$$
\item[(b)] if $g$ is another outer envelope on $M$, majorized by $f$,
$$
g\le f
$$
then $g$ is majorized by $f^{\text{\BSQ}\text{\SQ}}$:
$$
g\le f^{\text{\BSQ}\text{\SQ}}
$$
 \eit
\eit
 }

\paragraph{Rectangles in $\mathcal O(M)$.}

We call a set $ E\subseteq \mathcal O(M)$, $1\in E$,
 \bit
\item[---] a {\it rectangle, generated by a locally bounded
function}\index{rectangle} $f:M\to[1;+\infty)$, if
$$
 E=f^\text{\BSQ}
$$

\item[---] a {\it rectangle, generated by a bounded set} $D\subseteq \mathcal
O(M)$, $1\in D$, if
$$
 E= D^{\text{\SQ}\text{\BSQ}}
$$

\item[---] a {\it rectangle} in $\mathcal O(M)$, if the following equivalent
conditions hold:
 \bit
\item[(i)] $E$ is a rectangle, generated by some locally bounded function
$f:M\to[1;+\infty)$
$$
 E=f^\text{\BSQ}
$$

\item[(ii)] $E$  is a rectangle, generated by  some bounded set $D\subseteq
\mathcal O(M)$, $1\in D$,
$$
 E= D^{\text{\SQ}\text{\BSQ}}
$$

\item[(iii)] $E$  is a rectangle, generated by itself:
$$
 E= E^{\text{\SQ}\text{\BSQ}}
$$
 \eit
 \eit
\bpr The equivalence of conditions (i), (ii), (iii) is proved in the same way
as in the case of outer envelopes. \epr

\bigskip

\centerline{\bf Properties of rectangles:}
 {\it
 \bit
\item[(i)] Every rectangle $E$ in $\mathcal O(M)$ is an absolutely convex
compact set in $\mathcal O(M)$.

\item[(ii)] For any bounded set $D\subseteq \mathcal O(M)$ the rectangle
$D^{\text{\SQ}\text{\BSQ}}$ is the smallest rectangle in $\mathcal O(M)$,
containing $D$:
 \bit
\item[(a)] $D^{\text{\SQ}\text{\BSQ}}$ is a rectangle in $\mathcal O(M)$,
containing $D$:
$$
 D\subseteq D^{\text{\SQ}\text{\BSQ}}
$$
\item[(b)] if $E$ is another rectangle in $\mathcal O(M)$, containing $D$,
$$
 D\subseteq E,
$$
then $E$ contains $D^{\text{\SQ}\text{\BSQ}}$:
$$
 D^{\text{\SQ}\text{\BSQ}}\subseteq E
$$
 \eit

\item[(iii)] The rectangles in $\mathcal O(M)$ form a fundamental system of
compact sets in $\mathcal O(M)$: every compact set $D$ in $\mathcal O(M)$ is
contained in some rectangle.

 \eit
 }

\btm\label{TH:f<->D} The formulas
 \begin{align}
 D &=f^\text{\BSQ}, & f &= D^\text{\SQ}
 \end{align}
establish a bijection between outer envelopes $f$ on $M$ and rectangles $D$ in
$\mathcal O(M)$. \etm
 \bpr
By definition of outer envelopes and rectangles, the operations $f\mapsto
f^\text{\BSQ}$ and $D\mapsto D^\text{\SQ}$ turn outer envelopes into rectangles
and rectangles into outer envelopes. Moreover, these operations are mutually
inverse on those two classes -- if $f$ is an outer envelope, then
$f^{\text{\BSQ}\text{\SQ}}=f$, i.e. the composition of the operations
$\text{\BSQ}$ and $\text{\SQ}$ gives back to the initial function:
$$
f\mapsto f^\text{\BSQ}\mapsto f^{\text{\BSQ}\text{\SQ}}=f
$$
Just like this the composition of the operations $\text{\SQ}$ and $\text{\BSQ}$
returns to the set $D$, if initially it was chosen as a rectangle:
$$
D\mapsto D^\text{\SQ}\mapsto D^{\text{\SQ}\text{\BSQ}}=D
$$
 \epr

\subsection{Lemma on polars}

Recall that {\it polar}\index{polar} of a set $A$ in a locally convex space $X$
is the set $A^\circ$ of linear continuous functionals $f:X\to\C$, bounded by 1
on $A$:
$$
A^\circ=\{f\in X^\star:\quad\sup_{x\in A}|f(x)|\le 1\}
$$
If $X$ is a stereotype space and $A$ a subset in the dual space $X^\star$, then
because of the equality $(X^\star)^\star=X$ it is convenient to consider the
polar $A^\circ\subseteq (X^\star)^\star$ as a subset in $X$. We denote this set
by $^\circ\kern-2pt A$ (and again call it {\it polar} of $A$):
$$
{^\circ\kern-2pt A}=\{x\in X:\quad\sup_{f\in A}|f(x)|\le 1\}
$$
An important observation for us is that if $A\subseteq X^\star$, and we take
its polar $^\circ\kern-2pt A$, and after that ``the polar of the polar''
$(^\circ\kern-2pt A)^\circ$ -- this set is called {\it bipolar}\index{bipolar}
of $A$ -- then it turns out that $(^\circ\kern-2pt A)^\circ$ is exactly the
closed absolutely convex hull of the set $A$ in $X^\star$ (i.e. the closure of
the set of linear combinations of the form $\sum_{i=1}^n\lambda_i\cdot a_i$,
where $a_i\in A$, $\sum_{i=1}^n|\lambda_i|\le 1$):
 \beq\label{bipolar}
(^\circ\kern-2pt A)^\circ=\cabsconv A
 \eeq
That is the essence the classical {\it theorem on bipolar}\index{theorem!on
bipolar} as applied to stereotype spaces.

In the special case when $A=D$ is a set in $\mathcal O(M)$, its polar $D^\circ$
is the set of analytical functionals $\alpha\in \mathcal O^\star(M)$, bounded
by 1 on $D$:
$$
 D^\circ=\{\alpha\in \mathcal O^\star(M):\quad
\sup_{u\in D}|\alpha(u)|\le 1\}
$$
On the contrary, if $A$ is a set of analytical functionals, $A\subseteq\mathcal
O^\star(M)$, then its polar $^\circ\kern-2pt A$ in $\mathcal O(M)$ is the set
of functions $u\in\mathcal O(M)$ on which all the functionals $\alpha\in A$ are
bounded by 1:
$$
{^\circ\kern-2pt A}=\{u\in {\mathcal O}(M):\quad \sup_{\alpha\in
A}|\alpha(u)|\le 1\}
$$

\blm[\bf on polars]\label{LM-o-polyarah} The passage to polars possesses the
following properties.
 \bit
\item[(a)] For every bounded set $D$ in $\mathcal O(M)$ containing the unit,
its outer envelope $D^\text{\SQ}$ is connected with its polar $D^\circ$ through
the identity
 \beq\label{frac_1_varPhi^square_x}
\frac{1}{D^\text{\SQ}(x)} =\max \{\lambda>0:\quad \lambda\cdot\delta^x\in
D^\circ\}
 \eeq

\item[(b)] For any locally bounded function $f:M\to[1;+\infty)$ the rectangle
$f^\text{\BSQ}$ is the polar of functionals $\frac{1}{f(x)}\cdot\delta^x$:
 \beq
f^\text{\BSQ} = {^\circ\kern-2pt\left\{\frac{1}{f(x)}\cdot\delta^x;\quad x\in
M\right\}}
 \eeq

\item[(c)] The polar of the rectangle $f^\text{\BSQ}$ is an absolutely convex
hull of functionals $\frac{1}{f(x)}\cdot\delta^x$:
 \beq\label{f^blacksquare^circ}
(f^\text{\BSQ})^\circ = \cabsconv \left\{\frac{1}{f(x)}\cdot\delta^x;\quad x\in
M\right\}
 \eeq
 \eit
 \elm
\bpr (a) For $\lambda>0$ we have:
$$
\lambda\cdot\delta^x\in D^\circ\quad\Longleftrightarrow\quad \sup_{u\in
D}|\lambda\cdot\delta^x(u)|\le 1 \quad\Longleftrightarrow\quad
 D^\text{\SQ}(x)=\sup_{u\in D}|\delta^x(u)|\le\frac{1}{\lambda}
\quad\Longleftrightarrow\quad \lambda\le \frac{1}{D^\text{\SQ}(x)}
$$
(b) is a reformulation of the definition of $f^\text{\BSQ}$:
 \begin{multline*}
u\in
f^\text{\BSQ}\quad\overset{\eqref{f^blacksquare}}{\Longleftrightarrow}\quad
\forall x\in M\quad |u(x)|=|\delta^x(u)|\le f(x)
\quad\Longleftrightarrow \\
\Longleftrightarrow\quad\sup_{x\in
M}\left|\frac{1}{f(x)}\cdot\delta^x(u)\right|\le 1\quad\Longleftrightarrow\quad
u\in {^\circ\left\{\frac{1}{f(x)}\cdot\delta^x;\quad x\in M\right\}}
 \end{multline*}
(c) follows from (b) and from the theorem on bipolar:
 $$
(f^\text{\BSQ})^\circ = \left({^\circ\left\{\frac{1}{f(x)}\cdot\delta^x;\quad
x\in M\right\}}\right)^\circ=\eqref{bipolar}=\cabsconv
\left\{\frac{1}{f(x)}\cdot\delta^x;\quad x\in M\right\}
 $$ \epr

\subsection{Inner envelopes on $M$ and rhombuses in $\mathcal O^\star(M)$}

In this subsection we shall study closed absolutely convex neighborhoods of
zero $\varDelta$ in $\mathcal O^\star(M)$, satisfying the following two
equivalent conditions:
 {\it
 \bit
\item[(A)]\label{opredelenie-Delta} the polar ${^\circ\kern-2pt \varDelta}$ of
$\varDelta$ contains the unit $1\in{\mathcal O}(M)$:
 \beq\label{uslovie-na-A-2}
1\in {^\circ\kern-2pt \varDelta}
 \eeq
\item[(B)] the value of all functionals $\alpha\in \varDelta$ on the unit
$1\in{\mathcal O}(M)$ does not exceed 1:
 \beq\label{uslovie-na-A-1}
\forall \alpha\in \varDelta\quad |\alpha(1)|\le 1
 \eeq
 \eit }\noindent
These conditions imply another one (which is not, however, equivalent to (A)
and (B)):
 {\it
 \bit
\item[(C)] a functional of the form $\lambda\cdot\delta^x$, where $\lambda>0$,
can belong to $\varDelta$ only if $\lambda\le 1$:
 \beq\label{uslovie-na-A}
\forall x\in M\quad\forall\lambda>0 \quad \Big(\lambda\cdot\delta^x\in
\varDelta\;\Longrightarrow\; \lambda\le 1\Big).
 \eeq
 \eit }
\bpr If $\lambda\cdot\delta^x\in \varDelta$, then $\forall u\in
{^\circ\varDelta}$\quad $|\lambda\cdot\delta^x(u)|\le 1$. In particular, for $u=1$
we have $|\lambda\cdot\delta^x(1)|=\lambda\cdot 1 \le 1$, i.e. $\lambda\le 1$.
 \epr

We say that a function $\ph:M\to(0,+\infty)$ is {\it locally separated from
zero}, if for each point $x\in M$ there exists a neighborhood $U\owns x$ such
that
$$
\inf_{y\in U} \ph(y)>0
$$

\paragraph{Operations $\protect\BLZ$ and $\protect\LZ$.}

 \bprop
If a function $\ph:M\to(0,1]$ is locally separated from zero, then the formula
 \beq\label{f^blacklozenge}
\ph^{\text{\BLZ}} :=\cabsconv\{\ph(x)\cdot\delta^x;\; x\in M\}
 \eeq
defines an absolutely convex neighborhood of zero in the space of analytical
functionals $\mathcal O^\star(M)$, satisfying
\eqref{uslovie-na-A-2}-\eqref{uslovie-na-A}. \eprop

\bpr The function $f(x)=\frac{1}{\ph(x)}$ is locally bounded and has values in
$[1;+\infty)$. Hence, by Lemma \ref{LM-o-polyarah},
$$
(f^\text{\BSQ})^\circ=\eqref{f^blacksquare^circ}=\cabsconv
\left\{\frac{1}{f(x)}\cdot\delta^x;\quad x\in M\right\}=\cabsconv
\left\{\ph(x)\cdot\delta^x;\quad x\in M\right\}
$$
This set is a closed absolutely convex neighborhood of zero in  $\mathcal
O^\star(M)$, since it is the polar of the compact set $f^\text{\BSQ}$ in
$\mathcal O(M)$. Besides this, since $f\ge 1$, the compact set $f^\text{\BSQ}$
contains the unit, so its polar $\cabsconv\left\{\ph(x)\cdot\delta^x;\quad x\in
M\right\}$ satisfies conditions (A),(B),(C) on page
\pageref{opredelenie-Delta}.
 \epr

 \bprop
For any closed absolutely convex neighborhood of zero $\varDelta$ in the space
of analytical functionals $\mathcal O^\star(M)$, satisfying
\eqref{uslovie-na-A-2}-\eqref{uslovie-na-A-1} the formula
 \beq\label{fi^lozenge}
\varDelta^{\text{\LZ}} (x):=\sup\{\lambda>0: \; \lambda\cdot\delta^x\in \varDelta\}
 \eeq
defines a continuous and locally separated from zero function
$\varDelta^{\text{\LZ}}:M\to(0;1]$.
 \eprop
\bpr Since $\varDelta$ is a neighborhood of zero in $\mathcal O^\star(M)$, its
polar $D={^\circ\kern-2pt \varDelta}$ is a compact in $\mathcal O^\star(M)$,
and from \eqref{uslovie-na-A-2} we have $1\in D$. By Proposition \ref{PROP-nepr-vnesh-ogib}, the outer envelope
$D^\text{\SQ}$ of this compact set is continuous and locally bounded. Hence,
the function
$$
\varDelta^{\text{\LZ}} (x):=\sup\{\lambda>0: \; \lambda\cdot\delta^x\in \varDelta=
D^\circ\}=\eqref{frac_1_varPhi^square_x}=\frac{1}{D^\text{\SQ}(x)}
$$
is continuous and locally separated from zero.
 \epr
The following properties are proved similarly with \eqref{sv-0}, \eqref{sv-1}
and \eqref{sv-2}.

\bigskip

\centerline{\bf Properties of the operations ${\text{\BLZ}}$ and ${\text{\LZ}}$:}
 \begin{align}\label{sv-0-int}
\ph&\le \psi\quad\Longrightarrow\quad \ph^{\text{\BLZ}}\subseteq
\psi^{\text{\BLZ}}, & \varDelta &\subseteq \varGamma\quad\Longrightarrow\quad
 \varDelta^{\text{\LZ}}\le \varGamma^{\text{\LZ}}
 \end{align}
 \begin{align}\label{sv-1-int}
\ph& \le (\ph^{\text{\BLZ}})^{\text{\LZ}}, & (
\varDelta^{\text{\LZ}})^{\text{\BLZ}}&\subseteq \varDelta
 \end{align}
 \begin{align}\label{sv-2-int}
((\ph^{\text{\BLZ}})^{\text{\LZ}})^{\text{\BLZ}} &=\ph^{\text{\BLZ}}, & ((
\varDelta^{\text{\LZ}})^{\text{\BLZ}})^{\text{\LZ}} &= \varDelta^{\text{\LZ}}
 \end{align}

\paragraph{Inner envelopes on $M$.}

Let us use the following supplementary notations:
 \begin{align}\label{sv-3-1}
\ph^{{\text{\BLZ}}\kern-0.5pt{\text{\LZ}}} &:=(\ph^{\text{\BLZ}})^{\text{\LZ}} &
 \varDelta^{{\text{\LZ}}\kern-0.5pt{\text{\BLZ}}} &:=( \varDelta^{\text{\LZ}})^{\text{\BLZ}}
 \end{align}
Formulas \eqref{sv-0}, \eqref{sv-1}, \eqref{sv-2} imply
 \begin{align}\label{sv-0-0-1}
\ph&\le \psi\quad\Longrightarrow\quad \ph^{{\text{\BLZ}}\kern-0.5pt{\text{\LZ}}}\le
\psi^{{\text{\BLZ}}\kern-0.5pt{\text{\LZ}}}, & \varDelta &\subseteq
\varGamma\quad\Longrightarrow\quad
 \varDelta^{{\text{\LZ}}\kern-0.5pt{\text{\BLZ}}}\subseteq \varGamma^{{\text{\LZ}}\kern-0.5pt{\text{\BLZ}}}
 \end{align}
 \begin{align}\label{sv-4-1}
\ph &\le \ph^{{\text{\BLZ}}\kern-0.5pt{\text{\LZ}}}, &
\varDelta^{{\text{\LZ}}\kern-0.5pt{\text{\BLZ}}} &\subseteq
 \varDelta
 \end{align}
 \begin{align}\label{sv-6-1}
(\ph^{{\text{\BLZ}}\kern-0.5pt{\text{\LZ}}})^{{\text{\BLZ}}\kern-0.5pt{\text{\LZ}}}
&=\ph^{{\text{\BLZ}}\kern-0.5pt{\text{\LZ}}} & (
\varDelta^{{\text{\LZ}}\kern-0.5pt{\text{\BLZ}}})^{{\text{\LZ}}\kern-0.5pt{\text{\BLZ}}}
&= \varDelta^{{\text{\LZ}}\kern-0.5pt{\text{\BLZ}}}
 \end{align}

Let us call a locally separated from zero function $\psi:M\to(0;1]$
 \bit
\item[---] an {\it inner envelope\index{inner envelope} for an absolutely
convex neighborhood of zero $\varDelta$} in $\mathcal O^\star(M)$,
$1\in{^\circ\varDelta}$, if
$$
\psi=\varDelta^{\text{\LZ}}
$$

\item[---] an {\it inner envelope for a locally separated from zero function}
$\ph:M\to(0;1]$, if
$$
\psi=\ph^{{\text{\BLZ}}\kern-0.5pt{\text{\LZ}}}
$$

\item[---] an {\it inner envelope on $M$}, if it satisfies the following
equivalent conditions:
 \bit
\item[(i)] $\psi$ is an inner envelope for some absolutely convex neighborhood
of zero $\varDelta$ in $\mathcal O^\star(M)$, $1\in{^\circ\varDelta}$,
$$
\psi=\varDelta^{\text{\LZ}}
$$

\item[(ii)] $\psi$ is an inner envelope for some locally separated from zero
function $\ph:M\to(0;1]$.
$$
\psi=\ph^{{\text{\BLZ}}\kern-0.5pt{\text{\LZ}}}
$$

\item[(iii)] $\psi$ is an inner envelope for itself:
$$
\psi^{{\text{\BLZ}}\kern-0.5pt{\text{\LZ}}} = \psi
$$
 \eit
 \eit
The equivalence of the conditions (i), (ii), (iii) is proved in the same way as
for outer envelopes.

\bigskip

\centerline{\bf Properties of the inner envelopes:}
 {\it
 \bit
\item[(i)] Every inner envelope $\psi$ on $M$ is a continuous function on $M$.

\item[(ii)] For any locally separated from zero function $\ph:M\to(0;1]$ its
inner envelope $\ph^{{\text{\BLZ}}\kern-0.5pt{\text{\LZ}}}$ is the least inner
envelope on $M$, majorizing $\ph$:
 \bit
\item[(a)] $\ph^{{\text{\BLZ}}\kern-0.5pt{\text{\LZ}}}$ is an inner envelope on
$M$, majorizing $\ph$:
$$
\ph \le \ph^{{\text{\BLZ}}\kern-0.5pt{\text{\LZ}}}
$$
\item[(b)] if $\psi$ is another inner envelope on $M$, majorizing $\ph$,
$$
\ph\le \psi
$$
then $\psi$ majorizes $\ph^{{\text{\BLZ}}\kern-0.5pt{\text{\LZ}}}$:
$$
\ph^{{\text{\BLZ}}\kern-0.5pt{\text{\LZ}}}\le \psi
$$
 \eit
\eit
 }

\paragraph{Rhombuses in $\mathcal O^\star(M)$}

We call a set $\varGamma\subseteq \mathcal O^\star(M)$,
 \bit
\item[---] a {\it rhombus, generated by a locally separated from zero
function}\index{rhombus} $\ph:M\to(0;1]$, if
$$
\varGamma=\ph^{\text{\BLZ}}
$$

\item[---] a {\it rhombus, generated by an absolutely convex neighborhood of
zero} $\varDelta\subseteq \mathcal O^\star(M)$, $1\in{^\circ\varDelta}$, if
$$
\varGamma=\varDelta^{{\text{\LZ}}\kern-0.5pt{\text{\BLZ}}}
$$

\item[---] a {\it rhombus} in $\mathcal O^\star(M)$, if it satisfies the
following equivalent conditions:
 \bit
\item[(i)] $\varGamma$ is a rhombus, generated by some locally separated from
zero function $\ph:M\to(0;1]$,
$$
\varGamma=\ph^{\text{\BLZ}}
$$

\item[(ii)] $\varGamma$  is a rhombus, generated by some absolutely convex
neighborhood of zero $\varDelta\subseteq \mathcal O^\star(M)$,
$1\in{^\circ\varDelta}$,
$$
\varGamma=\varDelta^{{\text{\LZ}}\kern-0.5pt{\text{\BLZ}}}
$$

\item[(iii)] the rhombus generated by $\varGamma$, coincides with $\varGamma$:
$$
\varGamma=\varGamma^{{\text{\LZ}}\kern-0.5pt{\text{\BLZ}}}
$$
 \eit
 \eit
\bpr The equivalence of conditions (i), (ii), (iii) is proved in the same way
as in the case of rectangles. \epr

\bigskip

\centerline{\bf Properties of rhombuses:}
 {\it
 \bit
\item[(i)] Every rhombus $\varDelta$ in $\mathcal O^\star(M)$ is a closed
absolutely convex neighborhood of zero in $\mathcal O^\star(M)$.

\item[(ii)] For any neighborhood of zero $\varDelta\subseteq \mathcal
O^\star(M)$ the rhombus $\varDelta^{{\text{\LZ}}\kern-0.5pt{\text{\BLZ}}}$ is the
greatest rhombus in $\mathcal O^\star(M)$, contained in $\varDelta$:
 \bit
\item[(a)] $\varDelta^{{\text{\LZ}}\kern-0.5pt{\text{\BLZ}}}$ is a rhombus in
$\mathcal O^\star(M)$, contained in $\varDelta$:
$$
\varDelta^{{\text{\LZ}}\kern-0.5pt{\text{\BLZ}}}\subseteq \varDelta
$$
\item[(b)] if $\varGamma$ is another rhombus in $\mathcal O^\star(M)$,
contained in $\varDelta$,
$$
\varGamma\subseteq \varDelta,
$$
then $\varGamma$ is contained in
$\varDelta^{{\text{\LZ}}\kern-0.5pt{\text{\BLZ}}}$:
$$
\varGamma\subseteq \varDelta^{{\text{\LZ}}\kern-0.5pt{\text{\BLZ}}}
$$
 \eit

\item[(iii)] Rhombuses in $\mathcal O^\star(M)$ form a fundamental system of
neighborhoods of zero in $\mathcal O^\star(M)$: every neighborhood of zero
$\varGamma$ in $\mathcal O^\star(M)$ contains some rhombus.

 \eit
 }

By analogy with Theorem \ref{TH:f<->D} one can prove

\btm The formulas
 \begin{align}
\varDelta &=\ph^{\text{\BLZ}}, & \ph &=\varDelta^{\text{\LZ}}
 \end{align}
establish a bijection between inner envelopes $\ph$ on $M$ and rhombuses
$\varDelta$ in $\mathcal O^\star(M)$. \etm

\subsection{Duality between rectangles and rhombuses}

Lemma \ref{LM-o-polyarah} implies

\btm\label{TH:formuly-dlya-pryamoug-i-rombov} The following equalities hold
\begin{align}
 \label{f^blacksquare^circ=frac-1-f^blacklozenge}
(f^\text{\BSQ})^\circ &=\left(\frac{1}{f}\right)^{\text{\BLZ}}, &
{^\circ\kern-2pt(\ph^{\text{\BLZ}})} &=\left(\frac{1}{\ph}\right)^\text{\BSQ},
\\
 \label{D^circ^lozenge=frac-1-D^square}
(D^\circ)^{\text{\LZ}} &=\frac{1}{D^\text{\SQ}}\;\;, &
\left({^\circ\kern-2pt\varDelta}\right)^\text{\SQ}
&=\frac{1}{\varDelta^{\text{\LZ}}}\;\;,
\\
\label{1/f^blacksquare^square=(1/f)^blacklozenge^lozenge}
 \frac{1}{f^{\text{\BSQ}\text{\SQ}}}
&=\left(\frac{1}{f}\right)^{{\text{\BLZ}}\kern-0.5pt{\text{\LZ}}}, &
\frac{1}{\ph^{{\text{\BLZ}}\kern-0.5pt{\text{\LZ}}}}
&=\left(\frac{1}{\ph}\right)^{\text{\BSQ}\text{\SQ}}, \\
 \label{D^SQ^BSQ^circ=D^circ^lozenge^blacklozenge}
\Big(D^{\text{\SQ}\text{\BSQ}}\Big)^\circ
&=\Big(D^\circ\Big)^{{\text{\LZ}}\kern-0.5pt{\text{\BLZ}}}, &
{^\circ\kern-0.5pt\Big(\varDelta^{{\text{\LZ}}\kern-0.5pt{\text{\BLZ}}}\Big)}
&=\Big({^\circ\kern-2pt\varDelta}\Big)^{\text{\SQ}\text{\BSQ}},
\end{align}
where $f:M\to[1;+\infty)$ is an arbitrary locally bounded function,
$\ph:M\to(0;1]$ an arbitrary locally separated from zero function, $D$ an
arbitrary absolutely convex compact set in $\mathcal O(M)$, $\varDelta$ an
arbitrary closed absolutely convex neighborhood of zero in $\mathcal
O^\star(M)$. \etm

\bpr 1. The first formula in \eqref{f^blacksquare^circ=frac-1-f^blacklozenge}
follows from \eqref{f^blacksquare^circ}:
$$
(f^\text{\BSQ})^\circ =\eqref{f^blacksquare^circ}= \cabsconv
\left\{\frac{1}{f(x)}\cdot\delta^x;\quad x\in M\right\}=\eqref{f^blacklozenge}=
\left(\frac{1}{f}\right)^{\text{\BLZ}}
$$
After that the substitution $f=\frac{1}{\ph}$ gives the second formula:
$$
\left(\frac{1}{\ph}\right)^{\text{\BSQ}\ \circ}
=\ph^{\text{\BLZ}}\qquad\Longrightarrow\qquad
\left(\frac{1}{\ph}\right)^{\text{\BSQ}}={^\circ\left(\left(\frac{1}{\ph}\right)^{\text{\BSQ}\
\circ}\right)} ={^\circ\left(\ph^{\text{\BLZ}}\right)}
$$

2. The first formula in \eqref{D^circ^lozenge=frac-1-D^square} follows from
\eqref{frac_1_varPhi^square_x}:
$$
\frac{1}{D^\text{\SQ}(x)}=\eqref{frac_1_varPhi^square_x} =\max \{\lambda>0:\
\lambda\cdot\delta^x\in D^\circ\}=\eqref{fi^lozenge}=(D^\circ)^{\text{\LZ}}(x)
 \qquad\Longrightarrow\qquad
(D^\circ)^{\text{\LZ}} =\frac{1}{D^\text{\SQ}}
$$
Then the substitution $D={^\circ\varDelta}$ gives the second formula:
$$
D^\text{\SQ}=\frac{1}{(D^\circ)^{\text{\LZ}}}
 \qquad\Longrightarrow\qquad
({^\circ\varDelta})^\text{\SQ}=\frac{1}{(({^\circ\varDelta})^\circ)^{\text{\LZ}}}=\frac{1}{\varDelta^{\text{\LZ}}}
$$

3. Now the first formula in
\eqref{1/f^blacksquare^square=(1/f)^blacklozenge^lozenge} follows from the
first formula in \eqref{D^circ^lozenge=frac-1-D^square} and from the first
formula in \eqref{f^blacksquare^circ=frac-1-f^blacklozenge} by substitution
$D=f^{\text{\BSQ}}$:
$$
\frac{1}{D^\text{\SQ}}=(D^\circ)^{\text{\LZ}}
 \qquad\Longrightarrow\qquad
\frac{1}{f^\text{\BSQ\SQ}}=(f^{\text{\BSQ}\
\circ})^{\text{\LZ}}=\eqref{f^blacksquare^circ=frac-1-f^blacklozenge}=\left(\frac{1}{f}\right)^{{\text{\BLZ}}{\text{\LZ}}}
$$
The second formula in \eqref{1/f^blacksquare^square=(1/f)^blacklozenge^lozenge}
follows from the second formula in \eqref{D^circ^lozenge=frac-1-D^square} and
the second formula in \eqref{f^blacksquare^circ=frac-1-f^blacklozenge} after the
substitution $\varDelta=\ph^{\text{\BLZ}}$:
$$
\frac{1}{\varDelta^{\text{\LZ}}}=\left({^\circ\kern-2pt\varDelta}\right)^\text{\SQ}
 \qquad\Longrightarrow\qquad
\frac{1}{\ph^{{\text{\BLZ}}{\text{\LZ}}}}=\left({^\circ\kern-2pt(\ph^{\text{\BLZ}})}\right)^\text{\SQ}=
\eqref{f^blacksquare^circ=frac-1-f^blacklozenge}=\left(\frac{1}{\ph}\right)^{\text{\BSQ}\text{\SQ}}
$$

4. The first formula in \eqref{D^SQ^BSQ^circ=D^circ^lozenge^blacklozenge}
follows from the first formula in \eqref{D^circ^lozenge=frac-1-D^square} and
from the first formula in \eqref{f^blacksquare^circ=frac-1-f^blacklozenge}:
$$
(D^\circ)^{\text{\LZ}} =\frac{1}{D^\text{\SQ}}
 \qquad\Longrightarrow\qquad
(D^\circ)^{{\text{\LZ}}{\text{\BLZ}}}=\left(\frac{1}{D^\text{\SQ}}\right)^{\text{\BLZ}}
=\eqref{f^blacksquare^circ=frac-1-f^blacklozenge}=(D^\text{\SQ\BSQ})^\circ
$$
Finally the second formula in \eqref{D^SQ^BSQ^circ=D^circ^lozenge^blacklozenge}
follows from the second formula in \eqref{D^circ^lozenge=frac-1-D^square} and
from the second formula in \eqref{f^blacksquare^circ=frac-1-f^blacklozenge}:
$$
\left({^\circ\kern-2pt\varDelta}\right)^\text{\SQ}=\frac{1}{\varDelta^{\text{\LZ}}}
 \qquad\Longrightarrow\qquad
\left({^\circ\kern-2pt\varDelta}\right)^\text{\SQ\BSQ}=\left(\frac{1}{\varDelta^{\text{\LZ}}}\right)^\text{\BSQ}
=\eqref{f^blacksquare^circ=frac-1-f^blacklozenge}={^\circ\kern-0.5pt\Big(\varDelta^{{\text{\LZ}}\kern-0.5pt{\text{\BLZ}}}\Big)}
$$
 \epr

Theorem \ref{TH:formuly-dlya-pryamoug-i-rombov} implies two important
propositions.

\btm\label{TH-biject-vnesh-i-vnutr-ogib} The passage to the inverse function
\begin{align}
f&=\frac{1}{\ph}, & \ph &=\frac{1}{f}
\end{align}
establish a bijection between the outer envelopes $f$ and the inner envelopes
$\ph$ on $M$.
 \etm
\bpr If $\ph$ is an inner envelope, then $\ph^{{\text{\BLZ}}{\text{\LZ}}}=\ph$,
hence
$\left(\frac{1}{\ph}\right)^\text{\BSQ\SQ}=\eqref{1/f^blacksquare^square=(1/f)^blacklozenge^lozenge}=
\frac{1}{\ph^{{\text{\BLZ}}\kern-0.5pt{\text{\LZ}}}}=\frac{1}{\ph}$, i.e.
$\frac{1}{\ph}$ is an outer envelope. On the contrary, if $f$ is an outer
envelope, then $f^\text{\BSQ\SQ}=f$, therefore
$\left(\frac{1}{f}\right)^{{\text{\BLZ}}\kern-0.5pt{\text{\LZ}}}=\eqref{1/f^blacksquare^square=(1/f)^blacklozenge^lozenge}=
\frac{1}{f^\text{\BSQ\kern-0.5pt\SQ}}=\frac{1}{f}$, i.e. $\frac{1}{f}$ is an inner
envelope.
 \epr

\btm\label{TH-biject-pryam-i-romb} The passage to polar
\begin{align}
D&={^\circ\kern-2pt\varDelta}, & \varDelta &=D^\circ
\end{align}
establishes a bijection between rectangles $D$ in $\mathcal O(M)$ and rhombuses
$\varDelta$ in $\mathcal O^\star(M)$.
 \etm
 \bpr
If $\varDelta$ is a rhombus, then
$\varDelta^{{\text{\LZ}}\kern-0.5pt{\text{\BLZ}}}=\varDelta$, hence
$({^\circ\kern-2pt\varDelta})^\text{\SQ\BSQ}=\eqref{D^SQ^BSQ^circ=D^circ^lozenge^blacklozenge}=
{^\circ\kern-2pt(\varDelta^{{\text{\LZ}}\kern-0.5pt{\text{\BLZ}}})}={^\circ\kern-2pt\varDelta}$,
i.e. ${^\circ\kern-2pt\varDelta}$ is a rectangle. On the contrary, if $D$ is a
rectangle, then $D^\text{\BSQ\SQ}=D$, hence
$(D^\circ)^{{\text{\LZ}}\kern-0.5pt{\text{\BLZ}}}=\eqref{D^SQ^BSQ^circ=D^circ^lozenge^blacklozenge}=
(D^\text{\BSQ\SQ})^\circ=D^\circ$, i.e. $D^\circ$ is a rhombus. Since the
passage to polar is a bijection between closed absolutely convex sets, it is a
bijection between rectangles and rhombuses.
 \epr

Theorems \ref{TH-biject-vnesh-i-vnutr-ogib}, \ref{TH-biject-pryam-i-romb} and
\ref{TH:formuly-dlya-pryamoug-i-rombov} imply

 \btm
The following equalities hold:
\begin{align}
\left( (f^\text{\BSQ})^\circ\right)^{\text{\LZ}} &=\frac{1}{f}, &
\left({^\circ\kern-2pt(\ph^{\text{\BLZ}})}\right)^\text{\SQ} &=\frac{1}{\ph} \\
\left(\frac{1}{(D^\circ)^{\text{\LZ}}}\right)^\text{\BSQ} &=D, &
\left(\frac{1}{\left({^\circ\kern-2pt\varDelta}\right)^\text{\SQ}}\right)^{\text{\BLZ}}
&=\varDelta.
\end{align}
where $f:M\to[1;+\infty)$ is an arbitrary outer envelope, $\ph:M\to(0;1]$ an
arbitrary inner envelope, $D$ an arbitrary rectangle in $\mathcal O(M)$,
$\varDelta$ an arbitrary rhombus in $\mathcal O^\star(M)$. \etm
 \bpr
If $f$ is an outer envelope, then by Theorem
\ref{TH-biject-vnesh-i-vnutr-ogib}, $\frac{1}{f}$ is an inner envelope, so
$$
\left(
(f^\text{\BSQ})^\circ\right)^{\text{\LZ}}=\eqref{f^blacksquare^circ=frac-1-f^blacklozenge}=
\left(\frac{1}{f}\right)^{{\text{\BLZ}}\kern-0.5pt{\text{\LZ}}}=\frac{1}{f}
$$
The other formulas are proved by analogy.
 \epr

\section{Stein groups and Hopf algebras connected to them}
\label{SEC:stein-groups}

\subsection{Stein groups, linear groups and algebraic groups}\label{SUBSEC-lin-groups}

A complex Lie group $G$ is called a {\it Stein group}\index{Stein!group}, if
$G$ is a Stein manifold \cite{Grauert-Remmert}. By the Matsushima-Morimoto
Theorem \cite[XIII.5.9]{Neeb}, for complex groups this is equivalent to the
condition of holomorphic separability we mentioned at the page
\pageref{holom-separability}:
$$
\forall x\ne y\in G\quad \exists u\in {\mathcal O}(G)\quad u(x)\ne u(y)
$$
{\it Dimension} of a Stein group is its dimension as a complex manifold.

Special cases of Stein groups are {\it linear complex groups}\index{linear
group}. They are defined as complex Lie groups which can be represented as
closed complex Lie subgroups in the general linear group $\GL_n(\C)$. In other
words, a complex group $G$ is linear, if it is isomorphic to some closed
complex subgroup $H$ in $\GL_n(\C)$ (i.e. there is an isomorphism of groups
$\ph:G\to H$ which is at the same time a biholomorphic mapping).

Even more narrow class are {\it complex affine algebraic
groups}\index{algebraic group}. These are subgroups $H$ in $\GL_n(\C)$, which
are at the same time algebraic submanifolds. This means that the group $H$ must
be a common set of zeroes for some finite set of polynomials $u_1,...,u_k$ on
$\GL_n(\C)$ (by polynomial here we can understand a polynomial of matrix
elements):
$$
H=\{x\in \GL_n(\C):\ u_1(x)=...=u_k(x)=0\}
$$
If a complex group $G$ is isomorphic to some algebraic group $H$ (i.e. there
exists an isomorphism of groups $\ph:G\to H$ which is at the same time a
biholomorphic mapping), then $G$ is also considered an algebraic group, since
the algebraic operations on $G$ are regular mappings with respect to the
structure of algebraic manifold inherited from $H$.

A Stein group $G$ is called {\it compactly generated}\index{compactly generated
group}, if it has a generating compact set, i.e. a compact set $K\subseteq G$
such that
$$
G=\bigcup_{n\in\N}K^n,\qquad K^n=\underbrace{K\cdot...\cdot K}_{\text{$n$
factors}}
$$

Let us note some examples.

\bex {\it Complex torus} we have mentioned in \ref{rectangles-buses}, i.e. the
quotient group $\C/(\Z+i\Z)$, is an example of a complex group, which is not a
Stein group. \eex

\bex Every {\it discrete group} $G$ is a Stein group (of zero dimension).
Compact sets in $G$ are nothing more than finite sets, hence $G$ will be
compactly generated if and only if it is finitely generated. Thus, say, a free
group with infinite set of generators can be considered as an example a Stein
group which is not compactly generated.
 \eex

\bex A discrete group $G$ is algebraic if and only if it is {\it finite}. In
this case it can be represented as a group of transformations of the space
$\C^n$, where $n=\card G$ is the number of elements of $G$. For this $\C^n$
must be represented as the space of functions from $G$ into $\C$:
$$
x\in\C^G\quad\Longleftrightarrow\quad x:G\to\C,
$$
Then the imbedding of $G$ into $\GL(\C^G)$ (the group of a nondegenerate linear
transformations of the space $\C^G$) is defined by the formula
$$
\ph:G\to \GL(\C^G):\qquad \ph(g)(x)(h)=x(h\cdot g),\qquad g,h\in G,\quad x\in
\C^G
$$
\eex

\bex The general linear group $\GL_n(\C)$, i.e. the group of a nondegenerate
linear transformations of the space $\C^n$, is an algebraic group (of dimension
$n^2$). Certainly, $\GL_n(\C)$ is compactly generated. As a corollary, {\it
every linear group is compactly generated}.

\eex

\bex The additive group of complex numbers $\C$ is a complex algebraic group
(of dimension 1), since it can be embedded into $\GL(2,\C)$ by formula
$$
\ph(\lambda)=\begin{pmatrix}1 & \lambda \\ 0 & 1 \end{pmatrix},\qquad
\lambda\in\C
$$
\eex

\bex\label{EX:Z} The additive group $\Z$ of integers is a complex linear group
(of dimension 0), since it can be embedded into $\GL(2,\C)$ by the same formula
$$
\ph(n)=\begin{pmatrix}1 & n \\ 0 & 1 \end{pmatrix},\qquad n\in\Z
$$
But the difference with $\C$ is that $\Z$ is not algebraic: neither in this
embedding into $\GL_n(\C)$, nor in any other one, $\Z$ is a common set for a
family of polynomials (this is a result of the fact that $\Z$ is discrete and
infinite). \eex

\bex\label{EX:C^x} A multiplicative group $\C^\times$ of nonzero complex
numbers
$$
\C^\times=\C\setminus\{0\}
$$
can be represented as general linear group (of nondegenerate transformations of
the space $\C$),
$$
\C^\times\cong\GL_1(\C),
$$
Thus $\C^\times$ is a complex algebraic group (of dimension 1). We call this
group {\it complex circle}\index{complex circle}. \eex

\bex Consider the action of group $\C$ on itself by exponents:
$$
\ph:\C\to\Aut(\C),\qquad \ph(a)(x)=x\cdot e^a
$$
The semidirect product of $\C$ and $\C$ with respect to this action, i.e. a
group $\C\ltimes\C$, coinciding with the Cartesian product $\C\times\C$, but
endowed with a more complicated multiplication
$$
(a,x)\cdot(b,y):=(a+b,x\cdot e^b+y)
$$
is a (connected) linear complex group, since it can be embedded into
$\GL_3(\C)$ by the homomorphism
$$
(x,a)\mapsto \begin{pmatrix} e^a & 0 & 0 \\ x & 1 & 0 \\ 0 & 0 &
e^{ia}\end{pmatrix}
$$
But $\C\ltimes\C$ is not algebraic, since its center
$$
Z(\C\ltimes\C)=\{(2\pi i n,0);\; n\in\Z\}
$$
is an infinite discrete subgroup (this does not happen with algebraic groups).
 \eex

\subsection{Hopf algebras
$\mathcal O(G)$, $\mathcal O^\star(G)$, $\mathcal R(G)$, $\mathcal
R^\star(G)$}\label{R-O}

The Hopf algebras $\C^G$ and $\C_G$ we were talking about in
\ref{SEC:ster-algebry-Hopfa}\ref{C^G-i-C_G}, are interesting not so in
themselves but more as guiding examples for various similar constructions
arising in situations when the initial group $G$ is endowed with some
supplementary structure and functions on $G$ preserve this structure. In
particular, in \cite[examples 10.24-10.27]{Akbarov} the author noted standard
examples of algebras of functions and functionals, which are stereotype Hopf
algebras. If we add Hopf algebras $\C^G$ and $\C_G$ to those examples, we
obtain the following list:
 \begin{center}
 \begin{tabular}{|l|l|l|}
 \hline
  && \\
  class of groups & algebra of functions & algebra of functionals \\
  && \\
 \hline
  "pure" groups &
  \begin{tabular}{l}
  algebra $\C^G$ \phantom{$1^{1^{1^1}}$}\\
  of all functions on $G$
  \\
 \end{tabular}
   &
  \begin{tabular}{l}
  algebra $\C_G$ \phantom{$1^{1^{1^1}}$}\\
  of point charges on $G$
  \\
 \end{tabular}
   \\
 \hline
 algebraic groups &
  \begin{tabular}{l}
  algebra ${\mathcal R}(G)$ \phantom{$1^{1^{1^1}}$}\\
  of polynomials on $G$
  \\
 \end{tabular}
   &
  \begin{tabular}{l}
  algebra ${\mathcal R}^\star(G)$ \phantom{$1^{1^{1^1}}$}\\
  of currents of degree 0 on $G$
  \\
 \end{tabular}
 \\
 \hline
 Stein groups &
 \begin{tabular}{l}
  algebra ${\mathcal O}(G)$ \phantom{$1^{1^{1^1}}$}\\
  of holomorphic functions on $G$
  \\
 \end{tabular}
   &
  \begin{tabular}{l}
  algebra ${\mathcal O}^\star(G)$ \phantom{$1^{1^{1^1}}$}\\
  of analytical \\ functionals on $G$
  \\
 \end{tabular}
 \\
 \hline
 Lie groups &
  \begin{tabular}{l}
  algebra ${\mathcal E}(G)$ \phantom{$1^{1^{1^1}}$}\\
  of smooth functions on $G$
  \\
 \end{tabular}
   &
  \begin{tabular}{l}
  algebra ${\mathcal E}^\star(G)$ \phantom{$1^{1^{1^1}}$}\\
  of distributions on $G$
  \\
 \end{tabular}
 \\
 \hline
 locally compact groups &
   \begin{tabular}{l}
 algebra ${\mathcal C}(G)$ \phantom{$1^{1^{1^1}}$}\\
  of continuous functions on $G$
  \\
 \end{tabular}
   &
  \begin{tabular}{l}
  algebra ${\mathcal C}^\star(G)$ \phantom{$1^{1^{1^1}}$}\\
  of Radon measures on $G$
  \\
 \end{tabular}
 \\
 \hline
 \end{tabular}
 \end{center}
In \cite{Akbarov} the last four examples were mentioned without proof, so we
see fit to explain here why those algebras are indeed stereotype Hopf algebras.
Stein groups and algebraic groups will be case studies for us.

\paragraph{Hopf algebras $\mathcal O(G)$ and $\mathcal O^\star(G)$ on a Stein group $G$.}

If $G$ is a Stein group, then the proposition that algebra ${\mathcal O}(G)$ of
holomorphic functions on $G$ (with the usual topology of uniform convergence on
compact sets in $G$) is an injective stereotype Hopf algebra, is proved exactly
like this was done for $\C^G$. A significant aspect in those reasonings is the
isomorphism of functors connecting Cartesian product of groups $\times$ with
the corresponding tensor products of functional spaces $\odot$ and
$\circledast$ \cite[Theorem 8.13]{Akbarov},
 \beq\label{O(GxH)=O(G)-o-O(H)}
{\mathcal O}(G\times H)\stackrel{\rho_{G,H}}{\cong} {\mathcal O}(G)\odot
{\mathcal O}(H)\stackrel{@^{-1}}{\cong} {\mathcal O}(G)\circledast {\mathcal
O}(H)
 \eeq
-- here $\rho_{G,H}$ is defined by identity analogous to
\eqref{u-boxdot-v->u-odot-v}:
 \beq\label{u-boxdot-v->u-odot-v-R-O}
\rho_{G,H}(u\boxdot v)=u\odot v,\qquad u\in{\mathcal O}(G),\quad v\in {\mathcal
O}(H)
 \eeq
and the function $u\boxdot v$ is again defined by formula \eqref{g-h-na-S-T}.

After defining $\rho_{G,H}$, the multiplication and comultiplication in
${\mathcal O}(G)$ are defined initially on (or with values in) the space
${\mathcal O}(G\times G)$ of functions on the Cartesian square $G\times G$, and
then the passage to tensor square is carried out with the help of the
isomorphism  $\rho_{G,G}$:
 \begin{align}
\label{umn-v-O(G)} \text{multiplication:} & & &
\mu=\widetilde{\mu}\circ\rho_{G,G}:{\mathcal O}(G)\odot {\mathcal O}(G)\to
{\mathcal O}(G),
 & &
\widetilde{\mu}(v)(t)=v(t,t) \\
\label{edin-v-O(G)} \text{unit:} & & & \iota:\C\to {\mathcal O}(G),& &
\iota(\lambda)(t)=\lambda \\
 \text{comultiplication:} & & &
\label{koumn-v-O(G)} \varkappa=\rho_{G,G}\circ\widetilde{\varkappa}: {\mathcal
O}(G)\to {\mathcal O}(G)\odot {\mathcal O}(G),
& & \widetilde{\varkappa}(u)(s,t)=u(s\cdot t) \\
\label{koed-v-O(G)}\text{counit:} & & &  \e:{\mathcal O}(G)\to\C,& &
\e(u)=u(1_G) \\
\label{antipod-v-O(G)} \text{antipode:} & & & \sigma:{\mathcal
O}(G)\to{\mathcal O}(G),& & \sigma(u)(t)=u(t^{-1})
 \end{align}
This can be illustrated by the following picture:
 \beq\label{str-morph-v-O(G)}
 \begin{diagram}
 \node[3]{{\mathcal O}(G\times G)}
 \arrow{se,l}{\widetilde{\mu}}
 \arrow[2]{s,r}{\rho_{G,G}}
 \\
 \node{\C}\arrow{e,t}{\iota} \node{{\mathcal O}(G)}
 \arrow{ne,l}{\widetilde{\varkappa}}
 \arrow{se,r}{\varkappa}
 \node[2]{{\mathcal O}(G)}\arrow{e,t}{\e}
 \node{\C}
 \\
 \node[3]{{\mathcal O}(G)\odot {\mathcal O}(G)}
 \arrow{ne,r}{\mu}
 \end{diagram}
 \eeq
The fact that this defines an injective stereotype Hopf algebra is proved
literally like Theorem \ref{TH-C^G-algebra-Hopfa}. And again, like in the case
of $\C^G$, the Hopf algebra ${\mathcal O}(G)$ becomes a projective Hopf algebra
(hence, a rigid Hopf algebra in the sense of definition
\ref{SEC:ster-algebry-Hopfa}\ref{ster-alg-Hopf}), because of the second
equality in \eqref{O(GxH)=O(G)-o-O(H)}.

\btm\label{TH-O} For any Stein group $G$
 \bit
\item[---] the algebra $\mathcal O(G)$ of holomorphic functions on $G$ is a
rigid stereotype Hopf algebra with respect to algebraic operations, defined by
formulas \eqref{umn-v-O(G)}-\eqref{antipod-v-O(G)};

\item[---] its dual algebra $\mathcal O^\star(G)$ of analytical functionals on
$G$ is a rigid stereotype Hopf algebra with respect to the dual algebraic
operations.
 \eit
If in addition the group $G$ is compactly generated, then $\mathcal O(G)$ is a
nuclear Hopf-Fr\'echet algebra, and $\mathcal O^\star(G)$ a nuclear
Hopf-Brauner algebra.
 \etm

\paragraph{Hopf algebras $\mathcal R(G)$ and $\mathcal R^\star(G)$ on an affine algebraic group $G$.}

Let $G$ be an affine algebraic group, ${\mathcal R}(G)$ the algebra of
polynomials on $G$ (with the strongest locally convex topology). Like in the
previous cases, the identity
 $$
\rho_{G,H}(u\boxdot v)=u\odot v,\qquad u\in{\mathcal R}(G),\quad v\in {\mathcal
R}(H)
 $$
(again, $u\boxdot v$ is defined by formula \eqref{g-h-na-S-T}) defines an
isomorphism of functors, connecting the Cartesian product of groups $\times$
and tensor products of functional spaces $\odot$ and $\circledast$
\cite[Theorem 8.16]{Akbarov}:
 \beq\label{R(G-x-H)=R(G)-o-R(H)}
{\mathcal R}(G\times H)\stackrel{\rho_{G,H}}{\cong} {\mathcal R}(G)\odot
{\mathcal R}(H)\stackrel{@^{-1}}{\cong} {\mathcal R}(G)\circledast {\mathcal
R}(H)
 \eeq
Those isomorphisms then define algebraic operations on ${\mathcal R}(G)$ by
formulas, analogous to \eqref{umn-v-O(G)}-\eqref{antipod-v-O(G)}. As a result
we come to the following

\btm\label{TH-R} For any affine algebraic group $G$
 \bit
\item[---] the algebra $\mathcal R(G)$ of polynoimals on $G$ is a rigid
stereotype Hopf-Brauner algebra;

\item[---] the dual algebra $\mathcal R^\star(G)$ of currents of degree 0 on
$G$ is an injective stereotype Hopf-Fr\'echet algebra.
 \eit
 \etm

\paragraph{Convolutions in $\mathcal R^\star(G)$ and $\mathcal O^\star(G)$.}

For some further calculations it is useful to record the formulas defining
convolution in the algebras of functionals $\mathcal R^\star(G)$ and $\mathcal
O^\star(G)$. This definition is anticipated by formulas for shift, antipode and
convolution between a functional and a function:
 \begin{align}\label{DEF:u-cdot-a}
(u\cdot a) (x)&:=u (a\cdot x), & (a\cdot u) (x)&:=u (x\cdot a), & & u\in
\mathcal O(G), a,x\in G
 \\ \label{DEF:alpha-cdot-a}
(\alpha \cdot a)(u)&:=\alpha (a\cdot u), & (a\cdot \alpha) (u)&:=\alpha (u\cdot
a), & &\alpha \in \mathcal O^\star(G), u\in \mathcal O(G), a\in G
 \\ \label{DEF:tilde-u}
 \tilde {u}(x)&:=u(x^{-1}), & \tilde {\alpha}(u)&:=\alpha
 (\tilde{u}),& &\alpha \in \mathcal O^\star(G), u\in \mathcal O(G), x\in G
 \\ \label{DEF:alpha-*-u}
(\alpha * u)(x)&:=\alpha (\widetilde{x\cdot u}), &
 (u *\alpha)(x)&:=\alpha (\widetilde{u\cdot x}) & &\alpha \in
\mathcal O^\star(G), u\in \mathcal O(G), x\in G
 \end{align}
Then the convolution of functionals is defined by formula
 \begin{align}\label{DEF:svertka}
 \alpha * \beta (u)&:=\alpha(\widetilde{\beta * \tilde{u}})=
\tilde{\alpha} (\beta * \tilde{u})=\beta (\tilde{\alpha}*u) & &\alpha, \beta
\in \mathcal O^\star(G), u\in \mathcal O(G)
 \end{align}
In particular, the convolution with the delta-functional is a shift:
 \begin{align}\label{DEF:svertka-s-delta^a}
& \delta^a * \beta=a\cdot\beta && \beta*\delta^a=\beta\cdot a && \beta \in
\mathcal O^\star(G), a\in G
 \end{align}

\subsection{Examples}\label{Examples-of-Stein-groups}

It is useful to illustrate the latter theorems by several examples. For this we
choose the main examples of Abelian Stein groups -- $\Z$, $\C^\times$, $\C$ --
apart from everything else, these examples will be useful below in
\ref{SEC-hol-duality=gen-pontryagin} and \ref{SEC:R_q(C^X-x-C)}.

\paragraph{Algebras ${\mathcal O}(\Z)$ and ${\mathcal O}^\star(\Z)$.}

We have already noted in Example \ref{EX:Z} that the group $\Z$ of integers can
be considered as a complex group (of dimension 0). Since $\Z$ is discrete, each
function on $\Z$ is automatically holomorphic, so the algebra ${\mathcal
O}(\Z)$ formally coincides with the algebra $\C^{\Z}$, and the algebra
${\mathcal O}^\star(\Z)$ with the algebra $\C_{\Z}$:
 $$
{\mathcal O}(\Z)=\C^{\Z},\qquad {\mathcal O}^\star(\Z)=\C_{\Z}
 $$
As a corollary the structure of these algebras is described by formulas
\eqref{umnozhenie-v-C^G-infty}-\eqref{umnozhenie-na-bazise-v-C_G-infty}: the
characteristic functions of singletons
 \beq\label{DEF:1_n,n-in-Z}
1_n(m)=\begin{cases}0,& m=n\\ 1,& m=n\end{cases},\qquad m\in\Z,\; n\in\Z
 \eeq
form a basis in the stereotype space ${\mathcal O}(\Z)=\C^{\Z}$, and
delta-functionals
$$
\delta^k(u)=u(k),\qquad u\in{\mathcal O}(\Z)
$$
a dual (algebraic) basis in ${\mathcal O}^\star(\Z)=\C_{\Z}$:
 $$
\langle 1_n,\delta^k\rangle=\begin{cases}0,& n\ne k\\ 1,& n=k\end{cases},
 $$
It is convenient to represent elements of ${\mathcal O}(\Z)$ and ${\mathcal
O}^\star (\Z)$ in the form of series (which converge in these spaces)
 \begin{align}
 \label{razlozhenie-u-v-Z}
u &\in {\mathcal O}(\Z)=\C^\Z & & \Longleftrightarrow & u&=\sum_{n\in\Z}
u(n)\cdot 1_n, & u(n) & =\delta^n(u),
 \\
 \label{razlozhenie-alpha-v-Z}
\alpha &\in {\mathcal O}^\star(\Z)=\C_\Z & & \Longleftrightarrow & \alpha
&=\sum_{n\in\Z} \alpha_n\cdot \delta^n, & \alpha_n & =\alpha(1_n),
 \end{align}
where the action of $\alpha$ on $u$ is described by formula
 $$
\langle u,\alpha\rangle=\sum_{n\in\Z} u(n)\cdot\alpha_n
 $$
The operations of multiplication in ${\mathcal O}(\Z)=\C^{\Z}$ and in
${\mathcal O}^\star (\Z)=\C_{\Z}$ are represented by series:
 \beq\label{umnozhenie-v-O(Z)-i-O^star(Z)}
u\cdot v=\sum_{n\in\Z} u(n)\cdot v(n)\cdot 1_n, \qquad \alpha
* \beta=\sum_{k\in\Z} \left(\sum_{i\in\Z}\alpha_i\cdot \beta_{k-i}\right)\cdot \delta^k,
 \eeq
(in the first case this is the coordinate-wise multiplication, and in the
second case the multiplication of power series).

\bprop The algebra ${\mathcal O}(\Z)=\C^{\Z}$ of functions on $\Z$ is a nuclear
Hopf-Fr\'echet algebra with the topology generated by seminorms
 \beq\label{polunormy-v-H(Z)}
||u||_N=\sum_{|n|\le N} |u(n)|,\qquad N\in\N.
 \eeq
and with algebraic operations defined on basis elements $1_k$ by formulas
 \begin{align}
\label{umn-v-C^Z} & 1_m\cdot 1_n=\begin{cases}1_m,& m=n
\\ 0,& m\ne n
\end{cases} && 1_{{\mathcal R}^\star(\C^\times)}=\sum_{n\in\Z} 1_n \\
\label{koumn-v-C^Z}  &\varkappa(1_n)=\sum_{m\in\Z} 1_m\odot 1_{n-m} & &
\e(1_n)=\begin{cases}1,& n=0\\ 0,& n\ne 0\end{cases} \\
\label{antipod-v-C^Z} & \sigma(1_n)=1_{-n} &&
 \end{align}
 \eprop

\bprop The algebra ${\mathcal O}^\star(\Z)=\C_{\Z}$ of point charges on $\Z$ is
a nuclear Hopf-Brauner algebra with the topology generated by seminorms
 \beq\label{|alpha|_r-v-Z}
|||\alpha|||_r=\sum_{n\in\Z} r_n\cdot |\alpha_n|\qquad (r_n\ge 0)
 \eeq
and algebraic operations defined on basis monomials $\delta^k$ by formulas
 \begin{align}
\label{umn-v-C_Z} & \delta^k*\delta^l=\delta^{k+l} && 1_{{\mathcal R}(\C^\times)}=\delta^0 \\
\label{koumn-v-C_Z}  &\varkappa(\delta^k)= \delta^k\circledast \delta^k & &
\e(z^k)=1 \\
\label{antipod-v-C_C_Z} & \sigma(\delta^k)=\delta^{-k} &&
 \end{align}
 \eprop
\bpr It can be not obvious here that seminorms \eqref{|alpha|_r-v-Z} indeed
define the topology in ${\mathcal O}^\star(\Z)=\C_{\Z}$. The auxiliary
statement used here will be also useful for us below in Lemma
\ref{LM-||alpha||_C-v-Z}, so we formulate it separately: \epr

\blm\label{LM-polunormy-v-H*(Z)} If $p$ is a continuous seminorm on ${\mathcal
O}^\star(\Z)$, and $r_n=p(\delta^n)$, then $p$ is majorized by seminorm
\eqref{|alpha|_r-v-Z}:
 \beq\label{p(alpha)-le-|||alpha|||_r-O*(Z)}
p(\alpha)\le |||\alpha|||_r
 \eeq
 \elm
\bpr
 \beq\label{tsep-p-le-|||.|||_r}
p(\alpha)=p\left(\sum_{n\in\Z} \alpha_n\cdot\delta^n\right)\le \sum_{n\in\Z}
|\alpha_n|\cdot p(\delta^n)=\sum_{n\in\Z} |\alpha_n|\cdot r_n=|||\alpha|||_r
 \eeq
 \epr

\paragraph{Algebras ${\mathcal R}(\C^\times)$, ${\mathcal
R}^\star(\C^\times)$, ${\mathcal O}(\C^\times)$, ${\mathcal
O}^\star(\C^\times)$.}\label{SUBSEC:algebry-na-C^x}

In Example \ref{EX:C^x} we denoted by $\C^\times$ the multiplicative group of
nonzero complex numbers:
$$
\C^\times:=\C\setminus\{0\}.
$$
(the multiplication in $\C^\times$ is the usual multiplication of complex
numbers). We call this group {\it complex circle}.

The algebra ${\mathcal R}(\C^\times)$ of polynomials on $\C^\times$ consists of
Laurent polynomials,  i.e. of function of the form
$$
u=\sum_{n\in\Z} u_n\cdot z^n
$$
where $z^n$ are monomials on $\C^\times$:
 \beq\label{DEF:z^n-in-R(C^x)}
z^n(x):=x^n,\qquad x\in \C^\times,\qquad n\in\Z
 \eeq
and almost all $u_n\in\C$ vanish:
$$
\card\{n\in\Z:\ u_n\ne 0\}<\infty
$$

And the algebra ${\mathcal R}^\star(\C^\times)$ of currents on $\C^\times$
consists of functionals
$$
\alpha=\sum_{k\in\Z} \alpha_k\cdot \zeta^k
$$
where $\zeta_k$ is the functional of finding the $k$-th Laurent coefficient:
 \beq\label{DEF:zeta_n}
\zeta_k(u)=\frac{1}{2\pi i}\int_{|z|=1} \frac{u(z)}{z^{k+1}}\d z=\int_0^1
e^{-2\pi i k t}u(e^{2\pi i t})\d t,\qquad u\in{\mathcal R}(\C^\times)
 \eeq
and $\alpha_k\in\C$ is an arbitrary sequence.

Monomials $\zeta_k$ and $z^n$ act on each other by formula
 \beq\label{zeta^k(z^n)}
\langle z^n,\zeta_k\rangle=\begin{cases}1,& n=k \\ 0,& n\ne k \end{cases},
 \eeq
so the action of a current $\alpha$ on a polynomial $u$ is described by formula
 $$
\langle u,\alpha\rangle=\sum_{n\in\Z} u_n\cdot\alpha_n
 $$
and the operations of multiplication in ${\mathcal R}(\C^\times)$ and in
${\mathcal R}^\star (\C^\times)$ are represented by series as follows:
 \beq\label{umnozhenie-v-R(C-x)-i-R^star(C-x)}
u\cdot v=\sum_{n\in\Z} \left(\sum_{i\in\Z} u_i\cdot v_{n-i}\right)\cdot z^n,
\qquad \alpha * \beta=\sum_{k\in\Z} \alpha_k\cdot \beta_k\cdot \zeta_k,
 \eeq
(in the first case this is the multiplication of series, and in the second case
the coordinate-wise multiplication).

 \bprop The mapping
$$
u\in{\mathcal R}(\C^\times)\mapsto \{u_k;\;k\in\Z\}\in\C_{\Z}
$$
is an isomorphism of nuclear Hopf-Brauner algebras
  \beq\label{R(C^times)=C_Z}
{\mathcal R}(\C^\times)\cong \C_{\Z}
  \eeq
 \eprop

\bprop The algebra ${\mathcal R}(\C^\times)$ of polynomials on the complex
circle $\C^\times$ is a nuclear Hopf-Fr\'echet algebra with the topology
generated by seminorms
 \beq\label{|||u|||_r-v-R(C^x)}
|||u|||_r=\sum_{n\in\Z} r_n\cdot |u_n|\qquad (r_n\ge 0)
 \eeq
and algebraic operations defined on monomials $z^k$ by formulas
 \begin{align}
\label{umn-v-R(C^times)} & z^k\cdot z^l=z^{k+l} && 1_{{\mathcal R}(\C^\times)}=z^0 \\
\label{koumn-v-R(C^times)}  &\varkappa(z^k)= z^k\odot z^k & &
\e(z^k)=1 \\
\label{antipode-in-R(C^times)} & \sigma(z^k)=z^{-k} &&
 \end{align}
 \eprop

 \bprop The mapping
$$
\alpha\in{\mathcal R}^\star(\C^\times)\mapsto \{\alpha_k;\;k\in\Z\}\in\C^{\Z}
$$
is an isomorphism of nuclear Hopf-Fr\'echet algebras
  \beq\label{R^star(C^times)=C^Z}
{\mathcal R}^\star(\C^\times)\cong \C^{\Z}
  \eeq
 \eprop

\bprop The algebra ${\mathcal R}^\star(\C^\times)$ of currents on the complex
circle $\C^\times$ is a nuclear Hopf-Fr\'echet algebra with the topology
generated by seminorms
 \beq\label{polunormy-v-R*(C^x)}
||\alpha||_N=\sum_{|n|\le N} |\alpha_n|\qquad (N\in\N)
 \eeq
and algebraic operations defined on basis elements $\zeta_k$ by formulas
 \begin{align}
\label{umn-v-R^star(C^times)} & \zeta_k*\zeta_l=\begin{cases}\zeta_l,& k=l
\\ 0,& k\ne l
\end{cases} && 1_{{\mathcal R}^\star(\C^\times)}=\sum_{n\in\Z} \zeta_n \\
\label{koumn-v-R^star(C^times)}  &\varkappa(\zeta_k)=\sum_{l\in\Z}
\zeta_l\circledast \zeta_{k-l} & &
\e(\zeta_k)=\begin{cases}1,& k=0\\ 0,& k\ne 0\end{cases} \\
\label{antipode-in-R^star(C^times)} & \sigma(\zeta_k)=\zeta_{-k} &&
 \end{align}
 \eprop

As usual, by symbol ${\mathcal O}(\C^\times)$ we denote the algebra of
holomorphic functions on the complex circle $\C^\times$, and by ${\mathcal
O}^\star(\C^\times)$ its dual algebra of analytic functionals on $\C^\times$.

It is useful to represent the elements of algebras ${\mathcal O}(\C^\times)$
and ${\mathcal O}^\star (\C^\times)$ by series
 \begin{align}
 \label{razlozhenie-u-v-C-x}
 u &\in {\mathcal O}(\C^\times) & &
\Longleftrightarrow & u&=\sum_{n\in\Z} u_n\cdot z^n, & u_n\in\C:\quad & \forall
C>0\quad \sum_{n\in\Z}|u_n|\cdot C^{|n|}<\infty & \Big(u_n & =\zeta_n(u)\Big)
 \\
 \label{razlozhenie-alpha-v-C-x}
\alpha &\in {\mathcal O}^\star(\C^\times) & & \Longleftrightarrow & \alpha
&=\sum_{n\in\Z} \alpha_n\cdot \zeta_n, & \alpha_n\in\C:\quad & \exists C>0\quad
\forall n\in\N\quad |\alpha_n|\le C^{|n|} & \Big(\alpha_n & =\alpha(z^n)\Big)
 \end{align}
Like in the case of ${\mathcal R}(\C^\times)$ and ${\mathcal R}^\star
(\C^\times)$, the action of $\alpha$ on $u$ is described by formula
 $$
\langle u,\alpha\rangle=\sum_{n\in\Z} u_n\cdot\alpha_n,
 $$
and the operations of multiplication in ${\mathcal O}(\C^\times)$ and
${\mathcal O}^\star (\C^\times)$ can be written as the usual multiplication of
series in the first case, and a coordinate-wise multiplication in the second
case:
 \beq\label{umnozhenie-v-O(C-x)-i-O^star(C-x)}
u\cdot v=\sum_{n\in\Z} \left(\sum_{i\in\Z} u_i\cdot v_{n-i}\right)\cdot z^n,
\qquad \alpha * \beta=\sum_{n\in\Z} \alpha_n\cdot \beta_n\cdot \zeta_n,
 \eeq

\bprop The algebra ${\mathcal O}(\C^\times)$ of holomorphic functions on the
complex circle $\C^\times$ is a nuclear Hopf-Fr\'echet algebra with the
topology generated by seminorms
 \beq\label{polunormy-v-H(C-x)}
||u||_C=\sum_{n\in\Z} |u_n|\cdot C^{|n|},\qquad C\ge 1.
 \eeq
and algebraic operations defined on monomials $z^k$ by the same formulas
\eqref{umn-v-R(C^times)}-\eqref{antipode-in-R(C^times)} as in the case of
${\mathcal R}(\C^\times)$:
 \begin{align*}
 & z^k\cdot z^l=z^{k+l} && 1_{{\mathcal O}(\C^\times)}=z^0 \\
  &\varkappa(z^k)= z^k\odot z^k & &
\e(z^k)=1 \\
 & \sigma(z^k)=z^{-k} &&
 \end{align*}
 \eprop
\bpr Every usual seminorm $|u|_K=\max_{x\in K}|u(x)|$, where $K$ is a compact
set in $\C^\times$, is majorized by some seminorm $||u||_C$, namely the one
with $C=\max_{x\in K}\max\{|x|,\frac{1}{|x|}\}$:
$$
|u|_K=\max_{x\in K}|u(x)|=\max_{x\in K}\left|\sum_{n\in\Z} u_n\cdot
x^n\right|\le\max_{x\in K}\sum_{n\in\Z} |u_n|\cdot |x^n|\le\sum_{n\in\Z}
|u_n|\cdot C^{|n|}=||u||_C
$$
On the contrary, for every number $C\ge 1$ we can take a compact set
$K=\{t\in\C:\;\frac{1}{C+1}\le |t|\le C+1\}$, and then from the Cauchy formulas
for Laurent coefficients
$$
|u_n|\le
|u|_K\cdot\min\left\{(C+1)^n;\frac{1}{(C+1)^n}\right\}=|u|_K\cdot(C+1)^{-|n|}=
\frac{|u|_K}{(C+1)^{|n|}}
$$
we have that $||u||_C$ is majorized by $|u|_K$ (with coefficient $(C+1)^2$):
$$
||u||_C=\sum_{n\in\Z} |u_n|\cdot C^{|n|}\le \sum_{n\in\Z}
\frac{|u|_K}{(C+1)^{|n|}}\cdot C^{|n|}\le |u|_K\cdot\sum_{n\in\Z}
\left(\frac{C}{C+1}\right)^{|n|}=(C+1)^2\cdot|u|_K
$$
\epr

\bprop\label{PROP:O*(C^x)-Hopf} The algebra ${\mathcal O}^\star(\C^\times)$ of
analytical functionals on the complex circle $\C^\times$ is a nuclear
Hopf-Brauner algebra with the topology generated by seminorms
 \beq\label{|alpha|_r-v-C-x}
|||\alpha|||_r=\sup_{u\in E_r}|\alpha(u)|= \sum_{n\in\Z} r_n\cdot
|\alpha_n|\qquad \left(r_n\ge 0: \forall C>0\qquad \sum_{n\in\Z} r_n\cdot
C^{|n|}<\infty\right)
 \eeq
and algebraic operations defined on basis elements $\zeta_k$ by the same
formulas \eqref{umn-v-R^star(C^times)}-\eqref{antipode-in-R^star(C^times)} as
in the case of ${\mathcal R}^\star(\C^\times)$:
 \begin{align*}
 & \zeta_k*\zeta_l=\begin{cases}\zeta_l,& k=l
\\ 0,& k\ne l
\end{cases} && 1_{{\mathcal O}^\star(\C^\times)}=\sum_{n\in\Z} \zeta_n \\
  &\varkappa(\zeta_k)=\sum_{l\in\Z} \zeta_l\circledast
\zeta_{k-l} & &
\e(\zeta_k)=\begin{cases}1,& k=0\\ 0,& k\ne 0\end{cases} \\
 & \sigma(\zeta_k)=\zeta_{-k} &&
 \end{align*}
 \eprop
\bpr For every sequence of nonnegative numbers $r_n\ge 0$ satisfying the
condition
 \beq\label{forall-R>0-sum_r_n-R^n<infty-x}
\forall C>0\qquad \sum_{n\in\Z} r_n\cdot C^{|n|}<\infty
 \eeq
the set
 \beq\label{E_r-v-H(C-x)}
E_r=\{u\in {\mathcal O}(\C^\times):\; \forall n\in\Z\quad |u_n|\le r_n\}
 \eeq
is compact in ${\mathcal O}(\C^\times)$ since it is closed and is contained in
the rectangle $f^\text{\BSQ}$, where
$$
f(t)=\sum_{n\in\Z} r_n\cdot |t|^n,\qquad t\in\C^\times
$$
Hence, $E_r$ generates a continuous seminorm $\alpha\mapsto \max_{u\in
E_r}|\langle u,\alpha\rangle|$ on ${\mathcal O}^\star(\C^\times)$. But this is
exactly the seminorm from \eqref{|alpha|_r-v-C-x}:
$$
\max_{u\in E_r}|\langle u,\alpha\rangle|=\max_{|u|\le r_n} \left|\sum_{n\in\Z}
u_n\cdot\alpha_n \right|=\sum_{n\in\Z} r_n\cdot|\alpha_n|=|||\alpha|||_r
$$
It remains to verify that seminorms $|||\cdot|||_r$ indeed generate the
topology of the space ${\mathcal O}^\star(\C^\times)$. This follows from: \epr

\blm\label{LM-polunormy-v-H*(C-x)-1} If $p$ is a continuous seminorm on
${\mathcal O}^\star(\C^\times)$, then the family of numbers
$$
r_n=p(\zeta_n)
$$
satisfy condition \eqref{forall-R>0-sum_r_n-R^n<infty-x}, and $p$ is majorized
by the seminorm $|||\cdot|||_r$:
 \beq\label{p(alpha)-le-|||alpha|||_r-O^(C^x)}
p(\alpha)\le |||\alpha|||_r
 \eeq
 \elm
\bpr The set
$$
D=\{u\in{\mathcal O}(\C^\times):\; \sup_{\alpha\in {\mathcal
O}^\star(\C^\times):\; p(\alpha)\le 1}|\alpha(u)|\le 1\}
$$
is compact in ${\mathcal O}(\C^\times)$ generating the seminorm $p$:
$$
p(\alpha)=\sup_{u\in D}|\alpha(u)|
$$
Therefore,
$$
r_n=p(\zeta_n)=\sup_{u\in D} |\zeta_n(u)|=\sup_{u\in D} |u_n|
$$
For any $C>0$ we have:
$$
\infty>\sup_{u\in D}||u||_{C+1}=\sup_{u\in D}\sum_{n\in\Z} |u_n|\cdot
(C+1)^{|n|}\ge \sup_{u\in D}\sup_n |u_n|\cdot (C+1)^{|n|}=\sup_n \Big(r_n\cdot
(C+1)^{|n|}\Big)
$$
$$
\Downarrow
$$
$$
\exists M>0\quad \forall n\in\Z\qquad r_n\le \frac{M}{(C+1)^{|n|}}
$$
$$
\Downarrow
$$
$$
\sum_{n\in\Z} r_n\cdot C^{|n|} \le \sum_{n\in\Z} \frac{M}{(C+1)^{|n|}}\cdot
C^{|n|}<\infty
$$
i.e. $r_n$ indeed satisfy \eqref{forall-R>0-sum_r_n-R^n<infty-x}. The formula
\eqref{p(alpha)-le-|||alpha|||_r-O^(C^x)} is proved by a sequence of
inequalities analogous to \eqref{tsep-p-le-|||.|||_r}. \epr

\paragraph{The chain
${\mathcal R}(\C)\subset{\mathcal O}(\C)\subset{\mathcal
O}^\star(\C)\subset{\mathcal R}^\star(\C)$.}

By symbol ${\mathcal R}(\C)$ we denote the usual algebra of polynomials on the
complex plane $\C$. Let $t^k$ denote the monomial of degree $k\in\N$ on $\C$:
 \beq\label{def-t^k}
t^k(x):=x^k,\qquad x\in\C,\; k\in\N
 \eeq
Every polynomial $u\in {\mathcal R}(\C)$ is uniquely represented by the series
(with finite number of nonzero terms)
 \beq\label{razlozhenie-u-v-R(C)}
u=\sum_{k\in\N} u_k\cdot t^k,\qquad u_k\in\C: \quad \card\{k\in\N: \; u_k\ne
0\}<\infty
 \eeq
so the monomials $t^k$ form an algebraic basis in the space ${\mathcal R}(\C)$.
The multiplication in ${\mathcal R}(\C)$ is the usual multiplication of
polynomials
$$
u\cdot v=\sum_{k\in\N} \left(\sum_{i=0}^n u_{k-i}\cdot v_i\right)\cdot t^k
$$
and the topology in ${\mathcal R}(\C)$ is defined as the strongest locally
convex topology. This implies

 \bprop The mapping
$$
u\in{\mathcal R}(\C)\mapsto \{u_k;\;k\in\N\}\in\C_{\N}
$$
is an isomorphism of topological vector spaces
  \beq\label{R(C)=C_N}
{\mathcal R}(\C)\cong \C_{\N}
  \eeq
 \eprop

\brem We can interpret formula \eqref{R(C)=C_N} as isomorphism of algebras, if
the multiplication in $\C_{\N}$ is defined by the same formula
\eqref{umnozhenie-na-bazise-v-C_G-infty} as for algebra $\C_G$ of point charges
on the group $G$ (defined in \ref{SEC:ster-algebry-Hopfa}\ref{C^G-i-C_G}), we
should only remember here that the set $\N$ is not a group, but a monoid with
respect to the additive operation used for this purpose.

On the other hand formula \eqref{R(C)=C_N} is not an isomorphism of Hopf
algebras, in particular, since $\C_{\N}$ is not a Hopf algebra at all with
respect to the operations defined in
\ref{SEC:ster-algebry-Hopfa}\ref{C^G-i-C_G}: for $x\in\N$ the inverse element
$x^{-1}=-x$ does not exist when $x\ne 0$, so the antipode here cannot be
defined by equality \eqref{antipode-in-C_G}.\erem

\bprop The algebra ${\mathcal R}(\C)$ of polynomials on the complex plane $\C$
is a nuclear Hopf-Brauner algebra with the topology generated by seminorms
 \beq
\norm{u}_r=\sum_{k\in\N}r_k\cdot |u_k|,\qquad (r_k\ge 0)
 \eeq
and algebraic operations defined on monomials $t^k$ by formulas
 \begin{align}
\label{umn-v-R(C)} & t^k\cdot t^l=t^{k+l} && 1_{{\mathcal R}(\C)}= t^0 \\
\label{koumn-v-R(C)}  &\varkappa(t^k)=\sum_{i=0}^k\begin{pmatrix}k
\\ i\end{pmatrix}\cdot t^{k-i}\odot t^i & &
\e(t^k)=\begin{cases}1,& k=0\\ 0,& k>0\end{cases} \\
\label{antipode-in-R(C)} & \sigma(t^k)=(-1)^k\cdot t^k &&
 \end{align}
 \eprop
\bpr The algebraic operations which we did not find yet -- comultiplication,
counit and antipode -- are computed by formulas
\eqref{koumn-v-O(G)}-\eqref{antipod-v-O(G)}. For instance, comultiplication:
 $$
\widetilde{\varkappa}(t^k)(x,y)=t^k(x+y)=(x+y)^k=\sum_{i=0}^k\begin{pmatrix}k
\\ i\end{pmatrix}\cdot x^{k-i}\cdot y^i=\sum_{i=0}^k\begin{pmatrix}k
\\ i\end{pmatrix}\cdot t^{k-i}\boxdot t^i(x,y)
 $$
 $$
 \Downarrow
 $$
 $$
 \widetilde{\varkappa}(t^k)= \sum_{i=0}^k\begin{pmatrix}k
\\ i\end{pmatrix}\cdot t^{k-i}\boxdot t^i
 $$
 $$
 \Downarrow
 $$
 $$
\varkappa(t^k)=\rho_{\C,\C}(\widetilde{\varkappa}(t^k))=
\sum_{i=0}^k\begin{pmatrix}k
\\ i\end{pmatrix}\cdot t^{k-i}\odot t^i
 $$
\epr

Following the terminology of \cite{Akbarov}, we call a {\it current of degree
$0$}\index{current}, or a {\it current} on $\C$ an arbitrary linear functional
$\alpha:{\mathcal R}(\C)\to\C$ on the space of polynomials ${\mathcal R}(\C)$
(every such a functional is automatically continuous). Formula \eqref{R(C)=C_N}
implies that every compact set in the space of polynomials ${\mathcal R}(\C)$
is finite-dimensional, hence it is contained in a convex hull of a finite set
of basis monomials $t^k$. Therefore the topology of the space ${\mathcal
R}^\star(\C)$ of currents on $\C$ (which is formally defined as the topology of
uniform convergence on compact sets in ${\mathcal R}(\C)$) can be defined as
the topology of convergence on monomials $t^k$, i.e. the topology generated by
seminorms
$$
\norm{\alpha}_N=\sum_{k\in N} \abs{\alpha(t^k)},
$$
where $N$ is an arbitrary finite set in $\N$.

The functional of taking derivative of the $k$-th derivative in the point $0$
is a typical example of a current:
 \beq\label{def-tau^k}
\tau^k(u)=\left(\frac{\d^k}{\d x^k}u(x)\right)\Bigg|_{x=0},\qquad u\in{\mathcal
R}(\C)
 \eeq
By the Taylor theorem, these functionals are connected with the coefficients
$u_k$ in the decomposition \eqref{razlozhenie-u-v-R(C)} of a polynomial
$u\in{\mathcal R}(\C)$ through the formula
 \beq
u_k=\frac{1}{k!}\cdot\tau^k(u)
 \eeq
Therefore the action of a current $\alpha\in {\mathcal R}^\star(\C)$ on a
polynomial $u\in{\mathcal R}(\C)$ can be written by formula
 \beq \label{deistvie-alpha-na-u-v-R(C)}
\alpha(u)=\alpha\left(\sum_{k\in\N}u_k\cdot t^k\right)=\sum_{k\in\N}u_k\cdot
\alpha(t^k)=\sum_{k\in\N}\frac{1}{k!}\cdot\tau^k(u)\cdot
\alpha(t^k)=\left(\sum_{k\in\N}\frac{1}{k!}\cdot
\alpha(t^k)\cdot\tau^k\right)(u)
 \eeq
This means that $\alpha$ is decomposed into a series in terms of $\tau^k$:
 \beq \label{razlozhenie-alpha-v-R(C)}
\alpha=\sum_{k\in\N} \alpha_k\cdot \tau^k,\qquad
\alpha_k=\frac{1}{k!}\cdot\alpha(t^k)
 \eeq
We can deduce from this that currents $\tau^k$ form a basis in the topological
vector spaces ${\mathcal R}^\star(\C)$: every functional $\alpha\in{\mathcal
R}^\star(\C)$ can be uniquely represented by a converging in ${\mathcal
R}^\star(\C)$ series \eqref{razlozhenie-alpha-v-R(C)}, where coefficients
$\alpha_k\in\C$ continuously depend on $\alpha\in{\mathcal R}^\star(\C)$.

From \eqref{deistvie-alpha-na-u-v-R(C)} it follows that the action of the
current $\alpha$ on a polynomial $u$ is defined by formula
$$
\alpha(u)=\sum_{k\in\N}\underbrace{\frac{1}{k!}\cdot\tau^k(u)}_{u_k}\cdot
\underbrace{\frac{1}{k!}\cdot \alpha(t^k)}_{\alpha_k}\cdot k!=\sum_{k\in\N}
u_k\cdot\alpha_k\cdot k!
$$
and the basis currents $\tau^k$ act on monomials $t^k$ by formula
$$
\tau^k(t^n)=\begin{cases}0,& n\ne k\\ n!,& n=k\end{cases},
$$
This means that the system $\tau^k$ {\it is not a dual basis} for $t^k$: it
differs from the dual basis by the coefficients $k!$. However, since ${\mathcal
R}^\star(\C)\cong \C^{\N}$ (this follows from \eqref{R(C)=C_N}) and, on the
other hand, in $\C^{\N}$ any two bases are isomorphic (Theorem
\ref{isomorph-bazisov}), we have

\bprop The mapping
$$
\alpha\in{\mathcal R}^\star(\C)\mapsto \{\alpha_k;\; k\in\N\}\in\C^{\N}
$$
is an isomorphism of topological vector spaces:
 \beq\label{R*(C)=C^N}
{\mathcal R}^\star(\C)\cong \C^{\N}
 \eeq
 \eprop
This isomorphism is not however, an isomorphism of algebras, since the
multiplication in ${\mathcal R}^\star(\C)$ is not coordinate-wise (i.e. not by
formula \eqref{umnozhenie-v-C^G-infty-po-bazisu}, as it could defined in
$\C^{\N}$ by analogy with the case of
\ref{SEC:ster-algebry-Hopfa}\ref{C^G-i-C_G}), but as power series:
 \beq\label{umnozhenie-v-R^star(C)}
\alpha * \beta=\sum_{k\in\N} \left(\sum_{i=0}^k \alpha_{k-i}\cdot
\beta_i\right)\cdot \tau^k,
 \eeq
This follows from formula \eqref{umn-v-R*(C)} below:

\bprop The algebra ${\mathcal R}^\star(\C)$ of currents of zero degree on the
complex plane $\C$ is a nuclear Hopf-Fr\'echet algebra with the topology
generated by seminorms
 \beq
\norm{\alpha}_K=\sum_{k=0}^K |\alpha_k|,\qquad (K\in\N)
 \eeq
and algebraic operations defined on basis elements $\tau^k$ by formulas
 \begin{align}
\label{umn-v-R*(C)} & \tau^k*\tau^l=\tau^{k+l} && 1_{{\mathcal R}^\star(\C)}=\tau^0 \\
\label{koumn-v-R*(C)}  &\varkappa(\tau^k)=\sum_{i=0}^k\begin{pmatrix}k
\\ i\end{pmatrix}\cdot \tau^{k-i}\circledast\tau^i & &
\e(\tau^k)=\begin{cases}1,& k=0\\ 0,& k>0\end{cases} \\
\label{antipode-in-R*(C)} & \sigma(\tau^k)=(-1)^k\cdot \tau^k &&
 \end{align}
 \eprop
\bpr There are many ways to prove these formulas, for instance, to prove
\eqref{umn-v-R*(C)} one can use formula \eqref{koumn-v-R(C)} for
comultiplication in ${\mathcal R}(\C)$:
 \begin{multline*}
(\tau^k*\tau^l)(t^m)=(\tau^k\circledast\tau^l)(\varkappa(t^m))=\eqref{koumn-v-R(C)}=
(\tau^k\circledast\tau^l)\left(\sum_{i=0}^m \begin{pmatrix}m \\ i \end{pmatrix}
\cdot t^{m-i}\odot t^i \right)=\\=\sum_{i=0}^m \begin{pmatrix}m \\ i
\end{pmatrix}\cdot \tau^k(t^{m-i})\cdot\tau^l(t^i)=\left\{\begin{matrix} \begin{pmatrix}m \\
i \end{pmatrix}\cdot (m-l)!\cdot l!,& m=k+l \\ 0,& m\ne
k+l\end{matrix}\right\}=\left\{\begin{matrix} m!,& m=k+l \\ 0,& m\ne
k+l\end{matrix}\right\}=\tau^{k+l}(t^m).
 \end{multline*}
(on each monomial $t^m$ the action of functionals $\tau^k*\tau^l$ and
$\tau^{k+l}$ coincide, hence they coincide themselves). \epr

Let us consider now algebra ${\mathcal O}(\C)$ of entire functions and its dial
algebra ${\mathcal O}^\star(\C)$ of analytical functionals on complex plane
$\C$. Again, let $t^k$ denote the monomial of degree $k\in\N$ on $\C$, and
$\tau^k$ the functional of taking $k$-th derivative in the point 0:
$$
t^k(z)=z^k\qquad \tau^k(u)=\left(\frac{\d^k}{\d
z^k}u(z)\right)\Bigg|_{z=0}\qquad \Big(x\in\C,\qquad u\in{\mathcal O}(\C)\Big)
$$
Then it is convenient to represent elements ${\mathcal O}(\C)$ and ${\mathcal
O}^\star (\C)$ as (converging in these spaces) series
 \begin{align}
\label{razlozhenie-u-v-C} u &\in {\mathcal O}(\C) & & \Longleftrightarrow &
u&=\sum_{n=0}^\infty u_n\cdot t^n, & u_n=\frac{1}{n!}\tau^n(u)\in\C:\quad &
\forall C>0\quad \sum_{n=0}^\infty|u_n|\cdot C^n<\infty
 \\
\label{razlozhenie-alpha-v-C} \alpha &\in {\mathcal O}^\star(\C) & &
\Longleftrightarrow & \alpha &=\sum_{n=0}^\infty \alpha_n\cdot \tau^n,  &
\alpha_n=\frac{1}{n!}\alpha(t^n)\in\C:\quad & \exists M,C>0\quad \forall
n\in\N\quad |\alpha_n|\le M\cdot \frac{C^n}{n!}
 \end{align}
The action of an analytical functional $\alpha$ on an entire function $u$ is
described by formula
$$
\langle u,\alpha\rangle=\sum_{n=0}^\infty u_n\cdot\alpha_n\cdot n!
$$
and the multiplications in ${\mathcal O}(\C)$ and in ${\mathcal O}^\star (\C)$
are defined by the same formulas as for the usual power series:
 \beq\label{umnozhenie-v-O(C)-i-O^star(C)}
u\cdot v=\sum_{n=0}^\infty \left(\sum_{i=0}^n u_i\cdot v_{n-i}\right)\cdot t^n,
\qquad \alpha * \beta=\sum_{k\in\N} \left(\sum_{i=0}^k \alpha_i\cdot
\beta_{k-i}\right)\cdot \tau^k,
 \eeq

\bprop The algebra ${\mathcal O}(\C)$ of entire functions on complex plane $\C$
is a nuclear Hopf-Fr\'echet algebra with the topology generated by seminorms
 \beq\label{polunormy-v-H(C)}
\norm{u}_C=\sum_{k\in\N} |u_k|\cdot C^k\qquad (C\ge 1)
 \eeq
and algebraic operations defined on monomials $t^k$ by the same formulas
\eqref{umn-v-R(C)}-\eqref{antipode-in-R(C)} as for ${\mathcal R}(\C)$:
 \begin{align*}
& \mu(t^k\odot t^l)=t^{k+l} && 1_{{\mathcal
R}(\C)}= t^0 \\
&\varkappa(t^k)=\sum_{i=0}^k\begin{pmatrix}k
\\ i\end{pmatrix}\cdot t^{k-i}\odot t^i & &
\e(t^k)=\begin{cases}1,& k=0\\ 0,& k>0\end{cases} \\
& \sigma(t^k)=(-1)^k\cdot t^k &&
 \end{align*}
 \eprop
\bpr The formulas for algebraic operations are proved similarly with the same
formulas for ${\mathcal R}(\C)$, so it remains only to explain, why the
topology is generated by seminorms \eqref{polunormy-v-H(C)}. This is a subject
of mathematical folklore (see e.g. \cite{Pirkovskii}): clearly, every usual
seminorm $|u|_K=\max_{z\in K}|u(z)|$, where $K$ is a compact set in $\C$, is
majorized by some seminorm $||u||_C$, namely by the one with $C=\max_{z\in
K}|z|$:
$$
|u|_K=\max_{z\in K}|u(z)|=\max_{z\in K}\left|\sum_{n=0}^\infty u_n\cdot
z^n\right|\le\max_{z\in K}\sum_{n=0}^\infty |u_n|\cdot
|z^n|\le\sum_{n=0}^\infty |u_n|\cdot C^n=||u||_C
$$
On the contrary, for every $C>0$ we can take a compact set $K=\{z\in\C:\;
|z|\le C+1\}$, and then from formulas for Cauchy coefficients
$$
|u_n|\le \frac{|u|_K}{(C+1)^n}
$$
it follows that $||u||_C$ is subordinated to $|u|_K$ (with the constant $C+1$):
$$
||u||_C=\sum_{n=0}^\infty |u_n|\cdot C^n\le \sum_{n=0}^\infty
\frac{|u|_K}{(C+1)^n}\cdot C^n\le |u|_K\cdot\sum_{n=0}^\infty
\left(\frac{C}{C+1}\right)^n=|u|_K\cdot
\frac{1}{1-\frac{C}{C+1}}=(C+1)\cdot|u|_K
$$
\epr

\bprop\label{PROP:O*(C)-Hopf} The algebra ${\mathcal O}^\star(\C)$ of
analytical functionals on complex plane $\C$ is a nuclear Hopf-Brauner algebra
with the topology generated by seminorms
 \beq\label{|alpha|_r-v-C}
|||\alpha|||_r=\sum_{k\in\N} r_k\cdot |\alpha_k|\cdot k!\qquad\left(r_k\ge 0:\
\forall C>0\ \sum_{k\in\N} r_k\cdot C^k<\infty \right)
 \eeq
and algebraic operations defined on basis elements $\tau^k$ by the same
formulas \eqref{umn-v-R*(C)}-\eqref{antipode-in-R*(C)} as in the case of
${\mathcal R}^\star(\C)$:
 \begin{align*}
& \mu(\tau^k\circledast\tau^l)=\tau^{k+l} && 1_{{\mathcal
R}^\star(\C)}=\tau^0 \\
&\varkappa(\tau^k)=\sum_{i=0}^k\begin{pmatrix}k
\\ i\end{pmatrix}\cdot \tau^{k-i}\circledast\tau^i & &
\e(\tau^k)=\begin{cases}1,& k=0\\ 0,& k>0\end{cases} \\
& \sigma(\tau^k)=(-1)^k\cdot \tau^k &&
 \end{align*}
 \eprop
 \bpr
Here again the formulas for algebraic operations are proved similarly with the
case of ${\mathcal R}^\star(\C)$, so we need only to explain, why the topology
is generated by seminorms \eqref{|alpha|_r-v-C}. Note first that for any
sequence of non-negative numbers $r_k\ge 0$, satisfying the condition
 \beq\label{forall-R>0-sum_r_n-R^n<infty}
\forall C>0\qquad \sum_{k\in\N} r_k\cdot C^k<\infty
 \eeq
the set
 \beq\label{E_r-v-H(C)}
E_r=\{u\in {\mathcal O}(\C):\; \forall k\in\N\quad |u_k|\le r_k\}
 \eeq
is compact in ${\mathcal O}(\C)$, since it is closed and is contained in the
rectangle $f^\text{\BSQ}$, where
$$
f(x)=\sum_{k\in\N} r_k\cdot |x|^k,\qquad x\in\C
$$
Hence, $E_r$ generates a continuous seminorm $\alpha\mapsto \max_{u\in
E_r}|\langle u,\alpha\rangle|$ on ${\mathcal O}^\star(\C)$. But this is exactly
seminorm \eqref{|alpha|_r-v-C}:
$$
\max_{u\in E_r}|\langle u,\alpha\rangle| =\max_{|u|\le r_k} \left|\sum_{k\in\N}
u_k\cdot\alpha_k\cdot k! \right|=\sum_{k\in\N} r_k\cdot|\alpha_k|\cdot
k!=|||\alpha|||_r
$$
It remains to verify that the seminorms $|||\cdot|||_r$ indeed generate the
topology of the space ${\mathcal O}^\star(\C^\times)$. This follows from: \epr

\blm\label{LM-polunormy-v-H*(C)} If $p$ is a continuous seminorm on ${\mathcal
O}^\star(\C)$, then the numbers $r_k=\frac{1}{k!}p(\tau^k)$ satisfy the
condition \eqref{forall-R>0-sum_r_n-R^n<infty}, and $p$ is majorized by the
seminorm $|||\cdot|||_r$:
 \beq\label{p(alpha)-le-|||alpha|||_r}
p(\alpha)\le |||\alpha|||_r
 \eeq
\elm
 \bpr The set
$$
D=\{u\in{\mathcal O}(\C):\; \sup_{\alpha\in {\mathcal O}^\star(\C):\;
p(\alpha)\le 1}|\alpha(u)|\le 1\}
$$
is compact in ${\mathcal O}(\C)$ and generates the seminorm $p$:
$$
p(\alpha)=\sup_{u\in D}|\alpha(u)|
$$
Hence
$$
r_k=\frac{1}{k!}p(\tau^k)=\frac{1}{k!}\sup_{u\in D} |\tau^k(u)|=\sup_{u\in D}
|u_k|
$$
For each $C>0$ we have:
$$
\infty>\sup_{u\in D}||u||_{C+1}=\sup_{u\in D}\sum_{k\in\N} |u_k|\cdot
(C+1)^k\ge \sup_{u\in D}\sup_k |u_k|\cdot (C+1)^k=\sup_k \Big(r_k\cdot
(C+1)^k\Big)
$$
$$
\Downarrow
$$
$$
\exists M>0\quad \forall k\in\N\qquad r_k\le \frac{M}{(C+1)^k}
$$
$$
\Downarrow
$$
$$
\sum_{k\in\N} r_k\cdot C^k \le \sum_{k\in\N} \frac{M}{(C+1)^k}\cdot C^k<\infty
$$
I.e., $r_k$ indeed satisfy condition \eqref{forall-R>0-sum_r_n-R^n<infty}.
Formula \eqref{p(alpha)-le-|||alpha|||_r} is proved by the chain of
inequalities analogous to \eqref{tsep-p-le-|||.|||_r}.
 \epr

\bprop The mappings
$$
t^k\mapsto t^k\mapsto \tau^k\mapsto\tau^k \qquad (k\in\N)
$$
define the chain of homomorphisms of rigid Hopf algebras:
 \beq\label{R->O->O*->R*}
{\mathcal R}(\C)\to {\mathcal O}(\C)\to {\mathcal O}^\star(\C)\to {\mathcal
R}^\star(\C)
 \eeq
\eprop

\section{Functions of exponential type on a Stein group} \label{SEC-O_exp(G)}

\subsection{Semicharacters and inverse semicharacters on Stein groups}

Let $G$ be a Stein group. Then
 \bit
\item[---] a locally bounded function $f:G\to[1,+\infty)$ is called a {\it
semicharacter}\index{semicharacter}, if it satisfies the following so called
{\it submultiplicativity inequality}:
 \beq
f(x\cdot y)\le f(x)\cdot f(y),\qquad x,y\in G
 \eeq

\item[---] a function $\ph:G\to(0;1]$, locally separated from zero, is called
{\it inverse semicharacter}\index{inverse semicharacter}, if it satisfies the
inverse inequality:
 \beq
\ph(x)\cdot \ph(y)\le\ph(x\cdot y),\qquad x,y\in G
 \eeq
 \eit
Clearly, if $f:G\to[1,+\infty)$ is a semicharacter, then the inverse function
 \beq\label{poluharakter-obratnyi-poluharakter}
\ph(x)=\frac{1}{f(x)}
 \eeq
is an inverse semicharacter, and vice versa.

\bigskip

\centerline{\bf Properties of semicharacters and inverse
semicharacters:}\label{PR-f+g}

 \bit
\item[(i)] The set of all semicharacters on $G$ is closed under the following
operations:
 \bit
\item[---] multiplication by a sufficiently big constant: $C\cdot f$ $\quad
(C\ge 1)$,

\item[---] multiplication: $f\cdot g$,

\item[---] addition: $f+g$,

\item[---] taking maximum: $\max\{f,g\}$.
 \eit

\item[(ii)] The set of all inverse semicharacters on $G$ is closed under the
following operations:
 \bit
\item[---] multiplication by a sufficiently small constant: $C\cdot \ph$ $\quad
(C\le 1)$,

\item[---] multiplication: $\ph\cdot\psi$,

\item[---] taking half of harmonic mean: $\frac{\ph\cdot\psi}{\ph+\psi}$

\item[---] taking minimum: $ \min\{\ph,\psi\}.$
 \eit
 \eit
\bpr Having in mind the duality between semicharacters and inverse
semicharacters reflected by formula \eqref{poluharakter-obratnyi-poluharakter},
we can consider only the case of semicharacters.

If $f$ is a semicharacter on $G$ and $C\ge 1$, then
$$
C\cdot f(x\cdot y)\le C\cdot f(x)\cdot f(y)\le \big(C\cdot f(x)\big)\cdot
\big(C\cdot f(y)\big)
$$
If $f$ and $g$ are semicharacters on $G$, then considering their multiplication
we have:
 \begin{multline*}
(f\cdot g)(x\cdot y)=f(x\cdot y)\cdot g(x\cdot y)\le f(x)\cdot f(y)\cdot
g(x)\cdot g(y)=\\=f(x)\cdot g(x)\cdot f(y)\cdot g(y) =(f\cdot g)(x)\cdot
(f\cdot g)(y)
 \end{multline*}
For their sum we have:
 \begin{multline*}
(f+g)(x\cdot y)=f(x\cdot y)+g(x\cdot y)\le f(x)\cdot f(y)+ g(x)\cdot g(y)\le\\
\le f(x)\cdot f(y)+g(x)\cdot f(y)+f(x)\cdot g(y)+ g(x)\cdot
g(y)=\Big(f(x)+g(x)\Big)\cdot \Big(f(y)+g(y)\Big) =(f+g)(x)\cdot (f+g)(y)
 \end{multline*}
The proof for maximum $\max\{f,g\}$ is based on the following evident
inequality:
 \beq\label{max}
\max\{a\cdot b, c\cdot d\}\le \max\{a, c\}\cdot \max\{b, d\}\qquad (a,b,c,d>0)
 \eeq
It implies:
 \begin{multline*}
\max\{f,g\}(x\cdot y)=\max\{f(x\cdot y), g(x\cdot y)\}\le \max\{f(x)\cdot f(y),
g(x)\cdot g(y)\}\le\eqref{max}\le\\ \le \max\{f(x), g(x)\}\cdot
\max\{f(y),g(y)\}= \max\{f, g\}(x)\cdot \max\{f,g\}(y)
 \end{multline*}
 \epr

\bex All the submultiplicative matrix norms (see \cite{Horn-Johnson}), for
instance
 \beq\label{matrichnye-normy}
||x||=\sum_{i,j=1}^n|x_{i,j}|,\qquad ||x||=\sqrt{\sum_{i,j=1}^n|x_{i,j}|^2}
 \eeq
are semicharacters on $\GL_n(\C)$. From the properties of semicharacters it
follows that the function of the form
 \beq\label{r_N(x)=||x||^N-cdot-||x^-1||^N}
r^N_C(x)=C\cdot\max\{||x||; ||x^{-1}||\}^N,\qquad C\ge 1,\; N\in\N
 \eeq
are again semicharacters on $\GL_n(\C)$.  \eex

\bprop\label{PROP-stroenie-poluharakterov-na-GL_n} For any submultiplicative
matrix norm $||\cdot||$ on $\GL_n(\C)$ the semicharacters of the form
\eqref{r_N(x)=||x||^N-cdot-||x^-1||^N} majorize all other semicharacters on
$\GL_n(\C)$. \eprop

\bpr Note from the very beginning that it is sufficient to consider the case
when $||\cdot||$ is the Euclid norm on the algebra $\M_n(\C)$ of all complex
matrices $n\times n$:
$$
||x||=\sup_{\xi\in\C^n: ||\xi||\le 1}||x(\xi)||,\qquad\text{where}\quad
||\xi||=\sqrt{\sum_{i=1}^n |\xi_i|^2},
$$
-- since any other norm on $\M_n(\C)$ majorizes the Euclid seminorm up to a
constant multiplier, this will prove our proposition.

Consider the set
$$
K:=\big\{x\in\GL_n(\C):\; \max\{||x||; ||x^{-1}||\}\le 2 \big\}=
\big\{x\in\GL_n(\C):\;\forall\xi\in\C^n\quad
\frac{1}{2}\cdot||\xi||\le||x(\xi)||\le 2\cdot||\xi|| \big\}
$$
It is closed and bounded in the algebra of matrices $\M_n(\C)$, hence it is
compact. This is a generating compact in $\GL_n(\C)$
 \beq\label{GL_n(C)=bigcup}
\GL_n(\C)=\bigcup_{m=1}^\infty K^m,\qquad K^m=\underbrace{K\cdot...\cdot
K}_{\text{$m$ factors}},
 \eeq
In addition,
 \beq\label{K^m:=max-||x||^m;||x^-1||^m-le-2^m}
K^m=\big\{x\in\GL_n(\C):\; \max\{||x||; ||x^{-1}||\}^m\le 2^m \big\}
 \eeq
Indeed, if $x\in K^m$, then $x=x_1\cdot...\cdot x_m$, where $||x_i||\le 2$ and
$||x_i^{-1}||\le 2$, hence
$$
||x||\le \prod_{i=1}^m ||x_i||\le 2^m,\qquad ||x^{-1}||\le \prod_{i=1}^m
||x_i^{-1}||\le 2^m
$$
On the contrary, suppose $\max\{||x||^m; ||x^{-1}||^m\}\le 2^m$. Consider the
polar decomposition: $x=r\cdot u$, where $r$ is a positive definite Hermitian,
and $u$ a unitary matrices. Let us decompose $r$ into a product $r=v\cdot
d\cdot v^{-1}$, where $v$ is unitary, and $d$ a diagonal matrices:
$$
d=\begin{pmatrix}d_1 & 0 & ... & 0 \\ 0 & d_2 & ... & 0 \\ ... & ... & ... &
... \\ 0 & 0 & ... & d_n \\
 \end{pmatrix},\qquad d_i\ge 0
$$
The $m$-th root
$$
\root{m}\of{d}=\begin{pmatrix}\root{m}\of{d_1} & 0 & ... & 0 \\ 0 &
\root{m}\of{d_2} & ... & 0
\\ ... & ... & ... &
... \\ 0 & 0 & ... & \root{m}\of{d_n} \\
 \end{pmatrix}
$$
has the following properties:
$$
\max_{i=1,...,n} d_i=||d||=||r||=||x||\le 2^m\quad\Longrightarrow\quad
||\root{m}\of{d}||=\max_{i=1,...,n} \root{m}\of{d_i}\le 2,
$$
$$
\max_{i=1,...,n} d_i^{-1}=||d^{-1}||=||r^{-1}||=||x^{-1}||\le
2^m\quad\Longrightarrow\quad ||\root{m}\of{d^{-1}}||=\max_{i=1,...,n}
\root{m}\of{d_i^{-1}}\le 2
$$
As a corollary, the matrix
$$
y=v\cdot \root{m}\of{d}\cdot v^{-1}
$$
belongs to $K$:
$$
\max\{||y||; ||y^{-1}||\}=\max\{||\root{m}\of{d}||;
||\root{m}\of{d}^{-1}||\}\le 2 \quad\Longrightarrow\quad y\in K
$$
Hence the matrix $y\cdot u$ also belongs to $K$, and we obtain
$$
x=v\cdot d\cdot v^{-1}\cdot u=v\cdot (\root{m}\of{d})^{m-1}\cdot v^{-1}\cdot
v\cdot (\root{m}\of{d})\cdot v^{-1}\cdot u= \underbrace{y^{m-1}}_{\scriptsize
\begin{matrix}\text{\rotatebox{-90}{$\in$}} \\ K^{m-1}\end{matrix}}\cdot \underbrace{y\cdot v}_{\scriptsize
\begin{matrix}\text{\rotatebox{-90}{$\in$}} \\ K\end{matrix}}\in
K^m
$$
We have proved formula \eqref{K^m:=max-||x||^m;||x^-1||^m-le-2^m}. Now let $f$
be an arbitrary semicharacter on $\GL_n(\C)$. Put
$$
C=\sup_{x\in K}f(x),\qquad N\ge \log_2 C
$$
and let us show that
 \beq
\forall x\in\GL_n(\C)\qquad f(x)\le r^N_C(x)
 \eeq
Take $x\in\GL_n(\C)$. From \eqref{GL_n(C)=bigcup} it follows that $x\in
K^m\setminus K^{m-1}$ for some $m\in\N$. By formula
\eqref{K^m:=max-||x||^m;||x^-1||^m-le-2^m} we have:
$$
2^{m-1}<\max\{||x||; ||x^{-1}||\}\le 2^m
$$
$$
\Downarrow
$$
$$
m-1<\log_2 \Big(\max\{||x||; ||x^{-1}||\}\Big)
$$
$$
\Downarrow
$$
 \begin{multline*}
f(x)\le\sup_{y\in K}f(y)^m\le C^m\le C\cdot C^{m-1}<C\cdot C^{\log_2
\Big(\max\{||x||; ||x^{-1}||\}\Big)}=\\ =C\cdot \Big(\max\{||x||;
||x^{-1}||\}\Big)^{\log_2 C} \le C\cdot \Big(\max\{||x||;
||x^{-1}||\}\Big)^N=r^N_C(x)
 \end{multline*}
\epr

In the special case when $n=1$ we have:

\bcor On the complex circle $\C^\times$ the semicharacters of the form
$$
r_C^N(t)=C\cdot \max\{|t|;|t|^{-1}\}^N,\qquad C\ge 1,\; N\in\N
$$
majorize all other semicharacters. \ecor

\bprop\label{PROP-poluh-porozhd-K} If $G$ is a compactly generated Stein group
and $K$ a compact neighborhood of identity in $G$, generating $G$,
 $$
G=\bigcup_{n=1}^\infty K^n,\qquad K^n=\underbrace{K\cdot...\cdot K}_{\text{$n$
factors}},
 $$
then for any $C\ge 1$ the rule
 \beq\label{poluh-opred-K}
h_C(x)=C^n\quad\Longleftrightarrow\quad x\in K^n\setminus K^{n-1}
 \eeq
defines a semicharacter $h_C$ on $G$. Such semicharacters form a fundamental
system among all semicharacters on $G$: every semicharacter $f$ on $G$ is
majorized by some semicharacter $h_C$
$$
f(x)\le h_C(x),\qquad x\in G,
$$
-- for this the constant $C$ should be chosen such that
 \beq\label{C-ge-max-t-in-K-f(t)}
C\ge\max_{t\in K}f(t)
 \eeq
 \eprop
\bpr The local boundedness of $h_C$ is obvious, so we need to check the
submultiplicativity inequality. Take $x,y\in G$ and choose $k,l\in\N$ such that
$$
x\in K^k\setminus K^{k-1},\qquad y\in K^l\setminus K^{l-1}.
$$
Then $x\cdot y\in K^{k+l}$, so
$$
h_C(x\cdot y)\le C^{k+l}=C^k\cdot C^l=h_C(x)\cdot h_C(y)
$$

If now $f$ is an arbitrary semicharacter, and $C$ satisfies the condition
\eqref{C-ge-max-t-in-K-f(t)}, then
$$
x\in K^n\setminus K^{n-1} \quad\Longrightarrow\quad f(x)\le \Big(\max_{t\in K}
f(t)\Big)^n=C^n=h_C(x)
$$
\epr

\subsection{Submultiplicative rhombuses and dually submultiplicative rectangles}

Let us introduce the following definitions:

 \bit
\item[---] a closed absolutely convex neighborhood of zero $\varDelta$ in
$\mathcal O^\star(G)$ is said to be {\it
submultiplicative}\label{DEF:submult-neighb-zero}\index{submultiplicative!neighborhood
of zero}, if for any functionals $\alpha,\beta$ from $\varDelta$ their
convolution $\alpha*\beta$ also belong to $\varDelta$:
$$
\forall\alpha,\beta\in\varDelta\qquad \alpha*\beta\in\varDelta
$$
shortly this is expressed by the following inclusion:
$$
\varDelta * \varDelta\subseteq \varDelta
$$

\item[---] an absolutely closed compact set $D$ in $\mathcal O(G)$ is said to
be {\it dually submultiplicative}\index{dually submultiplicative!compact set},
if its polar $D^\circ$ is a submultiplicative neighborhood of zero:
$$
 D^\circ * D^\circ\subseteq D^\circ
$$
\eit

\blm \label{f-obr-poluh->F-dual-subm}$\phantom{.}$
 \bit
\item[(a)] If $\ph:G\to(0;1]$ is an inverse semicharacter on $G$, then its
rhombus $\ph^{\text{\BLZ}}$ is a (closed, absolutely convex, and)
submultiplicative neighborhood of zero in $\mathcal O^\star(G)$.

\item[(b)] If $\varDelta\subseteq \mathcal O^\star(G)$ is a closed absolutely
convex and submultiplicative neighborhood of zero in $\mathcal O^\star(G)$,
then its inner envelope $\varDelta^{\text{\LZ}}$ is an inverse semicharacter on
$G$. \eit \elm
 \bpr
1. Let us denote
$$
\e_x=\ph(x)\cdot\delta^x\in {\mathcal O}^\star(G),
$$
then
 \beq
\ph^{\text{\BLZ}} = \cabsconv \left\{\ph(x)\cdot\delta^x;\quad x\in
M\right\}=\cabsconv \left\{\e_x;\quad x\in M\right\}
 \eeq
and the inclusion
$$
\ph^{\text{\BLZ}} * \ph^{\text{\BLZ}}\subseteq \ph^{\text{\BLZ}}
$$
is verified in three steps. First we should note that
$$
 \forall x,y\in G\qquad \e_x*\e_y\in \ph^{\text{\BLZ}}
$$
Indeed,
 $$
 \e_x*\e_y=\ph(x)\cdot\delta^x * \ph(y)\cdot\delta^y=
 \ph(x)\cdot \ph(y)\cdot\delta^{x\cdot y}=
 \underbrace{\frac{\ph(x)\cdot \ph(y)}{\ph(x\cdot y)}}_{\scriptsize \begin{matrix}\text{\rotatebox{-90}{$\le$}} \\ 1\end{matrix}}
 \cdot\e_{x\cdot y}\in
 \cabsconv \left\{\e_z;\quad z\in M\right\}=\ph^{\text{\BLZ}}
 $$
Second, we take finite absolutely convex combinations of functionals $\e_x$:
 \beq\label{abs-vyp-komb}
\alpha=\sum_{i=1}^k \lambda_i\cdot \e_{x_i},\quad \beta=\sum_{j=1}^l \mu_j\cdot
\e_{y_j}\qquad \left( \sum_{i=1}^k |\lambda_i|\le 1,\quad \sum_{j=1}^l
|\mu_j|\le 1 \right)
 \eeq
For them we have:
 $$
\alpha * \beta=\sum_{i=1}^k \lambda_i\cdot \e_{x_i} * \sum_{j=1}^l \mu_j\cdot
\e_{y_j}=\sum_{1\le i\le k,1\le j\le l}
 \underbrace{\lambda_i\cdot\mu_j}_{\scriptsize
\begin{matrix}\uparrow \\ \sum_{1\le i\le k,1\le j\le l} |\lambda_i\cdot\mu_j| \\ \text{\rotatebox{-90}{$=$}}
\\ \sum_{i=1}^k |\lambda_i|\cdot \sum_{j=1}^l |\mu_j| \\
\text{\rotatebox{-90}{$\le$}}\\ 1
 \end{matrix}}\cdot
\underbrace{\e_{x_i} *\e_{y_j}}_{\scriptsize
\begin{matrix}\text{\rotatebox{-90}{$\in$}} \\ \\ \ph^{\text{\BLZ}} \end{matrix}}\in
 \underbrace{\ph^{\text{\BLZ}}}_{\scriptsize\begin{matrix} \text{absolutely} \\
\text{convex} \\
\text{set}\end{matrix}}
 $$
And, third, for arbitrary $\alpha,\beta\in \ph^{\text{\BLZ}}$ the inclusion
$\alpha*\beta\in \ph^{\text{\BLZ}}$ becomes a corollary of two facts: that
functionals \eqref{abs-vyp-komb} are dense in the set $\ph^{\text{\BLZ}} =
\cabsconv\left\{\e_x;\ x\in M\right\}$, and that the convolution $*$ is
continuous.

2. We use here formula \eqref{fi^lozenge}:
 \begin{multline*}
\varDelta * \varDelta\subseteq \varDelta\quad\Longrightarrow\quad \forall
x,y\in G\qquad \underbrace{\varDelta^{\text{\LZ}}(x)\cdot\delta^x}_{\scriptsize
\begin{matrix}\text{\rotatebox{-90}{$\in$}} \\ \varDelta, \\ \text{by \eqref{fi^lozenge}}\end{matrix}} *
\underbrace{\varDelta^{\text{\LZ}}(y)\cdot\delta^y}_{\scriptsize
\begin{matrix}\text{\rotatebox{-90}{$\in$}} \\ \varDelta, \\ \text{by \eqref{fi^lozenge}}\end{matrix}}
=\varDelta^{\text{\LZ}}(x)\cdot\varDelta^{\text{\LZ}}(y)\cdot\delta^{x\cdot y} \in
 \varDelta \quad\Longrightarrow\\ \Longrightarrow\quad
\varDelta^{\text{\LZ}}(x)\cdot\varDelta^{\text{\LZ}}(y)\le \max \{\lambda>0:\;
\lambda\cdot\delta^{x\cdot y}\in \varDelta
\}=\eqref{fi^lozenge}=\varDelta^{\text{\LZ}}(x\cdot y)
\quad\Longrightarrow\quad \varDelta^{\text{\LZ}}(x)\cdot\varDelta^{\text{\LZ}}(y)\le
\varDelta^{\text{\LZ}}(x\cdot y)
 \end{multline*}
 \epr

\blm \label{f-poluh->F-subm}$\phantom{.}$
 \bit
\item[(a)] If $f:G\to[1;\infty)$ is a semicharacter on $G$, then its rectangle
$f^\text{\BSQ}$ is dually submultiplicative.

\item[(b)] If $D\subseteq \mathcal O(G)$ is a dually submultiplicative
absolutely convex compacts set, then its outer envelope $D^\text{\SQ}$ is a
semicharacter on $G$. \eit \elm
 \bpr
1. If $f:G\to[1;\infty)]$ is a semicharacter, then $\frac{1}{f}$ is an inverse
semicharacter, therefore by Lemma \ref{f-obr-poluh->F-dual-subm} (a) the
rhombus
$$
\left(\frac{1}{f}\right)^{\text{\BLZ}}=\eqref{f^blacksquare^circ=frac-1-f^blacklozenge}=
(f^\text{\BSQ})^\circ
$$
is a submultiplicative neighborhood of zero. This means that $f^\text{\BSQ}$ is
dually submultiplicative.

2. If $D\subseteq \mathcal O(G)$ is a dually submultiplicative absolutely
convex compact set, then its polar $D^\circ\subseteq \mathcal O^\star(G)$ is a
submultiplicative neighborhood of zero in $\mathcal O^\star(G)$, hence by Lemma
\ref{f-obr-poluh->F-dual-subm} (b) the inner envelope
$$
(D^\circ)^{\text{\LZ}}=\eqref{D^circ^lozenge=frac-1-D^square}=\frac{1}{D^\text{\SQ}}
$$
is an inverse semicharacter. Therefore, $D^\text{\SQ}$ must be a semicharacter.
 \epr

Lemmas \ref{f-obr-poluh->F-dual-subm} and \ref{f-poluh->F-subm} together with
formulas $\varDelta=\varDelta^{{\text{\LZ}}{\text{\BLZ}}}$ and $D=D^\text{\SQ\BSQ}$
for rhombuses and rectangles give the following

 \btm\label{TH-prayam-poluharakter}$\phantom{.}$
 \bit
\item[(a)] A rhombus $\varDelta$ in $\mathcal O^\star(G)$ is submultiplicative
if and only if its inner envelope $\varDelta^{\text{\LZ}}$ is an inverse
semicharacter on $G$.

\item[(b)] A rectangle $D$ in $\mathcal O(G)$ is dually submultiplicative if
and only if its outer envelope $D^\text{\SQ}$ is a semicharacter on $G$. \eit
\etm

The following result shows that the submultiplicative rhombuses form a
fundamental system among all submultiplicative closed absolutely convex
neighborhood of zero in $\mathcal O^\star(G)$:

\btm\label{subm-in-trubk} $\phantom{.}$
 \bit
\item[(a)] Every closed absolutely convex neighborhood of zero $\varDelta$ in
$\mathcal O^\star(G)$ contains some submultiplicative rhombus, namely
$\varDelta^{{\text{\LZ}}\kern-0.5pt{\text{\BLZ}}}$.

\item[(b)] Every dually submultiplicative absolutely convex compact set $D$ in
$\mathcal O(G)$ is contained in some dually submultiplicative rectangle, namely
in $D^{\text{\SQ}\text{\BSQ}}$.
 \eit
 \etm
 \bpr
1. If $\varDelta$ is a closed absolutely convex neighborhood of zero in
$\mathcal O^\star(G)$, then by Lemma \ref{f-obr-poluh->F-dual-subm}(b),
$\varDelta^{\text{\LZ}}$ is an inverse semicharacter, hence by Lemma
\ref{f-obr-poluh->F-dual-subm}(a), the rhombus
$\varDelta^{{\text{\LZ}}\kern-0.5pt{\text{\BLZ}}}$ is submultiplicative.

2. If $D$ is a dually submultiplicative absolutely convex compact set in
$\mathcal O(G)$, then its polar $D^\circ$ is a closed absolutely convex
submultiplicative neighborhood of zero in $\mathcal O^\star(G)$, and by what we
have already proved $(D^\circ)^{{\text{\LZ}}\kern-0.5pt{\text{\BLZ}}}$ is a
submultiplicative rhombus. I.e. the set
$(D^{\text{\SQ}\text{\BSQ}})^\circ=(D^\circ)^{{\text{\LZ}}\kern-0.5pt{\text{\BLZ}}}$
is submultiplicative. Therefore the rectangle $D^{\text{\SQ}\text{\BSQ}}$ is
dually submultiplicative.
 \epr

In accordance with the definitions of \ref{rectangles-buses}, we call a
function $f$ on $G$ an {\it enveloping semicharacter}\index{enveloping
semicharacter}, if it is an outer envelope and at the same time a semicharacter
on $G$.

\btm\label{f_N} If $G$ is a compactly generated Stein group, then the systems
of all semicharacters, all enveloping semicharacters and all dual
submultiplicative rectangles in $G$ contain countable cofinal subsystems:
 \bit
\item[--] there exists a sequence $h_N$ of semicharacters on $G$ such that
every semicharacter $g$ is majorized by some semicharacter $h_N$:
$$
g(x)\le h_N(x),\quad x\in G
$$

\item[--] there exists a sequence $f_N$ of enveloping semicharacters on $G$
such that every enveloping semicharacter $g$ is majorized by some semicharacter
$f_N$:
$$
g(x)\le f_N(x),\quad x\in G
$$

\item[--] there exists a sequence $E_N$ of dually submupltiplicative rectangles
in $G$ such that every dually submultiplicative rectangle $D$ in $G$ is
contained in some $E_N$:
 $$
 D\subseteq E_N
 $$
 \eit
\etm \bpr This follows from Proposition \ref{PROP-poluh-porozhd-K}: the
semicharacters $h_N$, $N\in\N$, defined in \eqref{poluh-opred-K}, are the
sequence we look for. The sequences $E_N$ and $f_N$ are defined as follows:
 \beq\label{f_N^black}
 E_N=\{u\in\mathcal O(G):\; \max_{x\in K^n}|u(x)|\le N^n\}=(h_N)^\text{\BSQ}
 \eeq
 \beq
f_N=(E_N)^\text{\SQ}=(h_N)^{\text{\BSQ}\text{\SQ}}
 \eeq
($E_N$ are dually submultiplicative rectangles by Lemma
\ref{f-poluh->F-subm}(a), and $f_N$ are enveloping semicharacters by Lemma
\ref{f-poluh->F-subm}(b)).

If now $D$  is a dually submultiplicative rectangle, then by Lemma
\ref{f-poluh->F-subm}(b), its outer envelope $D^\text{\SQ}$ is a semicharacter,
hence by Proposition \ref{PROP-poluh-porozhd-K}, there is $N\in\N$ such that
$$
D^\text{\SQ}\le h_N
$$
This means that $D$ is contained in some $E_N$:
$$
D=(D^\text{\SQ})^\text{\BSQ}\subseteq (h_N)^\text{\BSQ}=E_N
$$
This in its turn implies that
$$
D^\text{\SQ}\le (E_N)^\text{\SQ}=f_N,
$$
-- so every enveloping semicharacter (always being of the form $D^\text{\SQ}$
by Theorem \ref{TH-prayam-poluharakter}) is majorized by some $f_N$.
 \epr

\subsection{Holomorphic functions of exponential type}

\paragraph{Algebra $\mathcal O_{\exp} (G)$ of holomorphic functions of
exponential type.}\label{SUBSEC:algebra-O_exp(G)}

A holomorphic function $u\in\mathcal O(G)$ on a compactly generated Stein group
$G$ we call a {\it function of exponential type}\index{function of exponential
type}, if it is bounded by some semicharacter:
$$
|u(x)|\le f(x),\qquad x\in G\quad \Big(f(x\cdot y)\le f(x)\cdot f(y)\Big)
$$
The set of all holomorphic functions of exponential type on $G$ is denoted by
${\mathcal O}_{\exp}(G)$. It is a subspace in ${\mathcal O}(G)$ and, by Theorem
\ref{subm-in-trubk}, ${\mathcal O}_{\exp}(G)$ can be considered as the union of
all dually submultiplicative rectangles in ${\mathcal O}(G)$:
$$
{\mathcal O}_{\exp}(G)=\bigcup_{\scriptsize\begin{matrix}\text{$D$ is a
dually}\\ \text{submultiplicative}\\ \text{rectangle in $\mathcal
O(G)$}\end{matrix}}
 D=\bigcup_{\text{$f$ is a semicharacter in $G$}} f^\text{\BSQ}
$$
or, what is the same, the union of all subspaces of the form $\C D$, where $D$
is a dually submultiplicative rectangle in $\mathcal O(G)$:
 \beq\label{O-exp}
{\mathcal O}_{\exp}(G)=\bigcup_{\scriptsize\begin{matrix}\text{$D$ is a
dually}\\ \text{submultiplicative}\\ \text{rectangle in $\mathcal
O(G)$}\end{matrix}} \C D
 \eeq
This equality allows to endow ${\mathcal O}_{\exp}(G)$ with the natural
topology -- the topology of injective (locally convex) limit of Smith spaces
$\C D$:
 \beq\label{O-exp-top}
{\mathcal O}_{\exp}(G)=\underset{\scriptsize\begin{matrix}\text{$D$ is a
dually}\\ \text{submultiplicative}\\ \text{rectangle in $\mathcal
O(G)$}\end{matrix}}{\underset{\rightarrow}{\lim}} \C D
 \eeq
From Theorem \ref{f_N} it follows that in this limit the system of all dually
submultiplicative rectangles can be repalced by some countable subsystem:
 \beq\label{O-exp-top-E_N}
{\mathcal O}_{\exp}(G)=\underset{\scriptsize N\to
\infty}{\underset{\longrightarrow}{\lim}} \C E_N
 \eeq
Together with Theorem \ref{Brauner-K_n} this gives the following fact.

\btm\label{TH:O_exp--Brauner} The space ${\mathcal O}_{\exp}(G)$ of the
functions of exponential growth on a compactly generated Stein group $G$ is a
Brauner space. \etm

\bcor\label{O_exp(G):comp-subseteq-f^BSQ} If $G$ is a compactly generated Stein
group, then every bounded set $D$ in ${\mathcal O}_{\exp}(G)$ is contained in
some rectangle of the form $f^\text{\BSQ}$, where $f$ is some semicharacter on
$G$:
$$
D\subseteq f^\text{\BSQ}
$$
 \ecor
\bpr By Proposition \ref{Brauner-crit}, $D$ is contained in one of the compact
sets $E_N$. This set, being a dually submultiplicative rectangle, by Theorem
\ref{TH-prayam-poluharakter} has the form $f_N^\text{\BSQ}$ for some
semicharacter $f_N$, namely for $f_N=E_N^\text{\SQ}$. \epr

\btm  The space ${\mathcal O}_{\exp}(G)$ of the functions of exponential growth
on a compactly generated Stein group $G$ is a projective stereotype algebra
with respect to usual pointwise multiplication of functions. \etm

\bpr Note that if two functions $u,v\in{\mathcal O}(G)$ are bounded by
semicharacters $f$ and $g$, then their multiplication $u\cdot v$ is bounded by
the semicharacter $f\cdot g$:
$$
u\in f^\text{\BSQ},\ v\in g^\text{\BSQ}\quad\Longrightarrow\quad u\cdot v\in
(f\cdot g)^\text{\BSQ}
$$
In other words, the multiplication $(u,v)\mapsto u\cdot v$ in the space
${\mathcal O}(G)$ turns any compact of the form $f^\text{\BSQ}\times
g^\text{\BSQ}$ (where $f$ and $g$ are semicharacters) into the rectangle
$(f\cdot g)^\text{\BSQ}$.
$$
(u,v)\in f^\text{\BSQ}\times g^\text{\BSQ}\quad\mapsto\quad u\cdot v\in (f\cdot
g)^\text{\BSQ}
$$
Since this operation is continuous in ${\mathcal O}(G)$, it turns
$f^\text{\BSQ}\times g^\text{\BSQ}$ into $(f\cdot g)^\text{\BSQ}$ continuously.
On the other hand, $(f\cdot g)^\text{\BSQ}$, being a dually submultiplicative
rectangle, is continuously included into the space ${\mathcal O}_{\exp}(G)$.
Thus we obtain a continuous mapping
$$
(u,v)\in f^\text{\BSQ}\times g^\text{\BSQ}\quad\mapsto\quad u\cdot v\in
{\mathcal O}_{\exp}(G)
$$
If now $D$ and $E$ are arbitrary compact sets in the Brauner space ${\mathcal
O}_{\exp}(G)$, then by Corollary \ref{O_exp(G):comp-subseteq-f^BSQ}, they are
contained in some compact sets of the form $f^\text{\BSQ}$ and $g^\text{\BSQ}$,
where $f$ and $g$ are semicharacters. Hence,
$$
D\times E\subseteq f^\text{\BSQ}\times f^\text{\BSQ}
$$
From this we can conclude that the operation of multiplication is continuous
from $D\times E$ into ${\mathcal O}_{\exp}(G)$:
$$
(u,v)\in D\times E\subseteq f^\text{\BSQ}\times g^\text{\BSQ}\quad\mapsto\quad
u\cdot v\in {\mathcal O}_{\exp}(G)
$$
This is true for arbitrary compact sets $D$ and $E$ in ${\mathcal
O}_{\exp}(G)$. So we can apply \cite[Theorem 5.24]{Akbarov}: since ${\mathcal
O}_{\exp}(G)$, being a Brauner space, is co-complete (i.e. its stereotype dual
space is complete), by \cite[Theorem 2.4]{Akbarov} it is saturated. This
implies that any bilinear form on ${\mathcal O}_{\exp}(G)$, if continuous on
compact sets of the form $D\times E$, will be continuous ${\mathcal
O}_{\exp}(G)$ in the sense of definition \cite[\S\, 5.6]{Akbarov}. Being
applied to the operation $(u,v)\mapsto u\cdot v$ this means that the
multiplication is continuous in the sense of conditions of Proposition
\ref{PROP:proj-algebra}. Hence, ${\mathcal O}_{\exp}(G)$ is a projective
stereotype algebra. \epr

Let us note the following obvious fact:

\btm\label{TH-ogranichenie} If $H$ is a closed subgroup in a Stein group $G$,
then the restriction of any holomorphic function of exponential type on $G$ to
$H$ is a holomorphic function of exponential type:
 \beq
u\in {\mathcal O}_{\exp}(G)\quad\Longrightarrow\quad u|_H\in {\mathcal
O}_{\exp}(H)
 \eeq
 \etm

\paragraph{Algebra $\mathcal O_{\exp}^\star(G)$ of exponential analytic functionals.}

We shall call the elements of the dual stereotype space ${\mathcal
O}_{\exp}^\star(G)$, i.e. linear continuous functionals on ${\mathcal
O}_{\exp}(G)$  {\it exponential analytic functionals}\index{exponential
functional} on the group $G$. The space ${\mathcal O}_{\exp}^\star(G)$ is
called the {\it space of exponential analytic functionals} on $G$.

From Theorem \ref{TH:O_exp--Brauner} we have

\btm The space ${\mathcal O}_{\exp}(G)$ of exponential functionals on a
compactly generated Stein group $G$ is a Fr\'echet space. \etm

\btm\label{TH:O_exp^*(G)-algebra} The space ${\mathcal O}_{\exp}^\star(G)$ of
exponential functionals on a compactly generated Stein group $G$ is a
projective stereotype algebra with respect to usual convolution of functionals
defined by Formulas \eqref{DEF:u-cdot-a}-\eqref{DEF:svertka}:
$$
(\alpha,\beta)\quad\mapsto\quad \alpha*\beta
$$
\etm

\bpr 1. Take a function $u\in {\mathcal O}_{\exp}(G)$ and note that for any
point $s\in G$ its shift $u\cdot s$ (defined by \eqref{DEF:u-cdot-a}) is again
a function from ${\mathcal O}_{\exp}(G)$:
 $$
\forall s\in G\quad \forall u\in{\mathcal O}_{\exp}(G)\quad u\cdot s\in
{\mathcal O}_{\exp}(G)
 $$
Indeed, if we take a semicharacter $f$ majorizing $u$, we obtain:
 $$
\forall t\in G \quad |(u\cdot s)(t)|=|u(s\cdot t)|\le f(s\cdot t)\le f(s)\cdot
f(t)
 $$
i.e.,
 \beq\label{u-in-praym=>us-in-s-pryam}
u\in f^\text{\BSQ}\quad\Longrightarrow\quad u\cdot s\in f(s)\cdot f^\text{\BSQ}
 \eeq

2. Let us denote the mapping $s\mapsto u\cdot s$ by $\widehat{u}$,
$$
\widehat{u}:G\to {\mathcal O}_{\exp}(G)\quad\Big|\quad \widehat{u}(s):=u\cdot
s,\qquad s\in G
$$
and show that it is continuous. Let $s_i$ be a sequence of points in $G$,
tending to a point $s$:
$$
s_i\overset{G}{\underset{i\to\infty}{\longrightarrow}} s
$$
Then the sequence of holomorphic functions $\widehat{u}(s_i)\in {\mathcal
O}(G)$ tends to the holomorphic function $\widehat{u}(s)\in {\mathcal O}(G)$
uniformly on every compact set $K\subseteq G$, i.e. in the space ${\mathcal
O}(G)$:
$$
\widehat{u}(s_i)\overset{{\mathcal O}(G)}
{\underset{i\to\infty}{\longrightarrow}} \widehat{u}(s)
$$
At the same time, from \eqref{u-in-praym=>us-in-s-pryam} it follows that all
those functions are bounded by the semicharacter $C\cdot f$, where
$C=\max\{\sup_{i}f(s_i), f(s)\}$ (this value is finite since the sequence $s_i$
with its limit $s$ forms a compact set), so it lies in the rectangle generated
by the semicharacter $C\cdot f$:
$$
\widehat{u}(s_i)=u\cdot s_i\in C\cdot f^\text{\BSQ},\qquad
\widehat{u}(s)=u\cdot s\in C\cdot f^\text{\BSQ}
$$
In other words, $\widehat{u}(s_i)$ tends to $\widehat{u}(s)$ in the compact set
$C\cdot f^\text{\BSQ}$
$$
\widehat{u}(s_i)=u\cdot s_i\overset{C\cdot f^\text{\BSQ}}
{\underset{i\to\infty}{\longrightarrow}} u\cdot s=\widehat{u}(s)
$$
Hence, $\widehat{u}(s_i)$ tends to $\widehat{u}(s)$ in the space ${\mathcal
O}_{\exp}(G)$:
$$
\widehat{u}(s_i)\overset{{\mathcal O}_{\exp}(G)}
{\underset{i\to\infty}{\longrightarrow}} \widehat{u}(s)
$$

3. The continuity of the mapping $\widehat{u}:G\to{\mathcal O}_{\exp}(G)$
implies that for any functional $\beta\in {\mathcal O}_{\exp}^\star (G)$ the
function $\beta\circ \widehat{u}:G\to\C$ is holomorphic. We can use the Morera
theorem to prove this: consider a closed oriented hypersurface $\varGamma$ in
$G$ with a sufficiently small diameter so that integral over $\varGamma$ of
every holomorphic function vanishes, and show that the integral of $\beta\circ
\widehat{u}$ also vanishes:
 \beq\label{morera-1}
\int_{\varGamma} (\beta\circ \widehat{u})(s) \d s=0
 \eeq
Indeed, take a net of functionals $\{\beta_i; i\to\infty\}\subset {\mathcal
O}_{\exp}^\star (G)$ approximating $\beta$ in ${\mathcal O}_{\exp}^\star (G)$,
and having a form of linear combinations of delta-functionals:
$$
\beta_i=\sum_{k}\lambda^k_i\cdot\delta^{a^k_i},\qquad \beta_i\overset{{\mathcal
O}_{\exp}^\star (G)}{\underset{i\to\infty}{\longrightarrow}}\beta
$$
Then we obtain the following: since $\widehat{u}:G\to{\mathcal O}_{\exp}(G)$ is
continuous,
$$
\beta\circ \widehat{u}\overset{{\mathcal C}(G)} {\underset{\infty\gets
i}{\longleftarrow}}\beta_i\circ \widehat{u}
$$
From this it follows that for any Radon measure $\alpha\in {\mathcal C}(G)$
$$
\alpha(\beta\circ \widehat{u}) \underset{\infty\gets i}{\longleftarrow}
\alpha(\beta_i\circ \widehat{u})
$$
In particular for the functional of integration over the hypersurface
$\varGamma$ we have
 \begin{multline*}
\int_{\varGamma} (\beta\circ \widehat{u})(s) \d s \underset{\infty\gets
i}{\longleftarrow} \int_{\varGamma} (\beta_i\circ \widehat{u})(s)\d s=
\int_{\varGamma} \left(\sum_{k}\lambda^k_i\cdot\delta^{a^k_i}\circ
\widehat{u}\right)(s)\d s=\\= \sum_{k}\lambda^k_i\cdot\int_{\varGamma}
\left(\delta^{a^k_i}\circ \widehat{u}\right)(s)\d s=
\sum_{k}\lambda^k_i\cdot\int_{\varGamma} \delta^{a^k_i}(\widehat{u}(s))\d s=\\=
 \sum_{k}\lambda^k_i\cdot\int_{\varGamma}
\widehat{u}(s)(a^k_i)\d s= \sum_{k}\lambda^k_i\cdot\underbrace{\int_{\varGamma}
u(s\cdot a^k_i)\d s}_{\scriptsize\begin{matrix}\| \\ 0, \\ \text{since $u$ is
holomorphic} \end{matrix}}=0
 \end{multline*}
So indeed \eqref{morera-1} is true.

4. We have showed that for any functional $\beta\in {\mathcal O}_{\exp}^\star
(G)$ the function $\beta\circ \widehat{u}:G\to\C$ is holomorphic. Let us show
now that this function is of exponential type:
 \beq\label{beta-circ-widetilde-u-in-mathcal O-exp-G}
\forall u\in {\mathcal O}_{\exp}(G)\quad \forall \beta\in {\mathcal
O}_{\exp}^\star (G)\qquad \beta\circ \widehat{u}\in {\mathcal O}_{\exp}(G)
 \eeq
Indeed, since the functional $\beta\in {\mathcal O}_{\exp}^\star (G)$ is
bounded on the compact set $f^\text{\BSQ}\subseteq {\mathcal O}_{\exp}(G)$, it
is a bounded functional on the Banach representation of the Smith space $\C
f^\text{\BSQ}$, i.e. the following inequality holds:
 \beq
\forall v\in\C f^\text{\BSQ}\qquad |\beta(v)|\le M\cdot
\norm{v}_{f^\text{\BSQ}}
 \eeq
where
$$
M=\norm{\beta}_{(f^\text{\BSQ})^\circ}:=\max_{v\in
f^\text{\BSQ}}|\beta(v)|,\qquad \norm{v}_{f^\text{\BSQ}}:=\inf\{\lambda>0:\
v\in \lambda\cdot f^\text{\BSQ}\}
$$
Now from formula \eqref{u-in-praym=>us-in-s-pryam} we have
 \begin{multline*}
\forall s\in G\quad u\cdot s\in f(s)\cdot f^\text{\BSQ}
\quad\Longrightarrow\quad \norm{u\cdot s}_{f^\text{\BSQ}}:=\inf\{\lambda>0:\
u\cdot s\in \lambda\cdot f^\text{\BSQ}\}\le f(s) \quad\Longrightarrow\\
\Longrightarrow\quad |(\beta\circ \widehat{u})(s)|=|\beta(u\cdot s)|\le
M\cdot\norm{u\cdot s}_{f^\text{\BSQ}}\le M\cdot f(s)
 \end{multline*}
So the function $\beta\circ \widehat{u}$ is bounded by the semicharacter
$M\cdot f=\norm{\beta}_{(f^\text{\BSQ})^\circ}\cdot f$:
 \beq\label{beta-circ-widehat-u-in-norm-beta-cdot-f^BSQ}
u\in f^\text{\BSQ}\quad\Longrightarrow\quad \beta\circ \widehat{u}\in
\norm{\beta}_{(f^\text{\BSQ})^\circ}\cdot f^\text{\BSQ}
 \eeq

5. Now let us note that
$$
(\beta\circ\widehat{u})(s)=\beta(u\cdot
s)=\eqref{DEF:alpha-*-u}=(u*\widetilde{\beta})(s),\qquad s\in G
$$
In other words, we proved that the function
$u*\widetilde{\beta}=\beta\circ\widehat{u}$ belongs to the space ${\mathcal
O}_{\exp}(G)$, so for any functional $\alpha\in {\mathcal O}_{\exp}^\star(G)$
the convolution
$$
\alpha * \beta (u)=\eqref{DEF:svertka}=\alpha(u* \tilde{\beta})
$$
is defined. It remains to verify that the operation $(\alpha,\beta)\mapsto
\alpha*\beta$ is a continuous bilinear form, i.e. satisfies the conditions (i)
and (ii) of Proposition \ref{PROP:proj-algebra}.

Let $\alpha_i$ be a net tending to zero in ${\mathcal O}_{\exp}(G)$, and $B$ a
compact set in ${\mathcal O}_{\exp}^\star(G)$. Let us consider an arbitrary
compact set $D$ in ${\mathcal O}_{\exp}(G)$. By Corollary
\ref{O_exp(G):comp-subseteq-f^BSQ}, it must be contained in a rectangle
$f^\text{\BSQ}$, where $f$ is a semicharacter:
$$
D\subseteq f^\text{\BSQ}
$$
On the other hand, on the Smith space $\C f^\text{\BSQ}$ the norms of
functionals $\beta\in B$ are bounded:
 $$
\sup_{\beta\in B}\norm{v}_{f^\text{\BSQ}}=M<\infty
 $$
So by virtue of \eqref{beta-circ-widehat-u-in-norm-beta-cdot-f^BSQ} we obtain
that
$$
\forall u\in D\quad \forall \beta\in B\qquad u*\widetilde{\beta}=\beta\circ
\widehat{u}\in \norm{\beta}_{(f^\text{\BSQ})^\circ}\cdot f^\text{\BSQ}\subseteq
M\cdot f^\text{\BSQ}
$$
Thus the set $\{u*\widetilde{\beta};\ u\in D,\ \beta\in B\}$ is contained in
the compact set $M\cdot f^\text{\BSQ}$ in the space ${\mathcal
O}_{\exp}^\star(G)$. Hence, the net of functionals $\alpha_i$ tends to zero
uniformly on this set:
 $$
(\alpha_i*\beta)(u)=\alpha_i(u*\widetilde{\beta})\overset{\C}{\underset{i\to\infty}{\rightrightarrows}}
0\qquad \Big(u\in D,\ \beta\in B\Big)
 $$
This is true for any compact set $D$, so the net $\alpha_i*\beta$ tends to zero
in the space ${\mathcal O}_{\exp}^\star(G)$ uniformly by $\beta\in B$:
 $$
\alpha_i*\beta\overset{{\mathcal
O}_{\exp}^\star(G)}{\underset{i\to\infty}{\rightrightarrows}} 0\qquad
\Big(\beta\in B\Big)
 $$
We can omit the case of the inverse sequence of multipliers due to the identity
$$
\widetilde{\alpha*\beta}=\widetilde{\beta}*\widetilde{\alpha}
$$
 \epr

\subsection{Examples.}

\paragraph{Finite groups.} \label{EX-finite-groups} As we have noted in
\ref{SEC:stein-groups}\ref{SUBSEC-lin-groups}, every finite group $G$,
considered as a zero-dimensional complex manifold, is a linear complex Lie
group on which every function is holomorphic. On the other hand, every function
$u:G\to\C$ is bounded (since its set of values is finite), so $u$ must be a
function of exponential type. Thus, algebras ${\mathcal O}_{\exp}(G)$ and
${\mathcal O}(G)$ in this case coincide with each other and are equal to the
algebra $\C^G$ of all functions on $G$ (with the pointwise algebraic operations
and the topology of pointwise convergence):
$$
{\mathcal O}_{\exp}(G)={\mathcal O}(G)=\C^G
$$

\paragraph{Groups $\C^n$.} For the case of complex plane $\C$ our definition,
certainly, coincide with the classical one -- functions of exponential type on
$\C$ are entire functions $u\in{\mathcal O}(\C)$ growing not faster than the
exponential:
 $$
\exists A>0:\quad |u(x)|\le A^{|x|}\qquad (x\in\C).
 $$
According to the classical theorems on the growth of entire functions
\cite{Levin,Shabat-1}, this is equivalent to the condition that the derivatives
of $u$ in a fixed point, say in zero, grow not faster than the exponential:
 $$
\exists B>0:\quad |u^{(k)}(0)|\le B^k\qquad (k\in\N).
 $$

The same is true for several variables: functions of exponential type on $\C^n$
according to our definition will be exactly the functions $u\in{\mathcal
O}(\C^n)$ satisfying the condition
 \beq\label{exp-type-for-C}
\exists A>0:\quad |u(x)|\le A^{|x|}\qquad (x\in\C^n),
 \eeq
which turns out to be equivalent to the condition
 \beq\label{exp-type-for-C-u'}
\exists B>0:\quad |u^{(k)}(0)|\le B^{|k|}\qquad (k\in\N^n),
 \eeq
where the factorial $k!$ and the modulus $|k|$ of a multiindex
$k=(k_1,...,k_n)\in\N^n$ are defined by the equalities
$$
k!:=k_1\cdot ...\cdot k_n,\qquad |k|:=k_1+...+k_n
$$
The equivalence of \eqref{exp-type-for-C} to our definition follows from
Proposition \ref{PROP-poluh-porozhd-K}, and the equivalence between
\eqref{exp-type-for-C} and \eqref{exp-type-for-C-u'} is proved similarly with
the case of one variable: the implication \eqref{exp-type-for-C-u'}
$\Rightarrow$ \eqref{exp-type-for-C} is obvious, and the inverse implication
\eqref{exp-type-for-C} $\Rightarrow$ \eqref{exp-type-for-C-u'} is a corollary
of the Cauchy inequalities (see \cite{Shabat}) for the coefficients
$c_k=\frac{u^{(k)}(0)}{k!}$ of the Taylor series for the function $u$:
$$
\forall R>0\qquad |c_k|\le \frac{\max\limits_{|x|\le
R}|u(x)|}{R^{|k|}}\le\eqref{exp-type-for-C}\le \frac{A^R}{R^{|k|}}
$$
$$
\Downarrow
$$
$$
|c_k|\le
\min_{R>0}\frac{A^R}{R^{|k|}}=\frac{A^R}{R^{|k|}}\Bigg|_{R=\frac{|k|}{\log
A}}=\frac{A^{\frac{|k|}{\log A}}}{|k|^{|k|}}\cdot (\log
A)^{|k|}=\frac{1}{|k|^{|k|}}\cdot B^{|k|},\qquad B=A^{\frac{1}{\log A}}\cdot
\log A
$$
$$
\Downarrow
$$
$$
|u^{(k)}(0)|=k!\cdot |c_k|\le \frac{k!}{|k|^{|k|}}\cdot B^{|k|}\le B^{|k|}
$$

\paragraph{Groups $\GL_n(\C)$.}\label{EX-mnogochleny-na-GL=f-exp-rosta}

On the group $\GL_n(\C)$ the functions of exponential type are exactly
polynomials, i.e. the functions of the form
 \beq\label{mnogochleny-na-GL}
u(x)=\frac{P(x)}{(\det x)^m},\qquad x\in \GL_n(\C),\qquad m\in\Z_+
 \eeq
where $P$ is a usual polynomial on matrix elements of the matrix $x$, and $\det
x$ is its determinant. Thus, the equality holds:
 \beq\label{O(GL)=H_exp(GL)}
{\mathcal O}_{\exp}(\GL_n(\C))={\mathcal R}(\GL_n(\C))
 \eeq
\bpr 1. First, let us prove the inclusion ${\mathcal R}(\GL_n(\C))\subseteq
{\mathcal O}_{\exp}(\GL_n(\C))$. Note for this that every matrix element
$$
x\mapsto x_{k,l}
$$
is a function of exponential type, since it is bounded by, for instance the
first of the matrix norms in \eqref{matrichnye-normy}:
$$
|x_{k,l}|\le \sum_{i,j=1}^n|x_{i,j}|=||x||
$$
As a corollary, every polynomial $x\mapsto P(x)$ of matrix elements is also a
function of exponential type (since functions of exponential type form an
algebra).

Further, every power of the determinant, $x\mapsto (\det x)^{-m}$, is
multiplicative
$$
(\det (x\cdot y))^{-m}=(\det x)^{-m}\cdot (\det y)^{-m}
$$
so its absolute value is also multiplicative:
$$
\left|(\det (x\cdot y))^{-m}\right|=\left|(\det x)^{-m}\right| \cdot
\left|(\det y)^{-m}\right|
$$
Therefore, the function $x\mapsto \left|(\det x)^{-m}\right|$ is a
semicharacter on $\GL_n(\C)$, and $x\mapsto (\det x)^{-m}$ a function of
exponential type on $\GL_n(\C)$.

Now, if we multiply two functions of exponential type $x\mapsto P(x)$ and
$x\mapsto (\det x)^{-m}$, we obtain a function of exponential type $x\mapsto
P(x)/ (\det x)^m$.

2. Let us prove the inverse inclusion: ${\mathcal O}_{\exp}(\GL_n(\C))\subseteq
{\mathcal R}(\GL_n(\C))$. If $u$ is a holomorphic function of exponential type
on $\GL_n(\C)$, then by Proposition \ref{PROP-stroenie-poluharakterov-na-GL_n},
it must be bounded by some semicharacter of the form
\eqref{r_N(x)=||x||^N-cdot-||x^-1||^N}. In particular, for some $C\ge 1$ and
$N\in\N$
 \beq\label{|u(t)|-le-max-||x||^N-||x^-1||^N}
|u(x)| \le C\cdot\max\{||x||; ||x^{-1}||\}^N, \qquad
||x||=\sum_{i,j=1}^n|x_{i,j}|
 \eeq
The elements of the inverse matrix $x^{-1}$ are obtained from $x$ as
complementary minors, divided by the determinant, so we can treat them as
polynomials (of degree $n-1$) of $x_{i,j}$ and $(\det x)^{-1}$. This means that
one can estimate the righthand side of \eqref{|u(t)|-le-max-||x||^N-||x^-1||^N}
by a polynomial (of degree $N(n-1)$) of $|x_{i,j}|$ and $|\det x|^{-1}$ with
nonnegative coefficients:
$$
|u(x)| \le C\cdot\max\{||x||; ||x^{-1}||\}^N\le C\cdot
P\Big(\{|x_{i,j}|\}_{1\le i,j\le n},\; |\det x|^{-1}\Big)
$$
Hence if we multiply $u$ by $(\det x)^{N(n-1)}$, we obtain a holomorphic
function on $\GL_n(\C)$ bounded by a polynomial of $|x_{i,j}|$:
$$
\left|u(x)\cdot (\det x)^{N(n-1)}\right| \le C\cdot Q\Big(\{|x_{i,j}|\}_{1\le
i,j\le n}\Big)
$$
Such a function is locally bounded in the points of the analytical set $\det
x=0$, so by the Riemann extension theorem (see \cite{Shabat}), it can be
extended to a holomorphic function to the manifold $\M_n(\C)$ (of all complex
matrices). Therefore, we can consider $u(x)\cdot (\det x)^{N(n-1)}$ as a
holomorphic function on $\M_n(\C)=\C^{n^2}$. Since it has a polynomial growth,
it must be a polynomial $q$ of matrix elements $x_{i,j}$:
$$
u(x)\cdot (\det x)^{N(n-1)}=q(x)
$$
Hence
$$
u(x)=\frac{q(x)}{(\det x)^{N(n-1)}},
$$
i.e. $u\in {\mathcal R}(\GL_n(\C))$.
 \epr

\bcor On the complex circle $\C^\times$ the functions of exponential type are
exactly the Laurent polynomials:
$$
u(t)=\sum_{|n|\le N} u_n\cdot t^n,\qquad N\in\N
$$
 \ecor

\subsection{Injection $\flat_G:{\mathcal O}_{\exp}(G)\to {\mathcal
O}(G)$}

We need to denote the injection of ${\mathcal O}_{\exp}(G)$ into ${\mathcal
O}(G)$ by some symbol. Let us use for this the symbol $\flat_G$:
 \beq\label{H_exp->H}
\flat_G:{\mathcal O}_{\exp}(G)\to {\mathcal O}(G)
 \eeq
This mapping is always injective, a homomorphism of algebras, and by definition
of topology in ${\mathcal O}_{\exp}(G)$, always continuous. Below in Theorem
\ref{TH-H_exp-plotno-v-H} we shall see that if $G$ is a linear group, then
$\flat_G$ has dense image in ${\mathcal O}(G)$.

From equality \eqref{O(GL)=H_exp(GL)} and Theorem \ref{TH-ogranichenie} it
follows that if $G$ is an arbitrary linear group (with a given representation
as a closed subgroup in $\GL_n(\C)$), then every function $G$ that can be
extended to a polynomial on $\GL_n(\C)$, is a function of exponential type.
Thus, we have the following chain of inclusions:
 \beq\label{O-subset-H_exp-subset-H}
{\mathcal R}(G)\subseteq {\mathcal O}_{\exp}(G)\subseteq {\mathcal O}(G)
 \eeq
(here ${\mathcal R}(G)$ denotes the functions that can be extended to
polynomials on $\GL_n(\C)$ -- we need this specification since $G$ is not
necessarily an algebraic group).

\btm\label{TH-H_exp-plotno-v-H} If $G$ is a linear complex group, then the
algebra ${\mathcal O}_{\exp}(G)$ of holomorphic functions of exponential type
on $G$ is dense in the algebra ${\mathcal O}(G)$ of all holomorphic functions
on $G$. \etm
 \bpr
Let $\ph:G\to\GL_n(\C)$ be a holomorphic embedding as a closed subgroup. By one
of the corollaries from the Cartan theorem \cite[11.5.2]{Taylor}, every
holomorphic function $v\in {\mathcal O}(G)$ can be extended to a holomorphic
function $u\in {\mathcal O}(\GL_n(\C))$. Let us approximate $u$ uniformly on
compact set by polynomials $u_i\in {\mathcal R}(\GL_n(\C))$. By
\eqref{O(GL)=H_exp(GL)}, all polynomials $u_i$ are functions of exponential
type on $\GL_n(\C)$, hence their restrictions $u_i|_G$ are functions of
exponential type on $G$. Thus, $v$ is approximated by functions of exponential
type $u_i|_G$ uniformly on compact sets in $G$.
 \epr

\subsection{Nuclearity of the spaces ${\mathcal O}_{\exp}(G)$ and ${\mathcal O}^\star_{\exp}(G)$}

\btm\label{TH-nuclear} For any compactly generated Stein group $G$ the space
${\mathcal O}_{\exp}(G)$ is a nuclear Brauner space, and its dual space
${\mathcal O}^\star_{\exp}(G)$ a nuclear Fr\'echet space. \etm

We premise the proof of this fact by two lemmas. The first of them is true for
arbitrary complex manifold and is proved by the same technique that is applied
for the proof of nuclearity of  ${\mathcal O}(\C)$ (see \cite[6.4.2]{Pietsch}):

\blm\label{LM-ver-mera} If $M$ is a complex manifold, and $K$ and $L$ are two
compact sets in $M$, such that $K$ is strictly contained in $L$,
$$
K\subseteq \Int L
$$
then there exists a constant $C\ge 0$ and a probability measure $\mu$ on $L$
such that for any $u\in{\mathcal O}(G)$ we have
 \beq\label{ver-mera}
|u|_K\le C\cdot\int_L |\alpha(u)|\ \mu(\d \alpha)
 \eeq
As a corollary, the operator $u|_L\mapsto u|_K$ of restriction is absolutely
summing, and its quasinorm of absolute summing is not greater than $C$: for any
$u_1,...,u_l\in{\mathcal O}(M)$
 $$
\sum_{i=1}^l |u_i|_K\le C\cdot \int_L \sum_{i=1}^l |\alpha(u_i)|\; \mu(\d
\alpha)\le C\cdot \sup_{\alpha\in\cabsconv\left(\delta^L\right)} \sum_{i=1}^l
|\alpha(u_i)|
 $$
 \elm

\brem Here $\cabsconv\left(\delta^L\right)$ means the universal compact set in
the Smith space $\mathcal C^\star(L)$ dual to the Banach space $\mathcal C(L)$
of continuous functions on $L$, or what is the same, the unit ball in the
Banach space $\mathcal M(L)$ of Radon measures on $L$. We use this notation,
because it is convenient to denote by $\delta^L=\{\delta^x;\ x\in L\}$ the set
of all delta-functionals on $\mathcal C(L)$ -- then the polar $B^\circ$ of the
unit ball $B$ in $\mathcal C(L)$ coincides with the absolutely convex hull of
$\delta^L$:
$$
B^\circ=\cabsconv\left(\delta^L\right)=\{\alpha\in M(L):\ ||\alpha||\le 1\}
$$
(the closure with respect to topology of $\mathcal C^\star(L)$). \erem

\blm For any generating compact neighborhood of identity $K$ in $G$
$$
G=\bigcup_{n=1}^\infty K^n
$$
there are constants $C\ge 0$, $\lambda\ge 0$ such that for any $l\in\N$ and for
arbitrary $u_1,...,u_l\in{\mathcal O}(G)$, $n\in\N$ the following inequality
holds:
 \beq\label{quasinor-abs-sum}
\sum_{i=1}^l |u_i|_{K^n}\le
C\cdot\lambda^{n-1}\cdot\sup_{\alpha\in\cabsconv\left(\delta^{K^{2n+1}}\right)}
\sum_{i=1}^l |\alpha(u_i)|
 \eeq
 \elm
\bpr 1. The set $U=\Int K$ is an open neighborhood of identity in $G$, so the
system of shifts $\{x\cdot U;\; x\in K^2\}$ is an open covering of the compact
set $K^2$. Let us choose a finite subcovering, i.e. a finite set $F\subseteq
K^2$ such that
$$
K^2\subseteq \bigcup_{x\in F} x\cdot U= F\cdot U
$$
Then we obtain:
 \beq\label{K^n-subseteq-F^n-1-cdot K}
K^n\subseteq F^{n-1}\cdot K\subseteq F^{n-1}\cdot K^2\subseteq K^{2n+1},\qquad
n\in\N
 \eeq
Here is the proof: first we should note that
$$
F\subseteq K^2\quad \Longrightarrow\quad F^n\subseteq K^{2n},\qquad n\in\N
$$
After that in the sequence \eqref{K^n-subseteq-F^n-1-cdot K} it is sufficient
to check only the first inclusion:
 \beq\label{K^n-subseteq-F^n-1-cdot K-1}
K^n\subseteq F^{n-1}\cdot K,\qquad n\in\N
 \eeq
This is done by the induction: for $n=2$ we have
$$
K^2\subseteq F\cdot U\subseteq F\cdot K
$$
and, if \eqref{K^n-subseteq-F^n-1-cdot K-1} is true for some $n$, then for
$n+1$ we have:
$$
K^{n+1}=K^n\cdot K\subseteq F^{n-1}\cdot K\cdot K\subseteq F^{n-1}\cdot
K^2\subseteq F^{n-1}\cdot F\cdot K=F^n\cdot K
$$

2. Since the compact set $K^2$ strictly contains $K$, by Lemma
\ref{LM-ver-mera} there are $C$, $\mu$ such that
 \beq\label{ver-mera-1}
|u|_K\le C\cdot\int_{K^2} |\alpha(u)|\mu(\d \alpha)
 \eeq
Under the shift by an element $x\in G$ this inequality takes form
 \beq\label{ver-mera-2}
|u|_{x\cdot K}\le C\cdot\int_{x\cdot K^2} |\alpha(u)|(x\cdot \mu)(\d \alpha)
 \eeq
where $x\cdot \mu$ is a shift of the measure $\mu$:
 $$
(x\cdot \mu)(u)=\mu(u\cdot x),\qquad (u\cdot x)(t)=u(x\cdot t)
 $$
This implies that if $E$ is an arbitrary finite set in $G$, then the sum
 $$
\nu=\frac{1}{\card E}\sum_{x\in E}x\cdot \mu
 $$
is a probability measure on the set $E\cdot K^2$, with the property
 \begin{multline*}
|u|_{E\cdot K}\le \sum_{x\in E} |u|_{x\cdot K}\le C\cdot\sum_{x\in E}
\int_{x\cdot K^2} |\alpha(u)|(x\cdot \mu)(\d \alpha)=\\= C\cdot\sum_{x\in E}
\int_{E\cdot K^2} |\alpha(u)|(x\cdot \mu)(\d \alpha)= C\cdot \int_{E\cdot K^2}
|\alpha(u)| \left(\sum_{x\in E} x\cdot \mu\right)(\d \alpha)=\\=
C\cdot\card(E)\cdot \int_{E\cdot K^2} |\alpha(u)|\; \nu(\d \alpha)
 \end{multline*}
From this we deduce that the restriction operator $u|_{E\cdot K^2}\mapsto
u|_{E\cdot K}$ is absolutely summing, and the quasinorm of absolute summing can
be estimated by the constant $C\cdot\card(E)$: for all $u_1,...,u_l\in{\mathcal
O}(G)$
 $$
\sum_{i=1}^l |u_i|_{E\cdot K}\le C\cdot\card(E)\cdot \int_{E\cdot K^2}
\sum_{i=1}^l |\alpha(u_i)|\; \nu(\d \alpha)\le C\cdot\card(E)\cdot
\sup_{\alpha\in\cabsconv\left(\delta^{E\cdot K^2}\right)} \sum_{i=1}^l
|\alpha(u_i)|
 $$
Now from formulas \eqref{K^n-subseteq-F^n-1-cdot K} we have:
 \begin{multline*}
\sum_{i=1}^l |u_i|_{K^n}\le \sum_{i=1}^l |u_i|_{F^{n-1}\cdot K}\le
C\cdot\card(F^{n-1})\cdot \sup_{\alpha\in\cabsconv\left(\delta^{F^{n-1}\cdot
K^2}\right)} \sum_{i=1}^l |\alpha(u_i)|\le\\ \le C\cdot\card(F^{n-1})\cdot
\sup_{\alpha\in\cabsconv\left(\delta^{K^{2n+1}}\right)} \sum_{i=1}^l
|\alpha(u_i)|\le C\cdot(\card(F))^{n-1}\cdot
\sup_{\alpha\in\cabsconv\left(\delta^{K^{2n+1}}\right)} \sum_{i=1}^l
|\alpha(u_i)|
 \end{multline*}
To obtain \eqref{quasinor-abs-sum} we can put $\lambda=\card F$.
 \epr

\bpr[Proof of Theorem \ref{TH-nuclear}] Consider the sequence of rectangles as
in Theorem \ref{f_N}:
$$
E_N=(f_N)^\text{\BSQ}
$$
By Theorem \ref{f_N}, their union covers the whole space ${\mathcal
O}_{\exp}(G)$
$$
{\mathcal O}_{\exp}(G)=\bigcup_{N=1}^\infty E_N
$$
and moreover any dually submultiplicative rectangle $D\subseteq {\mathcal
O}(G)$ is contained in some rectangle $E_N$:
$$
D\subseteq E_N
$$
From this we deduce that in Formulas \eqref{O-exp} and \eqref{O-exp-top} the
system of all dually submultiplicative rectangles $D$ can be replaced by the
system of rectangles $ E_N$:
$$ {\mathcal O}_{\exp}(G)=\bigcup_{N=1}^\infty \C
E_N=\underset{N\to\infty}{\underset{\rightarrow}{\lim}} \C E_N
$$
This means that ${\mathcal O}_{\exp}(G)$ is a Brauner space, and ${\mathcal
O}^\star_{\exp}(G)$ a Fr\'echet space. To prove that both spaces are nuclear we
can only note that ${\mathcal O}_{\exp}(G)$ is conuclear. Consider $\C E_N$ as
Banach spaces with unit balls $E_N$. By \eqref{f_N^black}, the seminorm on $\C
E_N$ is defined by equality
 $$
p_N(u)=\sup_{n\in\N}\frac{1}{N^n}|u|_{K^n}
 $$
Its unit ball is the set $ E_N$:
 \begin{multline*}
 E_N=\{u\in{\mathcal O}(G):\;\forall n\in\N\quad |u|_{K^n}\le
N^n\}=\left\{u\in{\mathcal O}(G):\;\forall n\in\N\quad
\frac{1}{N^n}|u|_{K^n}\le 1\right\}=\\=\left\{u\in{\mathcal
O}(G):\;p_N(u)=\sup_{n\in\N}\frac{1}{N^n}|u|_{K^n}\le 1\right\}
 \end{multline*}
To prove that the space ${\mathcal
O}_{\exp}(G)=\underset{N\to\infty}{\underset{\rightarrow}{\lim}} \C E_N$ is
conuclear, we need to verify that for any $N\in\N$ there exists $M\in\N$,
$M>N$, such that the inclusion mapping $\C E_N\to\C E_M$ is absolutely summing,
i.e. such that for some constant $L>0$, for any $l\in\N$ and for all
$u_1,...,u_l\in \C E_N$ we have:
 \beq\label{psi_N-psi_M-abs-summ}
\sum_{i=1}^l p_M(u_i)\le L\cdot\sup_{\alpha\in(
E_N)^\circ}\sum_{i=1}^l|\alpha(u_i)|
 \eeq
This is proved as follows. First we need to note that for all $n,N\in\N$ the
following inclusion holds
 \beq\label{varPsi_N-in-N^n-B_n}
 E_N=\{u\in{\mathcal O}(G):\;\forall n\in\N\quad |u|_{K^n}\le
N^n\}\subseteq \{u\in{\mathcal O}(G):\; |u|_{K^n}\le N^n\}=N^n\cdot
 {^\circ{\delta^{K^n}}}
 \eeq
It implies following chain:
 $$ E_N\subseteq N^n\cdot
 {^\circ\kern-2pt\left(\delta^{K^n}\right)}
 $$
 $$
 \Downarrow
 $$
 $$
( E_N)^\circ\supseteq (N^n\cdot
 {^\circ\kern-2pt\left(\delta^{K^n}\right)})^\circ=\frac{1}{N^n}\cdot
( {^\circ\kern-2pt\left(\delta^{K^n}\right)})^\circ=\frac{1}{N^n}\cdot
\cabsconv (\delta^{K^n})
 $$
 $$
 \Downarrow
 $$
 $$
 \cabsconv (\delta^{K^n})\subseteq N^n\cdot( E_N)^\circ
 $$
 $$
 \Downarrow
 $$
 $$
 \cabsconv (\delta^{K^{2n+1}})\subseteq N^{2n+1}\cdot( E_N)^\circ
 $$
 $$
 \Downarrow
 $$
 \begin{multline*}
\sum_{i=1}^l |u_i|_{K^n}\le\eqref{quasinor-abs-sum}\le C\cdot\lambda^{n-1}\cdot
\sup_{\alpha\in\cabsconv\left(\delta^{K^{2n+1}}\right)}\sum_{i=1}^l
|\alpha(u_i)|\le C\cdot\lambda^{n-1}\cdot\sup_{\alpha\in N^{2n+1}\cdot(
E_N)^\circ}\sum_{i=1}^l |\alpha(u_i)|=\\=C\cdot\lambda^{n-1}\cdot\sup_{\beta\in
( E_N)^\circ} \sum_{i=1}^l |N^{2n+1}\cdot\beta(u_i)|= C\cdot\lambda^{n-1}\cdot
N^{2n+1}\cdot\sup_{\beta\in ( E_N)^\circ}\sum_{i=1}^l |\beta(u_i)|
 \end{multline*}
 $$
 \Downarrow
 $$
 $$
 \forall M>0\qquad
\sum_{i=1}^l\frac{1}{M^n}|u_i|_{K^n}\le C\cdot\frac{\lambda^{n-1}\cdot
N^{2n+1}}{M^n}\cdot\sup_{\beta\in ( E_N)^\circ}\sum_{i=1}^l |\beta(u_i)|
 $$
 $$
 \Downarrow
 $$
 \begin{multline*}
\sum_{i=1}^l p_M(u_i)=\sum_{i=1}^l\sup_{n\in\N}\frac{1}{M^n}|u_i|_{K^n}\le
\sum_{i=1}^l\left(\sum_{n=1}^\infty\frac{1}{M^n}|u_i|_{K^n}\right)=
\sum_{n=1}^\infty\left(\sum_{i=1}^l\frac{1}{M^n}|u_i|_{K^n}\right)\le\\
\le \sum_{n=1}^\infty\left(C\cdot\frac{\lambda^{n-1}\cdot
N^{2n+1}}{M^n}\cdot\sup_{\beta\in ( E_N)^\circ}\sum_{i=1}^l |\beta(u_i)|\right)
= C\cdot\left(\sum_{n=1}^\infty\frac{\lambda^{n-1}\cdot N^{2n+1}}{M^n}\right)
\cdot\sup_{\beta\in ( E_N)^\circ}\sum_{i=1}^l |\beta(u_i)|
 \end{multline*}
If we now choose $M\in\N$ sufficiently large such that
 $$
C\cdot \sum_{n=1}^\infty\frac{\lambda^{n-1}\cdot N^{2n+1}}{M^n} \le 1
 $$
then the constant $L$ in \eqref{psi_N-psi_M-abs-summ} can be chosen as $1$:
 $$
 \sum_{i=1}^l p_M(u_i)\le
\underbrace{C\cdot \sum_{n=1}^\infty\frac{\lambda^{n-1}\cdot
N^{2n+1}}{M^n}}_{\scriptsize \begin{matrix} \text{\rotatebox{-90}{$\le$}} \\ 1
\end{matrix}}\cdot\sup_{\beta\in ( E_N)^\circ}\sum_{i=1}^l |\beta(u_i)|
 $$
\epr

\subsection{Holomorphic mappings of exponential type and tensor products
of the spaces ${\mathcal O}_{\exp}(G)$ and ${\mathcal O}^\star_{\exp}(G)$}

\btm\label{mathcal-O-exp-G-times-H-cong-mathcal-O-exp-G-mathcal-O-exp-H} Let
$G$ and $H$ be two compactly generated Stein groups. The formula
 \beq\label{rho_G,H}
\rho_{G,H}(u\boxdot v)=u\odot v,\qquad u\in{\mathcal O}_{\exp}(G),\quad
v\in{\mathcal O}_{\exp}(H)
 \eeq
(where $u\boxdot v$ is  the function from \eqref{g-h-na-S-T}) defines a linear
continuous mapping
$$
\rho_{G,H}:{\mathcal O}_{\exp}(G\times H)\to {\mathcal O}_{\exp}(G)\oslash
{\mathcal O}_{\exp}^\star(H)={\mathcal O}_{\exp}(G)\odot {\mathcal
O}_{\exp}(H),
$$
This mapping is an isomorphism of stereotype spaces and is natural by $G$ and
$H$, i.e. is an isomorphism of bifunctors from the category $\mathfrak{SG}$ of
Stein groups into the category $\mathfrak{Ste}$ of stereotype spaces:
$$
\Big((G;H)\mapsto {\mathcal O}(G\times H)\Big)\rightarrowtail \Big((G;H)\mapsto
{\mathcal O}(G)\odot {\mathcal O}(H)\Big)
$$
Equivalently this mapping is defined by formula
 \beq\label{widetilde-w-s-t-star}
\rho_{G,H}(w)(\beta)=\beta\circ\widehat{w},\qquad w\in {\mathcal
O}_{\exp}(G\times H),\quad \beta\in {\mathcal O}_{\exp}^\star (H)
 \eeq
where
 \beq\label{widetilde-w-s-t=w-s-t}
\widehat{w}:G\to {\mathcal O}_{\exp}(H)\quad\Big|\quad
\widehat{w}(s)(t)=w(s,t),\qquad s\in G,\; t\in H
 \eeq
\etm

 \bcor
The following isomorphisms of functors hold:
 \beq\label{tenz-pr-O-exp}
{\mathcal O}_{\exp}(G\times H)\cong {\mathcal O}_{\exp}(G)\odot {\mathcal
O}_{\exp}(H)\cong {\mathcal O}_{\exp}(G)\circledast {\mathcal O}_{\exp}(H)
 \eeq
 \beq\label{tenz-pr-O-exp-star}
{\mathcal O}_{\exp}^\star(G\times H)\cong {\mathcal O}_{\exp}^\star(G)\odot
{\mathcal O}_{\exp}^\star(H)\cong {\mathcal O}_{\exp}^\star(G)\circledast
{\mathcal O}_{\exp}^\star(H)
 \eeq
 \ecor

To prove Theorem
\ref{mathcal-O-exp-G-times-H-cong-mathcal-O-exp-G-mathcal-O-exp-H} we have to
recall the notion of injective tensor product $A\odot B$ of sets $A$ and $B$ in
stereotype spaces $X$ and $Y$. According to notations \cite[(7.27)]{Akbarov},
$A\odot B$ is defined as a subset in the space $X\odot Y=X\oslash Y^\star$ of
operators $\varphi:Y^\star\to X$ containing only those operators satisfying the
condition $\varphi(B^\circ)\subseteq A$. That is the sense of the following
formula we use in Lemma
\ref{(f-otimes-g)^blacksquare=f^blacksquare-odot-g^blacksquare} below:
$$
A\odot B=A\oslash B^\circ.
$$

\blm\label{(f-otimes-g)^blacksquare=f^blacksquare-odot-g^blacksquare} If
$g:G\to\R_+$ and $h:H\to\R_+$ are two semicharacters on $G$ and $H$, then
$g\boxdot h$ is a semicharacter on $G\times H$, and the mapping $\rho_{G,H}$
defined in \eqref{widetilde-w-s-t-star}-\eqref{widetilde-w-s-t=w-s-t} is a
homeomorphism between compact sets $(g\boxdot h)^\text{\BSQ}\subseteq {\mathcal
O}_{\exp}(G\times H)$ and $g^\text{\BSQ}\odot h^\text{\BSQ}\subseteq {\mathcal
O}_{\exp}(G)\odot {\mathcal O}_{\exp}(H)$:
 \beq
 (g\boxdot h)^\text{\BSQ}\cong g^\text{\BSQ}\odot h^\text{\BSQ}
 \eeq
 \elm
\bpr Here at the beginning we use reasonings similar to those we used in the
proof of Theorem \ref{TH:O_exp^*(G)-algebra}.

1. Note first that for any function $w\in(g\boxdot h)^\text{\BSQ}\subseteq
{\mathcal O}_{\exp}(G\times H)$ Formula \eqref{widetilde-w-s-t=w-s-t} defines
some mapping
 $$
\widehat{w}:G\to{\mathcal O}_{\exp}(H)
 $$
Indeed, since $w$ is holomorphic on $G\times H$, it is holomorphic with respect
to each of two variables, so when $s\in G$ is fixed, then the function
$\widehat{w}(s):H\to\C$ is also holomorphic. At the same time it is bounded by
the semicharacter $g(s)\cdot h$:
 \beq\label{widetilde-w-(s)-in-g(s)-cdot-h^blacksquare}
\forall s\in G\qquad \widehat{w}(s)\in g(s)\cdot h^\text{\BSQ}
 \eeq
since
 $$
w\in(g\boxdot h)^\text{\BSQ}\quad\Longrightarrow\quad
|\widehat{w}(s)(t)|=|w(s,t)|\le g(s)\cdot h(t)\quad\Longrightarrow\quad
\widehat{w}(s)\in g(s)\cdot h^\text{\BSQ}
 $$
Thus, $\widehat{w}(s)$ is always a holomorphic function of exponential type on
$H$, i.e. $\widehat{w}((s)\in {\mathcal O}_{\exp}(H)$

2. Let us show that the mapping $\widehat{w}:G\to{\mathcal O}_{\exp}(H)$ is
continuous. Let $s_i$ be a sequence of points in $G$ tending to a point $s$:
$$
s_i\overset{G}{\underset{i\to\infty}{\longrightarrow}} s
$$
Then the sequence of holomorphic functions $\widehat{w}(s_i)\in {\mathcal
O}(H)$ tends to the holomorphic function $\widehat{w}(s)\in {\mathcal O}(H)$
uniformly on each compact set $K\subseteq H$, i.e. in the space ${\mathcal
O}(H)$:
$$
\widehat{w}(s_i)\overset{{\mathcal O}(H)}
{\underset{i\to\infty}{\longrightarrow}} \widehat{w}(s)
$$
On the other hand, all functions are bounded by the semicharacter $C\cdot h$,
where $C=\max\{\sup_{i}g(s_i), g(s)\}$ is a finite number, since the sequence
$s_i$ converges and together with its limit is a compact set:
$$
|\widehat{w}(s_i)(t)|\le g(s_i)\cdot h(t)\le C\cdot h(t)
$$
Thus, the functions $\widehat{w}(s_i)$ and $\widehat{w}(s)$ lie in a rectangle
generated by the semicharacter $C\cdot h$:
$$
 \{\widehat{w}(s_i);\widehat{w}(s)\}\subseteq(C\cdot h)^\text{\BSQ}
$$
In other words, $\widehat{w}(s_i)$ tends to $\widehat{w}(s)$ on the compact set
$(C\cdot h)^\text{\BSQ}$
$$
\widehat{w}(s_i)\overset{(C\cdot h)^\text{\BSQ}}
{\underset{i\to\infty}{\longrightarrow}} \widehat{w}(s)
$$
so $\widehat{w}(s_i)$ tends to $\widehat{w}(s)$ in the space ${\mathcal
O}_{\exp}(H)$:
$$
\widehat{w}(s_i)\overset{{\mathcal O}_{\exp}(H)}
{\underset{i\to\infty}{\longrightarrow}} \widehat{w}(s)
$$

3. From continuity of the mapping $\widehat{w}:G\to{\mathcal O}_{\exp}(H)$ it
follows that for any functional $\beta\in {\mathcal O}_{\exp}^\star (H)$ the
function $\beta\circ \widehat{w}:G\to\C$ is holomorphic. To prove this one can
use the Morera theorem: consider a closed oriented hypersurface $\varGamma$ in
$G$ of dimension $n=\dim G$ with a sufficiently little diameter and show that
 \beq\label{morera}
\int_{\varGamma} (\beta\circ \widehat{w})(s) \d s=0
 \eeq
Indeed, take a net of functionals $\{\beta_i; i\to\infty\}\subset {\mathcal
O}_{\exp}^\star (H)$ which are linear combinations of delta-functionals, and
approximate $\beta$ in ${\mathcal O}_{\exp}^\star (H)$:
$$
\beta_i=\sum_{k}\lambda^k_i\cdot\delta^{a^k_i},\qquad \beta_i\overset{{\mathcal
O}_{\exp}^\star (H)}{\underset{i\to\infty}{\longrightarrow}}\beta
$$
Then we obtain the following. Since $\widehat{w}:G\to{\mathcal O}_{\exp}(H)$ is
continuous, we have
$$
\beta\circ \widehat{w}\overset{{\mathcal C}(G)} {\underset{\infty\gets
i}{\longleftarrow}}\beta_i\circ \widehat{w}
$$
This implies that for any Radon measure $\alpha\in {\mathcal C}(G)$
$$
\alpha(\beta\circ \widehat{w}) \underset{\infty\gets i}{\longleftarrow}
\alpha(\beta_i\circ \widehat{w})
$$
In particular, for the functional of integrating by our hypersurface
$\varGamma$ we obtain
 \begin{multline*}
\int_{\varGamma} (\beta\circ \widehat{w})(s) \d s \underset{\infty\gets
i}{\longleftarrow} \int_{\varGamma} (\beta_i\circ \widehat{w})(s)\d s=
\int_{\varGamma} \left(\sum_{k}\lambda^k_i\cdot\delta^{a^k_i}\circ
\widehat{w}\right)(s)\d s=\\= \sum_{k}\lambda^k_i\cdot\int_{\varGamma}
\left(\delta^{a^k_i}\circ \widehat{w}\right)(s)\d s=
\sum_{k}\lambda^k_i\cdot\int_{\varGamma} \delta^{a^k_i}(\widehat{w}(s))\d s=\\=
 \sum_{k}\lambda^k_i\cdot\int_{\varGamma}
\widehat{w}(s)(a^k_i)\d s= \sum_{k}\lambda^k_i\cdot\underbrace{\int_{\varGamma}
w(s,a^k_i)\d
s}_{\scriptsize\begin{matrix}\| \\ 0, \\ \text{since $w$ is holomorphic} \\
\text{with respect to the first variable}\end{matrix}}=0
 \end{multline*}
I.e., indeed \eqref{morera} is true.

4. We understood that for any functional $\beta\in {\mathcal O}_{\exp}^\star
(H)$ the function $\beta\circ \widehat{w}:G\to\C$ is holomorphic. Let us show
now that it is of exponential type:
 \beq\label{beta-circ-widetilde-w-in-mathcal O-exp-G}
\forall w\in {\mathcal O}_{\exp}(G\times H)\quad \forall \beta\in {\mathcal
O}_{\exp}^\star (H)\qquad \beta\circ \widehat{w}\in {\mathcal O}_{\exp}(G)
 \eeq
Indeed, since the functional $\beta\in {\mathcal O}_{\exp}^\star (H)$ is
bounded on the compact set $h^\text{\BSQ}\subseteq {\mathcal O}_{\exp}(H)$, it
must be a bounded functional on the Banach representation of the Smith space
$\C h^\text{\BSQ}$, i.e.
 \beq
\forall v\in\C h^\text{\BSQ}\qquad |\beta(v)|\le M\cdot
\norm{v}_{h^\text{\BSQ}}
 \eeq
where
$$
M=\norm{\beta}_{(h^\text{\BSQ})^\circ}:=\max_{v\in
h^\text{\BSQ}}|\beta(v)|,\qquad \norm{v}_{h^\text{\BSQ}}:=\inf\{\lambda>0:\
v\in \lambda\cdot h^\text{\BSQ}\}
$$
So from formula \eqref{widetilde-w-(s)-in-g(s)-cdot-h^blacksquare} we have:
 \begin{multline*}
\widehat{w}(s)\in g(s)\cdot h^\text{\BSQ} \quad\Longrightarrow\quad
\norm{\widehat{w}(s)}_{h^\text{\BSQ}}:=\inf\{\lambda>0:\ \widehat{w}(s)\in
\lambda\cdot h^\text{\BSQ}\}\le g(s) \quad\Longrightarrow\\
\Longrightarrow\quad |\beta(\widehat{w}(s))|\le M\cdot g(s)
 \end{multline*}
I.e. the function $\beta\circ \widehat{w}$ is bounded by the semicharacter
$M\cdot g$ :
 \beq\label{|beta-widetilde-w-s|-le-B-cdot-g-s}
\beta\circ \widehat{w}\in M\cdot g^\text{\BSQ}
 \eeq

5. We have proved \eqref{beta-circ-widetilde-w-in-mathcal O-exp-G}. Now let us
show that for any $w\in(g\boxdot h)^\text{\BSQ}\subseteq {\mathcal
O}_{\exp}(G\times H)$ the mapping
 \beq\label{beta-mapsto-beta-circ-widetilde-w-in-mathcal O-exp-G}
\beta\in {\mathcal O}_{\exp}^\star (H)\mapsto \rho_{G,H}(w)(\beta)=\beta\circ
\widehat{w}\in {\mathcal O}_{\exp}(G)
 \eeq
is continuous, i.e.
 \beq
 \rho_{G,H}(w)\in {\mathcal O}_{\exp}(G)\oslash {\mathcal O}_{\exp}^\star (H)
 \eeq
This follows from \eqref{|beta-widetilde-w-s|-le-B-cdot-g-s}: if $\beta_i$ is a
net tending to zero in ${\mathcal O}_{\exp}^\star (H)$, then
$$
\beta_i\circ\widehat{w}\in M_i\cdot g^\text{\BSQ},\qquad M_i=\max_{v\in
h^\text{\BSQ}}|\beta_i(v)|\underset{i\to\infty}{\longrightarrow} 0
$$
$$
\Downarrow
$$
$$
\beta_i\circ\widehat{w}\overset{g^\text{\BSQ}}{\underset{i\to\infty}{\longrightarrow}}
0
$$
$$ \Downarrow $$ $$ \beta_i\circ\widehat{w}\overset{{\mathcal
O}_{\exp}(G)}{\underset{i\to\infty}{\longrightarrow}} 0
$$

6. Now we need to verify that
 \beq
\rho_{G,H}(w)\in g^\text{\BSQ}\odot h^\text{\BSQ}=g^\text{\BSQ}\oslash
(h^\text{\BSQ})^\circ\subseteq {\mathcal O}_{\exp}(G)\oslash {\mathcal
O}_{\exp}^\star (H)
 \eeq
Or, in other words,
 $$
\rho_{G,H}(w)\left((h^\text{\BSQ})^\circ\right)\subseteq g^\text{\BSQ}
 $$
This follows from \eqref{widetilde-w-(s)-in-g(s)-cdot-h^blacksquare}:
 $$
 \forall s\in G\qquad\widehat{w}(s)\in g(s)\cdot h^\text{\BSQ}
 $$
 $$
 \Downarrow
 $$
 $$
 \forall s\in G\qquad \frac{1}{g(s)}\widehat{w}(s)\in h^\text{\BSQ}
 $$
 $$
 \Downarrow
 $$
 $$
\forall s\in G\qquad\forall \beta\in (h^\text{\BSQ})^\circ\qquad  1\ge
\left|\beta\left(\frac{1}{g(s)}\widehat{w}(s)\right)\right|=
\frac{1}{g(s)}|(\beta\circ\widehat{w})(s)|=
\frac{1}{g(s)}|(\rho_{G,H}(w)(\beta))(s)|
 $$
 $$
 \Downarrow
 $$
 $$
\forall s\in G\qquad\forall \beta\in (h^\text{\BSQ})^\circ\qquad
|\rho_{G,H}(w)(\beta)(s)|\le g(s)
 $$
$$
 \Downarrow
 $$
 $$
\forall \beta\in (h^\text{\BSQ})^\circ\qquad \rho_{G,H}(w)(\beta)\in
g^\text{\BSQ}
 $$
 $$
 \Downarrow
 $$
$$
\rho_{G,H}(w)\left((h^\text{\BSQ})^\circ\right)\subseteq g^\text{\BSQ}
 $$

7. Let us show that the mapping
 \beq
w\in (g\boxdot h)^\text{\BSQ} \quad\mapsto\quad \rho_{G,H}(w)\in
g^\text{\BSQ}\odot h^\text{\BSQ}
 \eeq
is injective. Consider functionals of the form
 \beq
\delta^{s,t}:{\mathcal O}_{\exp}(G)\oslash {\mathcal
O}_{\exp}^\star(H)\to\C\quad\Big|\quad \delta^{s,t}(\ph)=\ph(\delta^t)(s)
 \eeq
Now we have: if $w\ne 0$, then for some $s\in G$, $t\in H$ we have $w(s,t)\ne
0$, so
$$
\delta^{s,t}(\rho_{G,H}(w))=\rho_{G,H}(w)(\delta^t)(s)=
(\delta^t\circ\widehat{w})(s)=
\delta^t\Big(\widehat{w}(s)\Big)=\widehat{w}(s)(t)=w(s,t)\ne 0
$$
thus, $\rho_{G,H}(w)\ne 0$.

8. Similarly it turns out that the mapping
 \beq
w\in (g\boxdot h)^\text{\BSQ} \quad\mapsto\quad \rho_{G,H}(w)\in
g^\text{\BSQ}\odot h^\text{\BSQ}
 \eeq
is surjective: for any $\ph\in g^\text{\BSQ}\odot
h^\text{\BSQ}=g^\text{\BSQ}\oslash (h^\text{\BSQ})^\circ\subset{\mathcal
O}_{\exp}(G)\oslash {\mathcal O}_{\exp}^\star(H)$ we set
$$
w(s,t)=\delta^s\circledast\delta^t(\ph)=\ph(\delta^t)(s),\qquad s\in G,\ t\in H
$$
and then, first, $w$ is a holomorphic function on $G\times H$, since it is
holomorphic with respect to every variable: when $t\in H$ is fixed, the object
$\ph(\delta^t)$ is an element of the space ${\mathcal O}_{\exp}(G)$, i.e. a
holomorphic function (of exponential type) on $G$, hence $w(\cdot,t)$ is
holomorphic with respect to the first variable, and when $s\in G$ is fixed, the
mapping $\beta\in {\mathcal O}_{\exp}^\star(H) \mapsto
(\delta^s\circ\ph)(\beta)$ is a continuous functional on the space ${\mathcal
O}_{\exp}^\star(H)$, i.e. by stereotype duality, an element of the space
${\mathcal O}_{\exp}(H)$:
 $$
(\delta^s\circ\ph)(\beta)=\beta(v),\qquad v\in{\mathcal O}_{\exp}(H)
 $$
Thus
 $$
w(s,t)=\ph(\delta^t)(s)=(\delta^s\circ\ph)(\delta^t)=\delta^t(v)=v(t)
 $$
i.e. the function $w(s,\cdot)$ is holomorphic with respect to the second
variable.

Further, turning the chain of item 6 to reverse direction, we obtain:
 $$
\ph\in g^\text{\BSQ}\odot h^\text{\BSQ}
 $$
 $$
 \Downarrow
 $$
 $$
\ph\left((h^\text{\BSQ})^\circ\right)\subseteq g^\text{\BSQ}
 $$
 $$
 \Downarrow
 $$
 $$
\forall \beta\in (h^\text{\BSQ})^\circ\qquad \ph(\beta)\in g^\text{\BSQ}
 $$
 $$
 \Downarrow
 $$
 $$
\forall t\in H\quad \frac{1}{h(t)}\delta^t\in
(h^\text{\BSQ})^\circ\quad\Longrightarrow\quad
\ph\left(\frac{1}{h(t)}\delta^t\right)\in g^\text{\BSQ}
 $$
 $$
 \Downarrow
 $$
 $$
\forall t\in H\qquad  \ph(\delta^t)\in h(t)\cdot g^\text{\BSQ}
 $$
 $$
 \Downarrow
 $$
 $$
\forall s\in G\quad \forall t\in H\qquad  |w(s,t)|=|\ph(\delta^t)(s)|\le
h(t)\cdot g(s)=|(g\boxdot h)(s,t)|
 $$
 $$
 \Downarrow
 $$
 $$
w\in (g\boxdot h)^\text{\BSQ}
 $$
And finally it remains to note that $w$ is an inverse image of $\ph$ under the
mapping $w\mapsto \rho_{G,H}(w)$:
$$
\forall s,t\qquad \rho_{G,H}(w)(\delta^t)(s)=(\delta^t\circ\widehat{w})(s)=
\delta^t(\widehat{w}(s))=\widehat{w}(s)(t)=w(s,t)=\ph(\delta^t)(s)
$$
 $$
 \Downarrow
 $$
$$
\rho_{G,H}(w)=\ph
$$

9. Thus, we obtained that the mapping
 \beq
w\in (g\boxdot h)^\text{\BSQ} \quad\mapsto\quad \rho_{G,H}(w)\in
g^\text{\BSQ}\odot h^\text{\BSQ}
 \eeq
is bijective. It remains to show now that it is continuous in both directions.
This follows from the fact that both $(g\boxdot h)^\text{\BSQ}$ and
$g^\text{\BSQ}\odot h^\text{\BSQ}$ are compact sets. Since functionals
$\delta^s\circledast\delta^t$ separate the points of the compact set
$g^\text{\BSQ}\odot h^\text{\BSQ}$, the (Hausdorff) topology they generate on
$g^\text{\BSQ}\odot h^\text{\BSQ}$ coincides with the initial topology of
$g^\text{\BSQ}\odot h^\text{\BSQ}$: if
$$
\delta^s\circledast\delta^t(\ph_i)\underset{i\to\infty}{\longrightarrow}\delta^s\circledast\delta^t(\ph)
$$
for any $s\in G$, $t\in H$, then
$$
\ph_i\overset{g^\text{\BSQ}\odot
h^\text{\BSQ}}{\underset{i\to\infty}{\longrightarrow}}\ph
$$
From this we have that the mapping $w\mapsto \rho_{G,H}(w)$ is continuous in
forward direction:
$$
w_i\overset{(g\boxdot h)^\text{\BSQ}}{\underset{i\to\infty}{\longrightarrow}} w
$$
 $$
 \Downarrow
 $$
$$
\forall (s,t)\in G\times H\qquad
\delta^s\circledast\delta^t(\rho_{G,H}(w_i))=w_i(s,t)\underset{i\to\infty}{\longrightarrow}
w(s,t)=\delta^s\circledast\delta^t(\rho_{G,H}(w))
$$
 $$
 \Downarrow
 $$
$$
\rho_{G,H}(w_i)\overset{g^\text{\BSQ}\odot
h^\text{\BSQ}}{\underset{i\to\infty}{\longrightarrow}}\rho_{G,H}(w)
$$
Thus, the operation $w\mapsto \rho_{G,H}(w)$ is a continuous bijective mapping
of the compact set $(g\boxdot h)^\text{\BSQ}$ into the compact set
$g^\text{\BSQ}\odot h^\text{\BSQ}$. This means that $\rho_{G,H}$ is a
homeomorphism between $(g\boxdot h)^\text{\BSQ}$ and $g^\text{\BSQ}\odot
h^\text{\BSQ}$.
 \epr

\bpr[Proof of Theorem
\ref{mathcal-O-exp-G-times-H-cong-mathcal-O-exp-G-mathcal-O-exp-H}]  Let us
note from the very beginning that if  $f$ is a semicharacter on $G\times H$,
then the functions
 $$
g(s)=f(s,1_H),\qquad h(t)=f(1_G,t),\qquad s\in G,\; t\in H
 $$
are again semicharacters (as restrictions of $f$ on subgroups), and the
function $g\boxdot h$ is a semicharacter on $G\times H$, majorizing $f$:
 \beq\label{f-le-g-otimes-h}
f\le g\boxdot h
 \eeq
Indeed,
$$
f(s,t)=f( (s,1_H)\cdot (1_G,t) )\le f(s,1_H)\cdot f(1_G,t)=g(s)\cdot
h(t)=(g\boxdot h)(s,t)
$$

From this it follows that every function $w\in {\mathcal O}_{\exp}(G\times H)$
is contained in some compact set of the form $(g\boxdot h)^\text{\BSQ}$ (since
$w$ is always contained in compact set of the form $f^\text{\BSQ}$). At the
same time the object $\rho_{G,H}(w)$ is an element of the set
$g^\text{\BSQ}\odot h^\text{\BSQ}$, i.e. an element of the space ${\mathcal
O}_{\exp}(G)\oslash {\mathcal O}_{\exp}^\star(H)$.

Thus, Formulas \eqref{widetilde-w-s-t=w-s-t} and \eqref{widetilde-w-s-t-star}
correctly define a mapping
 \beq
w\in {\mathcal O}_{\exp}(G\times H) \quad\mapsto\quad \rho_{G,H}(w)\in
{\mathcal O}_{\exp}(G)\oslash {\mathcal O}_{\exp}^\star(H)
 \eeq
and we only need to check its bijectivity and continuity in both directions.

1. The injectivity of $\rho_{G,H}$ follows from its injectivity on compact sets
$(g\boxdot h)^\text{\BSQ}$ (and from the fact that the compact sets $(g\boxdot
h)^\text{\BSQ}$ and $g^\text{\BSQ}\odot h^\text{\BSQ}$ are injectively included
into the spaces ${\mathcal O}_{\exp}(G\times H)$ and ${\mathcal
O}_{\exp}(G)\oslash {\mathcal O}_{\exp}^\star(H)$).

2. The surjectivity of $\rho_{G,H}$ follows from the fact that it surjectively
maps compact sets $(g\boxdot h)^\text{\BSQ}$ into compact sets
$g^\text{\BSQ}\odot h^\text{\BSQ}$ (and from the fact that the compact sets
$g^\text{\BSQ}\odot h^\text{\BSQ}$ cover all the space ${\mathcal
O}_{\exp}(G)\oslash {\mathcal O}_{\exp}^\star(H)$).

3. The continuity follows from the fact that $\rho_{G,H}$ continuously maps
every compact set $(g\boxdot h)^\text{\BSQ}$ into the compact set
$g^\text{\BSQ}\odot h^\text{\BSQ}$, hence into the space ${\mathcal
O}_{\exp}(G)\oslash {\mathcal O}_{\exp}^\star(H)$. This means that $\rho_{G,H}$
is continuous on each compact set $K$ in the Brauner space ${\mathcal
O}_{\exp}(G\times H)$, thus it must be continuous on the whole space ${\mathcal
O}_{\exp}(G\times H)$.

4. The continuity in reverse direction is proved in the same way: since the
inverse mapping continuously turns every compact set $g^\text{\BSQ}\odot
h^\text{\BSQ}$ into a compact set $(g\boxdot h)^\text{\BSQ}$, it is continuous
on each compact set in the Brauner space ${\mathcal O}_{\exp}(G)\oslash
{\mathcal O}_{\exp}^\star(H)$. Thus, it is continuous on the whole space
${\mathcal O}_{\exp}(G)\oslash {\mathcal O}_{\exp}^\star(H)$.

5. We have proved that the mapping defined by formulas
\eqref{widetilde-w-s-t-star} - \eqref{widetilde-w-s-t=w-s-t} is an isomorphism
of the stereotype spaces: ${\mathcal O}_{\exp}(G\times H)\cong {\mathcal
O}_{\exp}(G)\oslash {\mathcal O}_{\exp}^\star(H)={\mathcal O}_{\exp}(G)\odot
{\mathcal O}_{\exp}(H)$ (hence the identities \eqref{tenz-pr-O-exp} hold). Let
us show that this mapping satisfies the identity \eqref{rho_G,H}: if
$u\in{\mathcal O}_{\exp}(G)$, $v\in{\mathcal O}_{\exp}(H)$, then for the
mapping $\widehat{u\boxdot v}:G\to{\mathcal O}_{\exp}(H)$ defined by formula
\eqref{widetilde-w-s-t=w-s-t} we have the following logic chain:
$$
\widehat{u\boxdot v}(s)(t)=(u\boxdot v)(s,t)=u(s)\cdot v(t)
$$
$$
\Downarrow
$$
$$
\widehat{u\boxdot v}(s)=u(s)\cdot v
$$
$$
\Downarrow
$$
$$
\forall\beta\in{\mathcal O}_{\exp}^\star(H)\qquad \beta(\widehat{u\boxdot
v}(s))=u(s)\cdot \beta(v)
$$
$$
\Downarrow
$$
$$
\forall\beta\in{\mathcal O}_{\exp}^\star(H)\qquad \rho_{G,H}(u\boxdot
v)(\beta)=\beta\circ\widehat{u\boxdot v}=\beta(v)\cdot
u=\eqref{x-odot-y}=(u\odot v)(\beta)
$$
$$
\Downarrow
$$
$$
\rho_{G,H}(u\boxdot v)=u\odot v
$$
It remains now to note that since the elements of the form $u\circledast v$
generate a dense subspace in ${\mathcal O}_{\exp}(G)\circledast {\mathcal
O}_{\exp}(H)$, by the proven above nuclearity of the spaces ${\mathcal
O}_{\exp}$, the corresponding elements of the form $u\odot v$ must generate a
dense subspace in ${\mathcal O}_{\exp}(G)\odot {\mathcal O}_{\exp}(H)$, and
elements $u\boxdot v$ a dense subspace in ${\mathcal O}_{\exp}(G\times H)$.
This implies that the property \eqref{rho_G,H} uniquely define the mapping
$\rho_{G,H}$.
 \epr

\subsection{Structure of Hopf algebras on ${\mathcal O}_{\exp}(G)$ and ${\mathcal
O}^\star_{\exp}(G)$}\label{proverka-Hopfa-dlya-H(G)}

In \ref{SEC:stein-groups}\ref{R-O} we have mentioned the standard trick,
allowing to prove that functional algebras of a given class on groups are Hopf
algebras -- a sufficient condition for this is a natural isomorphism between
the functional algebra of the Cartesian product of groups $\times$ and the
corresponding tensor product of their functional algebras. Theorem
\ref{mathcal-O-exp-G-times-H-cong-mathcal-O-exp-G-mathcal-O-exp-H},
establishing the natural isomorphism
 $$
{\mathcal O}_{\exp}(G\times H)\stackrel{\rho_{G,H}}{\cong} {\mathcal
O}_{\exp}(G)\odot {\mathcal O}_{\exp}(H)
 $$
allows now to make the same conclusion about algebras ${\mathcal O}_{\exp}(G)$:

\btm\label{O-exp-(G)-algebra-Hopfa} For any compactly generated Stein group $G$
 \bit
\item[--] the space ${\mathcal O}_{\exp}(G)$ of holomorphic functions of
exponential type on $G$ is a nuclear Hopf-Brauner algebra with respect to the
algebraic operations defined by formulas, analogous to
\eqref{umn-v-O(G)}-\eqref{antipod-v-O(G)};

\item[--] its stereotype dual space ${\mathcal O}_{\exp}^\star(G)$ is a nuclear
Hopf-Fr\'echet algebra with respect to dual algebraic operations.
 \eit
\etm

\section{Arens-Michael envelope and holomorphic reflexivity}\label{SEC:Arens-Michael}

\subsection{Submultiplicative seminorms and Arens-Michael algebras}

A seminorm $p:A\to\R_+$ on an algebra $A$ is said to be {\it
submultiplicative}\index{submultiplicative!seminorm}, if it satisfies the
following condition
$$
p(u\cdot v)\le p(u)\cdot p(v),\qquad u,v\in A
$$
This is equivalent to the fact that the unit ball of this seminorm
$$
U=\{u\in A:\; p(u)\le 1\}
$$
satisfies the condition
$$
U\cdot U\subseteq U
$$
(such sets in $A$ are also said to be {\it
submultiplicative}\index{submultiplicative!set} like submultiplicative
neighborhoods of zero, defined on page \pageref{DEF:submult-neighb-zero}).

A topological algebra $A$ is called an {\it Arens-Michael
algebra}\index{Arens-Michael!algebra}, if it is complete (as a topological
vector space) and satisfies the following equivalent conditions:
 \bit
\item[(i)] the topology of $A$ is generated by a system of submultiplicative
seminorms:

\item[(ii)] $A$ has a local base of submultiplicative closed absolutely convex
neighborhoods of zero.
 \eit

\bex\label{EX-polunormy-v-O(Z)} Seminorms \eqref{polunormy-v-H(Z)} generating
the topology on ${\mathcal O}(\Z)$ are submultiplicative:
 $$
||u||_N=\sum_{|n|\le N} |u(n)|,\qquad N\in\N,
 $$
hence ${\mathcal O}(\Z)$ is an Arens-Michael algebra.
 \eex \bpr
Indeed,
$$
\norm{u\cdot v}_N=\sum_{|n|\le N}|(u\cdot
v)(n)|=\eqref{umnozhenie-v-O(Z)-i-O^star(Z)}=\sum_{|n|\le N}|u(n)\cdot v(n)|\le
\left(\sum_{|n|\le N}|u(n)|\right)\cdot\left(\sum_{|n|\le
N}|v(n)|\right)=\norm{u}_N\cdot\norm{v}_N
$$
 \epr

\bex\label{EX-polunormy-v-O(C^x)} Seminorms \eqref{polunormy-v-H(C-x)}
generating the topology on ${\mathcal O}(\C^\times)$, are submultiplicative:
 $$
||u||_C=\sum_{n\in\Z} |u_n|\cdot C^{|n|},\qquad C\ge 1,
 $$
hence ${\mathcal O}(\C^\times)$ is an Arens-Michael algebra.
 \eex \bpr
Indeed,
 \begin{multline*}
||u\cdot v||_C=\sum_{n\in\Z} |(u\cdot v)_n|\cdot C^{|n|}
=\eqref{umnozhenie-v-O(C-x)-i-O^star(C-x)}=\sum_{n\in\Z} \left|\sum_{i\in\Z}
u_i\cdot v_{n-i}\right|\cdot C^{|n|}\le\\ \le \sum_{n\in\Z} \sum_{i\in\Z}
|u_i|\cdot |v_{n-i}|\cdot C^{|i|}\cdot C^{|n-i|}=\left(\sum_{k\in\Z} |u_k|\cdot
C^{|k|}\right)\cdot \left(\sum_{l\in\Z} |u_l|\cdot C^{|l|}\right)=||u||_C\cdot
||v||_C
 \end{multline*}\epr

\bex\label{EX-polunormy-v-H(C)} Seminorms \eqref{polunormy-v-H(C)}, generating
the topology on ${\mathcal O}(\C)$, are submultiplicative:
 $$
||u||_C=\sum_{n=0}^\infty |u_n|\cdot C^n,\qquad C\ge 0,
 $$
hence ${\mathcal O}(\C)$ is an Arens-Michael algebra.
 \eex \bpr
Indeed,
 \begin{multline*}
||u\cdot v||_C=\sum_{n=0}^\infty |(u\cdot v)_n|\cdot C^n
=\eqref{umnozhenie-v-O(C)-i-O^star(C)}=\sum_{n=0}^\infty \left|\sum_{i=0}^n
u_i\cdot v_{n-i}\right|\cdot C^n\le \sum_{n=0}^\infty \sum_{i=0}^n |u_i|\cdot
|v_{n-i}|\cdot C^n=\\=\left(\sum_{k\in\N} |u_k|\cdot C^k\right)\cdot
\left(\sum_{l=0}^\infty |u_l|\cdot C^l\right)=||u||_C\cdot ||v||_C
 \end{multline*}
 \epr

Of the following three propositions the first two are evident, and the third
one follows from the A.~Yu.~Pirkovskii theorem \ref{TH-Pirkovskii} which we
formulate below in this section:

\bprop The algebra ${\mathcal O}(M)$ of holomorphic functions on any complex
manifold $M$ is an Arens-Michael algebra. \eprop

\bprop The algebra ${\mathcal O}_{\exp}^\star(G)$ of exponential functionals on
every Stein group is an Arens-Michael algebra. \eprop

\bprop The algebra ${\mathcal R}(M)$ of polynomials on a complex affine
algebraic manifold $M$, endowed with the strongest locally convex topology, is
an Arens-Michael algebra if an only if the manifold $M$ is finite. \eprop

\subsection{Arens-Michael envelopes}

The {\it Arens-Michael envelope}\index{Arens-Michael!envelope} of a topological
algebra $A$ is a (continuous) homomorphism $\pi:A\to B$ of $A$ into an
Arens-Michael algebra $B$ such that for any  (continuous) homomorphism
$\rho:A\to C$ of $A$ into an arbitrary Arens-Michael algebra $C$ there is a
unique  (continuous) homomorphism $\sigma:B\to C$ such that the following
diagram is commutative:
$$
\begin{diagram}
\node{A} \arrow[2]{e,t}{\pi} \arrow{se,b}{\rho} \node[2]{B}
\arrow{sw,b,--}{\sigma}
\\
\node[2]{C}
\end{diagram}
$$
From this definition it is clear that if $\pi:A\to B$ and $\rho:A\to C$ are two
Arens-Michael envelopes of $A$, then the arising homomorphism $\sigma:B\to C$
becomes an isomorphism of topological algebras (due to the uniqueness of
$\sigma$). Hence the Arens-Michael envelope of $A$  is defined uniquely up to
an isomorphism, and as a corollary, we can introduce a special notation for
this construction:
$$
\heartsuit_A:A\to A^\heartsuit
$$
This should be understood as follows: if we have a homomorphism $\ph:A\to B$,
then the record $\ph=\heartsuit_A$ means that $\ph:A\to B$ is an Arens-Michael
envelope of the algebra $A$; on the other hand, if we have an algebra $B$, then
the record $B=A^\heartsuit$ means that there exists a homomorphism $\ph:A\to B$
which is an Arens-Michael envelope of the algebra $A$
--- in this case the algebra $B$ is also called an {\it Arens-Michael envelope of the algebra} $A$.

\bprop\label{kriterij-A-M} A topological algebra $B$ is an Arens-Michael
envelope for a topological algebra $A$ if and only if
 \bit
\item[(i)] $B$ is an Arens-Michael algebra, and

\item[(ii)] there exists a continuous homomorphism $\pi:A\to B$ such that
 \bit
\item[(a)] the image $\pi(A)$ of the algebra $A$ under the action of $\pi$ is
dense in $B$,

\item[(b)] for any continuous submultiplicative seminorm $p:A\to\R_+$ there is
a continuous submultiplicative seminorm $\widetilde{p}:B\to\R_+$, such that the
seminorm $\widetilde{p}\circ\pi:A\to\R_+$ majorizes the seminorm $p$:
$$
p(a)\le\widetilde{p}(\pi(a)),\qquad a\in A
$$ \eit
 \eit
 \eprop

The Arens-Michael envelope can be constructed directly: to any
submultiplicative neighborhood of zero $U$ in $A$ we can assign the closed
ideal $\Ker U$ in $A$ defined by the equality
$$
\Ker U=\bigcap_{\varepsilon>0}\varepsilon\cdot U
$$
and a quotient algebra
$$
A/\Ker U
$$
endowed with the topology of normed space with the unit ball $U+\Ker U$. Then
the completion $(A/\Ker U)^\blacktriangledown$ becomes a Banach algebra.
Following \eqref{X/U} we denote such algebra by $A/U$:
$$
A/U:=(A/\Ker U)^\blacktriangledown
$$
The family of such algebras (with different $U$) forms a projective system. The
Its limit is an Arens-Michael envelope $A^\heartsuit$:
 \beq\label{konstr-opr-A-M}
A^\heartsuit=\kern-30pt\underset{\scriptsize\begin{matrix}\text{$U$ is a submultiplicative}\\
\text{neighborhood of zero in
$A$}\end{matrix}}{\underset{\longleftarrow}{\lim}} \kern-30pt A/U.
 \eeq

The following propositions show that the operation of taking the Arens-Michael
envelope commutes with the passage to direct sums and quotient algebras.

\bprop The Arens-Michael envelope of a direct sum $A_1\oplus ...\oplus A_n$ of
finite family of topological algebras $A_1,...,A_n$ coincides with the direct
sum of the Arens-Michael envelopes of these algebras:
$$
(A_1\oplus ...\oplus A_n)^\heartsuit\cong A_1^\heartsuit\oplus ...\oplus
A_n^\heartsuit
$$
 \eprop

\bprop Let $\pi:A\to A^\heartsuit$ be an Arens-Michael envelope of the algebra
$A$, and let $I$ be a close d ideal in $A$. Then the Arens-Michael envelope of
the quotient algebra $A/I$ coincides with the completion of the quotient
algebra $A^\heartsuit/\overline{\pi(I)}$ over the closure $\overline{\pi(I)}$
in $A^\heartsuit$ of the image of ideal $I$ under the action of the mapping
$\pi$:
$$
(A/I)^\heartsuit\cong
\big(A^\heartsuit/\overline{\pi(I)}\big)^\blacktriangledown
$$
 \eprop

An important example of the Arens-Michael envelope was constructed by
A~.Yu.~Pirkovskii:

\btm[A~.Yu.~Pirkovskii, \cite{Pirkovskii}]\label{TH-Pirkovskii} The
Arens-Michael envelope of the algebra ${\mathcal R}(M)$ of polynomials on an
affine algebraic manifold $M$ coincides with the algebra ${\mathcal O}(M)$ of
holomorphic functions on $M$:
 \beq\label{R(M)^heartsuit=O(M)}
\Big({\mathcal R}(M)\Big)^\heartsuit\cong {\mathcal O}(M)
 \eeq
\etm

\subsection{The mapping $\flat_G^\star:{\mathcal O}^\star(G)\to {\mathcal O}_{\exp}^\star(G)$
is an Arens-Michael envelope}

\btm\label{TH-AM-O-star} For any Stein group $G$ the mapping
$$
\flat_G^\star:{\mathcal O}^\star(G)\to {\mathcal O}_{\exp}^\star(G)
$$
is an Arens-Michael envelope of the algebra ${\mathcal O}^\star(G)$:
 \beq\label{AM-O-star}
\Big({\mathcal O}^\star(G)\Big)^\heartsuit\cong {\mathcal O}^\star_{\exp}(G)
 \eeq
\etm

\bpr From the representation \eqref{O-exp-top} it follows that this space is a
projective limit of the Banach quotient algebras:
 \begin{multline*}
{\mathcal
O}^\star_{\exp}(G)=\eqref{O-exp-top}=\Bigg(\underset{\scriptsize\begin{matrix}\text{$D$
is a dually}\\ \text{submultiplicative}\\
\text{rectangle in $\mathcal
O(G)$}\end{matrix}}{\underset{\longrightarrow}{\lim}}\kern-15pt \C
D\Bigg)^\star\kern10pt=\kern-25pt \underset{\scriptsize\begin{matrix}\text{$D$ is a dually}\\
\text{submultiplicative}\\ \text{rectangle in $\mathcal
O(G)$}\end{matrix}}{\underset{\longleftarrow}{\lim}}\kern-25pt  (\C
D)^\star=\eqref{(C K)^star=X^star-K^circ}=\kern-25pt \underset{\scriptsize\begin{matrix}\text{$D$ is a dually}\\
\text{submultiplicative}\\ \text{rectangle in $\mathcal
O(G)$}\end{matrix}}{\underset{\longleftarrow}{\lim}}\kern-25pt  {\mathcal
O}^\star(G)/ D^\circ=\\=\kern-25pt
\underset{\scriptsize\begin{matrix}\text{$\varDelta$ is a submultiplicative}\\
\text{rhombus in $\mathcal
O^\star(G)$}\end{matrix}}{\underset{\longleftarrow}{\lim}}\kern-25pt  {\mathcal
O}^\star(G)/ \varDelta =(\text{Theorem \ref{subm-in-trubk}(a)})=\kern-25pt
\underset{\scriptsize\begin{matrix}\text{$U$ is a submultiplicative}\\
\text{neighborhood of zero in $\mathcal
O^\star(G)$}\end{matrix}}{\underset{\longleftarrow}{\lim}}\kern-25pt  {\mathcal
O}^\star(G)/ U=\eqref{konstr-opr-A-M}=\Big({\mathcal
O}^\star(G)\Big)^\heartsuit
 \end{multline*}
 \epr

\subsection{The mapping $\flat_G:{\mathcal O}_{\exp}(G)\to {\mathcal O}(G)$ is an Arens-Michael envelope for the
groups with algebraic connected component of identity}

\btm\label{TH-AM-O} Let $G$ be a compactly generated Stein group, such that the
connected component of identity $G_e$ is an algebraic group. Then the mapping
$$
\flat_G:{\mathcal O}_{\exp}(G)\to {\mathcal O}(G)
$$
is an Arens-Michael envelope of the algebra ${\mathcal O}_{\exp}(G)$:
 \beq\label{AM-O}
\Big({\mathcal O}_{\exp}(G)\Big)^\heartsuit\cong {\mathcal O}(G)
 \eeq
 \etm

We shall prove this theorem in several steps.

1. Let initially $G$ be a discrete group. Recall that in \eqref{1_x-infty} we
agreed to denote by $1_x$ the characteristic functions of singletons $\{x\}$ in
$G$:
 \beq\label{1_x}
1_x(y)=\begin{cases}1, & y=x \\ 0& y\ne x\end{cases}
 \eeq
(since $G$ is discrete, the function $1_x$ can be considered as element of both
algebras $\mathcal O(G)$ and ${\mathcal O}_{\exp}(G)$).

\blm\label{LM_1_x} The functions $\{1_x;\; x\in G\}$ form a basis in the
topological vector spaces ${\mathcal O}(G)$ and ${\mathcal O}_{\exp}(G)$: for
any function $u\in {\mathcal O}(G)$ ($u\in {\mathcal O}_{\exp}(G)$) the
following equality holds
 \beq\label{series-in-H}
u=\sum_{x\in G} u(x)\cdot 1_x
 \eeq
where the series converges in ${\mathcal O}(G)$ (${\mathcal O}_{\exp}(G)$), and
its coefficients continuously depend on $u\in {\mathcal O}(G)$ ($u\in {\mathcal
O}_{\exp}(G)$).\elm

\bpr For the space ${\mathcal O}(G)$ this is evident, since for the case of
discrete group $G$ this space coincides with the space $\C^G$ of all functions
on $G$. Let us prove this for ${\mathcal O}_{\exp}(G)$: if $u\in {\mathcal
O}_{\exp}(G)$, then taking a majorizing semicharacter $f:G\to\R_+$,
$$
|u(x)|\le f(x),\qquad x\in G
$$
we obtain that the partial sums of the series \eqref{series-in-H} are contained
in the rectangle $f^\text{\BSQ}$, so the series \eqref{series-in-H} converges
(not only in ${\mathcal O}(G)$, but) in ${\mathcal O}_{\exp}(G)$ as well. On
the other hand, every coefficient $u(x)$ continuously depend $u$, if $u$ runs
over the rectangle $f^\text{\BSQ}$. By the definition of topology in ${\mathcal
O}_{\exp}(G)$, this means that $u(x)$ continuously depend on $u$, when $u$ runs
over ${\mathcal O}_{\exp}(G)$.
 \epr

\blm\label{sum-q(1_x)-f(x)<infty} If $G$ is a discrete finitely generated
group, then for any continuous seminorm $q:{\mathcal O}_{\exp}(G)\to\R_+$ and
for any semicharacter $f:G\to[1;+\infty)$ the number family $\{f(x)\cdot
q(1_x);\; x\in G\}$ is summable:
$$
\sum_{x\in G} f(x)\cdot q(1_x)<\infty
$$
\elm
 \bpr
Let $T$ is an absolutely convex compact set in ${\mathcal O}_{\exp}^\star(G)$
corresponding to the seminorm $q$:
 $$
q(u)=\sup_{\alpha\in T}|\alpha(u)|
 $$
Every rectangle $f^\text{\BSQ}$ is a compact set in ${\mathcal
O}_{\exp}^\star(G)$, so
 \begin{multline*}
\infty>\sup_{u\in f^\text{\BSQ}}\sup_{\alpha\in T}
|\alpha(u)|=\eqref{series-in-H}=\sup_{u\in f^\text{\BSQ}}\sup_{\alpha\in T}
\Big|\alpha\Big(\sum_{x\in G} u(x)\cdot 1_x\Big)\Big|=\sup_{u\in
f^\text{\BSQ}}\sup_{\alpha\in T}  \Big|\sum_{x\in G} u(x)\cdot
\alpha(1_x)\Big|\ge\\ \ge \sup_{\alpha\in T} \Big|\sum_{x\in G}
\underbrace{\frac{f(x)\cdot\overline{\alpha(1_x)}}{|\alpha(1_x)|}}_{\scriptsize\begin{matrix}\uparrow\\
\text{one of the values of}\\ \text{$u\in f^\text{\BSQ}$}\end{matrix}}\cdot
\alpha(1_x)\Big|= \sup_{\alpha\in T} \sum_{x\in G} f(x)\cdot |\alpha(1_x)|\ge
\sup_{\alpha\in T} \sup_{x\in G} f(x)\cdot |\alpha(1_x)|=\\=\sup_{x\in G}
f(x)\cdot \sup_{\alpha\in T} |\alpha(1_x)|=\sup_{x\in G} f(x)\cdot q(1_x)
  \end{multline*}
So we have that for any semicharacter $f:G\to[1;+\infty)$
$$
\sup_{x\in G} f(x)\cdot q(1_x)<\infty
$$
Now take a finite set $K$, generating $G$,
$$
\bigcup_{n=1}^\infty K^n=G
$$
and define a semicharacter $g:G\to[1;+\infty)$ by formula
$$
g(x)=R^n\quad\Longleftrightarrow\quad x\in K^n\setminus K^{n-1}
$$
where $R$ is a number, greater than the cardinality of $K$:
$$
R>\card K.
$$
Since the product $g\cdot f$ is also a semicharacter, we have:
$$
\sup_{x\in G} \Big[ g(x)\cdot f(x)\cdot q(1_x)\Big]<\infty
$$
$$
\Downarrow
$$
$$
\exists C>0\qquad \forall x\in G \qquad f(x)\cdot q(1_x)\le \frac{C}{g(x)}
$$
$$
\Downarrow
$$
 \begin{multline*}
\sum_{x\in G} f(x)\cdot q(1_x)\le \sum_{x\in G} \frac{C}{g(x)}=
\sum_{n=1}^\infty \sum_{x\in K^n\setminus K^{n-1}} \frac{C}{g(x)}=
\sum_{n=1}^\infty \sum_{x\in K^n\setminus K^{n-1}} \frac{C}{R^n}\le
\sum_{n=1}^\infty  \frac{C\cdot \card (K^n)}{R^n}\le\\ \le \sum_{n=1}^\infty
\frac{C\cdot (\card K)^n}{R^n}=C\cdot \sum_{n=1}^\infty \left(\frac{\card
K}{R}\right)^n<\infty
  \end{multline*}
 \epr

If $q:{\mathcal O}_{\exp}(G)\to\R_+$ is a continuous seminorm on ${\mathcal
O}_{\exp}(G)$, then let us call its {\it support} the set
 \beq
\supp(q)=\{x\in G:\; q(1_x)\ne 0\}
 \eeq

\blm\label{LM-supp-q} If $G$ is a discrete finitely generated group, then for
any submultiplicative continuous seminorm $q:{\mathcal O}_{\exp}(G)\to\R_+$
 \bit
\item[(a)] its support $\supp(q)$ is a finite set:
$$
\card \supp(q)<\infty
$$
\item[(b)] for any point $x\in\supp(q)$ the value of the seminorm $q$ on any
function $1_x$ is less than 1:
$$
q(1_x)\ge 1
$$
 \eit
\elm \bpr Let us prove (b) first. If $x\in\supp(q)$, i.e. $q(1_x)>0$, then:
$$
1_x=1_x^2\quad\Longrightarrow\quad q(1_x)=q(1_x^2)\le
q(1_x)^2\quad\Longrightarrow\quad 1\le q(1_x)
$$
Now (a). Since the constant identity $f(x)=1$ is a semicharacter on $G$, by
Lemma \ref{sum-q(1_x)-f(x)<infty} the number family $\{q(1_x);\; x\in G\}$ is
summing:
$$
\sum_{x\in G} q(1_x)<\infty
$$
On the other hand, by the condition (b) we have already proved, all non-zero
terms in this series are bounded from below by 1. Hence there is a finite
number of them.
 \epr

2. Let us now pass to the case when $G$ is compactly generated Stein group,
whose connected component of identity $G_e$ is an algebraic group. Let
${\mathcal LC}_{\exp}(G)$ denote a subalgebra in ${\mathcal O}_{\exp}(G)$
consisting of locally constant functions:
$$
u\in{\mathcal LC}_{\exp}(G)\quad\Longleftrightarrow\quad u\in {\mathcal
O}_{\exp}(G)\quad\&\quad \Big(\forall x\in G\quad \exists\; \text{neighborhood}
\; U\owns x\quad \forall y\in U \quad u(x)=u(y)\Big)
$$

 \blm\label{H-exp(G/G_e)=LC-exp(G)}
Let $\pi:G\to G/G_e$ denote the quotient mapping. For any function
$v\in{\mathcal O}_{\exp}(G/G_e)$ the composition $v\circ\pi$ is a locally
constant function of exponential type on $G$, and the mapping
$$
v\mapsto v\circ\pi
$$
establishes an isomorphism of topological algebras:
$$
{\mathcal O}_{\exp}(G/G_e)\cong {\mathcal LC}_{\exp}(G)
$$
 \elm

Let now for any coset $K\in G/G_e$ and for any function $u\in{\mathcal
O}_{\exp}(G)$ the symbol $u_K$ denote the function coinciding with $u$ on the
set $K\subset G$ and vanishing outside of $K$:
 \beq\label{u_K}
u_K(x)=\begin{cases}u(x),& x\in K \\ 0,& x\notin K\end{cases}
 \eeq
The following proposition is proved just like Lemma \ref{LM_1_x}:

\blm\label{LM-sum-u_K} For any coset $K\in G/G_e$ and for any function
$u\in{\mathcal O}_{\exp}(G)$
 \bit
\item[(a)] the function $u_K$ belongs to algebra ${\mathcal O}_{\exp}(G)$,

\item[(b)] the mapping $u\in{\mathcal O}_{\exp}(G)\mapsto u_K\in{\mathcal
O}_{\exp}(G)$ is continuous, and

\item[(c)] the series $\sum_{K\in G/G_e} u_K$ converges in the space ${\mathcal
O}_{\exp}(G)$ to the function $u$:
$$
u=\sum_{K\in G/G_e} u_K
$$
 \eit
 \elm

For any coset $K\in G/G_e$ let us consider the operator of projection
$$
P_K:{\mathcal O}_{\exp}(G)\to {\mathcal O}_{\exp}(G),\qquad P_K(u)=u_K
$$
and let ${\mathcal O}_{\exp}(K)$ denote its image in the space ${\mathcal
O}_{\exp}(G)$:
$$
{\mathcal O}_{\exp}(K)=P_K\Big({\mathcal O}_{\exp}(G)\Big)
$$
Clearly, ${\mathcal O}_{\exp}(K)$ is a closed subspace in ${\mathcal
O}_{\exp}(G)$, so ${\mathcal O}_{\exp}(K)$ can be endowed with the topology
induced from ${\mathcal O}_{\exp}(G)$ (this will be the same as the topology of
an immediate subspace in ${\mathcal O}_{\exp}(G)$), and with respect to this
topology ${\mathcal O}_{\exp}(K)$ is a Brauner space.

Let in addition ${\mathcal O}(K)$ denote the usual algebra of holomorphic
functions on the complex manifold $K$.

\blm\label{LM-H-exp-K-AM=H-K} The inclusion ${\mathcal O}_{\exp}(K)\subseteq
{\mathcal O}(K)$ is an Arens-Michael envelope:
$$
{\mathcal O}_{\exp}(K)^\heartsuit={\mathcal O}(K)
$$
 \elm
 \bpr
We have to note first that it is sufficient to consider the case of $K=G_e$,
since shifts turn inclusions ${\mathcal O}_{\exp}(K)\subseteq {\mathcal O}(K)$
into inclusions ${\mathcal O}_{\exp}(G_e)\subseteq {\mathcal O}(G_e)$. For this
case our proposition becomes a corollary of the Pirkovskii Theorem
\ref{TH-Pirkovskii}: by assumption, $G_e$ is an algebraic group, so we can
consider the algebra ${\mathcal R}(G_e)$ of polynomials on $G_e$. Then the
train of thought is illustrated by the following diagram (where the horizontal
arrows mean inclusions):
$$
\begin{diagram}
\node{{\mathcal R}(G_e)} \arrow{e,t}{}\arrow{se,b}{\rho_{\mathcal R}}
\node{{\mathcal O}_{\exp}(G_e)}\arrow{e,t}{} \arrow{s,b}{\rho} \node{{\mathcal
O}(G_e)} \arrow{sw,b,--}{\widetilde{\rho}}
\\
\node[2]{B}
\end{diagram}
$$
If $\rho:{\mathcal O}_{\exp}(G_e)\to B$ is an arbitrary morphism into the
Arens-Michael algebra in $B$, then by Pirkovskii's Theorem there arises a
unique morphism $\rho_{\mathcal R}:{\mathcal R}(G_e)\to B$. Since ${\mathcal
R}(G_e)$ is dense in ${\mathcal O}_{\exp}(G_e)$ in the topology of ${\mathcal
O}(G_e)$, the morphism $\widetilde{\rho}$ extends the morphism $\rho$. And
since ${\mathcal O}_{\exp}(G_e)$ is dense in ${\mathcal O}(G_e)$, this
extension is unique.
 \epr

\bpr[Proof of Theorem \ref{TH-AM-O}] Let $G$ be an arbitrary compactly
generated Stein group with algebraic connected component of identity and
$p:{\mathcal O}_{\exp}(G)\to\R_+$ a submultiplicative seminorm. Its restriction
$p|_{{\mathcal LC}_{\exp}(G)}$ to the subalgebra ${\mathcal LC}_{\exp}(G)$
defines by Lemma \ref{H-exp(G/G_e)=LC-exp(G)} a continuous seminorm $q$ on
${\mathcal O}_{\exp}(G/G_e)$,
$$
q(v)=p(v\circ\pi),
$$
and $q$ will be submultiplicative, like $p$. Hence by Lemma \ref{LM-supp-q},
the support of $q$ is finite:
$$
\card \supp(q)<\infty
$$
Being applied to seminorm $p$ this mean that there exists a finite family of
cosets $\{K_1,...,K_n\}\subseteq G/G_e$, for which
$$
p(1_{K_i})\ne 0
$$
while for the others $K\in G/G_e$, $K\notin\{K_1,...,K_n\}$,
 \beq\label{p(1_K)=0}
p(1_K)=0
 \eeq
(here $1_K$ denotes the image of the constant identity $u(x)=1$ under the
projection \eqref{u_K}).

As a corollary, for any function $u\in {\mathcal O}_{\exp}(G)$ and for any
coset $K\notin \{K_1,...,K_n\}$ we have
 \beq\label{p(u_K)=0}
p(u_K)=0
 \eeq
(since $p(u_K)=p(1_K\cdot u_K)\le p(1_K)\cdot p(u_K)=0\cdot p(u_K)=0$).

Denote now by $P$ and $S$ the projections to the spaces consisting of functions
vanishing outside and inside $K_1\cup ...\cup K_n$:
$$
P(u)(x)=\left\{\begin{matrix}u(x),& x\in K_1\cup ...\cup K_n \\ 0,& x\notin
K_1\cup ...\cup K_n\end{matrix}\right\}=u_{K_1}+...+u_{K_n}
 $$
 $$
S(u)(x)=\left\{\begin{matrix}u(x),& x\notin K_1\cup ...\cup K_n \\ 0,& x\in
K_1\cup ...\cup K_n\end{matrix}\right\}=\sum_{K\notin \{K_1,...,K_n\}} u_K
$$
(the latter series converges in ${\mathcal O}_{\exp}(G)$, by Lemma
\ref{LM-sum-u_K}). From \eqref{p(u_K)=0} we have
 $$
\forall u\in {\mathcal O}_{\exp}(G)\qquad p(S(u))=p\left(\sum_{K\notin
\{K_1,...,K_n\}} u_K\right)\le \sum_{K\notin \{K_1,...,K_n\}}
p(u_K)=\eqref{p(u_K)=0}=\sum_{K\notin \{K_1,...,K_n\}} 0=0,
 $$
what implies
 \beq\label{p=p-circ-rho}
p=p\circ P
 \eeq
(on the one hand, $ p(u)=p(P(u)+S(u))\le p(P(u))+p(S(u))= p(P(u))$, and on the
other hand, $p(P(u))=p(u-S(u))\le p(u)+p(S(u))= p(u)$).

Denote by $p_i$ the (submultiplicative) seminorms on ${\mathcal
O}_{\exp}(K_i)$, induced by $p$,
$$
p_i(v)=p(v),\quad v\in {\mathcal O}_{\exp}(K_i)
$$
By Lemma \ref{LM-H-exp-K-AM=H-K} and Proposition \ref{kriterij-A-M}, these
seminorms are majorized by some seminorms $\widetilde{p_i}$ and ${\mathcal
O}(K_i)$:
$$
p_i(v)\le \widetilde{p_i}(v),\quad v\in {\mathcal O}_{\exp}(K_i)
$$
From \eqref{p=p-circ-rho} we have the following estimation:
$$
p(u)=p\Big(\sum_{i=1}^n u_{K_i}+\sum_{K\notin \{K_1,...,K_n\}}
u_{K}\Big)=\eqref{p=p-circ-rho}=p\Big(\sum_{i=1}^n u_{K_i}\Big)\le \sum_{i=1}^n
p_i(u_{K_i})\le \sum_{i=1}^n \widetilde{p_i}(u_{K_i})
$$
Thus, our initial seminorm $p$ on ${\mathcal O}_{\exp}(G)$ is majorized by the
seminorm $\sum_{i=1}^n \widetilde{p_i}$ on ${\mathcal O}(G)$. Again applying
Proposition \ref{kriterij-A-M}, we obtain that the inclusion ${\mathcal
O}_{\exp}(G)\subseteq {\mathcal O}(G)$ is an Arens-Michael envelope.
 \epr

\subsection{Holomorphic reflexivity}\label{holom-reflexivnost}

We can now declare the following a result of our considerations. For a
compactly generated Stein group $G$ with the algebraic connected component of
identity two algebras of those we considered above, namely,  ${\mathcal
O}^\star(G)$ and ${\mathcal O}_{\exp}(G)$, have the following curious property:
every such algebra $H$, being a rigid stereotype Hopf algebra, has an
Arens-Michael envelope $H^\heartsuit$, which also has a structure of rigid
stereotype Hopf algebra, and
 \bit{
\item[(i)] the natural homomorphism
$$
\heartsuit_H:H\to H^\heartsuit
$$
is a homomorphism of rigid Hopf algebras, and

\item[(ii)] the dual mapping
$$
(\heartsuit_H)^\star:(H^\heartsuit)^\star\to H^\star
$$
is an Arens-Michael envelope of the algebra $(H^\heartsuit)^\star$:
$$
(\heartsuit_H)^\star=\heartsuit_{(H^\heartsuit)^\star}
$$
 }\eit
Let us note the following in view of this:

\bprop\label{PROP:edinstvennost-Hopfa-na-H^heartsuit} For an arbitrary rigid
stereotype Hopf algebra $H$ the structure of rigid Hopf algebra on the
Arens-Michael envelope $H^\heartsuit$, satisfying conditions (i) and (ii), if
exists, is unique. \eprop
 \bpr
Note first that (i) and (ii) immediately imply
 \bit
\item[(iii)] the mappings $\heartsuit_H$ and $(\heartsuit_H)^\star$ are
bimorphisms of stereotype spaces (i.e. are injective and have dense image in
the range).
 \eit
Indeed, the mappings $\heartsuit_H:H\to H^\heartsuit$ and
$(\heartsuit_H)^\star:(H^\heartsuit)^\star\to H^\star$ are epimorphisms (i.e.
have dense image), since they are Arens-Michael envelopes. On the other hand
they are dual to each other, so they must be monomorphisms (i.e. are
injective).

This implies everything. First, the multiplication and the unit on
$H^\heartsuit$ are defined uniquely by the condition that $\heartsuit_H:H\to
H^\heartsuit$ is the Arens-Michael envelope of $H$. Consider the dual mapping
$(\heartsuit_H)^\star:(H^\heartsuit)^\star\to H^\star$. Like $\heartsuit_H$,
this must be a homomorphism of Hopf algebras. Hence, it is a homomorphism of
algebras, and at the same time an injective mapping, by virtue of (iii). This
means that the multiplication and the unit in $(H^\heartsuit)^\star$ are
defined uniquely since they are unduced from $H^\star$.

Thus, conditions (i) and (ii) impose rigid conditions on multiplication, unit,
comultiplication and counit in $H^\heartsuit$, and allow to define no more than
one structure of bialgebra on $H^\heartsuit$. On the other hand, we know that
antipode, if exists is also unique, so the structure of Hopf algebra on
$H^\heartsuit$ is also unique.
 \epr

It is convenient to sketch out conditions (i) and (ii) by the diagram,
 \beq\label{diagramma-refleksivnosti}
 \xymatrix @R=1.pc @C=1.pc
 {
 H
 & \ar@{|->}[r]^{\heartsuit} & &
 H^{\heartsuit}
 \\
 & & &
 \ar@{|->}[d]^{\star}
 \\
 \ar@{|->}[u]^{\star}
 & & &
 \\
 H^\star
 & &
 \ar@{|->}[l]_{\heartsuit}
 &
 (H^{\heartsuit})^\star
 }
 \eeq
with the following sense: first, in the corners of the square there are rigid
stereotype Hopf algebras, and the horizontal arrows (the Arens-Michael
operations $\heartsuit$) are their homomorphisms, and, second, the alternation
of the operations $\heartsuit$ and $\star$ (no matter which place you begin
with) at the fourth step returns back to the initial Hopf algebra (of course,
up to an isomorphism of functors).

The rigid stereotype Hopf algebras $H$, satisfying conditions (i) and (ii),
will be called {\it holomorphically reflexive}\index{holomorphically reflexive
Hopf algebra}, and the diagram \eqref{diagramma-refleksivnosti} for such
algebras the {\it reflexivity diagram}\index{reflexivity diagram}. The
justification of the term ``reflexivity'' in this case is the following: if we
denote by some symbol, say $\widehat{\ }\,\,$, the composition of operations
$\heartsuit$ and $\star$,
$$
\widehat H:=(H^\heartsuit)^\star
$$
and call such object a Hopf algebra {\it holomorphically dual to $H$}, then $H$
becomes naturally isomorphic to its second dual Hopf algebra:
 \beq\label{H-cong-(H^*)^*}
H\cong \widehat{\widehat H}
 \eeq
This is a corollary of Proposition
\eqref{PROP:edinstvennost-Hopfa-na-H^heartsuit}: since for holomorphically
reflexive Hopf algebras the passage $H\mapsto H^\heartsuit$ uniquely defines
the structure of Hopf algebra on $H^\heartsuit$, the isomorphism of algebras
$$
((H^\heartsuit)^\star)^\heartsuit\cong H^\star,
$$
postulating in axiom (ii), automatically must be an isomorphism of Hopf
algebras. The passage to the dual Hopf algebras exactly gives
\eqref{H-cong-(H^*)^*}.

Theorems \ref{TH-AM-O-star} and \ref{TH-AM-O} imply:

 \btm\label{dvoistvennost-H(G)...}
If $G$ is a Stein group with the algebraic connected component of identity,
then the algebras ${\mathcal O}^\star(G)$ and ${\mathcal O}_{\exp}(G)$ are
holomorphically reflexive, and the reflexivity diagram for them has the form:
 \beq\label{chetyrehugolnik-O-O*}
 \xymatrix @R=1.pc @C=2.pc
 {
 {\mathcal O}^\star(G)
 & \ar@{|->}[r]^{\heartsuit}_{\eqref{AM-O-star}} & &
 {\mathcal O}_{\exp}^\star(G)
 \\
 & & &
 \ar@{|->}[d]^{\star}
 \\
 \ar@{|->}[u]^{\star}
 & & &
 \\
 {\mathcal O}(G)
 & &
 \ar@{|->}[l]_{\heartsuit}^{\eqref{AM-O}}
 &
 {\mathcal O}_{\exp}(G)
 }
 \eeq
(the numbers under the horizontal arrows are references to the formulas in our
text).
 \etm

\bex For the group $\GL_n(\C)$ the reflexivity diagram
\eqref{chetyrehugolnik-O-O*} is as follows:
 \beq\label{chetyrehugolnik-O(GL)-O*(GL)}
 \xymatrix @R=1.pc @C=2.pc
 {
 {\mathcal O}^\star(\GL_n(\C))
 & \ar@{|->}[r]^{\heartsuit}_{\eqref{O(GL)=H_exp(GL)}} & &
 {\mathcal R}^\star(\GL_n(\C))
 \\
 & & &
 \ar@{|->}[d]^{\star}
 \\
 \ar@{|->}[u]^{\star}
 & & &
 \\
 {\mathcal O}(\GL_n(\C))
 & &
 \ar@{|->}[l]_{\heartsuit}^{\eqref{R(M)^heartsuit=O(M)}}
 &
 {\mathcal R}(\GL_n(\C))
 }
 \eeq
\eex

\section{Holomorphic reflexivity as a generalization of Pontryagin duality}
\label{SEC-hol-duality=gen-pontryagin}

\subsection{Pontryagin duality for compactly generated Stein groups}

The complex circle $\C^\times$ we were talking about in
\ref{SEC:stein-groups}\ref{SUBSEC-lin-groups}, occupies among all Abelian
compactly generated Stein groups the same place, as the usual ``real'' circle
$\T=\R / \Z$ among all locally compact Abelian groups (or among all Abelian
compactly generated real Lie groups), since for the Abelian compactly generated
Stein groups the following variant of Pontryagin's duality theory holds.

Let $G$ be an Abelian compactly generated Stein group. Let us call an arbitrary
holomorphic homomorphism of a Stein group $G$ into the complex circle
$$
\chi\in G^\bullet\quad\Longleftrightarrow\quad \chi:G\to\C^\times
$$
a {\it holomorphic character}\index{holomorphic character} on $G$. The set
$G^\bullet$ of all holomorphic characters on $G$ is a topological group with
respect to the pointwise multiplication and the topology of uniform convergence
on compact sets. The following theorem shows that the operation $G\mapsto
G^\bullet$ is analogous to the Pontryagin operation of passage to the dual
locally compact Abelian group:

\btm\label{PONT-AB} If $G$ is an Abelian compactly generated Stein group, then
its dual group $G^\bullet$ is again an Abelian compactly generated group and
the mapping
$$
\i_G:G\to G^{\bullet\bullet},\quad \i_G(x)(\chi)=\chi(x),\qquad x\in G,\;
\chi\in G^\bullet
$$
is an isomorphism of (topological groups and of) functors $G\mapsto G$ and
$G\mapsto G^{\bullet\bullet}$:
$$
G^{\bullet\bullet}\cong G
$$
\etm

In view of this fact we call $G^\bullet$ {\it dual complex
group}\index{dual!complex group} for the group $G$.

 \bpr
First we need to note that this is true for the special cases where $G=\C,
\;\C^\times,\; \Z$ and for the case of a finite Abelian group $G=F$. This is a
corollary of the following obvious formulas:
$$
\C^\bullet\cong \C,\quad (\C^\times)^\bullet\cong \Z,\quad \Z^\bullet\cong
\C^\times,\quad F^\bullet\cong F
$$
After that it remains to note that every compactly generated Stein group has
the form
$$
G\cong \C^l\times (\C^\times)^m\times \Z^n\times F\qquad (l,m,n\in\Z_+)
$$
so its dual group has the form
$$
G^\bullet\cong \C^l\times \Z^m\times (\C^\times)^n\times F\qquad (l,m,n\in\Z_+)
$$
So $G^\bullet$ is an Abelian compactly generated Stein group. The second dual
group $G^{\bullet\bullet}$ turns out to be isomorphic to $G$:
$$
G^{\bullet\bullet}\cong \C^l\times (\C^\times)^m\times \Z^n\times F\cong G
$$
 \epr

\subsection{Fourier transform as an Arens-Michael envelope}

If $G$ is an Abelian compactly generated Stein group, then every its
holomorphic character $\chi:G\to\C^\times$ is a holomorphic function on $G$. In
other words we can think of the dual complex group $G^\bullet$ as a subgroup in
the group of invertible elements of the algebra ${\mathcal O}(G)$ of
holomorphic functions on $G$:
$$
G^\bullet\subset{\mathcal O}(G)
$$
If we pass to dual objects by Theorem \ref{PONT-AB}, we obtain that the group
$G$ itself is included by the transformation $\i_G$ into the group of
invertible elements of the algebra ${\mathcal O}(G^\bullet)$ of holomorphic
functions on $G^\bullet$:
$$
\i_G:G\to G^{\bullet\bullet}\subset{\mathcal O}(G^\bullet).
$$
On the other hand, obviously, $G$ is included (through delta-functionals) into
algebra ${\mathcal O}^\star(G)$:
$$
\delta: G\to {\mathcal O}^\star(G)\qquad (x\mapsto\delta^x).
$$
By \cite[Theorem 10.12]{Akbarov} this implies that there exists a unique
homomorphism of stereotype algebras
$$
\sharp_G:{\mathcal O}^\star(G)\to {\mathcal O}(G^\bullet),
$$
such that the following diagram is commutative:
$$
\begin{diagram}
\node[2]{G} \arrow{se,t}{\i_G}\arrow{sw,t}{\delta}
\\
\node{{\mathcal O}^\star(G)}\arrow[2]{e,b,--}{\sharp_G} \node[2]{{\mathcal
O}(G^\bullet)}
\end{diagram}
$$
(this is the property of being group algebra for ${\mathcal O}^\star(G)$). It
is natural to call the homomorphism $\sharp_G:{\mathcal O}^\star(G)\to
{\mathcal O}(G^\bullet)$ the (inverse) {\it Fourier transform}\index{Fourier
transform} on the Stein group $G$, since it is defined by the same formula as
for the (inverse) Fourier transform for measures and distributions
\cite[31.2]{Hewitt-Ross}:
 \beq\label{alpha-chi=w-chi}
\overbrace{\alpha^\sharp(\chi)}^{\scriptsize \begin{matrix}
\text{value of the function $\alpha^\sharp\in \C^{G^\bullet}$}\\
\text{in the point $\chi\in G^\bullet$} \\ \downarrow \end{matrix}}\kern-35pt=
\kern-50pt\underbrace{\alpha(\chi)}_{\scriptsize \begin{matrix}\uparrow \\
\text{action of the functional $\alpha\in{\mathcal O}^\star(G)$}\\
\text{on the function $\chi\in G^\bullet\subseteq {\mathcal O}(G)$
}\end{matrix}} \kern-50pt ,\qquad \chi\in G^\bullet\qquad (\alpha\in {\mathcal
O}^\star(G),\quad w\in {\mathcal O}(G^\bullet))
 \eeq

\btm\label{O-exp^star=O} For any Abelian compactly generated Stein group $G$
its Fourier transform
$$
\sharp_G:{\mathcal O}^\star(G)\to {\mathcal O}(G^\bullet),
$$
is:
 \bit
\item[(a)] a homomorphism of rigid Hopf-Fr\'echet algebras, and

\item[(b)] an Arens-Michael envelope of the algebra ${\mathcal O}^\star(G)$.
 \eit
As a corollary the following isomorphisms of rigid Hopf-Fr\'echet algebras
hold:
 \beq\label{O_exp^*_G=O-G^bullet}
{\mathcal O}_{\exp}^\star(G)\cong \Big({\mathcal O}^\star(G)\Big)^\heartsuit
\cong {\mathcal O}(G^\bullet)
 \eeq
and the reflexivity diagram for $G$ takes the form
 \beq\label{chetyrehugolnik-O-O*-Abel}
 \xymatrix @R=1.pc @C=2.pc
 {
 {\mathcal O}^\star(G)
 &\ar@{|->}[rrr]^{\text{Fourier transform}}_{\eqref{alpha-chi=w-chi}} & &
 &&
 {\mathcal O}(G^\bullet)
 \\
 & & & &&
 \ar@{|->}[d]^{\star}
 \\
 \ar@{|->}[u]^{\star}
 & & & & &
 \\
 {\mathcal O}(G)
 & & & &
 \ar@{|->}[lll]_{\text{Fourier transform}}^{\eqref{alpha-chi=w-chi}} &
  {\mathcal O}^\star(G^\bullet)
 }
 \eeq
 \etm

Like Theorem \ref{PONT-AB}, this is proved by successive consideration of the
cases $G=\C, \;\C^\times,\; \Z$ and the case of arbitrary finite Abelian group
$G=F$. In the rest of this section up to the ``inclusion diagram'' we devote to
this.

\paragraph{Finite Abelian group.}\label{ex-finite}

As we told in \ref{SEC:stein-groups}\ref{SUBSEC-lin-groups}, every finite group
$G$ can be considered as a complex Lie group (of zero dimension), on which
every function is holomorphic. Moreover, in example
\ref{SEC-O_exp(G)}\ref{SUBSEC:algebra-O_exp(G)} we have noticed that every
function on $G$ has exponential type, so the algebras  ${\mathcal
O}_{\exp}(G)$, ${\mathcal O}(G)$ and $\C^G$ coincide:
$$
{\mathcal O}_{\exp}(G)={\mathcal O}(G)=\C^G
$$
If in addition $G$ is commutative, then the theorem \ref{O-exp^star=O} we
illustrate here turns into a formally more strong proposition:

\bprop If $G$ is a finite Abelian group, then the formula
\eqref{alpha-chi=w-chi} establishes an isomorphism of Hopf algebras:
 \beq\label{finite-groups}
{\mathcal O}_{\exp}^\star(G)={\mathcal O}^\star(G)=
\C_G\cong\C^{G^\bullet}={\mathcal O}(G^\bullet)= {\mathcal O}_{\exp}(G^\bullet)
 \eeq
 \eprop

\paragraph{Complex plane $\C$.}\label{ex-C} Let for every
$\lambda\in\C$ the symbol $\chi_\lambda$ denote a character on the group $\C$,
defined by formula:
$$
\chi_\lambda(t)=e^{\lambda\cdot t}
$$
The mapping $\lambda\in\C\mapsto\chi_\lambda\in\C^\bullet$ is an isomorphism of
complex groups
$$
\C\cong\C^\bullet
$$
and this isomorphism turns formula \eqref{alpha-chi=w-chi} into formula
 \beq\label{alpha-chi-lambda=w-lambda}
\alpha^\sharp(\lambda)=\alpha(\chi_\lambda),\qquad \lambda\in \C\qquad
(\alpha\in {\mathcal O}_{\exp}^\star(\C),\quad w\in {\mathcal O}(\C))
 \eeq
(we denote this isomorphism by the same symbol $\sharp$, although formally it
is a composition of mappings \eqref{alpha-chi=w-chi} and
$\lambda\mapsto\chi_\lambda$). As a result Theorem \ref{O-exp^star=O} being
applied to the group $G=\C$ is turned into

\bprop\label{PROP-H*(C)->H(C)} Formula \eqref{alpha-chi-lambda=w-lambda}
defines a homomorphism of stereotype Hopf algebras
$$
\sharp_{\C}:{\mathcal O}^\star(\C)\to {\mathcal O}(\C)
$$
which is an Arens-Michael envelope of the algebra ${\mathcal O}^\star(\C)$, and
establishes an isomorphism of Hopf-Fr\'echet algebras:
 \beq\label{O_exp*(C)=O(C)}
{\mathcal O}_{\exp}^\star(\C)\cong {\mathcal O}(\C)
 \eeq
\eprop

We shall need the following

\blm\label{LM-||alpha||_C-v-C} Seminorms of the form
 \beq\label{||alpha||_C-v-C}
||\alpha||_C=\sum_{k\in\N} |\alpha_k|\cdot C^k,\qquad C\ge 0,
 \eeq
(i.e. special case of seminorms \eqref{|alpha|_r-v-C}, when
$r_k=\frac{C^k}{k!}$) form a fundamental system in the set of all
submultiplicative continuous seminorms on ${\mathcal O}^\star(\C)$. \elm

\bpr As we had noted in \ref{SEC:stein-groups}\ref{Examples-of-Stein-groups},
the multiplication in ${\mathcal O}(\C)$ and in ${\mathcal O}^\star (\C)$ is
defined by the same formulas on series \eqref{umnozhenie-v-O(C)-i-O^star(C)}.
So we can say that the submultiplicativity of seminorms \eqref{||alpha||_C-v-C}
is already proven, since in Example \ref{EX-polunormy-v-H(C)} we had proven the
same fact for seminorms \eqref{polunormy-v-H(C)}, defined by the same formula
on the series.

Let us show that seminorms \eqref{||alpha||_C-v-C} form a fundamental system
among all submultiplicative continuous seminorms on ${\mathcal O}^\star(\C)$.
This is done like in Proposition \ref{PROP:O*(C)-Hopf}. Let $p$ be a
submultiplicative continuous seminorm:
$$
p(\alpha*\beta)\le p(\alpha)\cdot p(\beta)
$$
Put
$$
r_k=\frac{1}{k!}p(\zeta_k)
$$
Then
$$
(k+l)!\cdot r_{k+l}=p(\zeta_{k+l})=p(\zeta_k*\zeta_l)\le p(\zeta_k)\cdot
p(\zeta_l)=(k!\cdot r_k)\cdot (l!\cdot r_l)
$$
The sequence $A_k=r_k\cdot k!$ satisfies the recurrent inequality $A_{k+1}\le
A_k\cdot A_1$, which implies $A_k\le C^k$, for $C=A_1$. This in its turn
implies inequalities
 $$
r_k\le\frac{C^k}{k!}
 $$
Now using the same reasonings as in the proof of Proposition
\ref{PROP:O*(C)-Hopf}, we obtain:
$$
p(\alpha)\le\eqref{p(alpha)-le-|||alpha|||_r}\le |||\alpha|||_r=\sum_{k\in\N}
r_k\cdot |\alpha_k|\cdot k!\le \sum_{k\in\N} |\alpha_k|\cdot C^k=||\alpha||_C
$$ \epr

\bpr[Proof of Proposition \ref{PROP-H*(C)->H(C)}] Note from the very beginning
that the mapping $\sharp_{\C}:{\mathcal O}^\star(\C)\to {\mathcal O}(\C)$,
defined by Formula \eqref{alpha-chi-lambda=w-lambda} is continuous: by the
continuity of the mapping $\lambda\in\C\mapsto \chi_\lambda\in {\mathcal
O}(\C)$, every compact set $T$ in $\C$ is turned into a compact set
$\{\chi_\lambda;\;\lambda\in T\}$ in ${\mathcal O}(\C)$, so if the net of
functionals $\alpha_i$ tends to zero in ${\mathcal O}^\star(\C)$, then for any
compact set $T$ in $\C$ we have
$$
\alpha_i^\sharp(\lambda)=\alpha_i(\chi_\lambda)\underset{\lambda\in
T}{\rightrightarrows} 0,\qquad i\to\infty
$$
Thus, the functions $\alpha_i^\sharp$ tend to zero in ${\mathcal O}(\C)$.

Further, note that the mapping $\sharp_{\C}$ turns the functionals $\zeta_n$
into functions $z^n$:
$$
(\zeta_n)^\sharp(\lambda)=\zeta_n(\chi_\lambda)=\left(\frac{\d^n}{\d
t^n}e^{\lambda t}\right)\Bigg|_{t=0}=\lambda^n=z^n(\lambda)
$$
From this and from the continuity of $\sharp_{\C}$ it follows that the mapping
act on functionals $\alpha$ as the substitution of monomials $\zeta_n$ in the
decomposition \eqref{razlozhenie-alpha-v-C} by monomials $z^n$:
$$
\alpha^\sharp=\left(\sum_{n=0}^\infty \alpha_n\cdot
\zeta_n\right)^\sharp=\sum_{n=0}^\infty \alpha_n\cdot
(\zeta_n)^\sharp=\sum_{n=0}^\infty \alpha_n\cdot z^n
$$
This immediately implies the rest.

1. First, the mapping $\sharp_{\C}:{\mathcal O}^\star(\C)\to {\mathcal O}(\C)$
is an Arens-Michael envelope, since by Lemma \ref{LM-||alpha||_C-v-C}, every
submultiplicative continuous seminorm on ${\mathcal O}^\star(\C)$ is majorized
by a seminorm of the form \eqref{||alpha||_C-v-C}, which in its turn can be
extended by the mapping $\sharp_{\C}$ to a seminorm \eqref{polunormy-v-H(C)} on
${\mathcal O}(\C)$.

2. Second, the mapping $\sharp_{\C}:{\mathcal O}^\star(\C)\to {\mathcal O}(\C)$
is a homomorphism of algebras, since by formulas
\eqref{razlozhenie-u-v-C}-\eqref{razlozhenie-alpha-v-C} these algebras can be
considered as algebras of power series, where the multiplication is defined by
usual formulas for power series \eqref{umnozhenie-v-O(C)-i-O^star(C)}, and
$\sharp_{\C}$ will be just inclusion of one algebra into another, more wide,
algebra.

3. To prove that the mapping $\sharp_{\C}:{\mathcal O}^\star(\C)\to {\mathcal
O}(\C)$ is an isomorphism of coalgebras, let us note that the dual mapping
$$
(\sharp_{\C})^\star:{\mathcal O}^\star(\C)\to({\mathcal
O}^\star(\C))^\star={\mathcal O}(\C)^{\star\star}
$$
coincides with $\sharp_{\C}$ up to the isomorphism $\i_{{\mathcal
O}(\C)}:{\mathcal O}(\C)\cong {\mathcal O}(\C)^{\star\star}$:
 \beq\label{sharp_C^star=i_H(C)-circ-sharp_C}
\begin{diagram}
\node{{\mathcal O}^\star(\C)}
\arrow[2]{e,t}{(\sharp_{\C})^\star}\arrow{se,b}{\sharp_{\C}} \node[2]{{\mathcal
O}(\C)^{\star\star}}
\\
\node[2]{{\mathcal O}(\C)}\arrow{ne,b}{\i_{{\mathcal O}(\C)}}
\end{diagram}
 \eeq
This follows from the formula
 \beq\label{sharp_C-delta_x}
\sharp_{\C}(\delta^x)=\chi_x
 \eeq
Indeed,
$$
\sharp_{\C}(\delta^x)(\lambda)=\delta^x(\chi_\lambda)=\chi_\lambda(x)=e^{\lambda
x}=\chi_x(\lambda)
$$
Now we obtain:
$$
\sharp^\star(\delta^a)(\delta^b)=\delta^a(\sharp_{\C}(\delta^b))=\eqref{sharp_C-delta_x}
=\delta^a(\chi_b)=\sharp_{\C}(\delta^a)(b)=\delta^b(\sharp_{\C}(\delta^a))=
=\i_{{\mathcal O}(\C)}(\sharp_{\C}(\delta^a))(\delta^b)=(\i_{{\mathcal
O}(\C)}\circ\; \sharp)(\delta^a)(\delta^b)
$$
This is true for all $a,b\in\C$. On the other hand, the linear hull of
delta-functionals is dense in ${\mathcal O}^\star$, so
$$
(\sharp_{\C})^\star=\i_{{\mathcal O}(\C)}\circ\;\sharp_{\C}
$$
i.e. the diagram \eqref{sharp_C^star=i_H(C)-circ-sharp_C} is commutative. But
we have already proved that $\sharp_{\C}$ is a homomorphism of algebras, and
for $\i_{{\mathcal O}(\C)}$ this is obvious. Hence, $(\sharp_{\C})^\star$ is
also a homomorphism of algebras, and this means that $\sharp_{\C}$ is a
homomorphism of coalgebras.

4. Now it remains to prove that $\sharp_{\C}$ preserves antipode:
$$
\sigma_{{\mathcal O}(\C)}(\chi_\lambda)(x)=\chi_\lambda(-x)=e^{-\lambda
x}=(e^{\lambda x})^{-1}=\chi_\lambda(x)^{-1}
$$
$$
\Downarrow
$$
$$
\sigma_{{\mathcal O}(\C)}(\chi_\lambda)=\chi_\lambda^{-1}
$$
$$
\Downarrow
$$
$$
(\sigma_{{\mathcal O}^\star(\C)}(\alpha))^\sharp(\lambda)= (\sigma_{{\mathcal
O}^\star(\C)}(\alpha))(\chi_\lambda)= (\alpha\circ\sigma_{{\mathcal
O}(\C)})(\chi_\lambda)= \alpha(\sigma_{{\mathcal O}(\C)}(\chi_\lambda))=
\alpha(\chi_\lambda^{-1})=\alpha^\sharp(-\lambda)=\sigma_{{\mathcal
O}(\C)}(\alpha^\sharp)(\lambda)
$$
$$
\Downarrow
$$
$$
(\sigma_{{\mathcal O}^\star(\C)}(\alpha))^\sharp=\sigma_{{\mathcal
O}(\C)}(\alpha^\sharp)
$$
\epr

\paragraph{Complex circle $\C^\times$.}\label{ex-C-x} Let for any
$n\in\Z$ the symbol $z^n$ denote the character on the group $\C^\times$,
defined by formula
 $$
z^n(t)=t^n
 $$
The mapping $n\in\Z\mapsto z^n\in(\C^\times)^\bullet$ is a homomorphism of
complex groups
 $$
\Z\cong(\C^\times)^\bullet
 $$
and formula \eqref{alpha-chi=w-chi} under this isomorphism takes the form
 \beq\label{alpha-chi-lambda=w-lambda-C-x}
\alpha^\sharp(n)=\alpha(z^n),\qquad n\in\Z\qquad (\alpha\in {\mathcal
O}_{\exp}^\star(\C^\times))
 \eeq
(like in the previous example we denote this mapping by the same symbol
$\sharp$, although formally this is a composition of mappings
\eqref{alpha-chi=w-chi} and $n\mapsto z^n$). As a result Theorem
\ref{O-exp^star=O}, being applied to group $G=\C^\times$ is turned into

\bprop\label{PROP-H*(C-x)->H(Z)} Formula \eqref{alpha-chi-lambda=w-lambda-C-x}
defines a homomorphism of stereotype Hopf algebras
$$
\sharp_{\C^\times}:{\mathcal O}^\star(\C^\times)\to {\mathcal O}(\Z)
$$
which is an Arens-Michael envelope of the algebra ${\mathcal
O}^\star(\C^\times)$, and establishes an isomorphism of Hopf-Fr\'echet algebras
 \beq\label{O_exp*(C^x)=O(Z)}
{\mathcal O}_{\exp}^\star(\C^\times)\cong {\mathcal O}(\Z)=\C^\Z
 \eeq
\eprop

We shall need

\blm\label{LM-||alpha||_C-v-C-x} The seminorms of the form
 \beq\label{||alpha||_C-v-C-x}
||\alpha||_N=\sum_{|n|\le N}|\alpha_n|,\qquad N\in\N
 \eeq
-- i.e. the special case of seminorms \eqref{|alpha|_r-v-C-x}, when
$$
r_n=\begin{cases}1,& |n|\le N \\ 0,& |n|>N \end{cases}
$$
-- form a fundamental system in the set of all submultiplicative continuous
seminorms on ${\mathcal O}^\star(\C^\times)$.
 \elm
\bpr Submultiplicativity of seminorms \eqref{||alpha||_C-v-C-x} follows from
the formula for the operation of multiplication in ${\mathcal O}^\star
(\C^\times)$:
$$
||\alpha*\beta||_N=\sum_{|n|\le
N}|(\alpha*\beta)_n|=\eqref{umnozhenie-v-O(C-x)-i-O^star(C-x)}=\sum_{|n|\le
N}|\alpha_n\cdot\beta_n|\le \left(\sum_{|n|\le
N}|\alpha_n|\right)\cdot\left(\sum_{|n|\le
N}|\beta_n|\right)=||\alpha||_N\cdot||\beta||_N
$$
Let us show that seminorms \eqref{|alpha|_r-v-C-x} form a fundamental system
among all submultiplicative continuous seminorms on ${\mathcal
O}^\star(\C^\times)$. Let $p$ be a submultiplicative continuous seminorm:
$$
p(\alpha*\beta)\le p(\alpha)\cdot p(\beta)
$$
Put
$$
r_n=p(\zeta_n)
$$
Then
$$
r_n=p(\zeta_n)=p(\zeta_n*\zeta_n)\le p(\zeta_n)\cdot p(\zeta_n)=r_n^2
$$
i.e. $0\le r_n\le r_n^2$, hence $r_n\ge 1$, or $r_n=0$. But by Lemma
\ref{LM-polunormy-v-H*(C-x)-1}, the numbers $r_n$ must satisfy the condition
\eqref{forall-R>0-sum_r_n-R^n<infty-x}, which in its turn implies that $r_n\to
0$.  This is possible only if all those numbers, except maybe a finite
subfamily, vanishes:
$$
\exists N\in\N\quad \forall n\in\Z\quad |n|>N\quad\Longrightarrow\quad r_n=0
$$
Put $M=\max_{n}r_n$, then by Lemma \ref{LM-polunormy-v-H*(C-x)-1} we obtain:
$$
p(\alpha)\le |||\alpha|||_r=\sum_{n\in\Z} r_n\cdot |\alpha_n|= \sum_{|n|\le N}
r_n\cdot |\alpha_n|\le \sum_{|n|\le N} M\cdot|\alpha_n|=M\cdot ||\alpha||_N
$$ \epr

\bpr[Beginning of the proof of Proposition \ref{PROP-H*(C-x)->H(Z)}] Note that
the mapping $\sharp_{\C^\times}:{\mathcal O}^\star(\C^\times)\to {\mathcal
O}(\Z)$, defined by Formula \eqref{alpha-chi-lambda=w-lambda-C-x} is
continuous: if a net of functionals $\alpha_i$ tends to zero in ${\mathcal
O}^\star(\C^\times)$, then for any $n\in\Z$ we have
$$
\alpha_i^\sharp(n)=\alpha_i(z^n)\longrightarrow 0,\qquad i\to\infty
$$
This means that $\alpha_i^\sharp$ tends to zero in ${\mathcal O}(\Z)=\C^{\Z}$.

Further, let us note that the mapping $\sharp_{\C^\times}$ turns functionals
$\zeta_k$ into characteristic functions of singletons in $\Z$:
$$
(\zeta_k)^\sharp(n)=\zeta_k(z^n)=\eqref{zeta^k(z^n)}=\left\{\begin{matrix}1,& n=k \\
0,& n\ne k\end{matrix}\right\}=\eqref{DEF:1_n,n-in-Z}=1_k(n)
$$
From this and from the continuity of $\sharp_{\C^\times}$ it follows that this
mapping acts on functionals $\alpha$ as substitution of monomials $\zeta_n$ in
the decomposition \eqref{razlozhenie-alpha-v-C} by monomials $1_n$:
$$
\alpha^\sharp=\left(\sum_{n\in\Z} \alpha_n\cdot
\zeta_n\right)^\sharp=\sum_{n\in\Z} \alpha_n\cdot
(\zeta_n)^\sharp=\sum_{n\in\Z} \alpha_n\cdot 1_n
$$
This in its turn imply the most part of Proposition \ref{PROP-H*(C-x)->H(Z)}.

1. First, the mapping $\sharp_{\C^\times}:{\mathcal O}^\star(\C^\times)\to
{\mathcal O}(\Z)$ is an Arens-Michael envelope, so by Lemma
\ref{LM-||alpha||_C-v-C-x}, every submultiplicative continuous seminorm on
${\mathcal O}^\star(\C^\times)$ is majorized by a seminorm of the form
\eqref{||alpha||_C-v-C-x}, which in its turn can be extended by the mapping
$\sharp_{\C}$ to a continuous seminorm on ${\mathcal O}(\Z)=\C^{\Z}$.

2. Second, the mapping $\sharp_{\C^\times}:{\mathcal O}^\star(\C^\times)\to
{\mathcal O}(\Z)$ is a homomorphism of
algebras\label{homomorphism-sharp-C^times}, since by the formula for
multiplication \eqref{umnozhenie-v-O(C-x)-i-O^star(C-x)} in ${\mathcal
O}^\star(\C^\times)$, this mapping can be represented as a space of two-sided
sequences $\alpha_n$ with the coordinate-wise multiplication, which is included
in the wider space ${\mathcal O}(\Z)=\C^{\Z}$ of all two-sided sequences with
the coordinate-wise multiplication by the mapping $\sharp_{\C^\times}$.

3. we have to postpone to the next example (at page
\pageref{end-of-proof-PROP-H*(C-x)->H(Z)}) the proof of the fact that the
mapping $\sharp_{\C^\times}:{\mathcal O}^\star(\C^\times)\to {\mathcal O}(\Z)$
is an isomorphism of coalgebras.

4. Let us check that $\sharp_{\C^\times}$ preserves antipode:
$$
\sigma_{{\mathcal O}(\C^\times)}(z^n)(x)=z^n(x^{-1})=x^{-n}=z^{-n}(x)
$$
$$
\Downarrow
$$
$$
\sigma_{{\mathcal O}(\C^\times)}(z^n)=z^{-n}
$$
$$
\Downarrow
$$
$$
(\sigma_{{\mathcal O}^\star(\C^\times)}(\alpha))^\sharp(n)= (\sigma_{{\mathcal
O}^\star(\C^\times)}(\alpha))(z^n)= (\alpha\circ\sigma_{{\mathcal
O}(\C^\times)})(z^n)=
\alpha(\sigma_{\C^\times}(z^n))=\alpha(z^{-n})=\alpha^\sharp(-n)=\sigma_{{\mathcal
O}(\Z)}(\alpha^\sharp)(n)
$$
$$
\Downarrow
$$
$$
(\sigma_{{\mathcal O}^\star(\C^\times)}(\alpha))^\sharp=\sigma_{{\mathcal
O}(\Z)}(\alpha^\sharp)
$$
\epr

\paragraph{Group of integers $\Z$.}\label{ex-Z} Let for any $t\in\C^\times$ the
symbol $\chi_t$ denote a character on the group $\Z$, defined by formula
 \beq\label{chi_t(n)=t^n}
\chi_t(n)=t^n
 \eeq
The mapping $t\in\C^\times\mapsto\chi_t\in\Z^\bullet$ is an isomorphism of
complex group
 $$
\C^\times\cong\Z^\bullet
 $$
and Formula \eqref{alpha-chi=w-chi} under this isomorphism takes the form
 \beq\label{alpha-chi-lambda=w-lambda-Z}
\alpha^\sharp(t)=\alpha(\chi_t),\qquad t\in\C^\times\qquad (\alpha\in {\mathcal
O}_{\exp}^\star(\Z))
 \eeq
(again we denote this mapping by the same symbol $\sharp$, although it is a
composition of the mappings \eqref{alpha-chi=w-chi} and $t\mapsto \chi_t$). As
a result, Theorem \ref{O-exp^star=O}, being applied to group $G=\Z$, is turned
into

\bprop\label{PROP-H*(Z)->H(C-x)} Formula \eqref{alpha-chi-lambda=w-lambda-Z}
defines a homomorphism of rigid stereotype Hopf algebras
$$
\sharp_{\Z}:{\mathcal O}^\star(\Z)\to {\mathcal O}(\C^\times)
$$
which is an Arens-Michael envelope of the algebra ${\mathcal O}^\star(\Z)$, and
establishes an isomorphism of Hopf-Fr\'echet algebras
 \beq\label{O_exp*(Z)=O(C^x)}
{\mathcal O}_{\exp}^\star(\Z)\cong {\mathcal O}(\C^\times)
 \eeq
\eprop

We need

\blm\label{LM-||alpha||_C-v-Z} The seminorms of the form
 \beq\label{||alpha||_C-v-Z}
||\alpha||_C=\sum_{n\in\Z} |\alpha_n|\cdot C^{|n|},\qquad C\ge 1
 \eeq
(i.e. special case of seminorms \eqref{|alpha|_r-v-Z}, when $r_n=C^{|n|}$) form
a fundamental system in the set of all submultiplicative continuous seminorms
on ${\mathcal O}^\star(\Z)$.
 \elm
\bpr We had already noted the submultiplicativity of seminorms
\eqref{||alpha||_C-v-Z} in Example \ref{EX-polunormy-v-O(C^x)}. Let us show
that they form a fundamental system. Let $p$ be a submultiplicative continuous
seminorm on  on ${\mathcal O}^\star(\Z)$:
$$
p(\alpha*\beta)\le p(\alpha)\cdot p(\beta)
$$
Put $r_n=p(\delta^n)$. Then
$$
r_{k+l}=p(\delta^{k+l})=p(\delta^k*\delta^l)\le p(\delta^k)\cdot
p(\delta^l)=r_k\cdot r_l
$$
From this recurrent formula it follows that
$$
r_n\le M\cdot C^{|n|}
$$
where $M=r_0$, $C=\max\{r_1;r_{-1}\}$, and now by Lemma
\ref{LM-polunormy-v-H*(Z)} we obtain:
$$
p(\alpha)\le\eqref{p(alpha)-le-|||alpha|||_r-O*(Z)}\le
|||\alpha|||_r=\sum_{n\in\Z} r_n\cdot |\alpha_n|\le \sum_{n\in\Z} M\cdot
C^{|n|}\cdot |\alpha_n|=M\cdot ||\alpha||_C
$$ \epr

\bpr[Proof of Proposition \ref{PROP-H*(Z)->H(C-x)}] First of all, let us note
that the mapping $\sharp_{\Z}:{\mathcal O}^\star(\Z)\to {\mathcal
O}(\C^\times)$, defined by formula \eqref{alpha-chi-lambda=w-lambda-Z} is
continuous: by continuity of the mapping $t\in\C^\times\mapsto \chi_t\in
{\mathcal O}(\Z)$, every compact set $T$ in $\C^\times$ is turned into a
compact set $\{\chi_t;\; t\in T\}$ in ${\mathcal O}(\Z)$, so if a net of
functionals $\alpha_i$ turns to zero in ${\mathcal O}^\star(\Z)$, then for any
compact set $T$ in $\C^\times$ we have
$$
\alpha_i^\sharp(t)=\alpha_i(\chi_t)\underset{\lambda\in T}{\rightrightarrows}
0,\qquad i\to\infty
$$
This means that $\alpha_i^\sharp$ tends to zero in ${\mathcal O}(\C^\times)$.

Further, let us note that the mapping $\sharp_{\Z}$ turns functionals
$\delta^n$ into monomials $z^n$:
$$
(\delta^n)^\sharp(t)=\delta^n(\chi_t)=\chi_t(n)=t^n=z^n(t)
$$
This together with continuity of $\sharp_{\Z}$ implies that this mapping acts
on functionals $\alpha$ as substitution of monomials $\delta^n$ in the
decomposition \eqref{razlozhenie-alpha-v-Z} by monomials $z^n$:
$$
\alpha^\sharp=\left(\sum_{n\in\Z} \alpha_n\cdot
\delta^n\right)^\sharp=\sum_{n\in\Z} \alpha_n\cdot
(\delta^n)^\sharp=\sum_{n\in\Z} \alpha_n\cdot z^n
$$
This implies the rest.

1. First, the mapping $\sharp_{\Z}:{\mathcal O}^\star(\Z)\to {\mathcal
O}(\C^\times)$ is an Arens-Michael envelope, so by Lemma
\ref{LM-||alpha||_C-v-Z}, every submultiplicative seminorm on ${\mathcal
O}^\star(\Z)$ is majorized by a seminorm of the form \eqref{||alpha||_C-v-Z},
which in its turn can be extended by the mapping $\sharp_{\Z}$ to a seminorm
\eqref{polunormy-v-H(C-x)} on ${\mathcal O}(\C^\times)$.

2. Second, the mapping $\sharp_{\Z}:{\mathcal O}^\star(\Z)\to {\mathcal
O}(\C^\times)$ is a homomorphism of algebras, since by
\eqref{umnozhenie-v-O(Z)-i-O^star(Z)}-\eqref{umnozhenie-v-O(C-x)-i-O^star(C-x)}
these are algebras of power series with the usual multiplication of power
series, and $\sharp_{\C}$ is simply inclusion of one algebra into another,
wider algebra.

3. To prove that the mapping $\sharp_{\Z}:{\mathcal O}^\star(\Z)\to {\mathcal
O}(\C^\times)$ is an isomorphism of coalgebras, let us note that the dual
mapping
$$
(\sharp_{\Z})^\star:{\mathcal O}^\star(\C^\times)\to({\mathcal
O}^\star(\Z))^\star={\mathcal O}(\Z)^{\star\star}
$$
up to the isomorphism $\i_{{\mathcal O}(\Z)}:{\mathcal O}(\Z)\cong {\mathcal
O}(\Z)^{\star\star}$ coincides with $\sharp_{\C^\times}$:
 \beq\label{sharp_Z^star=i_H(Z)-circ-sharp_C-x}
\begin{diagram}
\node{{\mathcal O}^\star(\C^\times)}
\arrow[2]{e,t}{(\sharp_{\Z})^\star}\arrow{se,b}{\sharp_{\C^\times}}
\node[2]{{\mathcal O}(\Z)^{\star\star}}
\\
\node[2]{{\mathcal O}(\Z)}\arrow{ne,b}{\i_{{\mathcal O}(\Z)}}
\end{diagram}
 \eeq
This follows from formula
 \beq\label{sharp_Z-C-x-delta_x}
\delta^t(\sharp_{\Z}(\delta^n))=t^n=\delta^n(\sharp_{\C^\times}(\delta^t)),\qquad
t\in\C^\times,\quad n\in\Z
 \eeq
Indeed,
$$
\delta^t(\sharp_{\Z}(\delta^n))=\sharp_{\Z}(\delta^n)(t)=\delta^n(\chi_t)=\chi_t(n)=t^n
$$
and
$$
\delta^n(\sharp_{\C^\times}(\delta^t))=\sharp_{\C^\times}(\delta^t)(n)=\delta^t(z^n)=z^n(t)=t^n
$$
Now we obtain: for $t\in\C^\times$ and $n\in\Z$
$$
(\sharp_{\Z})^\star(\delta^t)(\delta^n)=\delta^t(\sharp_{\Z}(\delta^n))=\eqref{sharp_Z-C-x-delta_x}
=\delta^n(\sharp_{\C^\times}(\delta^t))=\i_{{\mathcal
O}(\Z)}(\sharp_{\C^\times}(\delta^t))(\delta^n)=(\i_{{\mathcal O}(\Z)}\circ\;
\sharp_{\C^\times})(\delta^t)(\delta^n)
$$
This is true for any $t\in\C^\times$ and for any $n\in\Z$. On the other hand,
delta-functionals generate a dense subspace in ${\mathcal O}^\star$, hence
$$
(\sharp_{\Z})^\star=\i_{{\mathcal O}(\Z)}\circ\;\sharp_{\C^\times}
$$
i.e. diagram \eqref{sharp_Z^star=i_H(Z)-circ-sharp_C-x} is commutative. Now we
can note that in the first part of proof of Proposition
\ref{PROP-H*(C-x)->H(Z)} (p.\pageref{homomorphism-sharp-C^times}) we have
already verified that $\sharp_{\C^\times}$ is a homomorphism of algebras.
Certainly, the mapping $\i_{{\mathcal O}(\Z)}$ is also a homomorphism, so we
conclude that $(\sharp_{\Z})^\star$ is a homomorphism of algebras, and this
means that $\sharp_{\Z}$ is a homomorphism of coalgebras.

4. Now it remains to check that $\sharp_{\Z}$ preserves antipode:
$$
\sigma_{{\mathcal
O}(\Z)}(\chi_t)(n)=\chi_t(-n)=t^{-n}=(t^{-1})^n=\chi_{t^{-1}}(n)
$$
$$
\Downarrow
$$
$$
\sigma_{{\mathcal O}(\Z)}(\chi_t)=\chi_{t^{-1}}
$$
$$
\Downarrow
$$
$$
(\sigma_{{\mathcal O}^\star(\Z)}(\alpha))^\sharp(t)= (\sigma_{{\mathcal
O}^\star(\Z)}(\alpha))(\chi_t)= (\alpha\circ\sigma_{{\mathcal O}(\Z)})(\chi_t)=
\alpha(\sigma_{{\mathcal O}(\Z)}(\chi_t))=
\alpha(\chi_{t^{-1}})=\alpha^\sharp(t^{-1})=\sigma_{{\mathcal
O}(\C^\times)}(\alpha^\sharp)(t)
$$
$$
\Downarrow
$$
$$
(\sigma_{{\mathcal O}^\star(\Z)}(\alpha))^\sharp=\sigma_{{\mathcal
O}(\C^\times)}(\alpha^\sharp)
$$
\epr

\bpr[End of the proof of Proposition
\ref{PROP-H*(C-x)->H(Z)}]\label{end-of-proof-PROP-H*(C-x)->H(Z)} In Proposition
\ref{PROP-H*(C-x)->H(Z)} it remains to proof that the mapping
$\sharp_{\C^\times}:{\mathcal O}^\star(\C^\times)\to {\mathcal O}(\Z)$ is an
isomorphism of coalgebras. Note that the dual mapping
$$
(\sharp_{\C^\times})^\star:{\mathcal O}^\star(\Z)\to({\mathcal
O}^\star(\C^\times))^\star={\mathcal O}(\C^\times)^{\star\star}
$$
coincides with $\sharp_{\Z}$ up to the isomorphism $\i_{{\mathcal
O}(\C^\times)}:{\mathcal O}(\C^\times)\cong {\mathcal
O}(\C^\times)^{\star\star}$:
 \beq\label{sharp_C-x^star=i_H(C-x)-circ-sharp_Z}
\begin{diagram}
\node{{\mathcal O}^\star(\Z)}
\arrow[2]{e,t}{(\sharp_{\C^\times})^\star}\arrow{se,b}{\sharp_{\Z}}
\node[2]{{\mathcal O}(\C^\times)^{\star\star}}
\\
\node[2]{{\mathcal O}(\C^\times)}\arrow{ne,b}{\i_{{\mathcal O}(\C^\times)}}
\end{diagram}
 \eeq
This also follows from \eqref{sharp_Z-C-x-delta_x}: if $t\in\C^\times$ and
$n\in\Z$
$$
(\sharp_{\C^\times})^\star(\delta^n)(\delta^t)=\delta^n(\sharp_{\C^\times}(\delta^t))=\eqref{sharp_Z-C-x-delta_x}
=\delta^t(\sharp_{\Z}(\delta^n))=\i_{{\mathcal
O}(\C^\times)}(\sharp_{\Z}(\delta^n))(\delta^t)=(\i_{{\mathcal O}(\Z)}\circ\;
\sharp_{\C^\times})(\delta^n)(\delta^t)
$$
This is true for all $t\in\C^\times$ and $n\in\Z$, and the delta-functionals
form a dense subspace in ${\mathcal O}^\star$, so
$$
(\sharp_{\C^\times})^\star=\i_{{\mathcal O}(\C^\times)}\circ\;\sharp_{\Z}
$$
i.e. Diagram \eqref{sharp_C-x^star=i_H(C-x)-circ-sharp_Z} is commutative. Now
we can note that in proof of Proposition \ref{PROP-H*(Z)->H(C-x)} we have
already verified that $\sharp_{\Z}$ is a homomorphism of algebras. And for
$\i_{{\mathcal O}(\C^\times)}$ this is obvious, so we obtain that
$(\sharp_{\C^\times})^\star$ is also a homomorphism of algebras, and this means
that $\sharp_{\C^\times}$ is a homomorphism of coalgebras. \epr

\paragraph{Proof of Theorem \ref{O-exp^star=O}} Now formulas
\eqref{finite-groups}, \eqref{O_exp*(C)=O(C)}, \eqref{O_exp*(C^x)=O(Z)},
\eqref{O_exp*(Z)=O(C^x)} prove isomorphism \eqref{O_exp^*_G=O-G^bullet} for
that cases $G=\C,\; \C^\times,\; \Z$ and for the case of finite group $G$, so
it remains just to apply formulas \eqref{tenz-pr-O-exp} and
\eqref{tenz-pr-O-exp-star}: let us decompose an Abelian compactly generated
Stein group $G$ into a direct product
$$
G=\C^l\times(\C^\times)^m\times\Z^n\times F
$$
where $F$ is a finite group. Then we have:
 \begin{multline*}
{\mathcal O}_{\exp}^\star(G)={\mathcal
O}_{\exp}^\star(\C^l\times(\C^\times)^m\times\Z^n\times
F)=\eqref{tenz-pr-O-exp}=\\= {\mathcal O}_{\exp}^\star(\C)^{\odot l}\odot
{\mathcal O}_{\exp}^\star(\C^\times)^{\odot m}\odot {\mathcal
O}_{\exp}^\star(\Z)^{\odot n}\odot {\mathcal O}_{\exp}^\star(F)=\\= {\mathcal
O}(\C^\bullet)^{\odot l}\odot {\mathcal O}((\C^\times)^\bullet)^{\odot m}\odot
{\mathcal O}(\Z^\bullet)^{\odot n}\odot {\mathcal
O}(F^\bullet)=\\=\eqref{tenz-pr-O-exp}= {\mathcal
O}((\C^\bullet)^l\times((\C^\times)^\bullet)^m\times(\Z^\bullet)^n\times
F^\bullet)= {\mathcal O} (G^\bullet)
 \end{multline*}

\subsection{Inclusion diagram}

\btm the following construction is a generalization of the Pontryagin duality
theory from the category of Abelian compactly generated Stein groups to the
category of compactly generated Stein groups with the algebraic connected
component of identity, {\sf
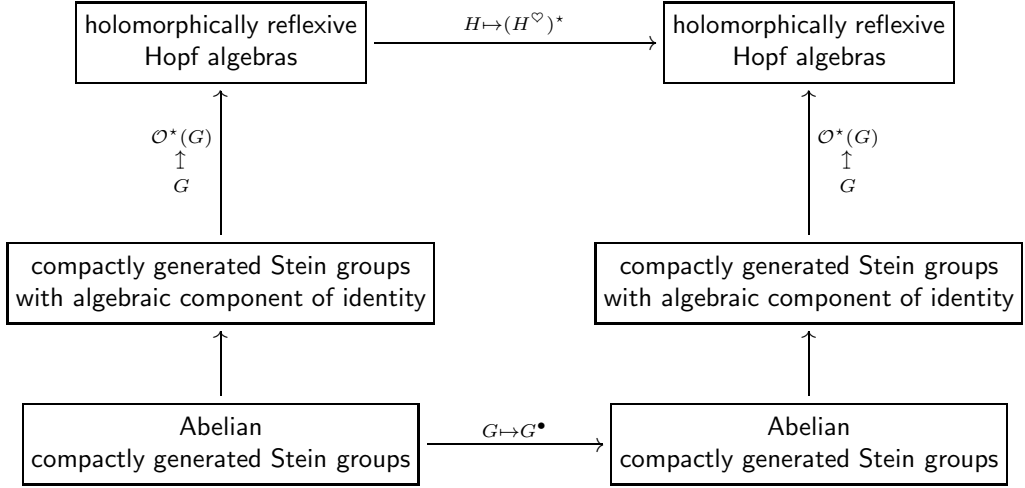
 $$
 \xymatrix
 {
 \boxed{\begin{matrix}
 \text{holomorphically reflexive}\\
 \text{Hopf algebras}
 \end{matrix}}
 \ar[rr]^{H\mapsto (H^\heartsuit)^\star} & &
 \boxed{\begin{matrix}
 \text{holomorphically reflexive}\\
 \text{Hopf algebras}
 \end{matrix}}
 \\ & & \\
 \boxed{\begin{matrix}
  \text{compactly generated Stein groups} \\
  \text{with algebraic component of identity}
 \end{matrix}} \ar[uu]^{\scriptsize\begin{matrix} {\mathcal O}^\star(G)\\
 \text{\rotatebox{90}{$\mapsto$}} \\ G\end{matrix}} & &
 \boxed{\begin{matrix}
  \text{compactly generated Stein groups} \\
  \text{with algebraic component of identity}
 \end{matrix}} \ar[uu]_{\scriptsize\begin{matrix} {\mathcal O}^\star(G)\\
 \text{\rotatebox{90}{$\mapsto$}} \\ G\end{matrix}} \\
 \boxed{\begin{matrix}
 \text{Abelian}\\
 \text{compactly generated Stein groups}
   \end{matrix}} \ar[u] \ar[rr]^{G\mapsto G^\bullet} & &
  \boxed{\begin{matrix}
 \text{Abelian}\\
 \text{compactly generated Stein groups}
\end{matrix}}\ar[u]
 }
 $$ }
and the commutativity of this diagram is established by the isomorphism of
functors
$$
{\mathcal O}^\star(G^\bullet)\cong \left(\Big({\mathcal
O}^\star(G)\Big)^\heartsuit\right)^\star.
$$
 \etm
\bpr In Theorem \ref{dvoistvennost-H(G)...} we have already showed that the
functor $G\mapsto {\mathcal O}^\star(G)$ acts from the category of compactly
generated Stein groups with algebraic component of identity into the category
of holomorphically reflexive rigid stereotype Hopf algebras. So we need only to
verify the commutativity of this ``categorical diagram''. This follows from
Theorem \ref{PONT-AB}: if $G$ is an Abelian compactly generated Stein group,
then the walk around the diagram gives the following objects:
$$
 \xymatrix
 {
 {\mathcal O}^\star (G) \ar@{|->}[r]^{\kern-40pt\heartsuit} & {\mathcal O}_{\exp}^\star (G)\cong
 \eqref{O_exp^*_G=O-G^bullet}\cong {\mathcal O} (G^\bullet) \ar@{|->}[r]^{\kern40pt\star} & {\mathcal O}^\star (G^\bullet) \\
 G \ar@{|->}[u]^{{\mathcal O}^\star} \ar@{|->}[rr]^{\bullet} & & G^\bullet
 \ar@{|->}[u]_{{\mathcal O}^\star}
 }
$$
\epr

\section{Appendix: holomorphic reflexivity of the quantum group `$az+b$'}\label{SEC:R_q(C^X-x-C)}

In this final section we, following our promises in Introduction, will show by
the example of the quantum group `$az+b$', that the holomorphic reflexivity
described above does not restrict itself on algebras of analytical functionals
${\mathcal O}^\star(G)$,  but lengthens into the theory of quantum groups.

\subsection{Quantum combinatorial formulas}

The theory of quantum groups has its own analog of elementary combinatorics,
used in the situations where the computations are applied to variables with the
following commutation law
 \beq\label{yx=qxy}
yx=qxy
 \eeq
where $q$ is a fixed number. For those computations in particular, the analogs
of usual binomial formulas are deduced. Some of them will be useful for us in
our constructions connected with `$az+b$', so we record them for further
references (we refer the reader to details in C.~Kassel's textbook
\cite{Kassel}).

For an arbitrary positive integer $n$ we put
 \beq\label{(n)!_q}
(n)_q:=1+q+...+q^{n-1}=\frac{q^n-1}{q-1},\qquad (n)!_q:=(1)_q (2)_q...(n)_q=
\frac{(q-1)(q^2-1)...(q^n-1)}{(q-1)^n}
 \eeq
The integer $(n)!_q$ will be called {\it quantum
factorial}\index{quantum!factorial} of the integer $n$. The {\it quantum
binomial coefficient}\index{quantum!binomial coefficient} is defined by formula
 \beq
\begin{pmatrix} n \\ k \end{pmatrix}_q:=\begin{cases}\frac{(n)!_q}{(k)!_q\cdot
(n-k)!_q}, & k\le n \\ 0, & k>n
\end{cases}
 \eeq

\btm[\bf quantum binomial formula] Let $x$ and $y$ be elements of an
associative algebra $A$ satisfying condition \eqref{yx=qxy}. Then for all
$n\in\N$
 \beq\label{q-binom}
(x+y)^n=\sum_{k=0}^n \begin{pmatrix} n \\ k \end{pmatrix}_q\cdot x^k\cdot
y^{n-k}
 \eeq
\etm

\btm[\bf quantum Chu-Wandermond formula] For all $l,m,n\in\N$
 \begin{multline}\label{C-W}
\begin{pmatrix} m+n \\ l \end{pmatrix}_q=\sum_{\max\{0,l-n\}\le i\le \min\{l,m\}}
q^{(m-i)\cdot(l-i)}\cdot
\begin{pmatrix} m \\ i \end{pmatrix}_q\cdot \begin{pmatrix} n \\ l-i \end{pmatrix}_q
=\\=\sum_{0\le i\le l} q^{(m-i)\cdot(l-i)}\cdot
\begin{pmatrix} m \\ i \end{pmatrix}_q\cdot \begin{pmatrix} n \\ l-i \end{pmatrix}_q
 \end{multline}
 \etm
\bpr\footnote{We give here the proof of the Chu-Wandermond formula just to draw
reader's attention to the limits of summing $\max\{0,l-n\}\le i\le
\min\{l,m\}$, which will be useful below.} If $l>m+n$, then for any $i=0,...,m$
we have $n<l-m\le l-i$, so $\begin{pmatrix} n
\\ l-i
\end{pmatrix}_q=0$. Thus, both sums in \eqref{C-W}
vanish. And the same happens with $\begin{pmatrix} m+n \\
l \end{pmatrix}_q$, so formula \eqref{C-W} is trivial.

Thus only the case of $l\le m+n$ is interesting. Consider the equality
$(x+y)^{m+n}=(x+y)^m(x+y)^n$. Removing the parenthesis in \eqref{q-binom} we
obtain
 \begin{multline*}
\sum_{l=0}^{m+n} \begin{pmatrix} m+n \\ l \end{pmatrix}_q\cdot x^l\cdot
y^{n-l}=\left(\sum_{i=0}^m \begin{pmatrix} m \\ i \end{pmatrix}_q\cdot x^i\cdot
y^{m-i}\right)\cdot\left(\sum_{j=0}^n \begin{pmatrix} n \\ j
\end{pmatrix}_q\cdot x^j\cdot y^{n-j}\right)=\\=
\sum_{i=0}^m\sum_{j=0}^n  \begin{pmatrix} m \\ i
\end{pmatrix}_q\cdot\begin{pmatrix} n \\ j \end{pmatrix}_q\cdot x^i\cdot
y^{m-i}\cdot  x^j\cdot y^{n-j}=\eqref{yx=qxy}= \sum_{i=0}^m\sum_{j=0}^n
\begin{pmatrix} m \\ i
\end{pmatrix}_q\cdot\begin{pmatrix} n \\ j \end{pmatrix}_q\cdot q^{(m-i)\cdot j}\cdot x^{i+j}\cdot
y^{m+n-i-j}
 \end{multline*}
Let us consider in the last sum the terms with indices $i$ and $j$ connected by
equality $i+j=l$. All those terms can be indexed by the parameter $i$, if we
express $j$ trough $i$ by formula $j=l-i$. We need only to note that to obtain
a bijection the index $i$ must vary in the following limits:
$$
\max\{0,l-n\}\le i\le \min\{l,m\}
$$
This follows from the restrictions on $i$ and $j$:
$$
\begin{cases}0\le i\le m \\ 0\le j\le n\end{cases}\quad\Longleftrightarrow\quad
\begin{cases}0\le i\le m \\ 0\le l-i\le n\end{cases}\quad\Longleftrightarrow\quad
\begin{cases}0\le i\le m \\ -n\le i-l\le 0\end{cases}\quad\Longleftrightarrow\quad
\begin{cases}0\le i\le m \\ l-n\le i\le l\end{cases}
$$
Now equating the coefficients at monomial $x^l\cdot y^{n-l}$, we obtain the
first equality in \eqref{C-W}:
$$
\begin{pmatrix} m+n \\ l \end{pmatrix}_q=\sum_{\max\{0,l-n\}\le i\le \min\{l,m\}}
q^{(m-i)\cdot(l-i)}\cdot
\begin{pmatrix} m \\ i \end{pmatrix}_q\cdot \begin{pmatrix} n \\ l-i \end{pmatrix}_q
$$
The second equality is evident, since for $i<l-n$ and $i>l$ we have
respectively $n<l-i$ or $l-i<0$, hence $\begin{pmatrix} n \\ l-i
\end{pmatrix}_q=0$, and the terms vanish:
$$
\sum_{\max\{0,l-n\}\le i\le \min\{l,m\}}= \underbrace{\sum_{0\le
i<\max\{0,l-n\}}}_{0}+\sum_{\max\{0,l-n\}\le i\le
\min\{l,m\}}+\underbrace{\sum_{\min\{l,m\}<i\le m}}_{0}=\sum_{0\le i\le m}
$$
\epr

\subsection{Hopf algebra of skew polynomials and similar constructions}

\paragraph{Tensor products $X\odot{\mathcal R}(\C)$, $X\circledast{\mathcal
R}(\C)$, $X\odot {\mathcal R}^\star(\C)$, $X\circledast {\mathcal
R}^\star(\C)$}

\btm\label{TH:R(C)-odot-X} Let $X$ be a stereotype space. Then
 \bit
 \item[---]
the elements of tensor product $X\odot {\mathcal R}(\C)$ are uniquely
represented as converging series
 \beq\label{predstavlenie-mnogochlenov-R(C)-odot-X}
u=\sum_{k\in\N} u_k\odot t^k,
 \eeq
with coefficients $u_k\in X$ continuously depending on $u\in X\odot {\mathcal
R}(\C)$,

 \item[---]
the elements of tensor product $X\circledast {\mathcal R}(\C)$  are uniquely
represented as converging series
 \beq\label{predstavlenie-mnogochlenov-R(C)-circledast-X}
u=\sum_{k\in\N} u_k\circledast t^k,
 \eeq
with coefficients $u_k\in X$ continuously depending on $u\in X\circledast
{\mathcal R}(\C)$,

 \item[---]
the elements of tensor product $X\odot{\mathcal R}^\star(\C)$ are uniquely
represented as converging series
 \beq\label{predstavlenie-mnogochlenov-R^star(C)-odot-X}
u=\sum_{k\in\N} u_k\odot\tau^k,
 \eeq
with coefficients $u_k\in X$ continuously depending on $u\in X\odot{\mathcal
R}^\star(\C)$,

 \item[---]
the elements of tensor product $X\circledast {\mathcal R}^\star(\C)$ are
uniquely represented as converging series
 \beq\label{predstavlenie-mnogochlenov-R^star(C)-circledast-X}
u=\sum_{k\in\N} u_k\circledast\tau^k,
 \eeq
with coefficients $u_k\in X$ continuously depending on $u\in X\circledast
{\mathcal R}^\star(\C)$,
 \eit
\etm

\paragraph{Algebras of skew polynomials $A\stackrel{\ph}{\odot}{\mathcal R}(\C)$
and skew power series $A\stackrel{\ph}{\circledast}{\mathcal R}^\star(\C)$}

\btm\label{TH-kosye-mnogochleny} Let $A$ be an injective stereotype algebra and
$\ph:A\to A$ its (continuous) automorphism. Then the formula
 \beq\label{umnozhenie-kosyh-mnogochlenov}
u\cdot v=\sum_{k\in\N} u_k\odot t^k\cdot\sum_{l\in\N} v_l\odot
t^l=\sum_{n\in\N} \sum_{i=0}^n
\underbrace{u_i\cdot\ph^i(v_{n-i})}_{\text{multiplication in $A$}}\odot t^n
 \eeq
defines an associative and continuous multiplication on the tensor product
$A\odot{\mathcal R}(\C)$, and turns it into an injective stereotype algebra
called {\rm algebra of skew polynomials (with respect to the automorphism
$\ph$) with coefficients in $A$} and is denoted by
$A\stackrel{\ph}{\odot}{\mathcal R}(\C)$. If in addition $A$ is a Brauner
algebra, then $A\stackrel{\ph}{\odot}{\mathcal R}(\C)$ is also a Brauner
algebra.
  \etm

 \btm\label{TH-kosye-step-ryady}
Let $A$ be a projective stereotype algebra and $\ph:A\to A$ its (continuous)
automorphism. Then the formula
 \beq\label{umnozhenie-kosyh-step-ryadov}
\alpha*\beta=\sum_{k\in\N} \alpha_k\circledast \tau^k * \sum_{l\in\N}
\beta_l\circledast \tau^l=\sum_{n\in\N} \sum_{i=0}^n
\underbrace{\alpha_i\cdot\ph^i(\beta_{n-i})}_{\text{multiplication in
$A$}}\circledast\tau^n
 \eeq
defines an associative and continuous multiplication on the tensor product
$A\circledast{\mathcal R}^\star(\C)$, and turns it into a projective stereotype
algebra called {\rm algebra of skew power series (with respect to the
automorphism $\ph$) with coefficients in $A$} and denoted by
$A\stackrel{\ph}{\circledast}{\mathcal R}^\star(\C)$. If in addition $A$ is a
Fr\'echet algebra, then $A\stackrel{\ph}{\circledast}{\mathcal R}^\star(\C)$ is
also a Fr\'echet algebra.
 \etm

\paragraph{Quantum pairs in stereotype Hopf algebras.}

Let $H$ ba an injective (resp., a projective) stereotype Hopf algebra and
suppose we have
 \bit
\item[--] a group-like central element $z$ in $H$

\item[--] a group-like central element $\omega$ in $H^\star$
 \eit
Suppose in addition that the operators $\M_{\omega}^\star$ and $\M_z^\star$
dual to the operators of multiplication by elements $\omega$ and $z$, act on
elements $z$ and $\omega$ as multiplication by a fixed number
$q\in\C\setminus\{0\}$:
 \beq\label{M_upsilon^star(z)=qz}
\M_{\omega}^\star(z)=q\cdot z,\qquad \M_z^\star(\omega)=q\cdot\omega
 \eeq
Then we call the pair $(z,\omega)$ a {\it quantum
pair}\index{quantum!pair}\label{quant-pair} in the Hopf algebra $H$ (with the
parameter $q$).

The key examples for us are the pairs of algebras $\langle {\mathcal
R}(\C^\times),{\mathcal R}^\star(\C^\times)\rangle$ and $\langle {\mathcal
O}(\C^\times),{\mathcal O}^\star(\C^\times)\rangle$ on the complex circle
$\C^\times$, we have considered in
\ref{SEC:stein-groups}\ref{SUBSEC:algebry-na-C^x}. There by symbol $z$ we
dented a monomial of power 1 on $\C^\times$:
$$
z(x):=x,\qquad x\in \C^\times
$$
Let us take a number $q\in\C^\times$ and consider the delta-functional
$\delta^q$ in the point $q$:
$$
\delta^q(u)=u(q)=\sum_{n\in\Z} u_n\cdot q^n,\qquad u\in{\mathcal R}(\C^\times)
$$
Certainly, this is a current on $\C^\times$, i.e. an element of the space
${\mathcal R}^\star(\C^\times)$. Its expansion in basis $\zeta_n$ has the form:
 \beq\label{delta^q=sum}
\delta^q=\sum_{n\in\Z} q^n\cdot \zeta_n
 \eeq

\bprop\label{quant-pair-z-d^q} Elements $(z,\delta^q)$ form a quantum pair in
the rigid stereotype Hopf algebras ${\mathcal R}(\C^\times)$ and ${\mathcal
O}(\C^\times)$ with the parameter $q$. \eprop

\bpr Let us prove this for ${\mathcal R}(\C^\times)$. We need to verify that
$z$ and $\delta^q$ are central and group-like elements. The first property is
trivial, since ${\mathcal R}(\C^\times)$ and ${\mathcal R}^\star(\C^\times)$
are commutative algebras. The fact that $z$ is a group-like element follows
from \eqref{koumn-v-R(C^times)}:
$$
\varkappa(z)=z\odot z
$$
And for $\delta^q$ this follows from multiplicativity of delta-functionals:
 \begin{multline*}
\langle u\odot v,\varkappa(\delta^q)\rangle= \langle u\cdot v,
\delta^q\rangle=(u\cdot v)(q)=u(q)\cdot v(q)=\langle u,\delta^q\rangle\cdot
\langle v,\delta^q\rangle=\langle u\odot v,\delta^q\circledast\delta^q\rangle
\\ \Longrightarrow\qquad
\varkappa(\delta^q)=\delta^q\circledast\delta^q.
 \end{multline*}
Finally the equalities \eqref{M_upsilon^star(z)=qz} are verified directly:
$$
\langle \M_{\delta^q}^\star z,\alpha\rangle=\langle  z,\delta^q*\alpha\rangle=
\langle
z,\sum_{m\in\Z}q^m\zeta_m*\sum_{n\in\Z}\alpha_n\zeta_n\rangle=\eqref{umnozhenie-v-R(C-x)-i-R^star(C-x)}=
\langle z,\sum_{n\in\Z}q^n\cdot \alpha_n\cdot\zeta_n\rangle=q\cdot\alpha_1=
q\cdot \langle z,\alpha\rangle
$$
and
 \begin{multline*}
\langle u,\M_z^\star \delta^q\rangle=\langle z\cdot u, {\delta^q}\rangle=
\left\langle \sum_{m\in\Z}u_m\cdot z^{m+1}, \sum_{n\in\Z}q^n\cdot
\zeta_n\right\rangle=\\=\left\langle \sum_{m\in\Z}u_m\cdot z^{m+1},
\sum_{l\in\Z}q^{l+1}\cdot \zeta_{l+1}\right\rangle=\sum_{m\in\Z}u_m\cdot
q^{m+1}=q\cdot \sum_{m\in\Z}u_m\cdot q^m=q\cdot \langle u,\delta^q\rangle
 \end{multline*}
 \epr

\paragraph{Hopf algebras $H\odot^z_\omega {\mathcal R}(\C)$ and $H^\star\circledast^\omega_z
{\mathcal R}^\star(\C)$.}

In the following theorem $\theta$ means the isomorphism of functors
\eqref{ABCD->ACBD}, and $(k)!_q$ the quantum factorial from \eqref{(n)!_q}:

\btm\label{TH-R-odot-H-i-R*-circledast-H*} Let $H$ be an injective stereotype
Hopf algebra and let $(z,\omega)$ be a quantum pair in $H$ with the parameter
$q\in\C^{\times}$. Then
 \bit
\item[(a)] the tensor product $H\odot {\mathcal R}(\C)$ has a unique structure
of injective Hopf algebra with the algebraic operations defined by formulas:
 \begin{align}
\label{umnozh-v-R(C)-odot-H} \text{multiplication:}&&& a\odot t^k\cdot b\odot
t^l=a\cdot (\M_{\omega}^\star)^k(b)\odot t^{k+l} &&\\
\label{edin-v-R(C)-odot-H} \text{unit:}&&& 1_{H\odot {\mathcal R}(\C)}=1_H\odot 1_{{\mathcal R}(\C)} &&\\
\label{koumn-v-R(C)-odot-H} \text{comultiplication:}&&& \varkappa(a\odot
t^k)=\sum_{i=0}^k\begin{pmatrix}k
\\ i\end{pmatrix}_q\cdot \theta\Big(\big(1_H\odot \M_z^i\big)(\varkappa_H(a))
\odot t^i\odot t^{k-i} \Big)=&& \\
\nonumber &&&\qquad =\sum_{(a)} \sum_{i=0}^k\begin{pmatrix}k
\\ i\end{pmatrix}_q\cdot a'\odot t^i\odot (z^i\cdot a'')\odot t^{k-i} &&\\
\label{koed-v-R(C)-odot-H} \text{counit:}&&& \e(a\odot
t^k)=\begin{cases}\e_H(a),& k=0\\ 0,&
k>0\end{cases} &&\\
\label{antipode-in-R(C)-odot-H} \text{antipode:}&&& \sigma(a\odot
t^k)=(-1)^k\cdot q^{-\frac{k(k+1)}{2}}\cdot z^{-k}\cdot
(\M_{\omega}^\star)^k(\sigma_H(a))\odot t^k &&
 \end{align}
$H\odot{\mathcal R}(\C)$ with such a structure of Hopf algebra is denoted by
$H\underset{\omega}{\stackrel{z}{\odot}} {\mathcal R}(\C)$; the common formula
for multiplication in this algebra has the form:
 \beq\label{obsh-umnozh-v-R(C)-odot-H}
u\cdot v=\left(\sum_{k\in\N}u_k\odot t^k\right)\cdot\left(\sum_{l\in\N}v_l\odot
t^l\right)=\sum_{m\in\N}\left(\sum_{k=0}^m u_k\cdot
(\M_{\omega}^\star)^k(v_{m-k})\right)\odot t^m
 \eeq

\item[(b)] the tensor product $H^\star\circledast {\mathcal R}^\star(\C)$ has a
unique structure of projective Hopf algebra with the algebraic operations
defined by formulas
 \begin{align}
\label{umn-v-R*(C)-*-H*} \text{multiplication:}&&& \alpha\circledast\tau^k
* \beta\circledast
 \tau^l=\alpha\cdot(\M_z^\star)^k(\beta) \circledast \tau^{k+l} &&\\
\label{edin-v-R*(C)-*-H*} \text{unit:}&&& 1_{H^\star\circledast{\mathcal
R}(\C)}=
1_{H^\star}\circledast 1_{{\mathcal R}^\star(\C)} &&\\
\label{koumn-v-R*(C)-*-H*}  \text{comultiplication:}&&&
\varkappa(\alpha\circledast \tau^k)=\sum_{i=0}^k\begin{pmatrix}k \\
i\end{pmatrix}_q\cdot
\theta\Big(\big(\id_{H^\star}\circledast\M_{\omega}^i\big)(\varkappa_{H^\star}(\alpha))
\circledast \tau^i \circledast \tau^{k-i}
\Big)= &&\\
\nonumber &&&\qquad =\sum_{(\alpha)} \sum_{i=0}^k\begin{pmatrix}k
\\ i\end{pmatrix}_q\cdot \alpha'\circledast\tau^i\circledast (\omega^i*\alpha'')\circledast \tau^{k-i} &&\\
\label{koed-v-R*(C)-*-H*}\text{counit:}&&& \e(\alpha\circledast\tau^k) =\begin{cases}\e_{H^\star}(\alpha),& k=0\\ 0,& k>0\end{cases} &&\\
\label{antipode-in-R*(C)-*-H*} \text{antipode:}&&&
\sigma(\alpha\circledast\tau^k)=(-1)^k\cdot q^{-\frac{k(k+1)}{2}}\cdot
\omega^{-k}* (\M_z^\star)^k(\sigma_{H^\star}(\alpha))\circledast\tau^k &&
 \end{align}
$H^\star\circledast{\mathcal R}^\star(\C)$ with such a structure of Hopf
algebra is denoted by
$H^\star\underset{z}{\stackrel{\omega}{\circledast}}{\mathcal R}^\star(\C)$;
the common formula for multiplication in this algebra has the form:
\beq\label{obsh-umnozh-v-R*(C)-*-H*}
\alpha*\beta=\left(\sum_{k\in\N}\alpha_k\odot
\tau^k\right)\cdot\left(\sum_{l\in\N}\beta_l\odot
\tau^l\right)=\sum_{m\in\N}\left(\sum_{k=0}^m \alpha_k*
(\M_z^\star)^k(\beta_{m-k})\right)\odot \tau^m
 \eeq

\item[(c)] the bilinear form
 \beq\label{<H-odot-R-C,H*-circledast-R*-C>}
\left\langle \sum_{k\in\N} u_k\odot t^k , \sum_{k\in\N} \alpha_k\circledast
\tau^k\right\rangle=\sum_{k\in\N}  \langle u_k,\alpha_k\rangle\cdot (k)!_q
 \eeq
turns $H\underset{\omega}{\stackrel{z}{\odot}}{\mathcal R}(\C)$ and
$H^\star\underset{z}{\stackrel{\omega}{\circledast}}{\mathcal R}^\star(\C)$
into a dual pair of stereotype Hopf algebras:
 \beq\label{(R-C-odot-H)*=R*-C-circledast-H*-Hopf}
\left(H\underset{\omega}{\stackrel{z}{\odot}}{\mathcal R}(\C)\right)^\star\cong
H^\star\underset{z}{\stackrel{\omega}{\circledast}}{\mathcal R}^\star(\C)
 \eeq
 \eit
 \etm

We divide the proof of this theorem in 7 lemmas. Some of them are evident, and
in those cases we omit the proof.

\blm The multiplication and the unit \eqref{umnozh-v-R(C)-odot-H},
\eqref{edin-v-R(C)-odot-H} endow $H\odot{\mathcal R}(\C)$ with the structure of
injective stereotype algebra, isomorphic to the algebra of skew polynomials
with coefficients in algebra $H$ and the generating automorphism
$$
\ph=\M_{\omega}^\star
$$
\elm

\blm The multiplication and the unit \eqref{umn-v-R*(C)-*-H*},
\eqref{edin-v-R*(C)-*-H*} endow $H^\star\circledast {\mathcal R}^\star(\C)$
with the structure of projective stereotype algebra, isomorphic to the algebra
of skew power series with coefficients in algebra $H^\star$ and the generating
automorphism
$$
\ph=\M_z^\star
$$
 \elm

\blm\label{LM:koalg-R(C)-odot-H} The bilinear form
\eqref{<H-odot-R-C,H*-circledast-R*-C>} turns the comultiplication
\eqref{koumn-v-R(C)-odot-H} into the multiplication \eqref{umn-v-R*(C)-*-H*},
and the counit \eqref{koed-v-R(C)-odot-H} into the counit
\eqref{edin-v-R*(C)-*-H*}:
 \beq\label{varkappa(u)->a*b}
\left\langle\varkappa (u),\alpha\circledast\beta\right\rangle= \left\langle
u,\alpha*\beta\right\rangle,\qquad \e(u)=\langle u, 1_{H^\star}\circledast
1_{{\mathcal R}^\star(\C)}\rangle
 \eeq
As a corollary, the comultiplication \eqref{koumn-v-R(C)-odot-H} and the counit
\eqref{koed-v-R(C)-odot-H} define the structure of injective stereotype
coalgebra on $H\odot{\mathcal R}(\C)$.\elm

\blm\label{LM:koalg-R*(C)-*-H*} The bilinear form
\eqref{<H-odot-R-C,H*-circledast-R*-C>} turns multiplication
\eqref{umnozh-v-R(C)-odot-H} into comultiplication \eqref{koumn-v-R*(C)-*-H*},
and unit \eqref{edin-v-R(C)-odot-H} into counit \eqref{koed-v-R*(C)-*-H*}:
 \beq\label{u.v->varkappa(a)}
\left\langle u\cdot v,\alpha\right\rangle= \left\langle u\odot
v,\varkappa(\alpha)\right\rangle,\qquad \langle 1_H\odot 1_{{\mathcal
R}(\C)},\alpha\rangle=\e(\alpha)
 \eeq
As a corollary, comultiplication \eqref{koumn-v-R*(C)-*-H*} and counit
\eqref{koed-v-R*(C)-*-H*} define the structure of projective stereotype
coalgebra on $H^\star\circledast{\mathcal R}^\star(\C)$. \elm

\bpr[Proof of Lemmas \ref{LM:koalg-R(C)-odot-H} and \ref{LM:koalg-R*(C)-*-H*}]
Because of the symmetry between formulas
\eqref{koumn-v-R(C)-odot-H}-\eqref{umn-v-R*(C)-*-H*} and
\eqref{umnozh-v-R(C)-odot-H}-\eqref{koumn-v-R*(C)-*-H*} it is sufficient to
prove \eqref{varkappa(u)->a*b}. For this we can take $u=a\odot t^k$. Then the
second equality becomes evident
$$
\e(a\odot t^k)=\left\{\begin{matrix}\e_H(a),& k=0\\ 0,&
k>0\end{matrix}\right\}=\langle a\odot t^k, 1_{H^\star}\circledast 1_{{\mathcal
R}^\star(\C)}\rangle
$$
and the first one is proved by the following chain:
 \begin{multline*}
\left\langle\varkappa (a\odot t^k),\alpha\circledast\beta\right\rangle=
\sum_{i=0}^k\begin{pmatrix} k \\ i\end{pmatrix}_q\cdot \left\langle
\theta\Big((1_H\odot \M_z^i)(\varkappa_H(a))\odot t^i\odot t^{k-i} \Big),
\sum_{l\in\N}\alpha_l\circledast\tau^l\circledast
\sum_{m\in\N}\beta_m\circledast\tau^m  \right\rangle=\\=
\sum_{l,m\in\N}\sum_{i=0}^k\begin{pmatrix} k \\ i\end{pmatrix}_q\cdot
\left\langle \theta\Big((1_H\odot \M_z^i)(\varkappa_H(a))\odot t^i\odot t^{k-i}
\Big), \alpha_l\circledast\tau^l\circledast \beta_m\circledast\tau^m
\right\rangle=\\= \sum_{l,m\in\N}\sum_{i=0}^k\begin{pmatrix} k \\
i\end{pmatrix}_q\cdot \left\langle (1_H\odot \M_z^i)(\varkappa_H(a))\odot
t^i\odot t^{k-i} , \theta\Big(\alpha_l\circledast\tau^l\circledast
\beta_m\circledast\tau^m\Big) \right\rangle=
 \\= \sum_{l,m\in\N}\sum_{i=0}^k\begin{pmatrix} k \\
i\end{pmatrix}_q\cdot \left\langle (1_H\odot \M_z^i)(\varkappa_H(a))\odot
t^i\odot t^{k-i} , \alpha_l\circledast \beta_m\circledast\tau^l\circledast
\tau^m \right\rangle=
 \\= \sum_{i=0}^k\begin{pmatrix} k \\
i\end{pmatrix}_q\cdot (i)!_q\cdot (k-i)!_q \left\langle (1_H\odot
\M_z^i)(\varkappa_H(a)), \alpha_i\circledast \beta_{k-i}\right\rangle=
 \\=
(k)!_q\cdot \sum_{i=0}^k  \left\langle \varkappa_H(a),
\alpha_i\circledast(\M_z^\star)^i(\beta_{k-i})\right\rangle= (k)!_q\cdot
\sum_{i=0}^k  \left\langle a,
\alpha_i*(\M_z^\star)^i(\beta_{k-i})\right\rangle=
 \\=
\left\langle a\odot t^k, \sum_{i=0}^k
\alpha_i*(\M_z^\star)^i(\beta_{k-i})\circledast\tau^k\right\rangle=
\left\langle a\odot t^k, \sum_{n\in\N}\left(\sum_{i=0}^n
\alpha_i*(\M_z^\star)^i(\beta_{n-i})\right)\circledast\tau^n\right\rangle=\\=
\left\langle a\odot t^k, \alpha*\beta\right\rangle
 \end{multline*}
\epr

\blm Comultiplication \eqref{koumn-v-R(C)-odot-H} and counit
\eqref{koed-v-R(C)-odot-H} are homomorphisms of injective stereotype algebras,
and as a corollary endow $H\odot {\mathcal R}(\C)$ with the structure of
injective stereotype bialgebra. \elm

\blm Comultiplication \eqref{koumn-v-R*(C)-*-H*} and counit
\eqref{koed-v-R*(C)-*-H*} are homomorphisms of projective stereotype algebras,
and as a corollary endow $H^\star\circledast {\mathcal R}^\star(\C)$ with the
structure of projective stereotype bialgebra.
  \elm

\bpr Again due to the symmetry of formulas it is sufficient here to prove the
first lemma. We will just check that the comultiplication
\eqref{koumn-v-R(C)-odot-H} is a homomorphism of algebras (and the reader is
supposed to check the identity for counit by analogy). Let us note the
following identities:
 \begin{align}
\label{k(1-odot-a-cdot-1-odot-b)} \varkappa(a\odot 1\cdot b\odot 1)&=\varkappa(a\odot 1)\cdot\varkappa(b\odot 1) \\
\label{k(t^k-odot-1-cdot-t^l-odot-1)} \varkappa(1\odot t^k\cdot 1\odot t^l)&=\varkappa(1\odot t^k)\cdot\varkappa(1\odot t^l) \\
\label{k(t^k-odot-1-cdot-1-odot-a)} \varkappa(1\odot t^k\cdot a\odot 1)&=\varkappa(1\odot t^k)\cdot\varkappa(a\odot 1) \\
\label{k(1-odot-a-cdot-t^l-odot-1)} \varkappa(a\odot 1\cdot 1\odot t^k) &=
\varkappa(a\odot 1)\cdot\varkappa(1\odot t^k)
 \end{align}
Indeed, for \eqref{k(1-odot-a-cdot-1-odot-b)} we have:
 \begin{multline*}
\varkappa(a\odot 1)\cdot\varkappa(b\odot 1)=\sum_{(a)} (a'\odot 1\odot a''\odot
1)\cdot \sum_{(b)} (b'\odot 1\odot b''\odot 1)=\\= \sum_{(a),(b)} a'b'\odot
1\odot a''b''\odot 1=\varkappa((a\cdot b)\odot 1)=\varkappa(a\odot 1\cdot
b\odot 1)
 \end{multline*}
For \eqref{k(t^k-odot-1-cdot-t^l-odot-1)}:
 \begin{multline*}
\varkappa(1\odot t^k)\cdot\varkappa(1\odot t^k)=\sum_{i=0}^k\begin{pmatrix}k
\\ i\end{pmatrix}_q\cdot 1\odot t^i\odot z^i\odot t^{k-i}\cdot
\sum_{j=0}^l\begin{pmatrix}l
\\ j\end{pmatrix}_q\cdot 1\odot t^j\odot z^j\odot t^{l-j}=\\=
\sum_{i=0}^k \sum_{j=0}^l \begin{pmatrix}k
\\ i\end{pmatrix}_q\cdot\begin{pmatrix}l
\\ j\end{pmatrix}_q\cdot \Big((1\odot t^i)\cdot(1\odot t^j)\Big) \odot\Big( (z^i\odot t^{k-i})\cdot
 (z^j\odot t^{l-j})\Big)=\eqref{umnozh-v-R(C)-odot-H}=\\=
\sum_{i=0}^k \sum_{j=0}^l \begin{pmatrix}k
\\ i\end{pmatrix}_q\cdot\begin{pmatrix}l
\\ j\end{pmatrix}_q\cdot \Big(1\odot t^{i+j}\Big)\odot\Big((z^i\cdot(\M_{\omega}^\star)^{k-i}(z^j)\odot t^{k+l-i-j}\Big)
 =\eqref{M_upsilon^star(z)=qz}
 =\\= \sum_{i=0}^k \sum_{j=0}^l \begin{pmatrix}k
\\ i\end{pmatrix}_q\cdot\begin{pmatrix}l
\\ j\end{pmatrix}_q\cdot q^{(k-i)j}\cdot 1\odot t^{i+j}\odot z^{i+j}
\odot  t^{k+l-i-j} =\\= {\scriptsize \begin{pmatrix} i+j=m \\ j=m-i \\ 0\le
m\le k+l
\\ \max\{0,m-l\}\le i\le \min\{k,m\}\end{pmatrix}}=\\=
 \sum_{m=0}^{k+l} \sum_{\max\{0,m-l\}\le i\le \min\{k,m\}}
\begin{pmatrix}k
\\ i\end{pmatrix}_q\cdot\begin{pmatrix}l
\\ m-i\end{pmatrix}_q\cdot q^{(k-i)(m-i)}
 \cdot 1\odot t^m\odot z^m\odot t^{k+l-m}=\eqref{C-W}=\\=
\sum_{m=0}^{k+l}\begin{pmatrix}k+l
\\ m\end{pmatrix}_q\cdot 1\odot t^m\odot z^m\odot t^{k+l-m}=
\varkappa(1\odot t^{k+l})=\varkappa(1\odot t^k\cdot 1\odot t^l)
 \end{multline*}
For \eqref{k(t^k-odot-1-cdot-1-odot-a)}:
 \begin{multline*}
\varkappa(1\odot t^k)\cdot\varkappa(a\odot 1)=\sum_{i=0}^k\begin{pmatrix}k
\\ i\end{pmatrix}_q\cdot 1\odot t^i\odot z^i\odot t^{k-i}\cdot \sum_{(a)} a'\odot 1\odot a''\odot
1=\\= \sum_{i=0}^k \sum_{(a)} \begin{pmatrix}k
\\ i\end{pmatrix}_q\cdot \Big( (1\odot t^i)\cdot (a'\odot 1)\Big)\odot \Big( (z^i\odot t^{k-i})\cdot (a''\odot
1)\Big)=\\=\sum_{i=0}^k\sum_{(a)} \begin{pmatrix}k
\\ i\end{pmatrix}_q\cdot (\M_{\omega}^\star)^i(a')\odot t^i \odot z^i(\M_{\omega}^\star)^{k-i}(a'')\odot
t^{k-i}=\eqref{(M_z*)^k}=\\=\varkappa( (\M_{\omega}^\star)^k(a)\odot
t^k)=\eqref{umnozh-v-R(C)-odot-H}= \varkappa(1\odot t^k \cdot a\odot 1)
 \end{multline*}
For \eqref{k(1-odot-a-cdot-t^l-odot-1)}:
 \begin{multline*}
\varkappa(a\odot 1)\cdot\varkappa(1\odot t^k)=\sum_{(a)} a'\odot 1\odot
a''\odot 1\cdot \sum_{i=0}^k\begin{pmatrix}k\\ i\end{pmatrix}_q\cdot 1\odot
t^i\odot z^i\odot t^{k-i}=\\= \sum_{(a)} \sum_{i=0}^k\begin{pmatrix}k\\
i\end{pmatrix}_q\cdot \Big((a'\odot 1)\cdot(1\odot t^i)\Big)\odot\Big((a''\odot
1)\cdot(z^i\odot t^{k-i})\Big)=\\=
\sum_{(a)} \sum_{i=0}^k\begin{pmatrix}k\\
i\end{pmatrix}_q\cdot a'\odot t^i\odot a'' z^i \odot
t^{k-i}=\eqref{koumn-v-R(C)-odot-H}= \varkappa(a \odot t^k)= \varkappa(a\odot 1
\cdot 1\odot t^k)
 \end{multline*}

From \eqref{k(1-odot-a-cdot-1-odot-b)}-\eqref{k(1-odot-a-cdot-t^l-odot-1)} it
follows that \eqref{koumn-v-R(C)-odot-H} is a homomorphism a algebras:
 \begin{multline*}
\varkappa(a\odot t^k \cdot b \odot t^l)=\varkappa(a\cdot
(\M_{\omega}^\star)^k(b) \odot
t^{k+l})=\varkappa(a\cdot(\M_{\omega}^\star)^k(b)\odot 1\cdot 1\odot
t^{k+l})=\eqref{k(1-odot-a-cdot-t^l-odot-1)}=\\=
\varkappa(a\cdot(\M_{\omega}^\star)^k(b)\odot 1)\cdot\varkappa(1\odot t^{k+l})
=\varkappa(a\odot 1\cdot (\M_{\omega}^\star)^k(b)\odot 1)\cdot \varkappa(1\odot
t^k\cdot 1\odot t^l)
=\\=\eqref{k(1-odot-a-cdot-1-odot-b)},\eqref{k(t^k-odot-1-cdot-t^l-odot-1)}=
\varkappa(a\odot 1)\cdot\varkappa((\M_{\omega}^\star)^k(b)\odot 1)\cdot
\varkappa(1\odot t^k)\cdot\varkappa(1\odot t^l)
=\\=\eqref{k(1-odot-a-cdot-t^l-odot-1)}= \varkappa(a\odot
1)\cdot\varkappa((\M_{\omega}^\star)^k(b)\odot 1\cdot 1\odot
t^k)\cdot\varkappa(1\odot t^l) =\\=
 \varkappa(a\odot
1)\cdot\varkappa(1\odot t^k\cdot b\odot 1)\cdot\varkappa(1\odot t^l) =\\=
\eqref{k(t^k-odot-1-cdot-1-odot-a)}=
 \varkappa(a\odot 1)\cdot\varkappa(1\odot t^k)\cdot\varkappa(b\odot 1)\cdot\varkappa(1\odot t^l) =
\eqref{k(1-odot-a-cdot-t^l-odot-1)}=\\=
 \varkappa(a\odot 1\cdot 1\odot t^k)\cdot\varkappa(b\odot 1\cdot 1\odot t^l) =
 \varkappa(a\odot t^k)\cdot\varkappa(b\odot t^l)
 \end{multline*}
\epr

\blm Formulas \eqref{antipode-in-R(C)-odot-H} and
\eqref{antipode-in-R*(C)-*-H*} define antipodes in bialgebras $H\odot {\mathcal
R}(\C)$ and $H^\star\circledast {\mathcal R}^\star(\C)$, dual to each other
with respect to the bilinear form \eqref{<H-odot-R-C,H*-circledast-R*-C>}:
 \beq\label{antipod*(H-o-R(C))=antipod(H*-*-R*(C))}
\langle\sigma(u),\alpha\rangle=\langle u,\sigma(\alpha)\rangle
 \eeq
  \elm

\bpr Let us show that formula \eqref{antipode-in-R(C)-odot-H} defines an
antipode in $H\odot {\mathcal R}(\C)$. First we need to verify that $\sigma$ is
an automorphism:
 \begin{multline*}
\sigma(a\odot t^k\cdot b\odot t^l)=\sigma(a\cdot(\M_{\omega}^\star)^k(b) \odot
t^{k+l})=\\= (-1)^{k+l}\cdot q^{-\frac{(k+l)(k+l+1)}{2}}\cdot z^{-k-l}\cdot
(\M_{\omega}^\star)^{k+l}\Bigg(\sigma_H\Big(a\cdot(\M_{\omega}^\star)^k(b)\Big)\Bigg)\odot
t^{k+l}=\\= (-1)^{k+l}\cdot q^{-\frac{(k+l)(k+l+1)}{2}}\cdot z^{-k-l}\cdot
(\M_{\omega}^\star)^{k+l}\Bigg(\sigma_H\Big((\M_{\omega}^\star)^k(b)\Big)\cdot\sigma_H(a)\Bigg)\odot
t^{k+l}=\\= (-1)^{k+l}\cdot q^{-\frac{(k+l)(k+l+1)}{2}}\cdot z^{-k-l}\cdot
(\M_{\omega}^\star)^l\Bigg(\Big((\M_{\omega}^\star)^k\circ\sigma_H\circ(\M_{\omega}^\star)^k\Big)(b)\Bigg)\cdot
(\M_{\omega}^\star)^{k+l}(\sigma_H(a))\odot t^{k+l}=\\=\eqref{M_z*-s-M_z*=s}=
(-1)^{k+l}\cdot q^{-\frac{k^2+l^2+k+l}{2}}\cdot q^{-kl}\cdot z^{-k-l}\cdot
(\M_{\omega}^\star)^l(\sigma_H(b))\cdot
(\M_{\omega}^\star)^{k+l}(\sigma_H(a))\odot t^{k+l}=\\= (-1)^{k+l}\cdot
q^{-\frac{k(k+1)+l(l+1)}{2}}\cdot z^{-l}\cdot
(\M_{\omega}^\star)^l(\sigma_H(b))\cdot
(\M_{\omega}^\star)^l\Big(z^{-k}(\M_{\omega}^\star)^k(\sigma_H(a))\Big)\odot
t^{k+l}=\\= (-1)^l\cdot q^{-\frac{l(l+1)}{2}}\cdot z^{-l}\cdot
(\M_{\omega}^\star)^l(\sigma_H(b))\odot t^l \cdot (-1)^k\cdot
q^{-\frac{k(k+1)}{2}}\cdot z^{-k}\cdot (\M_{\omega}^\star)^k(\sigma_H(a))\odot
t^k =\\=\sigma(b\odot t^l)\cdot\sigma(a\odot t^k)
 \end{multline*}
Now let us note that diagram \eqref{AX-antipode} becomes commutative if we put
there $a\odot 1$:
$$
 \xymatrix  @R=.3pc @C=.1pc
 {
 & \sum\limits_{(a)} a'\odot 1_{{\mathcal R}(\C)}\odot a''\odot 1_{{\mathcal R}(\C)} \ar@{|->}@/^4ex/[rr]^{\sigma\otimes 1_H} & & \sum\limits_{(a)} \sigma_H(a')\odot 1_{{\mathcal R}(\C)}\odot a''\odot 1_{{\mathcal R}(\C)} \ar@{|->}@/_3ex/[dr]_{\mu}& \\
 & & & & \sum\limits_{(a)} \sigma_H(a')\cdot a''\odot 1_{{\mathcal R}(\C)} \\
 & & & & \| \\
 a\odot 1 \ar@{|->}@/^2ex/[uuur]^{\varkappa}\ar@{|->}@/_2ex/[dddr]_{\varkappa}\ar@{|->}[rr]^{\e} & &  \e_H(a)\ar@{|->}[rr]^{\iota} & & \e_H(a)\cdot 1_{H}\odot 1_{{\mathcal R}(\C)}\\
 & & & & \| \\
 & & & & \sum\limits_{(a)} a'\cdot\sigma_H(a'')\odot 1_{{\mathcal R}(\C)} \\
 & \sum\limits_{(a)} a'\odot 1_{{\mathcal R}(\C)}\odot a''\odot 1_{{\mathcal R}(\C)}\ar@{|->}@/_4ex/[rr]_{1_H\otimes\sigma} & & \sum\limits_{(a)} a'\odot 1_{{\mathcal R}(\C)}\odot \sigma_H(a'')\odot 1_{{\mathcal R}(\C)} \ar@{|->}@/^3ex/[ur]^{\mu}&
 }
$$
or $1\odot t$:
$$
 \xymatrix  @R=.3pc @C=.1pc
 {
 & 1\odot 1\odot 1\odot t+1\odot t\odot z\odot 1 \ar@{|->}@/^4ex/[rr]^{\sigma\otimes 1_H} & &
 1\odot 1\odot 1\odot t-q^{-1} z^{-1}\odot t\odot z\odot 1 \ar@{|->}@/_3ex/[dr]_{\mu}& \\
 & & & & 1\odot t- q^{-1}q\cdot z^{-1}z\odot t \\
 & & & & \| \\
 1\odot t \ar@{|->}@/^2ex/[uuur]^{\varkappa}\ar@{|->}@/_2ex/[dddr]_{\varkappa}\ar@{|->}[rr]^{\e} & &  0\ar@{|->}[rr]^{\iota} & & 0\\
 & & & & \| \\
 & & & & -q^{-1} z^{-1}\odot t+q^{-1} \cdot z^{-1}\odot t  \\
 & 1\odot 1\odot 1\odot t+1\odot t\odot z\odot 1 \ar@{|->}@/_4ex/[rr]_{1_H\otimes\sigma} & &
 -q^{-1} 1\odot 1\odot z^{-1}\odot t+1\odot t\odot z^{-1}\odot 1 \ar@{|->}@/^3ex/[ur]^{\mu}&
 }
$$
By Lemma on antipode \ref{LM:antipod-x-y} this implies that the diagram
\eqref{AX-antipode} is commutative if we put there after various products of
$a\odot 1$ and $1\odot t$. In particular, if we put $a\odot t^k$. This means
that \ref{LM:antipod-x-y} is commutative (with arbitrary argument), and we
obtain that the mapping \eqref{antipode-in-R(C)-odot-H} is an antipode in
$H\odot {\mathcal R}(\C)$.

Further, by the symmetry of formulas, the mapping
\eqref{antipode-in-R*(C)-*-H*} is an antipode in $H^\star\circledast {\mathcal
R}^\star(\C)$.

It remains to verify that the bilinear form
\eqref{<H-odot-R-C,H*-circledast-R*-C>} turns the antipode
\eqref{antipode-in-R(C)-odot-H} into the antipode
\eqref{antipode-in-R*(C)-*-H*}. Clearly, formula
\eqref{antipod*(H-o-R(C))=antipod(H*-*-R*(C))} is equivalent to formula
$$
\langle\sigma(a\odot t^k),\alpha\circledast \tau^l\rangle= \langle a\odot
t^k,\sigma(\alpha\circledast \tau^l)\rangle
$$
For $k\ne l$ both sides here vanish, so we need check only the case $k=l$.
Indeed,
 \begin{multline*}
\langle\sigma(a\odot t^k),\alpha\circledast
\tau^k\rangle=\eqref{antipode-in-R(C)-odot-H}= \langle (-1)^k\cdot
q^{-\frac{k(k+1)}{2}}\cdot z^{-k}\cdot (\M_{\omega}^\star)^k(\sigma_H(a))\odot
t^k,\alpha\circledast \tau^k\rangle=\\= (-1)^k\cdot q^{-\frac{k(k+1)}{2}}\cdot
\langle z^{-k}\cdot (\M_{\omega}^\star)^k(\sigma_H(a)),\alpha \rangle\cdot
(k)!_q= (-1)^k\cdot q^{-\frac{k(k+1)}{2}}\cdot \langle
(\M_{\omega}^\star)^k(\sigma_H(a)),(\M_{z^{-1}}^\star)^k(\alpha) \rangle\cdot
(k)!_q=\\= (-1)^k\cdot q^{-\frac{k(k+1)}{2}}\cdot \langle
\sigma_H(a),\omega^k*(\M_{z^{-1}}^\star)^k(\alpha) \rangle\cdot (k)!_q =
(-1)^k\cdot q^{-\frac{k(k+1)}{2}}\cdot \langle
a,\sigma_H^\star\Big((\omega^k*(\M_{z^{-1}}^\star)^k(\alpha)\Big) \rangle\cdot
(k)!_q=\\= (-1)^k\cdot q^{-\frac{k(k+1)}{2}}\cdot \langle
a,\omega^{-k}*\sigma_H^\star\Big((\M_{z^{-1}}^\star)^k(\alpha)\Big)
\rangle\cdot (k)!_q=\eqref{M_z*-s-M_z*=s}=\\= (-1)^k\cdot
q^{-\frac{k(k+1)}{2}}\cdot \langle
a,\omega^{-k}*(\M_z^\star)^k(\sigma_H^\star(\alpha)) \rangle\cdot (k)!_q=
(-1)^k\cdot q^{-\frac{k(k+1)}{2}}\cdot \langle a\odot
t^k,\omega^{-k}*(\M_z^\star)^k(\sigma_H^\star(\alpha))\circledast \tau^k
\rangle=\\=\eqref{antipode-in-R*(C)-*-H*}=\langle a\odot
t^k,\sigma(\alpha\circledast \tau^k) \rangle
 \end{multline*}
\epr

\paragraph{Chains $H\odot^z_\omega {\mathcal R}(\C)\subset H\odot^z_\omega {\mathcal O}(\C)
\subset H\odot^z_\omega{\mathcal O}^\star(\C)\subset H\odot^z_\omega {\mathcal
R}^\star(\C)$ \\ and $H\circledast^z_\omega{\mathcal R}(\C)\subset
H\circledast^z_\omega{\mathcal O}(\C)\subset H\circledast^z_\omega {\mathcal
O}^\star(\C)\subset H\circledast^z_\omega{\mathcal R}^\star(\C)$.}

$\phantom{.}$

The same formulas and reasonings as we used in the proof of Theorem
\ref{TH-R-odot-H-i-R*-circledast-H*}, allow us to define, apart from
$H\underset{\omega}{\stackrel{z}{\odot}} {\mathcal R}(\C)$ and
$H^\star\underset{z}{\stackrel{\omega}{\circledast}}{\mathcal R}^\star(\C)$, a
series of similar stereotype Hopf algebras. These are algebras of
 \bit
\item[---] skew polynomials $H\underset{\omega}{\stackrel{z}{\odot}} {\mathcal
R}(\C)$ and $H\underset{\omega}{\stackrel{z}{\circledast}} {\mathcal R}(\C)$,

\item[---] skew entire functions $H\underset{\omega}{\stackrel{z}{\odot}}
{\mathcal O}(\C)$ and $H\underset{\omega}{\stackrel{z}{\circledast}} {\mathcal
O}(\C)$,

\item[---] skew analytic functionals $H\underset{\omega}{\stackrel{z}{\odot}}
{\mathcal O}^\star(\C)$ and $H\underset{\omega}{\stackrel{z}{\circledast}}
{\mathcal O}^\star(\C)$,

\item[---] skew power series $H\underset{\omega}{\stackrel{z}{\odot}} {\mathcal
R}^\star(\C)$ and $H\underset{\omega}{\stackrel{z}{\circledast}} {\mathcal
R}^\star(\C)$.
 \eit
Visually the connection between them can be illustrated by the following chains
of inclusions:
$$
H\underset{\omega}{\stackrel{z}{\odot}}{\mathcal R}(\C)\subset
H\underset{\omega}{\stackrel{z}{\odot}}{\mathcal O}(\C)\subset
H\underset{\omega}{\stackrel{z}{\odot}}{\mathcal O}^\star(\C)\subset
H\underset{\omega}{\stackrel{z}{\odot}}{\mathcal R}^\star(\C)
$$
and
$$
H\underset{\omega}{\stackrel{z}{\circledast}}{\mathcal R}(\C)\subset
H\underset{\omega}{\stackrel{z}{\circledast}}{\mathcal O}(\C)\subset
H\underset{\omega}{\stackrel{z}{\circledast}}{\mathcal O}^\star(\C)\subset
H\underset{\omega}{\stackrel{z}{\circledast}}{\mathcal R}^\star(\C)
$$
If $H$ is an injective Hopf algebra, then the upper chain is defined, if $H$ is
a projective Hopf algebra, then the lower chain is defined, and if $H$ is a
rigid stereotype Hopf algebra, then both chains are defined (and, certainly,
they coincide up to isomorphisms).

Theorem \ref{TH-R-odot-H-i-R*-circledast-H*} correctly define only the first
link in the first chain and the last link in the second chain. In all
conscience, to give accurate definition for all links we should formulate three
more analogous theorems.

To avoid those troubles we can either simply say that the other links are
defined by analogy (with replacing, if necessary, $\odot$ by $\circledast$, and
$\mathcal R$ by $\mathcal O$). Or we can unite all those four theorems (the one
already proven and three not yet formulated) into the following quite bulky
proposition:

\btm\label{TH-H-*-F-i-H*-o-F*} Let
 \bit
\item[---] $F$ denote one of the two Hopf algebras: ${\mathcal R}(\C)$ or
${\mathcal O}(\C)$,

\item[---] $H$ be an arbitrary injective stereotype Hopf algebra,

\item[---] $(z,\omega)$ be a quantum pair in $H$ with parameter
$q\in\C^{\times}$.
 \eit
Then\footnote{Here again $\theta$ is the isomorphism of functors
\eqref{ABCD->ACBD}, and $(k)!_q$ the quantum factorial defined in
\eqref{(n)!_q}.}
 \bit
\item[(a)] the tensor product $H\odot F$ possesses a unique structure of
injective Hopf algebra with algebraic operations, defined by formulas:
 \begin{align}
\label{umnozh-v-H-o-F} \text{multiplication:}&&& a\odot t^k\cdot b\odot
t^l=a\cdot (\M_{\omega}^\star)^k(b)\odot t^{k+l} &&\\
\label{edin-v-H-o-F} \text{unit:}&&& 1_{H\odot {\mathcal R}(\C)}=1_H\odot 1_{{\mathcal R}(\C)} &&\\
\label{koumn-v-H-o-F} \text{comultiplication:}&&& \varkappa(a\odot
t^k)=\sum_{i=0}^k\begin{pmatrix}k
\\ i\end{pmatrix}_q\cdot \theta\Big(\big(1_H\odot \M_z^i\big)(\varkappa_H(a))
\odot t^i\odot t^{k-i} \Big)=&& \\
\nonumber &&&\qquad =\sum_{(a)} \sum_{i=0}^k\begin{pmatrix}k
\\ i\end{pmatrix}_q\cdot a'\odot t^i\odot (z^i\cdot a'')\odot t^{k-i} &&\\
\label{koed-v-H-o-F} \text{counit:}&&& \e(a\odot t^k)=\begin{cases}\e_H(a),&
k=0\\ 0,&
k>0\end{cases} &&\\
\label{antipode-in-H-o-F} \text{antipode:}&&& \sigma(a\odot t^k)=(-1)^k\cdot
q^{-\frac{k(k+1)}{2}}\cdot z^{-k}\cdot (\M_{\omega}^\star)^k(\sigma_H(a))\odot
t^k &&
 \end{align}
$H\odot F$ with such a structure of Hopf algebra is denoted by
$H\underset{\omega}{\stackrel{z}{\odot}} F$;

\item[(b)] the tensor product $H^\star\circledast F^\star$ possesses a unique
structure of projective Hopf algebra with the algebraic operations, defined by
formulas:
 \begin{align}
\label{umn-v-H*-*-F*} \text{multiplication:}&&& \alpha\circledast\tau^k
* \beta\circledast
 \tau^l=\alpha\cdot(\M_z^\star)^k(\beta) \circledast \tau^{k+l} &&\\
\label{edin-v-H*-*-F*} \text{unit:}&&& 1_{H^\star\circledast{\mathcal R}(\C)}=
1_{H^\star}\circledast 1_{{\mathcal R}^\star(\C)} &&\\
\label{koumn-v-H*-*-F*}  \text{comultiplication:}&&&
\varkappa(\alpha\circledast \tau^k)=\sum_{i=0}^k\begin{pmatrix}k \\
i\end{pmatrix}_q\cdot
\theta\Big(\big(\id_{H^\star}\circledast\M_{\omega}^i\big)(\varkappa_{H^\star}(\alpha))
\circledast \tau^i \circledast \tau^{k-i}
\Big)= &&\\
\nonumber &&&\qquad =\sum_{(\alpha)} \sum_{i=0}^k\begin{pmatrix}k
\\ i\end{pmatrix}_q\cdot \alpha'\circledast\tau^i\circledast (\omega^i\cdot \alpha'')\circledast \tau^{k-i} &&\\
\label{koed-v-H*-*-F*}\text{counit:}&&& \e(\alpha\circledast\tau^k) =\begin{cases}\e_{H^\star}(\alpha),& k=0\\ 0,& k>0\end{cases} &&\\
\label{antipode-in-H*-*-F*} \text{antipode:}&&&
\sigma(\alpha\circledast\tau^k)=(-1)^k\cdot q^{-\frac{k(k+1)}{2}}\cdot
\omega^{-k}* (\M_z^\star)^k(\sigma_{H^\star}(\alpha))\circledast\tau^k &&
 \end{align}
$H^\star\circledast F^\star$ with such a structure of Hopf algebra is denoted
by $H^\star\underset{z}{\stackrel{\omega}{\circledast}} F^\star$;

\item[(c)] the bilinear form
 \beq\label{<H-o-F-C,H*-*-F*>}
\left\langle \sum_{k\in\N} u_k\odot t^k , \sum_{k\in\N} \alpha_k\circledast
\tau^k\right\rangle=\sum_{k\in\N}  \langle u_k,\alpha_k\rangle\cdot (k)!_q
 \eeq
turns $H\underset{\omega}{\stackrel{z}{\odot}} F$ and
$H^\star\underset{z}{\stackrel{\omega}{\circledast}}F^\star$ into dual pair of
stereotype Hopf algebras.
 \eit
If, under all other assumptions (except injectivity), $H$ is a projective
stereotype Hopf algebra, then
 \bit
\item[(a)] the tensor product $H\circledast F$ has a unique structure of
projective stereotype Hopf algebra with algebraic operations defined by
formulas \eqref{umnozh-v-H-o-F}-\eqref{antipode-in-H-o-F}, but with replacing
$\odot$ by $\circledast$; $H\circledast F$ with such a structure of Hopf
algebra is denoted by $H\underset{\omega}{\stackrel{z}{\circledast}}F$;

\item[(b)] the tensor product $H^\star\odot F^\star$ has a unique structure of
injective stereotype Hopf algebra with algebraic operations defined by formulas
\eqref{umn-v-H*-*-F*}-\eqref{antipode-in-H*-*-F*}, but with replacing
$\circledast$ by $\odot$; $H^\star\odot F^\star$ with such a structure of Hopf
algebra is denoted by $H^\star\underset{z}{\stackrel{\omega}{\odot}}F^\star$

\item[(c)] the bilinear form
 \beq\label{<H-*-F-C,H*-o-F*>}
\left\langle \sum_{k\in\N} u_k\circledast t^k , \sum_{k\in\N} \alpha_k\odot
\tau^k\right\rangle=\sum_{k\in\N}  \langle u_k,\alpha_k\rangle\cdot (k)!_q
 \eeq
turns $H\underset{\omega}{\stackrel{z}{\circledast}} F$ and
$H^\star\underset{z}{\stackrel{\omega}{\odot}}F^\star$ into a dual pair of
stereotype Hopf algebras.
 \eit
 \etm
As we have already told this is proved by analogy with Theorem
\ref{TH-R-odot-H-i-R*-circledast-H*}.

\bprop Let $H$ be an injective Hopf algebra, and $(z,\omega)$ a quantum pair in
$H$ with the parameter $q\in\C^{\times}$. Then the rules
$$
a\odot t^k\mapsto a\odot t^k\mapsto a\odot \tau^k\mapsto a\odot \tau^k \qquad
(k\in\N,\quad a\in H)
$$
uniquely define a chain of (continuous) homomorphisms of injective stereotype
Hopf algebras:
$$
H\underset{\omega}{\stackrel{z}{\odot}}{\mathcal R}(\C)\to
H\underset{\omega}{\stackrel{z}{\odot}}{\mathcal O}(\C)\to
H\underset{\omega}{\stackrel{z}{\odot}}{\mathcal O}(\C)^\star\to
H\underset{\omega}{\stackrel{z}{\odot}}{\mathcal R}(\C)^\star
$$
\eprop

\bprop Let $H$ be a projective Hopf algebra and $(z,\omega)$ a quantum pair in
$H$ with the parameter $q\in\C^{\times}$. Then the rules
$$
a\circledast t^k\mapsto a\circledast t^k\mapsto a\circledast \tau^k\mapsto
a\circledast \tau^k \qquad (k\in\N,\quad a\in H)
$$
uniquely define a chain of (continuous) homomorphisms of projective stereotype
Hopf algebras:
$$
H\underset{\omega}{\stackrel{z}{\circledast}}{\mathcal R}(\C)\to
H\underset{\omega}{\stackrel{z}{\circledast}}{\mathcal O}(\C)\to
H\underset{\omega}{\stackrel{z}{\circledast}}{\mathcal O}(\C)^\star\to
H\underset{\omega}{\stackrel{z}{\circledast}}{\mathcal R}(\C)^\star
$$
\eprop

\subsection{Quantum group $`az+b'={\mathcal
R}_q(\C^\times\ltimes\C)$}

Here we show that the quantum group `$az+b$' (defined in
\cite{Woronowicz,VanDaele,Wang,QSNG}) is a special case of the construction
described in Theorem \ref{TH-R-odot-H-i-R*-circledast-H*}.

\paragraph{The group $\C^\times\ltimes\C$ of affine transformation of a complex plane.}

The group of affine transformations of the complex plane, often denoted as
`$az+b$', from the algebraic point of view is a semidirect product
$\C^\times\ltimes\C$ of complex circle $\C^\times$ and complex plane $\C$,
where $\C^\times$ acts on $\C$ by usual multiplication. In other word,
$\C^\times\ltimes\C$ is a Cartesian product $\C^\times\times\C$ with algebraic
operations
 \begin{align*}
&\text{multiplication:} & (a,x)\cdot(b,y)&=(ab,xb+y) && (a,b\in\C^\times,\; x,y\in\C)\\
&\text{unit:} & 1_{\C^\times\ltimes\C}&=(1,0) && \\
&\text{inverse element:} & (a,x)^{-1}&=\left(\frac{1}{a},-\frac{x}{a}\right) &&
(a\in\C^\times,\; x\in\C)
 \end{align*}
Clearly, this is a connected Stein group. Moreover, $\C^\times\ltimes\C$ is an
algebraic group, since it can be represented as a linear group by matrices of
the form
$$
(a,x)=\begin{pmatrix}a & 0 \\ x & 1 \end{pmatrix}\qquad (a\in\C^\times,\;
x\in\C)
$$
(and multiplication, unit and inverse element become usual operations with
matrices).

\paragraph{Stereotype algebras ${\mathcal R}(\C^\times\ltimes\C)$ and ${\mathcal
R}^\star(\C^\times\ltimes\C)$.}\label{R(C^x-lx-C)-i-R*(C^x-lx-C)}

By symbol ${\mathcal R}(\C^\times\ltimes\C)$ we, as usual, denote the algebra
of polynomials on an algebraic group $\C^\times\ltimes\C$. According to the
general approach of \ref{SEC:stein-groups}\ref{R-O}, we endow the space
${\mathcal R}(\C^\times\ltimes\C)$ with the strongest locally convex topology.
The dual space of currents ${\mathcal R}^\star(\C^\times\ltimes\C)$ is an
algebra with respect to the usual convolution of functionals
\eqref{DEF:svertka}. Like any other dual space to a stereotype space, it is
endowed with the topology of uniform convergence on compact sets in ${\mathcal
R}(\C^\times\ltimes\C)$. In this case this is equivalent to the ${\mathcal
R}(\C^\times\ltimes\C)$-weak topology.

Recall that by $z^n$ and $t^k$ we denote basis monomials in the spaces
${\mathcal R}(\C^\times)$ and ${\mathcal R}(\C)$ (we defined them by formulas
\eqref{DEF:z^n-in-R(C^x)} and \eqref{def-t^k}). In accordance with the common
notation \eqref{g-h-na-S-T}, it is reasonable to denote the basis monomials in
the space of functions ${\mathcal R}^\star(\C^\times\ltimes\C)$ by symbol
$z^n\boxdot t^k$:
$$
(z^n\boxdot t^k)(a,x):=a^n\cdot x^k,\qquad a\in\C^\times,\; x\in\C
$$
Similarly, extending the old notations $\zeta_n$ and $\tau^k$ from
\eqref{DEF:zeta_n} and \eqref{def-tau^k}, we denote by $\zeta_n\boxedast\tau^k$
the functional on ${\mathcal R}(\C^\times\ltimes\C)$ of taking the $n$-th
coefficient of Laurent series with respect to the first variable and at the
same time the $k$-th derivative in the point $(1,0)$ with respect to the second
variable:
 \beq
\zeta_n\boxedast \tau^k(u)=\int_0^1 e^{-2\pi i n t}\frac{\d^k}{\d x^k}
u(x,e^{2\pi i t})\Big|_{x=0} \d t=\frac{\d^k}{\d x^k} \int_0^1 e^{-2\pi i n t}
u(x,e^{2\pi i t}) \d t \Big|_{x=0}
 \eeq

\bprop\label{bazis-v-R(C^x-lx-C)} 1) The functions $\{z^n\boxdot t^k;\,
n\in\Z,\, k\in\N\}$ form an algebraic basis in the space ${\mathcal
R}(\C^\times\ltimes\C)$ of polynomials on $\C^\times\ltimes\C$: every
polynomial $u\in {\mathcal R}(\C^\times\ltimes\C)$ is uniquely decomposed into
the series
 \beq\label{u-na-C-times-C^times}
u=\sum_{k\in\N, n\in\Z} u_{n,k}\cdot z^n\boxdot t^k, \qquad \card\{(n,k):\;
u_{n,k}\ne 0\}<\infty,
 \eeq
where the coefficients can be computed by formula
 \beq
u_{n,k}=\frac{1}{k!}\cdot\zeta_n\boxedast \tau^k(u)
 \eeq
The correspondence $u\leftrightarrow\{u_{n,k};\, n\in\Z,\,k\in\N\}$ establishes
an isomorphism of topological vector spaces
$$
{\mathcal R}(\C^\times\ltimes\C)\cong \C_{\Z\times\N}
$$
2) The functionals $\{\zeta_n\boxedast \tau^k ;\, n\in\Z,\, k\in\N\}$ form a
basis in the stereotype space ${\mathcal R}^\star(\C^\times\ltimes\C)$: every
functional $\alpha\in{\mathcal R}^\star(\C^\times\ltimes\C)$ is uniquely
decomposed into a (converging in the space ${\mathcal
R}^\star(\C^\times\ltimes\C)$) series
 \beq\label{potok-na-C-times-C^times}
\alpha=\sum_{k\in\N, n\in\Z}\alpha_{n,k}\cdot\zeta_n\boxedast \tau^k,
 \eeq
where the coefficients can be computed by formula
 \beq
\alpha_{n,k}=\frac{1}{k!}\cdot\alpha(z^n\boxdot t^k)
 \eeq
The correspondence $\alpha\leftrightarrow\{\alpha_{n,k};\, n\in\Z,k\in\N\}$
establishes an isomorphism of topological vector spaces
$$
{\mathcal R}^\star(\C^\times\ltimes\C)\cong \C^{\Z\times\N}
$$
3) The bases $\{z^n\boxdot t^k;\, n\in\Z,\, k\in\N\}$ and $\{\zeta_n\boxedast
\tau^k;\, n\in\Z,\, k\in\N\}$ are dual to each other up to the constant $k!$
 \beq
\langle z^m\boxdot t^k,\zeta_n*\tau^l\rangle=\langle
z^m,\zeta_n\rangle\cdot\langle t^k,\tau^l\rangle =\begin{cases}0, &
(m,k)\ne(n,l)
\\ k!,& (m,k)=(n,l)\end{cases}
 \eeq
and the action of functionals $\alpha\in{\mathcal R}^\star(\C^\times\ltimes\C)$
on polynomials $u\in {\mathcal R}(\C^\times\ltimes\C)$ is described by the
formula
 \beq
\langle u,\alpha\rangle=\sum_{n\in\Z, k\in\N}u_{n,k}\cdot\alpha_{n,k}\cdot k!
 \eeq
 \eprop

\brem The functional $\zeta_n\boxedast\tau^k$ can be represented as the
convolution of two components $\zeta_n$ and $\tau^k$ (and the order of their
multiplication becomes important here): if we denote by $Z_n$ the functional of
taking the $n$-th coefficient of Laurent series with respect to the first
variable in the point $(1,0)$,
 \beq
Z_n(u)=\zeta_n\boxedast\tau^0(u)=\frac{1}{2\pi i}\int_{|z|=1}
\frac{u(z,0)}{z^{n+1}}\d z=\int_0^1 e^{-2\pi i n t}u(e^{2\pi i t},0)\d t,\qquad
u\in {\mathcal R}(\C^\times\ltimes\C),
 \eeq
and by $T^k$ the functional of taking the $k$-th derivative with respect to the
second variable in the point $(1,0)$,
 \beq
T^k(u)=\zeta_0\boxedast\tau^k(u)=\frac{\d^k}{\d x^k}u(1,x)\Big|_{x=0},\qquad
u\in {\mathcal R}(\C^\times\ltimes\C),
 \eeq
then the following equalities hold:
 \beq
Z_n * T^k=\zeta_n\boxedast\tau^k,\qquad T^k* Z_n=Z_{n-k} *
T^k=\zeta_{n-k}\boxedast\tau^k
 \eeq
For the proof we can use formulas \eqref{DEF:svertka-s-delta^a}: first,
$$
\delta^a*T^k=a\cdot T^k
$$
$$
\Downarrow
$$
 \begin{multline}\label{tau^k*-delta-0-e^2pi-i-t}
(\delta^{(e^{2\pi i t},0)}*T^k)(u)=((e^{2\pi i t},0)\cdot
T^k)(u)=T^k(u\cdot(e^{2\pi i t},0))=\\=\frac{\d^k}{\d x^k}u\Big((e^{2\pi i
t},0)\cdot(1,x)\Big)\Big|_{x=0}=\frac{\d^k}{\d x^k}u(e^{2\pi i t},x)\Big|_{x=0}
 \end{multline}
$$
\Downarrow
$$
 \begin{multline*}
(Z_n*T^k)(u)=\left(\int_0^1 e^{-2\pi i n t}\cdot \delta^{(e^{2\pi i t},0)}\d
t*T^k\right)(u)=\left(\int_0^1 e^{-2\pi i n t}\cdot \delta^{(e^{2\pi i
t},0)*T^k}\d t\right)(u)=\\=\int_0^1 \left(e^{-2\pi i n t}\cdot\delta^{(e^{2\pi
i t},0)}*T^k\right)(u)\d t=\int_0^1 e^{-2\pi i n t}\cdot\frac{\d^k}{\d
x^k}\,u(e^{2\pi i t},x)\Big|_{x=0}\d t
 \end{multline*}
And, second,
$$
T^k*\delta^a=T^k\cdot a
$$
$$
\Downarrow
$$
 \begin{multline*}
(T^k*\delta^{(e^{2\pi i t},0)})(u)=(T^k\cdot(e^{2\pi i t},0))(u)=T^k((e^{2\pi i
t},0)\cdot u)=\frac{\d^k}{\d x^k}u\Big((1,x)\cdot(e^{2\pi i
t},0)\Big)\Big|_{x=0}=\\=\frac{\d^k}{\d x^k}u(e^{2\pi i t},x\cdot e^{2\pi i
t})\Big|_{x=0}=(e^{2\pi i t})^k\cdot\frac{\d^k}{\d y^k}u(e^{2\pi i
t},y)\Big|_{y=0}=\eqref{tau^k*-delta-0-e^2pi-i-t}=e^{2\pi i k
t}\cdot(\delta^{(e^{2\pi i t},0)}*T^k)(u)
 \end{multline*}
$$
\Downarrow
$$
 \beq\label{tau^k*delta^(e^2pi i t,0)}
T^k*\delta^{(e^{2\pi i t},0)}=e^{2\pi i k t}\cdot \delta^{(e^{2\pi i t},0)}*T^k
 \eeq
$$
\Downarrow
$$
 \begin{multline*}
T^k*Z_n=T^k*\int_0^1 e^{-2\pi i n t}\cdot \delta^{(e^{2\pi i t},0)}\d
t=\int_0^1 e^{-2\pi i n t}\cdot T^k*\delta^{(e^{2\pi i t},0)}\d
t=\eqref{tau^k*delta^(e^2pi i t,0)}=\\= \int_0^1 e^{-2\pi i n t}\cdot e^{2\pi i
k t}\cdot \delta^{(e^{2\pi i t},0)}*T^k\d t=\int_0^1 e^{-2\pi i (n-k) t}\cdot
\delta^{(e^{2\pi i t},0)}*T^k\d t=\\=\int_0^1 e^{-2\pi i (n-k) t}\cdot
\delta^{(e^{2\pi i t},0)}\d t*T^k=Z_{n-k}*T^k
 \end{multline*}
\erem

\paragraph{Stereotype algebras ${\mathcal O}(\C^\times\ltimes\C)$ and
${\mathcal O}^\star(\C^\times\ltimes\C)$.}

As always, we endow the algebra ${\mathcal O}(\C^\times\ltimes\C)$ of
holomorphic functions on $\C^\times\ltimes\C$ with the topology of uniform
convergence on compact sets in $\C^\times\ltimes\C$. Its dual algebra
${\mathcal O}^\star(\C^\times\ltimes\C)$ is endowed with the topology of
uniform convergence on compact sets in ${\mathcal O}(\C^\times\ltimes\C)$, and
the multiplication there is the usual convolution \eqref{DEF:svertka}.

\bprop\label{bazis-v-O(C^x-lx-C)} 1) The functions $\{z^n\boxdot t^k;\,
n\in\Z,\, k\in\N\}$ form a basis in the stereotype space ${\mathcal
O}(\C^\times\ltimes\C)$ of holomorphic functions on $\C^\times\ltimes\C$: every
function $u\in {\mathcal O}(\C^\times\ltimes\C)$ can be uniquely represented as
a sum of a (converging in ${\mathcal O}(\C^\times\ltimes\C)$) series
 \beq\label{u-v-O(C-times-C^times)}
u=\sum_{k\in\N, n\in\Z} u_{n,k}\cdot z^n\boxdot t^k, \qquad \forall C>0\quad
\sum_{k\in\N, n\in\Z} |u_{k,n}|\cdot C^{k+|n|}<\infty,
 \eeq
where the coefficients (continuously depend on $u$ and) can be computed by
formula
 \beq
u_{n,k}=\frac{1}{k!}\cdot\zeta_n\boxedast \tau^k(u)
 \eeq
The topology of ${\mathcal O}(\C^\times\ltimes\C)$ can be described by
seminorms:
 \beq\label{polunormy-v-O(C-rtimes-C-x)}
||u||_C=\sum_{k\in\N, n\in\Z} |u_{k,n}|\cdot C^{k+|n|},\qquad C\ge 1
 \eeq

2) The functionals $\{\zeta_n\boxedast \tau^k ;\, n\in\Z,\, k\in\N\}$ form a
basis in the space ${\mathcal O}^\star(\C^\times\ltimes\C)$: every functional
$\alpha\in{\mathcal O}^\star(\C^\times\ltimes\C)$ can be uniquely represented
as a sum of a (converging in ${\mathcal O}^\star(\C^\times\ltimes\C)$) series
 \beq\label{alpha-in-O*(C-times-C^times)}
\alpha=\sum_{k\in\N, n\in\Z}\alpha_{n,k}\cdot\zeta_n\boxedast \tau^k,
 \eeq
where the coefficients (continuously depend on $\alpha$ and) can be computed by
formula
 \beq
\alpha_{n,k}=\frac{1}{k!}\cdot\alpha(z^n\boxdot t^k)
 \eeq
The topology of the space ${\mathcal O}^\star(\C^\times\ltimes\C)$ can be
described by seminorms:
 \beq\label{polunormy-v-O*(C-rtimes-C-x)}
|||\alpha|||_r=\sum_{n\in\Z,k\in\N} r_{n,k}\cdot |\alpha_{n,k}|\cdot k!,\qquad
r_{n,k}\ge 0:\;\;\forall C\ge 1\;\; \sum_{n\in\Z,k\in\N} r_{k,n}\cdot
C^{k+|n|}<\infty
 \eeq

3) The bases $\{z^n\boxdot t^k;\, n\in\Z,\, k\in\N\}$ and $\{\zeta_n\boxedast
\tau^k;\, n\in\Z,\, k\in\N\}$ are dual to each other up to the constant $k!$:
 \beq
\langle z^m\boxdot t^k,\zeta_n*\tau^l\rangle=\langle
z^m,\zeta_n\rangle\cdot\langle t^k,\tau^l\rangle =\begin{cases}0, &
(m,k)\ne(n,l)
\\ k!,& (m,k)=(n,l)\end{cases}
 \eeq
and the action of functionals $\alpha\in{\mathcal O}^\star(\C^\times\ltimes\C)$
on functions $u\in {\mathcal O}(\C^\times\ltimes\C)$ is described by formula
 \beq
\langle u,\alpha\rangle=\sum_{n\in\Z, k\in\N}u_{n,k}\cdot\alpha_{n,k}\cdot k!
 \eeq
 \eprop

\paragraph{${\mathcal R}(\C^\times\ltimes\C)$,
${\mathcal R}^\star(\C^\times\ltimes\C)$, ${\mathcal O}(\C^\times\ltimes\C)$
and ${\mathcal O}^\star(\C^\times\ltimes\C)$ as Hopf algebras.}

The algebras ${\mathcal R}(\C^\times\ltimes\C)$ and ${\mathcal
O}(\C^\times\ltimes\C)$, being standard functional algebras on groups, are
endowed with the natural structure of Hopf algebras (we noted this general fact
in Theorems \ref{TH-R} and \ref{TH-O}). The following propositions describe the
structure of these Hopf algebras.

\bprop The algebra ${\mathcal R}(\C^\times\ltimes\C)$ (resp., algebra
${\mathcal O}(\C^\times\ltimes\C)$) is a nuclear Hopf-Brauner (resp.,
Hopf-Fr\'echet) algebra with the algebraic operations defined on basis elements
$z^n\boxdot t^k$ by formulas
 \begin{align}
\label{umn-v-R(C^x-x-C)} & z^m\boxdot t^k \cdot z^n\boxdot t^l=z^{m+n}\boxdot t^{k+l} && 1_{{\mathcal R}(\C^\times\ltimes\C)}=z^0t^0 \\
\label{koumn-v-R(C^x-x-C)}  &\varkappa(z^n\boxdot t^k)=\sum_{i=0}^k \begin{pmatrix}k\\
i\end{pmatrix}\cdot z^n t^i\odot  z^{n+i} t^{k-i} & &
\e(z^n\boxdot t^k)=\begin{cases}0,& (n,k)\ne(0,0)\\ 1,& (n,k)=(0,0) \end{cases} \\
\label{antipode-in-R(C^x-x-C)} & \sigma(z^n\boxdot t^k)=(-1)^k z^{-k-n}t^k &&
 \end{align}
As a corollary, the general formula of multiplication in ${\mathcal
R}(\C^\times\ltimes\C)$ is as follows:
 \beq\label{umn-u.v-R(C^x-x-C)}
u\cdot v=\sum_{n\in\Z, k\in\N} \Bigg(\sum_{j=0}^k \sum_{i\in\Z} u_{i,j}\cdot
v_{n-i,k-j}\Bigg)\cdot z^n\boxdot t^k.
 \eeq
 \eprop
\bpr The comultiplication, counit and antipode can be computed by formulas
\eqref{koumn-v-O(G)}-\eqref{antipod-v-O(G)}. For instance, the
comultiplication:
 \begin{multline*}
\widetilde{\varkappa}(z^n\boxdot t^k)\Big((a,x),(b,y)\Big)= (z^n\boxdot
t^k)\Big((a,x)\cdot(b,y)\Big)=(z^n\boxdot t^k)(ab,xb+y)=(ab)^n(xb+y)^k=\\=
a^nb^n\sum_{i=0}^k \begin{pmatrix}k \\ i\end{pmatrix} x^ib^iy^{k-i}=
\sum_{i=0}^k \begin{pmatrix}k \\ i\end{pmatrix} a^nx^i b^{n+i}y^{k-i}=
\sum_{i=0}^k \begin{pmatrix}k \\ i\end{pmatrix} (z^n\boxdot t^i)\boxdot
(z^{n+i}\boxdot t^{k-i})\Big((a,x),(b,y)\Big)
 \end{multline*}
$$
\Downarrow
$$
$$
\widetilde{\varkappa}(z^n\boxdot t^k)=\sum_{i=0}^k \begin{pmatrix}k \\
i\end{pmatrix}\cdot (z^n\boxdot t^i)\boxdot ( z^{n+i}\boxdot t^{k-i})
$$
$$
\Downarrow
$$
$$
\varkappa(z^n\boxdot t^k)=\sum_{i=0}^k \begin{pmatrix}k \\ i\end{pmatrix}\cdot
z^n\boxdot t^i\odot z^{n+i}\boxdot t^{k-i}
$$
 \epr

\bprop The algebra ${\mathcal R}^\star(\C^\times\ltimes\C)$ (resp., algebra
${\mathcal O}^\star(\C^\times\ltimes\C)$) is a nuclear Hopf-Fr\'echet (resp.,
Hopf-Brauner) algebra with the algebraic operations defined on basis elements
$\zeta_n\boxedast \tau^k$ by formulas
 \begin{align}
\label{umn-v-R*(C^x-x-C)} &
(\zeta_m\boxedast\tau^k)*(\zeta_n\boxedast\tau^l)=\begin{cases}\zeta_m\boxedast\tau^{k+l},&
m+k=n \\ 0, & m+k\ne n\end{cases} &&
1_{{\mathcal R}^\star(\C^\times\ltimes\C)}=\sum_{n\in\Z}\zeta_n\boxedast\tau^0 \\
\label{koumn-v-R*(C^x-x-C)} &\varkappa(\zeta_n\boxedast
\tau^k)=\sum_{m\in\Z}\sum_{i=0}^k\begin{pmatrix}k
\\ i\end{pmatrix}\cdot \zeta_m\boxedast\tau^i\circledast\zeta_{n-m}\boxedast\tau^{k-i} & &
\e(\zeta_n\boxedast \tau^k)=\begin{cases}1,& (n,k)=(0,0)\\ 0,& (n,k)\ne (0,0)\end{cases} \\
\label{antipode-in-R*(C^x-x-C)} & \sigma(\zeta^n\boxedast\tau^k)=(-1)^k\cdot
\zeta_{-n-k}\boxedast\tau^k &&
 \end{align}
As a corollary, the general formula of multiplication in ${\mathcal
R}^\star(\C^\times\ltimes\C)$ is as follows:
 \beq\label{umnozhenie-v-R_e^*}
\alpha*\beta=\sum_{n\in\Z, k\in\N}\Big(\sum_{j=0}^k
\alpha_{n,j}\cdot\beta_{n+j,k-j}\Big)\cdot\zeta_n\boxedast \tau^k.
 \eeq
 \eprop
\bpr All those formulas appear as duals of
\eqref{umn-v-R(C^x-x-C)}-\eqref{antipode-in-R(C^x-x-C)}. For example, the
comultiplication is computed as follows:
 \begin{multline*}
\langle u\odot v, \varkappa(\zeta_n\boxedast \tau^k)\rangle= \langle u\cdot v,
\zeta_n\boxedast \tau^k\rangle=\left\langle \sum_{r\in\Z, l\in\N}
\Bigg(\sum_{i=0}^l \sum_{m\in\Z} u_{m,i}\cdot v_{r-m,l-i}\Bigg)\cdot z^r t^l,
\zeta_n\boxedast \tau^k\right\rangle=\\=k!\cdot\sum_{i=0}^k \sum_{m\in\Z}
u_{m,i}\cdot v_{n-m,k-i}=k!\cdot\sum_{i=0}^k \sum_{m\in\Z} \frac{1}{i!}\langle
u, \zeta_m\boxedast\tau^i\rangle \cdot \frac{1}{(k-i)!} \langle v,
\zeta_{n-m}\boxedast\tau^{k-i}\rangle=\\= \sum_{m\in\Z} \sum_{i=0}^k
\begin{pmatrix}k \\ i\end{pmatrix}\cdot \langle u\odot v, \zeta_m\boxedast\tau^i\circledast
\zeta_{n-m}\boxedast\tau^{k-i}\rangle=\left\langle u\odot v, \sum_{m\in\Z}
\sum_{i=0}^k
\begin{pmatrix}k \\ i\end{pmatrix}\cdot  \zeta_m\boxedast\tau^i\circledast
\zeta_{n-m}\boxedast\tau^{k-i}\right\rangle
 \end{multline*}
$$
\Downarrow
$$
$$
\varkappa(\zeta_n\boxedast \tau^k)=\sum_{m\in\Z} \sum_{i=0}^k
\begin{pmatrix}k \\ i\end{pmatrix}\cdot  \zeta_m\boxedast\tau^i\circledast
\zeta_{n-m}\boxedast\tau^{k-i}
$$
\epr

\paragraph{Hopf algebras ${\mathcal R}_q(\C^\times\ltimes\C)$, ${\mathcal
R}_q^\star(\C^\times\ltimes\C)$, ${\mathcal O}_q(\C^\times\ltimes\C)$,
${\mathcal O}_q^\star(\C^\times\ltimes\C)$.}

The quantum group `$az+b$' can be defined as Hopf algebra ${\mathcal
R}(\C^\times\ltimes\C)$, where the algebraic operations are deformed in some
special way. We describe here this deformation, and together with the algebra
${\mathcal R}(\C^\times\ltimes\C)$ we shall consider the algebra ${\mathcal
O}(\C^\times\ltimes\C)$. Both these constructions will be useful below in
Theorem \ref{reflexivity-R_q(C^x-lx-C)}.

The constructions starts with the choice of a constant $q\in \C^\times$.

\bprop\label{PROP:def-R_q(C^x-lx-C)} On the stereotype space ${\mathcal
R}(\C^\times\ltimes\C)$ (respectively, ${\mathcal O}(\C^\times\ltimes\C)$)
there exists a unique structure of rigid stereotype Hopf algebra with the
algebraic operations defined on basis elements $z^n\boxdot t^k$ by formulas
 \begin{align}
\label{umn-v-R_q(C^x-x-C)} & z^m\boxdot t^k \cdot z^n\boxdot t^l=q^{kn}\cdot z^{m+n}\boxdot t^{k+l} && 1_{{\mathcal R}(\C^\times\ltimes\C)}=z^0t^0 \\
\label{koumn-v-R_q(C^x-x-C)}  &\varkappa(z^n\boxdot t^k)=\sum_{i=0}^k \begin{pmatrix}k\\
i\end{pmatrix}_q\cdot z^n\boxdot  t^i\odot  z^{n+i}\boxdot t^{k-i} & &
\e(z^n\boxdot t^k)=\begin{cases}0,& (n,k)\ne (0,0)\\ 1,& (n,k)=(0,0) \end{cases} \\
\label{antipode-in-R_q(C^x-x-C)} & \sigma(z^n\boxdot t^k)=(-1)^k\cdot
q^{-\frac{k(k+1)}{2}-kn}\cdot z^{-k-n}\boxdot t^k &&
 \end{align}
The space ${\mathcal R}(\C^\times\ltimes\C)$ (respectively, ${\mathcal
O}(\C^\times\ltimes\C)$) with such a structure of Hopf algebra is defined by
${\mathcal R}_q(\C^\times\ltimes\C)$ (respectively, by ${\mathcal
O}_q(\C^\times\ltimes\C)$). Besides this,
 \bit
\item[1)] the general formula of multiplication in ${\mathcal
R}_q(\C^\times\ltimes\C)$ (respectively, in ${\mathcal
O}_q(\C^\times\ltimes\C)$) has the form
 \beq\label{umn-u.v-R_q(C^x-x-C)}
u\cdot v=\sum_{n\in\Z, k\in\N} \Bigg(\sum_{i=0}^k \sum_{m\in\Z} q^{i(n-m)}\cdot
u_{m,i}\cdot v_{n-m,k-i}\Bigg)\cdot z^n\boxdot t^k.
 \eeq
\item[2)] The mapping
$$
z^n\boxdot t^k\mapsto z^n\odot t^k
$$
establishes an isomorphism between ${\mathcal R}_q(\C^\times\ltimes\C)$
(respectively, ${\mathcal O}_q(\C^\times\ltimes\C)$) and the Hopf algebra of
skew polynomials (entire functions) with coefficients in ${\mathcal
R}(\C^\times)$ (respectively, in ${\mathcal O}(\C^\times)$) with respect to the
quantum pair $(z,\delta^q)$ from Proposition \ref{quant-pair-z-d^q}:
 \beq\label{R_q(C^x-lx-C)-cong-R(C^x)-o-R(C)}
{\mathcal R}_q(\C^\times\ltimes\C)\cong{\mathcal R}(\C^\times)
\underset{\delta^q}{\stackrel{z}{\odot}} {\mathcal R}(\C)\qquad \Bigg({\mathcal
O}_q(\C^\times\ltimes\C)\cong{\mathcal O}(\C^\times)
\underset{\delta^q}{\stackrel{z}{\odot}} {\mathcal O}(\C)\Bigg)
 \eeq
\item[3)] For $q=1$ the Hopf algebra ${\mathcal R}_q(\C^\times\ltimes\C)$
(respectively, ${\mathcal O}_q(\C^\times\ltimes\C)$) turns into the Hopf
algebra ${\mathcal R}(\C^\times\ltimes\C)$ (respectively, ${\mathcal
O}(\C^\times\ltimes\C)$) with the structure of Hopf algebra described on page
\pageref{R(C^x-lx-C)-i-R*(C^x-lx-C)}:
 $$
{\mathcal R}(\C^\times\ltimes\C)={\mathcal R}_1(\C^\times\ltimes\C)\qquad
\Bigg({\mathcal O}(\C^\times\ltimes\C)={\mathcal O}_1(\C^\times\ltimes\C)
\Bigg)
 $$
\eit
 \eprop
\bpr Here everything starts from formulas
\eqref{R_q(C^x-lx-C)-cong-R(C^x)-o-R(C)}: the mapping $z^n\boxdot t^k\mapsto
z^n\odot t^k$ establishes an isomorphism of stereotype spaces
$$
{\mathcal R}_q(\C^\times\ltimes\C)\cong{\mathcal R}(\C^\times) \odot {\mathcal
R}(\C)
$$
(this is exactly the isomorphism, which for general case is described by the
identity \eqref{R(G-x-H)=R(G)-o-R(H)}). This isomorphism induces on ${\mathcal
R}(\C^\times\ltimes\C)$ the structure of rigid Hopf algebra from ${\mathcal
R}(\C^\times) \underset{\delta^q}{\stackrel{z}{\odot}} {\mathcal R}(\C)$ (where
this structure is defined by formulas
\eqref{umnozh-v-R(C)-odot-H}-\eqref{antipode-in-R(C)-odot-H}). In this
isomorphism the formulas
\eqref{umnozh-v-R(C)-odot-H}-\eqref{antipode-in-R(C)-odot-H} turn into formulas
\eqref{umn-v-R_q(C^x-x-C)}-\eqref{antipode-in-R_q(C^x-x-C)}. To derive them we
need to use the first formula in \eqref{M_upsilon^star(z)=qz}:
$$
\M_{\delta^q}^\star(z)=q\cdot z
$$
This implies
$$
(\M_{\delta^q}^\star)^k(z^n)=q^{kn}\cdot z^n
$$
Then, for instance, the formula of multiplication \eqref{umn-v-R_q(C^x-x-C)} is
derived from \eqref{umnozh-v-R(C)-odot-H} as follows:
 $$
z^m\boxdot t^k\cdot z^n\odot t^l=z^m\cdot (\M_{\delta^q}^\star)^k(z^n)\cdot
t^{k+l}=z^m\cdot q^{kn}\cdot z^n\odot t^{k+l}=q^{kn}\cdot z^{k+n}\odot t^{k+l}
 $$
The remaining formulas are deduced by analogy. The general formula for
multiplication \eqref{umn-u.v-R_q(C^x-x-C)} follows from
\eqref{umn-v-R_q(C^x-x-C)}:
 \begin{multline*}
u\cdot v=\left(\sum_{m\in\Z,k\in\N} u_{m,k}\cdot z^m\boxdot t^k\right)\cdot
\left(\sum_{n\in\Z,l\in\N} v_{n,l}\cdot z^n\boxdot
t^l\right)=\sum_{m\in\Z,k\in\N}\sum_{n\in\Z,l\in\N}u_{m,k}\cdot v_{n,l}\cdot
z^m\boxdot t^k\cdot z^n\boxdot
t^l=\\=\sum_{m\in\Z,k\in\N}\sum_{n\in\Z,l\in\N}u_{m,k}\cdot
v_{n,l}\cdot q^{kn}\cdot z^{m+n}t^{k+l}=\begin{pmatrix}n=r-m \\
l=s-k\end{pmatrix}=\\=\sum_{r\in\Z,s\in\N}\sum_{m\in\Z,k\in\N} q^{k(r-m)}\cdot
u_{m,k}\cdot v_{r-m,s-k}\cdot  z^{m+n}t^{k+l}
 \end{multline*}
It remains to add that for $q=1$ the formulas
\eqref{umn-v-R_q(C^x-x-C)}-\eqref{antipode-in-R_q(C^x-x-C)} turn into formulas
\eqref{umn-v-R(C^x-x-C)}-\eqref{antipode-in-R(C^x-x-C)}, so the Hopf algebra
${\mathcal R}_q(\C^\times\ltimes\C)$ (respectively, ${\mathcal
O}_q(\C^\times\ltimes\C)$) turns into Hopf algebra ${\mathcal
R}(\C^\times\ltimes\C)$ (respectively, ${\mathcal O}(\C^\times\ltimes\C)$).
\epr

\bprop On the stereotype space ${\mathcal R}^\star(\C^\times\ltimes\C)$
(respectively, ${\mathcal O}^\star(\C^\times\ltimes\C)$) there exists a unique
structure of rigid stereotype Hopf algebra with algebraic operations defined on
basis elements $\zeta_n\boxedast \tau^k$ by formulas
 \begin{align}
\label{umn-v-R_q*(C^x-x-C)} &
(\zeta_m\boxedast\tau^k)*(\zeta_n\boxedast\tau^l)=\begin{cases}\zeta_m\boxedast\tau^{k+l},&
m=n-k
\\ 0,& m\ne n-k\end{cases} &&
1_{{\mathcal R}^\star(\C^\times\ltimes\C)}=\sum_{n\in\Z}\zeta_n\boxedast\tau^0 \\
\label{koumn-v-R_q*(C^x-x-C)} &\varkappa(\zeta_n\boxedast
\tau^k)=\sum_{m\in\Z}\sum_{i=0}^k\begin{pmatrix}k
\\ i\end{pmatrix}_q\cdot q^{i(n-m)}\cdot \zeta_m\boxedast\tau^i\circledast\zeta_{n-m}\boxedast\tau^{k-i} & &
\e(\zeta_n\boxedast \tau^k)=\begin{cases}1,& k=0\\ 0,& k>0\end{cases} \\
\label{antipode-in-R_q*(C^x-x-C)} & \sigma(\zeta^n\boxedast\tau^k)=(-1)^k\cdot
q^{-\frac{k(k+1)}{2}+k(n+k)} \zeta_{-n-k}\boxedast\tau^k &&
 \end{align}
The space ${\mathcal R}^\star(\C^\times\ltimes\C)$ (respectively, ${\mathcal
O}^\star(\C^\times\ltimes\C)$) with this structure of Hopf algebra is denoted
by ${\mathcal R}_q^\star(\C^\times\ltimes\C)$ (respectively, ${\mathcal
O}_q^\star(\C^\times\ltimes\C)$). Besides this,
 \bit
\item[1)] the general formula for multiplication in ${\mathcal
R}_q^\star(\C^\times\ltimes\C)$ (respectively, in ${\mathcal
O}_q^\star(\C^\times\ltimes\C)$) has the form
 \beq\label{umnozhenie-a*b-R_q*(C^x-x-C)}
\alpha*\beta=\sum_{n\in\Z, k\in\N}\Big(\sum_{i=0}^k
\alpha_{n,i}\cdot\beta_{n+i,k-i}\Big)\cdot\zeta_n\boxedast \tau^k.
 \eeq
\item[2)]  the mapping
$$
\zeta_n\boxedast \tau^k\mapsto \zeta_n\circledast \tau^k
$$
establishes an isomorphism between ${\mathcal R}_q^\star(\C^\times\ltimes\C)$
(respectively, ${\mathcal O}_q^\star(\C^\times\ltimes\C)$) and the Hopf algebra
of skew power series (respectively, analytic functionals) with coefficients in
${\mathcal R}^\star(\C^\times)$ (respectively, in ${\mathcal
O}^\star(\C^\times)$) with respect to the quantum pair $(\delta^q,z)$ from
Proposition \ref{quant-pair-z-d^q}:
 \beq\label{R_q*(C^x-lx-C)-cong-R_q*(C^x)-o-R(C)}
{\mathcal R}_q^\star(\C^\times\ltimes\C)\cong{\mathcal R}^\star(\C^\times)
\underset{z}{\stackrel{\delta^q}{\odot}} {\mathcal R}^\star(\C),\qquad
{\mathcal O}_q^\star(\C^\times\ltimes\C)\cong{\mathcal O}^\star(\C^\times)
\underset{z}{\stackrel{\delta^q}{\odot}} {\mathcal O}^\star(\C)
 \eeq
\item[3)] for $q=1$ the Hopf algebra ${\mathcal R}_q^\star(\C^\times\ltimes\C)$
(respectively, ${\mathcal O}_q^\star(\C^\times\ltimes\C)$) turns into the Hopf
algebra ${\mathcal R}^\star(\C^\times\ltimes\C)$ (respectively, ${\mathcal
O}^\star(\C^\times\ltimes\C)$) with the structure of Hopf algebra defined on
page \pageref{R(C^x-lx-C)-i-R*(C^x-lx-C)}:
 $$
{\mathcal R}^\star(\C^\times\ltimes\C)={\mathcal
R}_1^\star(\C^\times\ltimes\C)\qquad \Bigg({\mathcal
O}^\star(\C^\times\ltimes\C)={\mathcal O}_1^\star(\C^\times\ltimes\C) \Bigg)
 $$
\eit
 \eprop

\bpr Here again everything is based on formulas
\eqref{R_q*(C^x-lx-C)-cong-R_q*(C^x)-o-R(C)}: the mapping $\zeta_n\boxedast
\tau^k\mapsto \zeta_n\circledast \tau^k$ establishes an isomorphism of
stereotype spaces
$$
{\mathcal R}_q^\star(\C^\times\ltimes\C)\cong{\mathcal R}^\star(\C^\times)
\circledast{\mathcal R}^\star(\C)
$$
(this is the dual isomorphism for \eqref{R(G-x-H)=R(G)-o-R(H)}). This
isomorphism induces on ${\mathcal R}^\star(\C^\times\ltimes\C)$ a structure of
rigid Hopf algebra from ${\mathcal R}^\star(\C^\times)
\underset{z}{\stackrel{\delta^q}{\circledast}}{\mathcal R}^\star(\C)$ (where
this structure is defined by formulas
\eqref{umn-v-R*(C)-*-H*}-\eqref{antipode-in-R*(C)-*-H*}). Under this
isomorphism formulas \eqref{umn-v-R*(C)-*-H*}-\eqref{antipode-in-R*(C)-*-H*}
turn into formulas
\eqref{umn-v-R_q*(C^x-x-C)}-\eqref{antipode-in-R_q*(C^x-x-C)}. To derive them
we can use two formulas:
 \beq\label{M_z*(zeta_n)=zeta_n-1}
\M_z^\star(\zeta_n)=\zeta_{n-1}, \qquad \delta^q*\zeta_n=q^n\cdot\zeta_n
 \eeq
The first of them (it will be useful in proving formulas for multiplication and
antipode), is deduced as follows:
 \begin{multline*}
\langle u,\M_z^\star(\zeta_n)\rangle=\langle z\cdot
u,\zeta_n\rangle=\left\langle z\cdot \sum_{m\in\Z}u_m\cdot
z^m,\zeta_n\right\rangle= \left\langle \sum_{m\in\Z}u_m\cdot
z^{m+1},\zeta_n\right\rangle=\\=\left\langle \sum_{m\in\Z}u_{m-1}\cdot
z^m,\zeta_n\right\rangle=u_{n-1}=\langle u,\zeta_{n-1}\rangle
 \end{multline*}
And the second one (it will be useful in proving the formula for
comultiplication) as follows:
$$
\delta^q*\zeta_n=\eqref{delta^q=sum}=\sum_{m\in\Z}q^m\cdot\zeta_m*\zeta_n=
\eqref{umn-v-R^star(C^times)}=q^n\cdot\zeta_n
$$
\epr

\paragraph{${\mathcal R}_q(\C^\times\ltimes\C)$ as an algebra with generators and defining relations.}

It remains to explain why the constructed algebra ${\mathcal
R}_q(\C^\times\ltimes\C)$ indeed can be identified with the quantum group
`$az+b$', i.e. with the Hopf algebra which is usually denoted as `$az+b$' (see
\cite{Woronowicz,VanDaele,Wang,QSNG}). Formally `$az+b$' is defined as the Hopf
algebra with three generators $t$, $z$, $z^{-1}$ and the defining relations
 \begin{align}
\label{tozhd-dlya-R} & t\cdot z=q\cdot z\cdot t && z\cdot z^{-1}=1 && 1=z^{-1}\cdot z \\
\label{k-R}&\varkappa(t)=t\otimes 1+z\otimes t, && \varkappa(z)=z\otimes z, &&
\varkappa(z^{-1})=z^{-1}\otimes z^{-1} \\
\label{e-R} &\e(t)=0, && \e(z)=1, && \e(z^{-1})=0 \\
\label{s-R} & \sigma(t)=-t\cdot z^{-1}, && \sigma(z)=z^{-1}, &&
\sigma(z^{-1})=z,
 \end{align}

\bprop\label{R_q(C^x-x-C)-porozhd-t-z} The mapping
$$
z\mapsto z\boxdot 1,\qquad z^{-1}\mapsto z^{-1}\boxdot 1,\qquad t\mapsto
1\boxdot t
$$
is uniquely extended to an isomorphism between Hopf algebras `$az+b$' and
${\mathcal R}_q(\C^\times\ltimes\C)$.
 \eprop
\bpr From \eqref{umn-v-R_q(C^x-x-C)} it follows that in ${\mathcal
R}_q(\C^\times\ltimes\C)$ the following identities hold
$$
1\boxdot t\cdot z\boxdot 1=q\cdot z\boxdot 1\cdot 1\boxdot t,\qquad z\boxdot
1\cdot z^{-1}\boxdot 1=1,\qquad 1=z^{-1}\boxdot 1\cdot z\boxdot 1
$$
in which one can recognize formulas \eqref{tozhd-dlya-R}, transformed by our
mapping. This means that our mapping can be uniquely extended to some
homomorphism of algebras $\ph:`az+b'\to{\mathcal R}_q(\C^\times\ltimes\C)$.
This homomorphism is a bijection since it turns the algebraic basis
$\{z^nt^k;\, n\in\Z,\,k\in\N\}$ in algebra `$az+b$' into the algebraic basis
$\{z^n\boxdot t^k;\, n\in\Z,\,k\in\N\}$ in algebra ${\mathcal
R}_q(\C^\times\ltimes\C)$.

Thus, $\ph$ is an isomorphism of algebras. Note then that $\ph$ preserves the
comultiplication on generators:
$$
(\ph\otimes\ph)(\varkappa(t))=(\ph\otimes\ph)(t\otimes 1+z\otimes t)=1\boxdot
t\otimes 1\boxdot 1+z\boxdot 1\otimes 1\boxdot t=\varkappa(1\boxdot
t)=\varkappa(\ph(t))
$$
and similarly,
$$
(\ph\otimes\ph)(\varkappa(z))=\varkappa(\ph(z)),\qquad
(\ph\otimes\ph)(\varkappa(z^{-1}))=\varkappa(\ph(z^{-1}))
$$
Since the comultiplication, like $\ph$, is a homomorphism of algebras, this
implies that the same formulas are true for the arguments of the form $z^nt^k$,
and thus for all elements of algebra `$az+b$'. We obtain that $\ph$ preserves
the comultiplication (on all elements). Similarly it is proved that $\ph$
preserves counit and antipode (here we need to use the fact that the counit is
a homomorphism, and antipode an antihomomorphism). Hence, $\ph$ is an
isomorphism of Hopf algebras.
 \epr

\subsection{Reflexivity of ${\mathcal R}_q(\C^\times\ltimes\C)$}

\paragraph{Reflexivity diagram for ${\mathcal R}_q(\C^\times\ltimes\C)$.}

\btm\label{reflexivity-R_q(C^x-lx-C)} For any $q\in\C^\times$ the rigid Hopf
algebra ${\mathcal R}_q(\C^\times\ltimes\C)$ is holomorphically reflexive, and
 \bit
\item[1)] for $|q|=1$ its reflexivity diagram has the form
 \beq\label{R_q(C-rtimes-C^times)-|q|=1}
 \xymatrix @R=1.pc @C=1.pc
 {
 {\mathcal R}_q(\C^\times\ltimes\C)  \cong  {\mathcal
R}(\C^\times)\underset{\delta^q}{\overset{z}{\odot}}{\mathcal R}(\C)
 & \ar@{|->}[r]^{\heartsuit} & &
 {\mathcal O}(\C^\times)\underset{\delta^q}{\overset{z}{\odot}}{\mathcal O}(\C)  \cong  {\mathcal O}_q(\C^\times\ltimes\C)
 \\
 & & &
 \ar@{|->}[d]^{\star}
 \\
  \ar@{|->}[u]^{\star}
 & & &
 \\
 {\mathcal R}_q^\star(\C^\times\ltimes\C)  \cong  {\mathcal
R}^\star(\C^\times)\underset{z}{\overset{\delta^q}{\odot}}{\mathcal
R}^\star(\C)
 & &
 \ar@{|->}[l]_{\heartsuit}
 &
 {\mathcal O}^\star(\C^\times)\underset{z}{\overset{\delta^q}{\odot}}{\mathcal
O}^\star(\C) \cong  {\mathcal O}_q^\star(\C^\times\ltimes\C)
 }
 \eeq
\item[2)] for $|q|\ne 1$ the form:
 \beq\label{R_q(C-rtimes-C^times)-|q|-ne-1}
 \xymatrix @R=1.pc @C=1.pc
 {
 {\mathcal R}_q(\C^\times\ltimes\C)  \cong
 {\mathcal R}(\C^\times)\underset{\delta^q}{\overset{z}{\odot}}{\mathcal R}(\C)
 & \ar@{|->}[r]^{\heartsuit} & &
 {\mathcal O}(\C^\times)\underset{\delta^q}{\overset{z}{\odot}}{\mathcal R}^\star(\C)
 \\
 & & &
 \ar@{|->}[d]^{\star}
 \\
 \ar@{|->}[u]^{\star}
 & & &
 \\
 {\mathcal R}_q^\star(\C^\times\ltimes\C)  \cong
 {\mathcal R}^\star(\C^\times)\underset{z}{\overset{\delta^q}{\odot}}{\mathcal R}^\star(\C)
 & &
 \ar@{|->}[l]_{\heartsuit}
 &
 {\mathcal O}^\star(\C^\times)\underset{z}{\overset{\delta^q}{\odot}}{\mathcal R}(\C)
 }
\eeq
 \eit
\etm

For symmetry one can note here that all the multipliers here -- ${\mathcal
R}(\C)$, ${\mathcal R}^\star(\C)$, ${\mathcal R}(\C^\times)$, ${\mathcal
R}^\star(\C^\times)$, ${\mathcal O}(\C^\times)$, ${\mathcal
O}^\star(\C^\times)$ -- are nuclear spaces and appear only in pair
Fr\'echet-Fr\'echet and Brauner-Brauner. So if in diagrams
\eqref{R_q(C-rtimes-C^times)-|q|=1}-\eqref{R_q(C-rtimes-C^times)-|q|-ne-1} all
(or some of) the injective tensor products $\odot$ are replaced by projective
tensor products $\circledast$, then we shall obtain isomorphic diagrams.

The rest of this section is devoted to the proof of Theorem
\ref{reflexivity-R_q(C^x-lx-C)}. We carry out it in three steps.

\paragraph{${\mathcal
R}_{q}(\C^\times\ltimes\C)\stackrel{\heartsuit}{\to}{\mathcal
O}_{q}(\C^\times\ltimes\C)$ for $|q|=1$.}

\bprop For $|q|=1$ the algebra ${\mathcal O}_{q}(\C^\times\ltimes\C)$ is an
Arens-Michael algebra, and the natural inclusion ${\mathcal
R}_{q}(\C^\times\ltimes\C)\subseteq {\mathcal O}_{q}(\C^\times\ltimes\C)$ is an
Arens-Michael envelope of the algebra ${\mathcal R}_{q}(\C^\times\ltimes\C)$:
$$
{\mathcal R}_{q}(\C^\times\ltimes\C)^\heartsuit={\mathcal
O}_{q}(\C^\times\ltimes\C)
$$
 \eprop
\bpr This follows from the fact that for $|q|=1$ the seminorms
\eqref{polunormy-v-O(C-rtimes-C-x)} on ${\mathcal R}_{q}(\C^\times\ltimes\C)$
are submultiplicative. \epr

\paragraph{${\mathcal R}(\C^\times)\odot^z_{\delta^q}{\mathcal R}(\C)
\stackrel{\heartsuit}{\to}{\mathcal O}(\C^\times)\odot^z_{\delta^q}{\mathcal
R}^\star(\C)$ for $|q|\ne 1$.}

Recall that the structure of algebras ${\mathcal R}(\C^\times)$ and ${\mathcal
O}(\C^\times)$ was discussed in
\ref{SEC:stein-groups}\ref{SUBSEC:algebry-na-C^x}. In particular, we have noted
there that the topology of ${\mathcal O}(\C^\times)$ is generated by seminorms
\eqref{polunormy-v-H(C-x)}:
 $$
||u||_C=\sum_{n\in\Z} |u_n|\cdot C^{|n|},\qquad C\ge 1.
 $$
In Example \ref{EX-polunormy-v-O(C^x)} noted also that these seminorms are
submultiplicative. Let us state two more their properties: first, obviously,
 \beq\label{|u|_C-le-|u|_D}
C\le D\quad\Longrightarrow\quad \norm{u}_C\le \norm{u}_D
 \eeq
and, second, for $|q|<1$ we have
 \beq
|q|^n\le |q|^{-|n|},\qquad n\in\Z,
 \eeq
This implies
 $$
\norm{\M_{\delta^q}^\star(u)}_C=\sum_{n\in\Z} |u_n|\cdot |q|^n\cdot C^{|n|}\le
\sum_{n\in\Z} |u_n|\cdot |q|^{-|n|}\cdot C^{|n|}=\sum_{n\in\Z} |u_n|\cdot
\left(\frac{C}{|q|}\right)^{|n|}=\norm{u}_{C/|q|}
 $$
i.e.,
 \beq\label{|ph_q^i(u)|_C-le-|u|_C/|q|^i}
\forall i\in\N \qquad \norm{(\M_{\delta^q}^\star)^i(u)}_C\le \norm{u}_{C/|q|^i}
\qquad (|q|<1)
 \eeq
Similarly, for $|q|>1$ we have
 \beq
|q|^n\le |q|^{|n|},\qquad n\in\Z,
 \eeq
so
 $$
\norm{\M_{\delta^q}^\star(u)}_C=\sum_{n\in\Z} |u_n|\cdot |q|^n\cdot C^{|n|}\le
\sum_{n\in\Z} |u_n|\cdot |q|^{|n|}\cdot C^{|n|}=\sum_{n\in\Z} |u_n|\cdot
(C\cdot|q|)^{|n|}=\norm{u}_{C\cdot|q|}
 $$
i.e.,
 \beq\label{|ph_q^i(u)|_C-le-|u|_C/|q|^i-|q|>1}
\forall i\in\N \qquad \norm{(\M_{\delta^q}^\star)^i(u)}_C\le \norm{u}_{C\cdot
|q|^i} \qquad (|q|>1)
 \eeq

In accordance with Proposition \ref{R_q(C^x-x-C)-porozhd-t-z}, let us identify
in calculations the symbols $1\boxdot t$ and $z\boxdot 1$ with the symbols $z$
and $t$:
$$
1\boxdot t\equiv t,\qquad z\boxdot 1\equiv z
$$

 \blm\label{LM-polunormy-na-R(C)-odot-C_Z}
If $|q|\ne 1$ and $r$ is a submultiplicative seminorm on ${\mathcal
R}(\C^\times)\underset{\delta^q}{\overset{z}{\odot}}{\mathcal R}(\C)$, then for
some $K\in\N$
 \beq\label{forall k-ge-K-r(t^k)=0}
\forall k>K\qquad r(t^k)=0
 \eeq
 \elm
\bpr 1. Suppose first that $|q|<1$ and take $K\in\N$ such that
$|q|^K<\frac{1}{r(z)\cdot r(z^{-1})}$. Then for any $k>K$ we obtain $|q|^k\cdot
r(z)\cdot r(z^{-1})<1$, and so
 \begin{multline*}
 r(t^k)=r(t^k\cdot z^l\cdot z^{-l} )= |q|^{lk}\cdot r( z^l\cdot t^k
\cdot  z^{-l})\le  |q|^{lk}\cdot r( z^l)\cdot r(t^k)\cdot r(z^{-l})\le\\ \le
|q|^{lk}\cdot r(z)^l\cdot r(t)^k\cdot r(z^{-1})^l= r(t)^k\cdot (|q|^k\cdot r(z)
\cdot r(z^{-1}))^l\underset{l\to\infty}{\longrightarrow} 0
  \end{multline*}
2. On the contrary, suppose $|q|>1$. Then we can take $K\in\N$ such that
$|q|^K>r(z)\cdot r(z^{-1})$, and for all $k>K$ we have $\frac{r(z)\cdot
r(z^{-1})}{|q|^k}<1$, so
 \begin{multline*}
 r(t^k)=r(t^k\cdot z^{-l}\cdot z^l )= |q|^{-lk}\cdot r( z^{-l}\cdot t^k
\cdot  z^l)\le  |q|^{-lk}\cdot r( z^l)\cdot r(t^k)\cdot r(z^{-l})\le\\ \le
|q|^{-lk}\cdot r(z)^l\cdot r(t)^k\cdot r(z^{-1})^l= r(t)^k\cdot
\left(\frac{r(z) \cdot
r(z^{-1})}{|q|^k}\right)^l\underset{l\to\infty}{\longrightarrow} 0
  \end{multline*}
\epr

\bprop\label{PROP:p-na-R_q(C-times-C^x)} 1) Suppose $|q|<1$, and $D\ge 1$ and
$K\in\N$ are such that
 \beq\label{C-cdot-|q|^K-ge-1}
D\cdot |q|^K\ge 1
 \eeq
Then the seminorm $p_{D,K}:{\mathcal
R}(\C^\times)\underset{\delta^q}{\overset{z}{\odot}}{\mathcal R}(\C)\to\R_+$,
defined by the equality
 \beq\label{p-na-R_q(C-times-C^x)}
p_{D,K}\left(\sum_{k\in\N} u_k\odot t^k\right)
=\sum_{k=0}^K\underbrace{||u_k||_{D\cdot
|q|^k}}_{\scriptsize\begin{matrix}\text{seminorm}\\
\text{\eqref{polunormy-v-H(C-x)}}\end{matrix}}=
\sum_{k=0}^K\sum_{n\in\Z}|u_{k,n}|\cdot \left(D\cdot |q|^k\right)^{|n|}
 \eeq
is submultiplicative. On the contrary, every submupltiplicative seminorm on
${\mathcal R}(\C^\times)\underset{\delta^q}{\overset{z}{\odot}}{\mathcal
R}(\C)$ is subordinate to a seminorm of the form \eqref{p-na-R_q(C-times-C^x)}.

2) If $|q|>1$, and $D\ge 1$ and $K\in\N$ are such that
 \beq\label{D/|q|^K-ge-1}
\frac{D}{|q|^K}\ge 1
 \eeq
then the seminorm $p_{D,K}:{\mathcal
R}(\C^\times)\underset{\delta^q}{\overset{z}{\odot}}{\mathcal R}(\C)\to\R_+$,
defined by the equality
 \beq\label{p-na-R_q(C-times-C^x)-|q|>1}
p_{D,K}\left(\sum_{k\in\N} u_k\odot t^k\right)
=\sum_{k=0}^K\underbrace{||u_k||_{\frac{D}{|q|^K}}}_{\scriptsize\begin{matrix}\text{seminorm}\\
\text{\eqref{polunormy-v-H(C-x)}}\end{matrix}}=
\sum_{k=0}^K\sum_{n\in\Z}|u_{k,n}|\cdot \left(\frac{D}{|q|^K}\right)^{|n|}
 \eeq
is submultiplicative. On the contrary, every submultiplicative seminorm on
${\mathcal R}(\C^\times)\underset{\delta^q}{\overset{z}{\odot}}{\mathcal
R}(\C)$ is subordinate to some seminorm of the form
\eqref{p-na-R_q(C-times-C^x)}. \eprop
 \bpr
Consider the case $|q|<1$. The submultiplicativity of the seminorm
\eqref{p-na-R_q(C-times-C^x)} is verified directly:
 \begin{multline*}
p_{D,K}(u\cdot v)=\sum_{k=0}^K \norm{(u\cdot v)_k}_{D\cdot
|q|^k}=\eqref{obsh-umnozh-v-R(C)-odot-H}=\sum_{k=0}^K \norm{\sum_{i=0}^k
u_i\cdot (\M_{\delta^q}^\star)^i(v_{k-i})}_{D\cdot |q|^k}=\\=
 \sum_{k=0}^K \sum_{i=0}^k \norm{u_i\cdot (\M_{\delta^q}^\star)^i(v_{k-i})}_{D\cdot |q|^k}\le
 \sum_{k=0}^K \sum_{i=0}^k \norm{u_i}_{D\cdot |q|^k}\cdot \norm{(\M_{\delta^q}^\star)^i(v_{k-i})}_{D\cdot |q|^k}
 \le \eqref{|u|_C-le-|u|_D},\eqref{|ph_q^i(u)|_C-le-|u|_C/|q|^i} \le\\ \le
 \sum_{k=0}^K \sum_{i=0}^k \norm{u_i}_{D\cdot |q|^i}\cdot \norm{v_{k-i}}_{D\cdot
 |q|^{k-i}}\le  \left(\sum_{i=0}^K \norm{u_i}_{D\cdot |q|^i}\right)\cdot
 \left(\sum_{j=0}^K \norm{v_j}_{D\cdot |q|^j}\right)=p_{D,K}(u)\cdot p_{D,K}(v)
 \end{multline*}
Then, let $r$ be a submultiplicative seminorm on ${\mathcal
R}(\C^\times)\underset{\delta^q}{\overset{z}{\odot}}{\mathcal R}(\C)$. Let us
choose by Lemma \ref{LM-polunormy-na-R(C)-odot-C_Z} an integer $K\in\N$ such
that \eqref{forall k-ge-K-r(t^k)=0} holds, and put
$$
L=\max_{0\le k\le K} r(t^k)
$$
Note that $r$ is a submultiplicative seminorm on the subalgebra $1\odot
{\mathcal R}(\C^\times)$, consisting of functions of the form $a\odot 1$,
$a\in{\mathcal R}(\C^\times)$, and isomorphic to ${\mathcal R}(\C^\times)$.
Hence, on this subalgebra $r$, must be subordinate to some seminorm
\eqref{polunormy-v-H(C-x)}:
 \beq
r(a\odot 1)\le M\cdot \norm{a}_C,\qquad a\in{\mathcal R}(\C^\times)
 \eeq
for some $C\ge 1$, $M>0$. Choose now $D\ge 1$ such that
 \beq\label{C-le-D|q|^K-le...le-D|q|-le-D}
C\le D\cdot |q|^K\le D\cdot |q|^{K-1}\le ...\le D\cdot |q|\le D
 \eeq
Then
 \begin{multline*}
r(u)=r\left(\sum_{k\in\N} u_k\odot t^k\right)\le \sum_{k\in\N} r(u_k)\cdot
r(t^k)=\sum_{k=0}^K r(u_k)\cdot r(t^k)\le \sum_{k=0}^K M\cdot \norm{u_k}_C\cdot
L=\\= L\cdot M\cdot \sum_{k=0}^K \norm{u_k}_{D\cdot|q|^k}= L\cdot M\cdot
p_{D,K}(u)
 \end{multline*}
Thus, $r$ is subordinate to some seminorm $p_{D,K}$.

The case $|q|>1$ is considered similarly, but instead of
\eqref{|ph_q^i(u)|_C-le-|u|_C/|q|^i} we have to apply here formula
\eqref{|ph_q^i(u)|_C-le-|u|_C/|q|^i-|q|>1}, and instead of
\eqref{C-le-D|q|^K-le...le-D|q|-le-D} the chain
 $$
C\le \frac{D}{|q|^K}\le \frac{D}{|q|^{K-1}}\le ...\le \frac{D}{|q|}\le D
 $$
 \epr

\bprop For $|q|\ne 1$ the formula
 \beq
\zeta_n\odot t^k\quad\mapsto\quad z_n\odot \tau^k
 \eeq
uniquely defines a homomorphism of Hopf algebras
$$
{\mathcal R}(\C^\times)\underset{\delta^q}{\overset{z}{\odot}}{\mathcal R}(\C)
\to {\mathcal O}(\C^\times)\underset{\delta^q}{\overset{z}{\odot}}{\mathcal
R}^\star(\C),
$$
which is an Arens-Michael envelope of the algebra ${\mathcal
R}(\C^\times)\underset{\delta^q}{\overset{z}{\odot}}{\mathcal R}(\C)$.
 \eprop
\bpr It remains to check here that ${\mathcal O}(\C^\times)\odot{\mathcal
R}^\star(\C)$ is a completion of ${\mathcal R}(\C^\times)\odot{\mathcal R}(\C)$
with respect to the seminorms \eqref{p-na-R_q(C-times-C^x)} or, depending on
$q$, seminorms \eqref{p-na-R_q(C-times-C^x)-|q|>1}.
 \epr

\paragraph{${\mathcal O}^\star(\C^\times)\odot_z^{\delta^q}{\mathcal R}(\C)
\stackrel{\heartsuit}{\to} {\mathcal
R}^\star(\C^\times)\odot_z^{\delta^q}{\mathcal R}^\star(\C)$ for arbitrary
$q$.}

Consider the space ${\mathcal O}^\star(\C^\times)$. From the first formula in
\eqref{M_z*(zeta_n)=zeta_n-1}
 $$
\M_z^\star(\zeta_n)=\zeta_{n-1},
 $$
one can deduce the identity:
 \beq\label{M_z*^k(alpha)}
(\M_z^\star)^k(\alpha)=(\M_z^\star)^k\left(\sum_{n\in\Z}\alpha_n\cdot
\zeta_n\right)=\sum_{n\in\Z}\alpha_n\cdot
\zeta_{n-k}=\sum_{m\in\Z}\alpha_{m+k}\cdot \zeta_m
 \eeq
Recall seminorms $||\cdot||_N$ on ${\mathcal O}^\star(\C^\times)$, defined by
formulas \eqref{||alpha||_C-v-C-x}:
$$
\norm{\alpha}_N=\sum_{|n|\le N}|\alpha_n|,\qquad N\in\N ,\qquad
\alpha=\sum_{n\in\Z}\alpha_n\cdot\zeta_n
$$
Note the following their properties:
 \beq\label{p_M(beta)-le-p_N(beta)}
M\le N\quad\Longrightarrow\quad ||\beta||_M\le ||\beta||_N
 \eeq
and
 \beq\label{p_M(psi(beta))-le-p_M+1(beta)}
\norm{(\M_z^\star)^i(\beta)}_M=\sum_{|n|\le M}
|(\M_z^\star)^i(\beta)_n|=\eqref{M_z*^k(alpha)}= \sum_{|n|\le M}
|\beta_{n+i}|\le \sum_{|n|\le M+i} |\beta_{n}|=\norm{\beta}_{M+i}
 \eeq

By Theorem \eqref{TH-R-odot-H-i-R*-circledast-H*}, the general formula for
multiplication in the algebra ${\mathcal
O}^\star(\C^\times)\underset{z}{\overset{\delta^q}{\odot}}{\mathcal R}(\C)$ has
the form
 \beq\label{umnozhenie-R(C)-psi-odot-O^star(C^times)}
\alpha\cdot\beta=\sum_{k\in\N}\left(\sum_{i=0}^k \alpha_i*
(\M_z^\star)^i(\beta_{k-i})\right)\odot t^k
 \eeq
and on basis elements looks as follows:
 \begin{multline}\label{umnozhenie-na-bazise-v-R(C)-psi-odot-O^star(C^times)}
\zeta_m\odot t^k\cdot \zeta_n\odot
t^l=\eqref{umnozh-v-R(C)-odot-H}=\zeta_m*(\M_z^\star)^k(\zeta_n)\odot
t^{k+l}=\eqref{M_z*(zeta_n)=zeta_n-1}=\\
=\zeta_m*\zeta_{n-k}\odot t^{k+l}=\begin{cases}\zeta_m\odot t^{k+l}, & m=n-k
\\ 0, & m\ne n-k
\end{cases}=\begin{cases}\zeta_m\odot t^{k+l}, &
m+k=n \\ 0, & m+k\ne n \end{cases}
 \end{multline}

\blm\label{LM-nepr-subm-polunorm-na-R(C)-psi-odot-O^star(C^times)} Every
continuous submultiplicative seminorm $p$ on the algebra ${\mathcal
O}^\star(\C^\times)\underset{z}{\overset{\delta^q}{\odot}}{\mathcal R}(\C)$
vanishes on almost all elements of the basis $\{ \zeta_n\odot t^k;\; n\in\Z,\;
k\in\N\}$ (i.e. on all but a finite subfamily):
 $$
\card\{(n,k)\in\Z\times\N:\ p\left(\zeta_n\odot t^k\right)\ne 0\}<\infty
 $$
 \elm
\bpr Put $p_{n,k}=p(\zeta_n\odot t^k)$ and note that
 \beq\label{p_k+m,n-le-p_k,m+n-cdot-p_m,n}
p_{m,k+l}\le p_{m,k}\cdot p_{m+k,l},\qquad k,m\in\N,\; n\in\Z
 \eeq
Indeed, from \eqref{umnozhenie-na-bazise-v-R(C)-psi-odot-O^star(C^times)} it
follows that
$$
\zeta_m\odot t^{k+l}=\zeta_m\odot t^k\cdot \zeta_{m+k}\odot t^l
$$
so
$$
p_{m,k+l}=p(\zeta_m\odot t^{k+l})\le p(\zeta_m\odot t^k)\cdot
p(\zeta_{m+k}\odot t^l)= p_{m,k}\cdot p_{m+k,l}
$$
Consider the sets
$$
S_k=\{n\in\Z:\; p_{n,k}\ne 0\}
$$
and note the following two things.

1. From the continuity and submupltiplicativity of a seminorm $p$ on ${\mathcal
R}(\C)\underset{z}{\overset{\delta^q}{\odot}}{\mathcal O}^\star(\C^\times)$ it
follows that the functional
$$
p_0(\alpha)=p(\alpha\odot 1)=p(\alpha\odot t^0),\qquad \alpha\in {\mathcal
O}^\star(\C^\times)
$$
is a continuous submultiplicative seminorm on ${\mathcal O}^\star(\C^\times)$:
$$
p_0(\alpha*\beta)=p((\alpha*\beta)\odot 1)=p\big((\alpha\odot 1)\cdot
(\beta\odot 1)\Big)\le p(\alpha\odot 1)\cdot p(\beta\odot 1)=p_0(\alpha)\cdot
p_0(\beta)
$$
This means, by Lemma \ref{LM-||alpha||_C-v-C-x}, that $p_0$ is subordinate to
some seminorm of the form \eqref{||alpha||_C-v-C-x}:
$$
p_0(\alpha)\le C\cdot \norm{\alpha}_N=C\cdot \sum_{|n|\le N}|\alpha_n|,
\qquad\alpha\in {\mathcal O}^\star(\C^\times)
$$
This in its turn implies that the set $S_0=\{n\in\Z:\;
p_0(\zeta_n)=p(1\odot\zeta_n)\ne 0\}$ must be finite:
$$
\card S_0<\infty
$$

2. From the inequalities \eqref{p_k+m,n-le-p_k,m+n-cdot-p_m,n} one can deduce
the following implications:
$$
\begin{cases}
p_{m,k+1}\le p_{m,k}\cdot p_{m+k,1} \\
p_{m,k+1}\le p_{m,k+1}\cdot p_{m+k+1,0}
\end{cases}
\quad \Longrightarrow\quad
\begin{cases}
S_{k+1}\subseteq S_k
\\
S_{k+1}\subseteq S_0-(k+1)
\end{cases}
 \quad \Longrightarrow\quad
S_{k+1}\subseteq S_k\cap (S_0-k-1)
 $$
(here $S-i$ means the shift of the set $S$ by $i$ units to left on the group
$\Z$). So we can conclude that $S_k$ form a tapering chain of finite (since
$S_0$ is finite) sets:
$$
S_0\supseteq S_1\supseteq ...\supseteq S_k\supseteq S_{k+1}\supseteq ...
$$
And at the minimum after the number $K=\max S_0-\min S_0$ this chain vanishes:
$$
S_0\supseteq S_1\supseteq ...\supseteq S_K\supseteq S_{K+1}=\varnothing
$$
(since $S_{K+1}\subseteq S_0\cap (S_0-K-1)=\varnothing$).  \epr

\bprop\label{PROP:stand-nepr-subm-polunorm-na-R(C)-psi-odot-O^star(C^times)}
For any $N\in\N$ the functional
 \beq\label{stand-nepr-subm-polunorm-na-R(C)-psi-odot-O^star(C^times)}
r_N(\alpha)= r_N\left(\sum_{k\in\N} \alpha_k\odot t^k\right) = \sum_{k=0}^N
\norm{\alpha_k}_{N-k}=\sum_{k=0}^N \sum_{|n|\le N-k}|\alpha_{n,k}|
 \eeq
is a continuous submultiplicative seminorm on ${\mathcal
R}(\C)\underset{z}{\overset{\delta^q}{\odot}}{\mathcal O}^\star(\C^\times)$.
Every continuous submultiplicative seminorm on ${\mathcal
R}(\C)\underset{z}{\overset{\delta^q}{\odot}}{\mathcal O}^\star(\C^\times)$ is
subordinated to some seminorm $r_N$.
 \eprop
\bpr The functional $r_N$ is a continuous seminorm because
$\alpha\mapsto\alpha_k$ are linear continuous mappings:
$$
r_N(\lambda\cdot\alpha+\beta)=\sum_{k=0}^N
\norm{\lambda\cdot\alpha_k+\beta_k}_{N-k}\le |\lambda|\cdot \sum_{k=0}^N
\norm{\alpha_k}_{N-k}+\sum_{k=0}^N \norm{\beta_k}_{N-k}=|\lambda|\cdot
r_N(\alpha)+r_N(\beta)
$$
Let us consider the submultiplicativity:
 \begin{multline*}
r_N(\alpha\cdot\beta)=\eqref{umnozhenie-R(C)-psi-odot-O^star(C^times)}=
r_N\left(\sum_{k\in\N} \sum_{i=0}^k\alpha_i*(\M_z^\star)^i(\beta_{k-i})\odot
t^k\right)= \sum_{k=0}^N
\norm{\sum_{i=0}^k\alpha_i*(\M_z^\star)^i(\beta_{k-i})}_{N-k}\le
\\ \le \sum_{k=0}^N \sum_{i=0}^k \norm{\alpha_i}_{N-k}\cdot
\norm{(\M_z^\star)^i(\beta_{k-i})}_{N-k}\le\eqref{p_M(psi(beta))-le-p_M+1(beta)}\le
\sum_{k=0}^N \sum_{i=0}^k \norm{\alpha_i}_{N-k}\cdot \norm{\beta_{k-i}}_{N-k+i}\le\\
\le \sum_{k=0}^N \sum_{i=0}^k \underbrace{\norm{\alpha_i}_{N-k}}_{\scriptsize
\begin{matrix}\text{\rotatebox{-90}{$\le$}} \\ \norm{\alpha_i}_{N-i}, \\ \text{since} \\ k\ge i, \\
\text{hence} \\ N-k\le N-i, \\ \text{and this allows} \\
\text{to apply \eqref{p_M(beta)-le-p_N(beta)}}
\end{matrix}}\cdot\norm{\beta_{k-i}}_{N-(k-i)}
 \le  \sum_{k=0}^N \sum_{i=0}^k
\norm{\alpha_i}_{N-i}\cdot\norm{\beta_{k-i}}_{N-(k-i)}\le\\ \le
\left(\sum_{i=0}^N \norm{\alpha_i}_{N-i}\right)\cdot \left(\sum_{j=0}^N
\norm{\beta_j}_{N-j}\right)=r_N(\alpha)\cdot r_N(\beta)
 \end{multline*}
It remains to check that every continuous submultiplicative seminorm $p$ on
${\mathcal R}(\C)\underset{z}{\overset{\delta^q}{\odot}}{\mathcal
O}^\star(\C^\times)$ is subordinate to some seminorm $r_N$. This follows from
Lemma \ref{LM-nepr-subm-polunorm-na-R(C)-psi-odot-O^star(C^times)}: since $p$
vanishes on almost all basis elements $t^k\odot\zeta_n$, we can find a number
$N\in\N$ such that
 $$
\{(n,k)\in\Z\times\N:\ p(\zeta_n\odot t^k)\ne 0\}\subseteq
\{(n,k)\in\Z\times\N:\ k\le N \ \& \ |n|\le N-k \}
 $$
Then we put
$$
C=\max\{ p(\zeta_n\odot t^k);\; (n,k):\; k\le N, \; |n|\le N-k\}
$$
and obtain:
 \begin{multline*}
p(\alpha)=p\left(\sum_{n\in\Z, k\in\N}\alpha_{n,k}\cdot \zeta_n\odot
t^k\right)\le\sum_{n\in\Z, k\in\N}|\alpha_{n,k}|\cdot p(\zeta_n\odot
t^k)=\sum_{k=0}^N\sum_{|n|\le N-k}|\alpha_{n,k}|\cdot p(\zeta_n\odot t^k)\le \\
\le C\cdot \sum_{k=0}^N\sum_{|n|\le N-k}|\alpha_{n,k}|=
\eqref{stand-nepr-subm-polunorm-na-R(C)-psi-odot-O^star(C^times)} =C\cdot
r_N(\alpha)
 \end{multline*}
\epr

\bprop The formula
 \beq
\zeta_n\odot t^k\quad\mapsto\quad \zeta_n\odot \tau^k
 \eeq
defines a homomorphism of Hopf algebras
$$
{\mathcal O}^\star(\C^\times)\underset{z}{\overset{\delta^q}{\odot}}{\mathcal
R}(\C) \to {\mathcal R}^\star(\C^\times)
\underset{z}{\overset{\delta^q}{\odot}}{\mathcal R}^\star(\C)
$$
which is an Arens-Michael envelope of the algebra ${\mathcal
O}^\star(\C^\times)\underset{z}{\overset{\delta^q}{\odot}}{\mathcal R}(\C)$.
 \eprop
\bpr It remains only to check that ${\mathcal R}^\star(\C^\times)
\odot{\mathcal R}^\star(\C)$ is a completion of ${\mathcal
O}^\star(\C^\times)\odot{\mathcal R}(\C)$ with respect to seminorms
\eqref{stand-nepr-subm-polunorm-na-R(C)-psi-odot-O^star(C^times)}.
 \epr

\vfill\eject

\printindex

\vfill\eject\addcontentsline{toc}{section}{Contents}

\tableofcontents

\end{document}